\documentclass[10pt,twoside,english,a4paper]{article}

\usepackage{graphicx}  
\usepackage{caption}
\usepackage{subcaption}
\usepackage[dvipsnames]{xcolor}
\usepackage{amsmath}
\usepackage{amsfonts}
\usepackage{amssymb}

\usepackage{subfiles} 
\usepackage{xspace} 

\usepackage[most]{tcolorbox}
\usepackage{amssymb}
\usepackage{mathtools}
\usepackage{stmaryrd}

\usepackage[numbers,sort&compress]{natbib}
\bibpunct{[}{]}{,}{n}{}{;}
\usepackage[]{defjs3}  
\usepackage{xcolor}
\usepackage{booktabs}
\usepackage{placeins}
\usepackage{comment}

\newcommand{\Cc}{\mathbb{C}_{\mathrm{c}}}
\newcommand{\Ce}{\mathbb{C}_{\mathrm{e}}}
\newcommand{\Cmicro}{\mathbb{C}_{\mathrm{micro}}}
\newcommand{\Cmacro}{\mathbb{C}_{\mathrm{macro}}}
\newcommand{\Lc}{L_{\mathrm{c}}}
\newcommand{\Bdis}{{\bP}}
\newcommand{\dis}{{P}}
\newcommand{\curl}{{\textrm{curl}}}
\newcommand{\Curl}{{\textrm{Curl}}}
\newcommand{\IL}{\mathbb{L}}
\newcommand{\Jm}{\mathbb{J}_{\mathrm{m}_1}}
\newcommand{\Jc}{\mathbb{J}_{\mathrm{c}_1}}
\newcommand{\Te}{\mathbb{J}_{\mathrm{m}_2}}
\newcommand{\Tc}{\mathbb{J}_{\mathrm{c}_2}}
\newcommand{\LJ}{\mathbb{L}_{\mathrm{J}}}

\usepackage{comment}
\usepackage{fancyhdr}
\usepackage{authblk}

\title{Static size-effects meet the dynamic scattering properties of finite-sized  mechanical metamaterials: a relaxed micromorphic study with parameter identification via two-stage static-dynamic optimization}

\author{Mohammad Sarhil$^{1,}$\thanks{Corresponding author: mohammad.sarhil@tu-dortmund.de} }
\author{Leonardo Andres Perez Ramirez$^{1}$}
\author{Max Jendrik Voss$^{1}$}
\author{Angela Madeo$^{1}$}

\affil{ $^1$Institute for Structural Mechanics and Dynamics, 
Technical University Dortmund, Dortmund, Germany}

\date{December, 2025}

\begin{document}

\maketitle

\begin{abstract}
Mechanical metamaterials exhibit size-effects when a few unit-cells are subjected to static loading because no clear micro-macro scale separation holds and the  characteristic length of the deformation becomes comparable to the unit-cell size. These size-effects typically manifest themselves as a strengthening of the response in a form summarized as "smaller is stiffer". Moreover, the dynamical behavior of mechanical metamaterials is very remarkable, featuring unique phenomena such as dispersive behavior and  band-gaps where elastic waves cannot propagate over specific frequency ranges. In these frequency ranges, the wavelength becomes gradually comparable to the unit-cell size, giving rise to microstructure related phenomena which become particularly visible in the reflection/transmission patterns where an incident wave hits the metamaterial's interfaces. This raises the question of whether the static size-effects and dynamic reflection/transmission patterns are correlated.

 In this work, we investigate the interaction of the static size-effects and the dynamic scattering response of mechanical metamaterials by employing the relaxed micromorphic model. We introduce a two-stage optimization procedure to identify the material parameters. In the first stage, the static material parameters are identified by exploiting the static size-effects through a least squares fitting procedure based  on the  total energy. The dynamic parameters are determined in the second stage by fitting the dispersion curves of the relaxed micromorphic model to those of the fully discretized microstructure. At this second stage, we assess the results obtained by fitting the dispersion curves in one and in two propagation directions, for both the relaxed micromorphic model (RMM) with curvature and its reduced counterpart (RRMM) without curvature. For a finite-size scattering problem, we examine the RMM formulation, which accounts for the strengthening response associated with size effects, and compare it with the RRMM formulation  at different excitation frequencies. The  results show that, especially at relatively high frequencies when the wavelength becomes small enough to interact with the unit-cell, that RMM with curvature provides a significantly better agreement with the reference microstructured solution than the RRMM. 
 
 Overall, this work demonstrates a strong interaction between the static size-effects and the dynamic scattering response of mechanical metamaterials especially in the dispersive frequency range and in the band-gap range. These results point out that static size-effects are related to complex scattering properties in the dynamical regime. 

\end{abstract}

\section*{keywords} 
mechanical metamaterials, size-effects, dispersion curves,  band-gap, enriched continua, relaxed micromorphic model. 

\section{Introduction} 
\label{introduction}

The advances in additive manufacturing technologies enable a new generation of artificially-made materials with optimized geometry on smaller scales, which can be specifically tailored to unlock exotic mechanical properties and address new applications \cite{BenauRitDalRazBer:2021:acm,AamOzkKazOnu:2023:mmam,XueLuqZhe:2023:ai3}. Mechanical metamaterials are composed of a periodic repetition of unit-cells which exhibit extraordinary properties surpassing classical bulk materials in both static and dynamic behavior \cite{CheMenZhi:2022:ada,RenDasTraNgoXie:2018:ama,SinMuk:2023:pmp,SurGauDu:2019:mma,GupBec:2025:amd,BerVitChrvan:2017:fmm}.   These non-conventional properties are inherited from the complex microstructure rather than the  intrinsic properties of the constituent materials. The unique static properties of these mechanical metamaterials span to cover negative Poisson's ratio (auxetics) \cite{HouSil:2015:mwn,AkaNogMatSatYanYam:2023:tpt}, high stiffness-to-density ratio \cite{BerWadMCM:2017:mma},  negative thermal expansion \cite{JiaWenZhoKaiXuj:2021:scf}, negative-compressibility transitions \cite{QuKKadWeg:2017:pmw,NicMot:2012:mmw}, among many others.  The dynamic behavior of mechanical metamaterials is even more interesting and can play a key role for wave control tools and gadgets for diverse engineering applications, such as sound, vibrations, seismic and shock absorbers \cite{KimMahOhOh:2023:mpa,ArjBarRobVan:2024:amf,NicXihHonFra:2024:aoe,TinRenZehZha:2021:boa,GaoZhaDen:2022:amf,WanXuDua:2023:rob}. These materials exhibit a range of  exotic features such as band-gaps \cite{XiShuBinXiaBoLia:2025:amm}, cloaking \cite{MarMal:2022:moa,GupNik:2025:cpom}, channeling \cite{ZhaXuCheFuYin:2024:dot}, and negative refraction \cite{ZanBiaSanPit:2022:mea} besides their high impact resistance (energy absorption) \cite{CheHeCheLuLuLi:2025:ido}. Nowadays, designing new specific metamaterials with the aid of artificial intelligence is attracting immense interest in the scientific community \cite{QiaKamChe:2025:agt,GuoYanYu:2021:aia}.

Modeling heterogeneous metamaterials is a challenging task due to the substantial computational cost arising from the multicomponent, multimaterial microstructures with highly complex geometries, making full-resolution simulation infeasible at scales that are pertinent for engineering applications. Classical periodic homogenization is widely used as an efficient procedure to reduce the computational cost, see for example \cite{Zoh:2004:hma,Zoh:2005:lna,GeeKouMatYvo:2017:hma,GerKouBre:2010:msc}. Homogenized models derived from variational methods for lattice metamaterials in the continuum limit (many unit-cells) are presented in \cite{UllAriAndOrt:2025:cmi,AriConOrt:2024:hac}.  However, the first-order homogenization procedures are only applicable when a clear micro-macro scale separation is satisfied, allowing the introduction of a representative volume element (RVE), which is well established in the literature, and the results are accurate sufficiently far from the boundaries. Metamaterials reveal size-effects when scale-separation does not hold which is the case when the characteristic length of the deformation does not span a  sufficiently large number of unit-cells \cite{UllAriAnd:2024:fas,AbaVazNew:2022:iom,Yang;2021:seo}. Dealing with the dynamic behavior of metamaterials via classical periodic homogenization is impossible (except in the low-frequency range) since their response carries both nonlocal and local mechanical mechanisms \cite{CheFleSep:2025:nmm}, leading to interesting dynamical phenomena such as band-gaps (or stopband), which are the frequency ranges in which  waves cannot propagate \cite{DalBal:2022:aro}.

Enriched continua provide a natural extension of classical homogenized models either by embedding at the macro-scale higher-order differential operators  such as in second-gradient theories \cite{Min:1964:msi,MinEsh:1968:ofsg,AskAif:2011:gei,Eri:1972:lto} (which include strain-gradient \cite{MinEsh:1968:ofsg,KhaNii:2020:asg,GohCagErd:2022:iid,AufBouBre:2010:sge}, stress-gradient \cite{ForSab:2012:mrc,HutKarFor:2020:kac} and couple stress \cite{AlAbdMah:2015:sdb,SkrEre:2020:ote,EskWanChe:2025:ewp} models) or by introducing additional kinematics independent of the macroscopic translational displacement field, such as in the Mindlin-Eringen micromorphic  \cite{Eri:1967:moc,EriSub:1964:nto,Min:1964:msi} and the Cosserat \cite{ForSab:1998:com,RugHaLak:2019:cel,Lake:2023:cse} theories. These additional macroscopic kinematics account for the deformation at the micro-scale. Second-gradient theories recover size-effects, however, they are unable to reproduce band-gaps since the micro-vibrational modes of the unit-cell must be independent of the macroscopic displacement field. However, the enriched continua with additional kinematics can capture both the size-dependent and micro-inertial phenomena.

In order to model the size-effects of mechanical metamaterials, enriched continua have been employed widely in the literature, see for example \cite{HosNii:2022:3sg,YanTimGiaMue:2023:esg,XiaXuetal:2025:ios,EskshaAkb:2024:use,DosRodFraSer:2024:enm,FinForHohGum:2020:dac,Abi:2019:rtp}.  The size-effect phenomenon is pronounced for bending of metamaterial beams with a few unit-cells for instance \cite{Bar:2025:sei,SarSchSchNef:2023:mts,Sarschschnef:2023:seo,ShoHasOli:2024:nmt,KhaVilJar:2018:msd} but it also appears as well in torsion \cite{DunWhe:2020:sea}, axial \cite{AmeRokPeeGeer:2018:sei,GueEre:2019:ipc} and shear loading \cite{DinHue:2021:iot}. It is worth mentioning that size-effects in metamaterials can  manifest as an increase or a decrease in relative stiffness by reducing the size of the metamaterial (fewer unit-cells) \cite{DunWhe:2020:sea,KirAmsHue:2023:otq} where the geometry of the unit-cell leads to stiffer or more compliant behavior. An increase in stiffness when reducing the size is referred to as a positive size-effect while decreasing stiffness is termed as a negative size-effect.  This is different from the size-effects observed in classical bulk materials where size-effects originate from the defects and arrangement of dislocations, cracks, grains, etc \cite{KeiWenHam:2024:sea}, where "smaller is stiffer" is the main trend because of  the better microstructural control during manufacturing of small specimens which cannot be maintained for larger specimens. This strengthening size-effect can be related as well to the dimensional constraints arising from the testing system \cite{ZhuBusDun:2008;mms}. Modeling the dynamical behavior  of metamaterials via enriched continua has not gained the same attention as the static  size-effects. Modeling the dynamical behavior of materials with band-gap through enriched continua theories has gained less attention than the static size-effects, see for example \cite{ZhoWeiLiTan:2017:moa,WanLiuSohLia:2022:obg,SemJirHor:2024:imm,WanLiu:2023:pbs,MisPoor:2016:gmb,GanDoGod:2025:dho,GanDoGod:2025:dho,WanBinLiuLia:2025:wld,RokAMePeeGee:2019:mch}. Dynamic homogenization models of locally resonant metamaterials into enriched continua accounting for micro-inertia are presented in \cite{SirKouGee:2016:hol,KucGeeKou:2026:ewa,LipDomGeeKou:2025:ecf,LiuSriGeeKou:2021:cho}, and are combined with model-order reduction techniques in \cite{RusKouFaiGee:2024:ard,LipDomGeeKou:2024:aem} and machine learning in \cite{HenMenDosGeeRok:2025:seg}.
 
The identification of static material parameters in enriched continua represents an active field of research in the mechanics community. Three main approaches can be found in literature for the determination of the unknown parameters. The first approach depends on asymptotic expansion procedures \cite{ThbAudLes:2026:asg,fiswagket:2021:mam,GuiSheGui:2024:tsa,AslAnkSeyIsm:2020:cpb,AbaBar:2021:ami,YanAbaTimMue:2020:dom} while the second approach is based on prior definitions of higher-order (mainly quadratic) Dirichlet boundary conditions \cite{ForSab:1998:com,AlaGanSad:2022:cch,LahGodGan:2022:sis,BerDeoGodPicGan:2017:cof,Hue:2017:hoa,Hut:2019:otmm}. The third approach enforces energy equivalence and fits the parameters by error minimization \cite{SheAbaBarBer:2021:iao,SarSchLewSchNef:2024:aca,AlaGanRedSad:2021:com}. However, several issues are noted in literature, such as the non-vanishing higher-order moduli for the homogeneous case and for the case of many unit-cells, although scale-separation holds, the first-order homogenization theory is not recovered. Up to now, no universal procedure has been agreed upon in the scientific community and the identification of the parameters appearing in enriched continua is still an open research topic.

 The relaxed micromorphic model (RMM), used in this work, employs a relaxed curvature in its formulation in terms of the $\Curl$ of a microdistortion field rather than its full gradient  as in the Mindlin-Eringen micromorphic theory \cite{Neff.2014}. The well-posedness of the static case was studied in \cite{Neff.2015}. This second-order microdistortion tensor is non-symmetric in general.  The use of a relaxed curvature leads to many advantages. The main two advantages are the relative simplicity compared to the full micromorphic model where the number of material parameters is notably less than in the full micromorphic theory with the associated fifth- and sixth-order tensors. Indeed, the curvature-related elastic tensor remains fourth-order as in classical elasticity theory.  The second advantage is the behavior of the relaxed micromorphic model as a two-scale micro-macro elasticity model due to the upper bound of the stiffness, which other generalized continua do not exhibit since they predict stiffness going to infinity .  Other advantages have been proved such as well-posedness for the crucial case of the symmetric stress tensor and the need only to give the tangential trace of the microdistortion field on the boundary and not the full field as in the full micromorphic theory. The relaxed micromorphic model has been used to model the dynamical behavior of different metamaterials in \cite{Aivaliotis.2020,Madeo.2015,dAgostino.2020,Demore.2022,Rizzi.2021,VosRizNefMad:MaL:2023} and for the static behavior in \cite{Neff.2020,SarSchSchNef:2023:mts,SarSchSchNef:2023:oti,SarSchLewSchNef:2024:aca}. 
 
 In this work, we investigate the interaction between the static size-effects and the dynamical behavior of mechanical metamaterials.   When the characteristic length of the deformation (e.g. wavelength) is comparable to the size of the unit-cell, which is the case in the band-gap range, the static strengthening response should be considered in the same way as when a few unit-cells are statically deformed (e.g. bending or shearing). For this purpose, the relaxed micromorphic model is employed which provides the framework that we can use to isolate the size-effects by considering or omitting  the curvature that accounts for the strengthened static response. In other words, we can switch between the relaxed micromorphic model (RMM) which considers the static strengthening phenomena and a reduced formulation of the relaxed  micromorphic model (RRMM) that omits the higher-order differential term (curvature) leading to an internal-variable model, which can capture the band-gaps but not the static size-effects.

The paper is organized as follows: In Section \ref{subsec:rmm_2_1} we introduce the relaxed micromorphic model and derive the associated strong forms together with the corresponding Dirichlet and Neumann boundary conditions. We present the modified tensors for the two-dimensional tetragonal case  in Section \ref{sec:tetr}. A brief description of the numerical implementation is given in Section \ref{sec:ni}. A two-stage static-dynamic procedure is proposed to identify the material parameters of the relaxed micromorphic model in Section \ref{sec:tl}. In the first stage presented in Section \ref{sec:se}, the static parameters are identified. To this end, a least-squares minimization of a cost function is presented which fits the total energy of the relaxed micromorphic model to that of the microstructured solution. To accelerate the static optimization procedure, we employ in Section \ref{sec:ann} an artificial neural network as a surrogate model for the FEM calculations. This surrogate provides an initial predictor, which is subsequently used as the starting point for the FEM-coupled optimization. In the second stage introduced in Section \ref{sec:idm}, the dynamic parameters are determined via least-squares fitting of the dispersion curves. A finite-size metamaterial example is used  in Section \ref{chapter:finitesizes} to study the differences of the two formulations with and without curvature terms. Finally, we conclude this work in Section \ref{sec:con}.

  \begin{figure}[h!]
\center
	\unitlength=1mm
	\begin{picture}(170,130)
	\put(0,5){\def\svgwidth{16.5 cm}{\small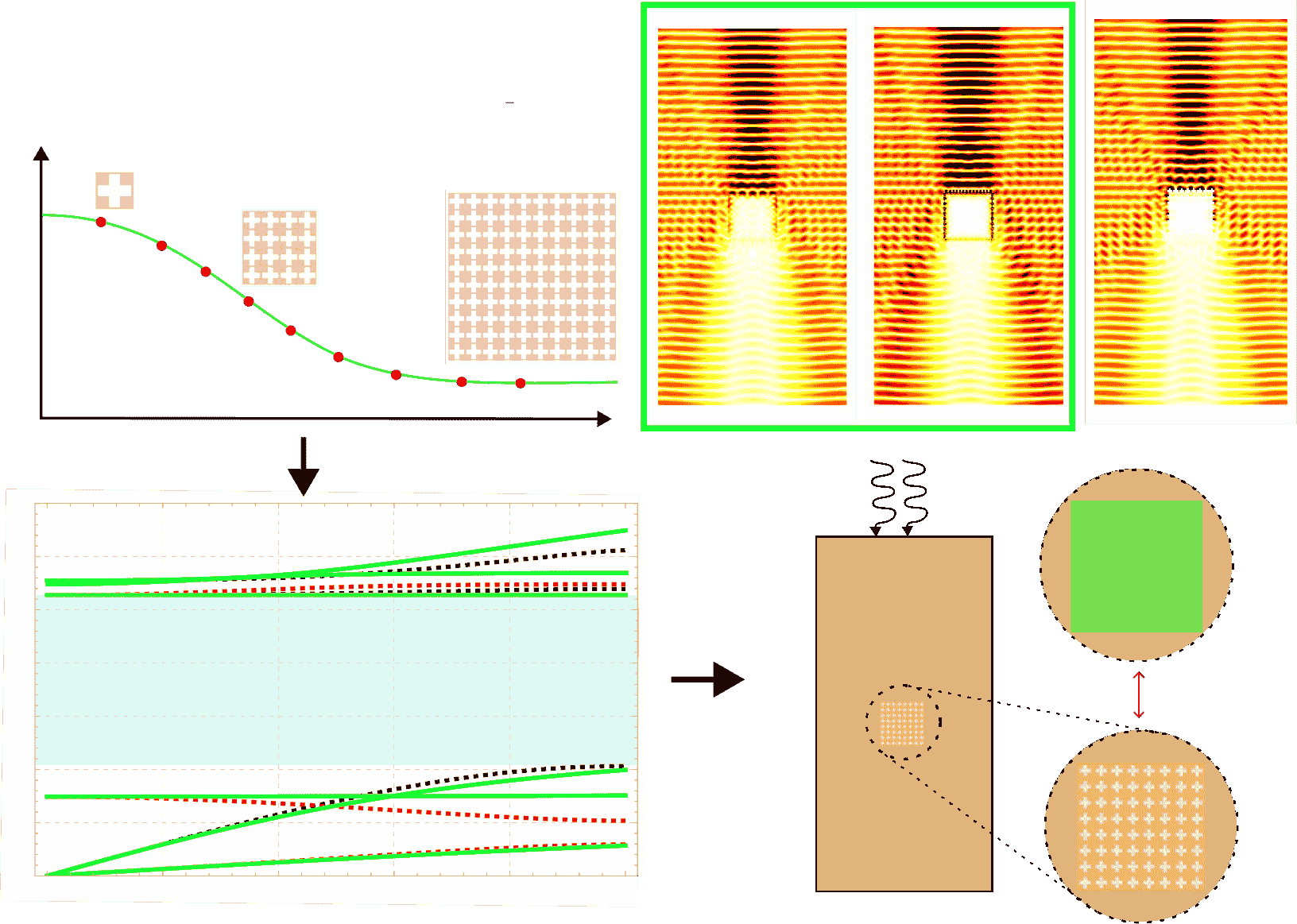}}
	\end{picture}
	\caption{Graphic abstract    }
	\label{Figure:summary}
\end{figure}

\section{The relaxed micromorphic model}
\label{subsec:rmm}

\subsection{Model formulation}
\label{subsec:rmm_2_1}

We consider a body $\Omega$ with a boundary $\partial \Omega$ over the time interval $[0,T]$, subjected to both dynamic and static loads. The relaxed micromorphic model, introduced in \cite{Neff.2014,Madeo.2015,Neff.2015}, is an enriched continuum that characterizes the kinematics of each material point by a displacement field $\bu \in \mathbb{R}^3$ and a  non-symmetric second-order micro-distortion field $\Bdis \in \mathbb{R}^{3\times3}$. These fields are defined by minimizing the action functional $\mathcal{A}$ which reads

\begin{align}
\mathcal{A}:=\iint\limits_{\Omega \times \left[0,T\right]} 
\mathcal{L} \left(\nabla \bu, \dot{\bu}, \nabla \dot{\bu}, \Bdis, \dot{\Bdis}, \Curl \Bdis,  \Curl \dot\Bdis  \right)
\, \textrm{d}x \, \textrm{d}t \, ,
\label{eq:act_func}
\end{align}

with the Lagrangian density given by 

\begin{align}
\mathcal{L}
\coloneqq&
K \left(\dot{\bu},\nabla \dot{\bu}, \dot{\Bdis} ,  \Curl \dot\Bdis \right) 
-
W \left(\nabla \bu, \Bdis, \Curl \Bdis\right)\,.
\label{eq:lagMicrom} 
\end{align}

The strain energy $W$ incorporates the $\Curl$ of the micro-distortion field, rather than  the full gradient used in  the classical Mindlin-Eringen micromorphic theory \cite{Min:1964:msi,EriSub:1964:nto,Eri:1968:mom}, and is given by the form

\begin{align}
W \left(\nabla \bu, \Bdis, \Curl \Bdis  \right)
=& 
\dfrac{1}{2} (   \langle \Ce \, \text{sym}\left(\nabla \bu -  \, \Bdis \right), \text{sym}\left(\nabla \bu -  \, \Bdis \right) \rangle
+  \langle \Cc \, \text{skew}\left(\nabla \bu -  \, \Bdis \right), \text{skew}\left(\nabla \bu -  \, \Bdis \right) \rangle 
\notag 
\\
&
+   \langle \Cmicro \, \text{sym}  \, \Bdis,\text{sym}  \, \Bdis \rangle
+ \langle \, \IL \, \Curl \Bdis , \Curl \Bdis \rangle  ) \, , 
\label{eq:strEneMic} 
\end{align}  

where $\Ce$ and $\Cmicro$ are 4th-order positive-definite elasticity tensors with both major and minor symmetries,   $\Cc$ is a 4th-order positive semi-definite rotational coupling tensor, $\IL$ is a 4th-order tensor acting on the non-symmetric curvature term $(\Curl \Bdis)$ and $\langle \bullet, \bullet \rangle$ indicates the scalar product of two tensors of the same order.

The kinetic energy $K$ in the relaxed micromorphic model is defined as

\begin{align}
K \left( \dot{\bu},\nabla \dot{\bu}, \dot{\Bdis}, \Curl \dot\Bdis \right) 
=&
\dfrac{1}{2} (  \rho \, \langle \dot{\bu},\dot{\bu} \rangle + 
\langle \Jm  \, \text{sym} \, \dot{\Bdis}, \text{sym} \, \dot{\Bdis} \rangle 
+  \langle \Jc \, \text{skew} \, \dot{\Bdis}, \text{skew} \, \dot{\Bdis} \rangle 
\notag
\\
&
+  \langle \Te \, \text{sym}\nabla \dot{\bu}, \text{sym}\nabla \dot{\bu} \rangle
+ \langle \Tc \, \text{skew}\nabla \dot{\bu}, \text{skew}\nabla \dot{\bu} \rangle 
\notag
\\
& 
+ 
\langle \LJ \Curl \dot\Bdis , \Curl \dot\Bdis \rangle  ) 
 \, , 
\label{eq:kinEneMic}
\end{align}

where $\rho$ is the  apparent density, $\Jm$, $\Jc$, $\Te$, $\Tc$ and $\LJ$ are 4th-order micro-inertia tensors.  The body is in equilibrium when the stationarity conditions of the action functional $\mathcal{A}$ with respect to the kinematic fields are satisfied; the first variation vanishes. The variation of the action functional with respect to the displacement, $\partial_\bu \mathcal{A} = 0$, yields  the following strong form of the generalized balance of linear momentum together with the associated Dirichlet and Neumann boundary conditions   

\begin{align}
&  \Div (\widetilde{\Bsigma} + \widehat{\Bsigma} ) - \rho\,\ddot{\bu} = \bzero 
\quad &\mathrm{in} \quad
&\Omega \, , 
\notag 
\\ 
&\bu = \overline{\bu}  
\quad &\mathrm{on} \quad
&\partial \Omega_\bu  \, , 
\notag
\\  
&\bt = ( \widetilde{\Bsigma} + \widehat{\Bsigma} ) \cdot \bn = \overline{\bt}   
\quad &\mathrm{on} \quad
&\partial \Omega_\Bsigma   \, ,
\label{eq:equiMic1}
\end{align}

satisfying $\partial \Omega_\bu \cap \partial \Omega_\Bsigma = \varnothing$ and $\partial \Omega_\bu \cup \partial \Omega_\Bsigma =  \partial \Omega$.  Similarly,  the variation of the action functional with respect to the micro-distortion field, $\partial_\Bdis \mathcal{A} = 0$, leads to the following strong form of the generalized balance of angular momentum together with the associated Dirichlet and Neumann boundary conditions

 \begin{align}
& \widetilde{\Bsigma} -  \overline{\Bsigma}  - \bS - \Curl ( \bbm + \overline{\bbm} ) =  \bzero
\quad  & \mathrm{in} \quad
&\Omega  \, ,
\notag
\\ 
& \Bdis_\bt = {\Bdis} \times \bn = \overline{\Bdis}_\bt 
\quad &\mathrm{on} \quad 
&\partial \Omega_\Bdis \, ,
\notag
\\ 
& \bT_\bbm = {\bbm} \times \bn = \overline{\bT}_{\bbm} 
\quad &\mathrm{on} \quad 
&\partial \Omega_\bbm \, , 
\label{eq:equiMic2}
\end{align} 
 
 satisfying $\partial \Omega_\Bdis \cap \partial \Omega_\bbm = \varnothing$ and $\partial \Omega_\Bdis \cup \partial \Omega_\bbm =  \partial \Omega$. Here, the cross product is performed on row-vectors of the associated second-order tensors. The previously introduced stress measures are defined as 

\begin{align}
& \widetilde{\Bsigma}
\coloneqq \Ce \,\text{sym}(\nabla \bu-\Bdis) + \Cc\,\text{skew}(\nabla \bu-\Bdis)
\, ,
\notag
\\ 
& \widehat{\Bsigma}
\coloneqq \Te \,\text{sym} \, \nabla\ddot{\bu} + \Tc \,\text{skew} \, \nabla\ddot{\bu}
\,,
\notag
\\
& \overline{\Bsigma}
\coloneqq \Jm\,\text{sym} \, \ddot{\Bdis} + \Jc \,\text{skew} \, \ddot{\Bdis}
\,,
\notag
\\
 & \bS
\coloneqq \mathbb{C}_{\rm micro}\, \text{sym} \Bdis
\,,
\notag
\\
 & \bbm
\coloneqq \, \IL \, \Curl  \Bdis 
\,,
\notag
 \\ 
 & \overline{\bbm}
\coloneqq \, \LJ \, \Curl  \ddot\Bdis 
\,. 
\label{eq:equiSigAll}
\end{align}

 \subsection{Material parameters in the two-dimensional case with tetragonal symmetry}
 \label{sec:tetr}

For the two-dimensional case, the micro-distortion field $\Bdis$ and its Curl take the form

\begin{align}
     \Bdis =   \begin{pmatrix}
        \dis_{11} &  \dis_{12} & 0\\ 
        \dis_{21} & \dis_{22} & 0\\ 
        0 & 0 & 0
        \end{pmatrix} 
        = \begin{pmatrix}
        (\Bdis^1)^T \\ 
        (\Bdis^2)^T\\ 
        \bzero
        \end{pmatrix} \,, 
       \quad 
    \Curl \Bdis =   \begin{pmatrix}
        0 & 0 & \curl^{2D} \Bdis^1\\ 
        0 & 0 & \curl^{2D} \Bdis^2\\ 
        0 & 0 & 0
        \end{pmatrix}\,, 
       \quad        
\end{align}
with $\curl^\textrm{2D} \Bdis^i= \dis_{i2,1}-\dis_{i1,2}$. Here, $\Bdis^i$ are the row-vectors of the tensor $\Bdis$. The reduced curvature measure, which is isotropic in 2D, can be written in a vector form as $\curl^\textrm{2D} \Bdis := (\curl^\textrm{2D} \Bdis^1, \curl^\textrm{2D} \Bdis^2)^T$ so that the curvature term simplifies to $\IL  \Curl \Bdis \rightarrow \mu \Lc^2 \, \Curl^\textrm{2D} \Bdis$ where $\mu > 0$ is a typical shear modulus and $\Lc > 0$ is the characteristic length parameter.  This applies as well to $\Curl \ddot \Bdis$ leading to $\LJ  \Curl \ddot \Bdis \rightarrow  \rho \, M  \Lc^4  \Curl^\textrm{2D} \ddot \Bdis$

  In this work, we consider a unit-cell with tetragonal symmetry. The tensors, which are acting on symmetric arguments, are written in Voigt notation as follows

\begin{align}
        \mathbb{C}_{\bullet}
        &= 
        \begin{pmatrix}
        \lambda_{\bullet} + 2 \mu_{\bullet}	& \lambda_{\bullet}  	& 0\\ 
        \lambda_{\bullet} 	& \lambda_{\bullet} + 2 \mu_{\bullet} & 0\\
        0 & 0 & \mu_{\bullet}^{*}
        \end{pmatrix}\,, \qquad &\bullet = {\textrm{e},\textrm{micro}}
        &
        \notag
        \\[2mm]
        \mathbb{J}_\star
        &=
        \rho \Lc^2
        \begin{pmatrix}
        \Lambda_\star + 2 M_\star & \Lambda_\star   & 0\\ 
        \Lambda_\star  &  \Lambda_\star + 2 M_\star  & 0\\ 
        0 & 0 &  M^{*}_\star\\ 
        \end{pmatrix}\,, \qquad &\star = {m_1,m_2} \,.
        &
        \label{eq:tensors}
\end{align}

The tensors $\Cc$, $\Jc$, and $\Tc$ are acting on skew-symmetric arguments. In the 2D case their contributions reduce to a single scalar parameter each; $\mu_\textrm{c}$, $\rho \Lc^2 M_{\textrm{c}1}$, and $\rho \Lc^2 M_{\textrm{c}_2}$,  where

\begin{align}
\IC \skew \bA = 2 \IC_{1212} \skew \bA  \quad  \textrm{for}  \quad \bA \in \mathbb{R}^{2\times2} \quad \textrm{and}  \quad \IC \in \mathbb{R}^{2\times2\times2\times2} \quad  \textrm{with} \quad  \IC_{ijkl} = - \IC_{ijlk} \,. 
\end{align}

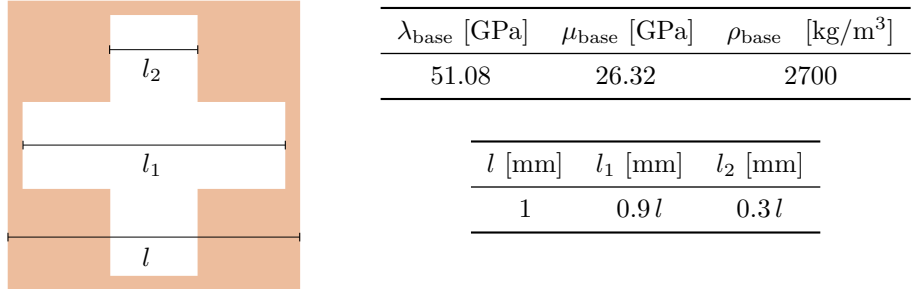
\begin{figure}[h!]
    \centering
    \begin{minipage}{.3\textwidth}
        \centering
       	\begin{picture}(3.5,4.4)
      \put(0,-1.8 cm){\includegraphics[width=4cm]{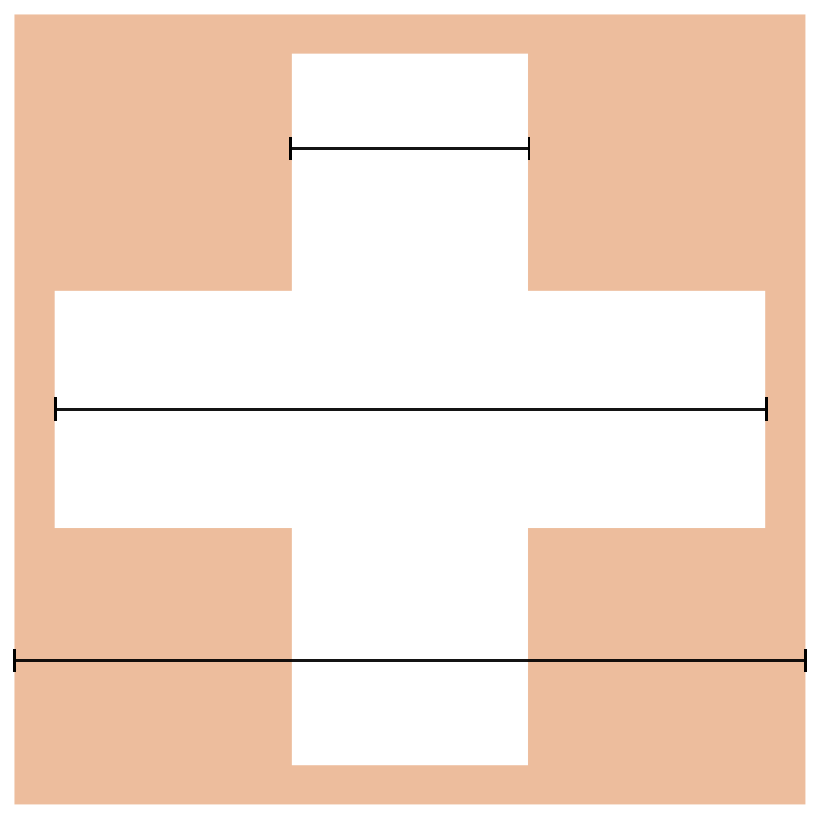}} 
      \put(1.838053286723413 cm,-0.150023909705 cm){\color[rgb]{0,0,0}\makebox(0,0)[lb]{\smash{$l_1$}}}%
      \put(1.838053286723413 cm,-1.40023909705 cm){\color[rgb]{0,0,0}\makebox(0,0)[lb]{\smash{$l$}}}%
       \put(1.838053286723413 cm,1.10990023909705 cm){\color[rgb]{0,0,0}\makebox(0,0)[lb]{\smash{$l_2$}}}%
        \end{picture}      
    \end{minipage}
    \hfill
    \begin{minipage}{.65\textwidth}
        \centering
        \begin{tabular}{ccc}
            \toprule
          $\lambda_\textrm{base}$ [GPa] & $\mu_\textrm{base}$ [GPa] & $\rho_\textrm{base}$ \, $[\textrm{kg}/\textrm{m}^3]$ \\
            \midrule
             $51.08 \,  $ & $26.32 \,  $ & $2700$ \\
            \bottomrule
        \end{tabular}
        
        \vspace{0.5 cm}
        \begin{tabular}{ccc}
            \toprule
            $l$ [mm] & $l_1$ [mm] & $l_2$ [mm]\\
            \midrule
           $1$  & $0.9 \, l$ & $0.3 \, l$ \\
            \bottomrule
            \vspace{0.5 cm}
        \end{tabular}
    \end{minipage}
            \caption{Unit-cell with the material and geometrical parameters.}
            \label{Fig;unit_cell}
\end{figure}

In Figure \ref{Fig;unit_cell}, we present the assumed unit-cell with the parameters of the base material (aluminum).

\FloatBarrier

\subsection{Numerical implementation} 
\label{sec:ni}

\begin{figure}[h!]
\center
	\unitlength=1mm
	\begin{picture}(50,40)
	\put(0,2){\def\svgwidth{4.5 cm}{\small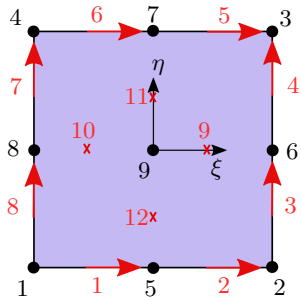}}
	\end{picture}
	\caption{The Q2NQ2 finite element. Black dots indicate the displacement degrees of freedom and the red arrows indicate  vectorial degrees of freedom representing the tangential projection of the microdistortion field on the element's vertices, see \cite{SchSarSchNef:2022:lhb}.}
	\label{Figure:T2NT2}
\end{figure}

The formulation of conforming finite elements for the relaxed micromorphic model was introduced in \cite{SchSarSchNef:2022:lhb,SarSchNefSch:2021:oat}. In this work, we employ the Q2NQ2 finite element shown in Figure \ref{Figure:T2NT2}. The Q2NQ2 discretizes the displacement using a second-order Lagrangian shape ansatz. The microdistortion field is discretized using the vectorial $H(\curl)$-conforming Nedelec formulation of second-order and first-type.  

The static optimization is performed in AceGen and AceFEM packages within the Mathematica software environment. For the dynamical parameter optimization, the analytic dispersion relations of the RMM are coded in Mathematica. The dynamical FEM simulations (frequency domain analysis) are performed on COMSOL Multiphysics which allows user-implemented weak formulations and offers an $H(\curl)$-conforming formulation. For the fully microstructured calculations, a classical T2 element with second-order Lagrangian formulation is used.

\section{Parameter identification; two levels static-dynamic optimization} 
\label{sec:tl}

We separate the identification procedure into two levels: a static stage and a dynamical stage.   The static parameters appearing in the strain energy $W$  will be determined based on static loading tests on a cluster of unit-cells where size-effects are pronounced. The dynamic parameters appearing in the kinetic energy $K$ will be identified from fitting  the dispersion curves of the relaxed micromorphic model with the ones obtained through Bloch–Floquet analysis of an infinite block of unit-cells. In this stage, the static parameters obtained from the static optimization are kept fixed.  This two-stage procedure allows to ensure uniqueness of the parameters' values, which is not the case when performing a purely dynamic identification.

\subsection{Identification of the static material parameters; size-effects} 
\label{sec:se}

Size-effects are a widely spread phenomenon in materials science.
They arise when the size of the microstructural components becomes comparable to the overall macroscopic specimen, preventing a clear scale separation. 
These effects are particularly pronounced in finite-sized metamaterials, where specimens may contain only a limited number of unit-cells. For this section, only static loads are considered with a very slow deformation rate so that dynamical effects are, in a first instance, excluded.  

 For the limiting case  when $\Lc \rightarrow \infty$, the micro-distortion field has to be a gradient field. With the so-called consistent boundary condition $\Bdis \cdot \Btau = \nabla \bu \cdot \Btau$  ($\Btau$ is the tangential vector on the boundary), the variational formulation leads to 

\begin{equation}
\Div \left( \Cmicro \sym \nabla \bu \right) = \bzero \,. 
\end{equation}

For the limiting case when the characteristic length $\Lc \rightarrow 0$, the generalized balances of linear and angular momentum are reduced to 

\begin{equation}
\label{eq:ce_relation}
\Div \left( \Cmacro \sym \nabla \bu \right) = \bzero \quad \textrm{with} \quad  \Cmacro := \Cmicro  (\Ce+\Cmicro)^{-1} \Ce \, 
\end{equation}

independently of the choice of the boundary conditions. Notice that the elasticity tensor $\Cmacro = (\Ce^{-1}+\Cmicro^{-1})^{-1}$ is a series sum of $\Ce$ and $\Cmicro$ and therefore it is softer than both.  

\begin{figure}[htpb]
\center
	\unitlength=1mm
	\begin{picture}(185,80)
	\put(0,2){\def\svgwidth{16.5 cm}{\small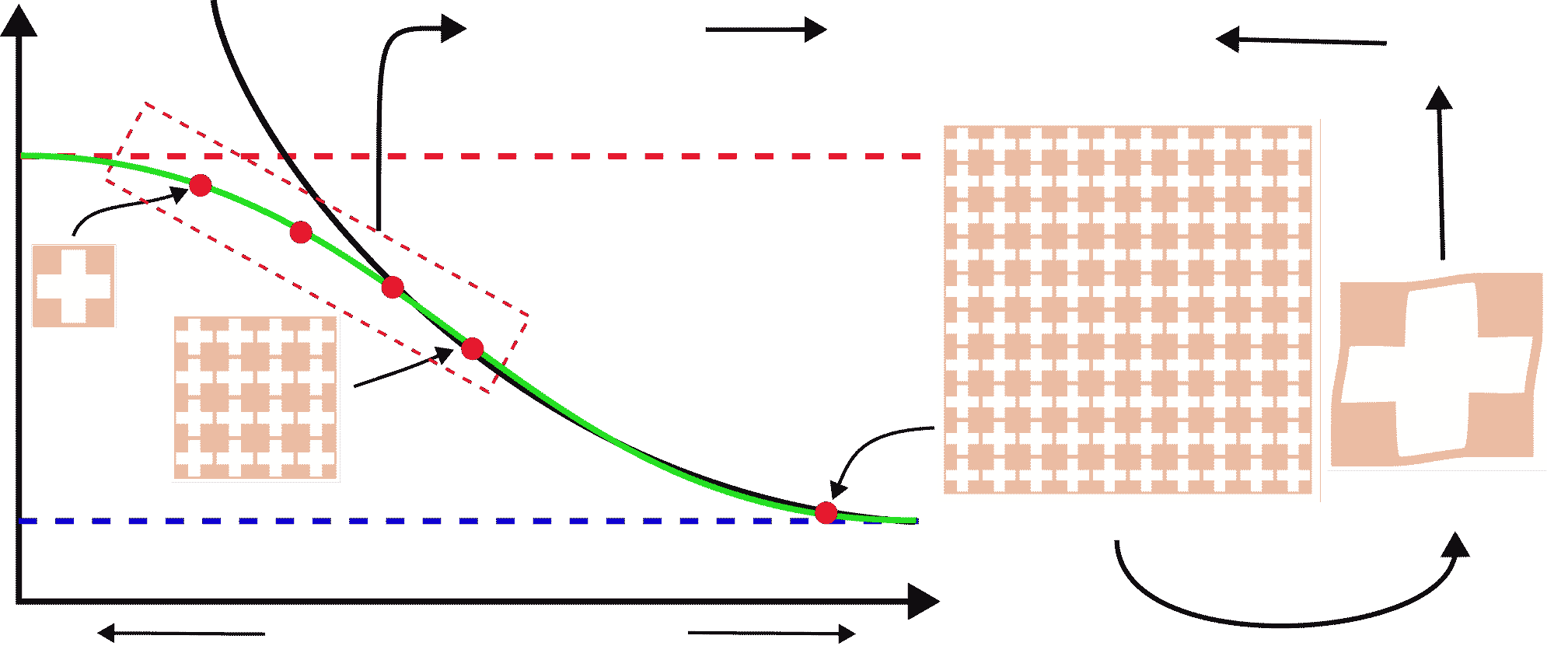}}
	\end{picture}
	\caption{The stiffness of the relaxed micromorphic model (RMM) is bounded from above and below.  Other generalized continua exhibit unbounded stiffness for small sizes.    }
	\label{Figure:static}
\end{figure}

The unique behavior of the relaxed micromorphic model in the static case is illustrated in Figure \ref{Figure:static}. For the limiting cases when the characteristic length $\Lc$ goes to $0$ and $\infty$, linear elasticity with elasticity tensors $\Cmacro$ and $\Cmicro$ are recovered, respectively. They bound the relaxed micromorphic model from below and above, differently from other enriched models which are typically only bounded from below. Thus, the relaxed micromorphic model behaves as a two-scale linear elasticity model and no other generalized continua exhibit such unique property.  The soft bound  represented by linear elasticity with $\Cmacro$ corresponds to the effective response of a large block of metamaterials with many unit-cells.  The stiff bound represented by linear elasticity with elasticity tensor $\Cmicro$ corresponds to the stiff response of small part of the material (e.g. a single unit-cell,  the homogeneous matrix or even stiffer).    

\subsubsection{Identification of $\Cmacro$} 

$\,$

The lower bound of the micromorphic model is given by linear elasticity with the effective elasticity tensor $\Cmacro$ in the limiting case when  $\Lc \rightarrow 0$. This corresponds to the soft response of an infinitely large block of unit-cells far away from the boundary where size-effects vanish.  The macroscopic response  is obtained by classical first-order homogenization theory built by enforcing periodic boundary conditions on a periodic representative volume element (here a unit-cell).  Note that the elasticity tensor $\Cmacro$ does not appear explicitly in the formulation of the RMM but is embedded in the identification of the elasticity tensor $\Ce$ through Equation (\ref{eq:ce_relation}).  The resulting Lam\'e coefficients of the effective elasticity tensor $\Cmacro$ read $\lambda_\textrm{macro} = 5.9$ [GPa], $\mu_\textrm{macro} = 0.627$ [GPa] and $ \mu^*_\textrm{macro} = 1.748$ [GPa].  These coefficients can also be equivalently identified by the slope of the acoustic dispersion curves near the origin in the low-frequency regime instead than by classical first-order homogenization.   

 \subsubsection{Identification of $\Ce$} 

$\,$ 

The elasticity tensor $\Ce$ has no direct physical interpretation like $\Cmicro$ and $\Cmacro$.  The identification of this tensor will follow through the relation 

\begin{equation}
\label{eq:reuss_like}
\Ce = \Cmicro   (\Cmicro-\Cmacro)^{-1}   \Cmacro\,, 
\end{equation}

which connects the tensors $\Ce$, $\Cmacro$ and $\Cmicro$. Thus, two of them are needed to identify the third  one.  Note that $\Cmicro$ has to be stiffer than $\Cmacro$ to keep $\Ce$ positive definite.   For tetragonal symmetry, the previous  relation turns into

\begin{equation}
\mu_\textrm{e} = \frac{\mu_\textrm{micro} \, \mu_\textrm{macro}}{\mu_\textrm{micro} - \mu_\textrm{macro}} , \quad \mu_\textrm{e}^* = \frac{\mu_\textrm{micro}^* \, \mu_\textrm{macro}^*}{\mu_\textrm{micro}^* - \mu_\textrm{macro}^*} , \quad \lambda_\textrm{e} + \mu_\textrm{e} = \frac{ (\lambda_\textrm{micro} + \mu_\textrm{micro}) \,  (\lambda_\textrm{macro} + \mu_\textrm{macro}) }{(\lambda_\textrm{micro} + \mu_\textrm{micro}) -  (\lambda_\textrm{macro} + \mu_\textrm{macro}) } \,. 
\end{equation}

These three correlations will be embedded in the optimization procedure to identify the tensor $\Cmicro$.

\subsubsection{Identification of $\Cmicro$, $\mu_\textrm{c}$ and $\mu \Lc^2$ \\} 
\label{sec:io}

We present an optimization algorithm to identify the rest of the unknown static parameters. The optimization is based on least-squares fitting of total energy following the concept presented in \cite{SarSchLewSchNef:2024:aca}. Since the primary focus of this paper is on modeling wave propagation for frequencies within and below the band gap (i.e., in the acoustic range), first-order affine deformation modes are considered to identify the remaining static parameters. Three affine deformation modes are considered which are pure shear,  pure deviatoric normal strain and equal biaxial strain. Any other linear first-order deformation mode can be expressed as a linear combination of these three fundamental modes.  

A cost function is introduced  as the sum of the squared differences between the total strain energy of the relaxed micromorphic model and that of the fully detailed microstructure with linear elasticity. It reads  

\begin{equation}
\label{Eq:w_minimization} 
r^2 = \min_{\Cmicro, \mu_c,  \, \mu \Lc^2}  \sum_{n}  \sum_{i} \lvert \lvert  \Pi^\textrm{het}_{\{i , n\}} (\bu, \IC ) -\Pi_{\{i , n\}} (\bu,\Bdis,\Cmicro, \mu_c,  \, \mu \Lc^2)  \rvert \rvert^2,     
\end{equation} 

 with the following energy definitions 
 
\begin{align}
\Pi^\textrm{het}_{\{i , n\}} (\bu, \mathbb{C}  )   =&  \int_{\B}  \Bvarepsilon : \mathbb{C} (\bx) :\Bvarepsilon \, \textrm{d} V \,, \notag  \\ 
   \Pi_{\{i , n\}} (\bu, \Bdis, \Cmicro, \mu_c,  \, \mu \Lc^2)  =&  \int_{\B}  W (\bu, \nabla \bu,  \Bdis, \Curl \Bdis, \Cmicro, \mu_c,  \, \mu \Lc^2)  \, \textrm{d} V \,,  
\end{align}   
 
where $\B =  [-\dfrac{n \, l}{2},\dfrac{n \, l}{2}] \times  [-\dfrac{n \, l}{2},\dfrac{n \, l}{2}]$, $n$ is the number of the unit-cells in the specimen and $i$ is the loading case following  $\bu = \overline{\bu}_i = \overline{\Bvarepsilon}_i \cdot \bx$ and $\Bdis \cdot \Btau = \nabla \overline{\bu}_i \cdot \Btau$ on  $\partial\B = \partial\B_u = \partial \B_\dis $ with $i=1,2,3$ and $\Btau$ is the tangential vector on the boundary $\partial \B$ . The consistent boundary condition is enforced on the whole boundary.  It imposes that the projection of the microdistortion onto the tangential plane of the boundary equals  the projection of the displacement gradient  onto the  tangential plane of the boundary.  The consistent boundary condition is required to recover linear elasticity with elasticity tensor $\Cmicro$ as rigorous upper bound. Imposing alternative boundary conditions on the microdistortion field $\Bdis$ preserves the size-effects in the relaxed micromorphic model but with a different upper bound.  A gradient descent algorithm is proposed to minimize the cost function.  To facilitate the optimization, we reformulate the minimization problem in Equation (\ref{Eq:w_minimization})  in terms of the increments of the unknowns 

\begin{equation}
\label{Eq:error}
\begin{aligned}
r^2=& \min_{
\Delta \Cmicro,
\Delta \mu_c,
\Delta   \mu \Lc^2  }  
  \sum_{n} \sum_{i}  \, \Big|\Big|  \Pi^\textrm{het}_{\{i , n\}} -  
 \Big( \Pi_{\{i , n\}}  +     \frac{\partial \Pi_{\{i , n\}}}{\partial \mu_\textrm{micro}}   \Delta \mu_\textrm{micro}  
  \\ &
   +   \frac{\partial \Pi_{\{i , n\}}}{\partial \mu_\textrm{micro}^*}   \Delta \mu_\textrm{micro}^* +  \frac{\partial \Pi_{\{i , n\}}}{\partial \lambda_\textrm{micro}}  \Delta \lambda_\textrm{micro}    + 
 \frac{\partial \Pi_{\{i , n\}}}{\partial \mu_c}  \Delta \mu_c  +  \frac{\partial \Pi_{\{i , n\}}}{\partial \mu \Lc^2}   \Delta   \mu \Lc^2  
  \Big)  \Big|\Big|^2 \,. 
  \end{aligned}
\end{equation}

 The incremental solution $\Delta\BLambda $ obtained from the least-squares minimization in Equation (\ref{Eq:error}) reads 

\begin{equation}
\Delta\BLambda = ( \bD^T \cdot \bD )^{-1} \cdot \bD^T  \cdot (\ba^\textrm{het} - \ba^\textrm{RMM})    \,, 
\end{equation}

with 

\begin{equation}
\ba^\textrm{het} = \left[\begin{array}{c}
\Pi^\textrm{het}_{\{1 , 1\}}  \\ 
\Pi^\textrm{het}_{\{1 , 2\}}  \\ 
\vdots \\
\Pi^\textrm{het}_{\{i_\textrm{max} , n_\textrm{max}\}}  \\ 
\end{array}
\right] \,, \quad 
\ba^\textrm{RMM}  = \left[\begin{array}{c}
\Pi_{\{1 , 1\}}  \\ 
\Pi_{\{1 , 2\}}  \\ 
\vdots \\
\Pi_{\{i_\textrm{max} , n_\textrm{max}\}}  \\ 
\end{array}
\right] \,, \quad  
\Delta \BLambda = 
\left[
\begin{array}{c}
\Delta \mu_\textrm{micro} \\
\Delta \mu^*_\textrm{micro} \\
\Delta \lambda_\textrm{micro} \\
\Delta \mu_c \\
\Delta \mu \Lc^2  
\end{array}
\right] \,, 
\end{equation}
and the partial derivative matrix is
\begin{equation}
\bD = 
\left[
\begin{array}{c c c c c}
  \frac{\partial \Pi_{\{1 , 1\}}}{\partial \mu_\textrm{micro}}  &   \frac{\partial \Pi_{\{1 , 1\}}}{\partial \mu_\textrm{micro}^*} & \frac{\partial \Pi_{\{1 , 1\}}}{\partial \lambda_\textrm{micro}} &  \frac{\partial \Pi_{\{1 , 1\}}}{\partial \mu_c}   &  \frac{\partial \Pi_{\{1 , 1\}}}{\partial \mu  \Lc^2}  \\[1 em] 
  \frac{\partial \Pi_{\{1 , 2\}}}{\partial \mu_\textrm{micro}}  &   \frac{\partial \Pi_{\{1 , 2\}}}{\partial \mu_\textrm{micro}^*} & \frac{\partial \Pi_{\{1 , 2\}}}{\partial \lambda_\textrm{micro}} &  \frac{\partial \Pi_{\{1 , 2\}}}{\partial \mu_c} &  \frac{\partial \Pi_{\{1 , 2\}}}{\partial \mu  \Lc^2}  \\
\vdots & \vdots & \vdots & \vdots & \vdots \\
  \frac{\partial \Pi_{\{i_\textrm{max} , n_\textrm{max}\}}}{\partial \mu_\textrm{micro}}  &   \frac{\partial \Pi_{\{i_\textrm{max} , n_\textrm{max}\}}}{\partial \mu_\textrm{micro}^*} & \frac{\partial \Pi_{\{i_\textrm{max} , n_\textrm{max}\}}}{\partial \lambda_\textrm{micro}} &  \frac{\partial \Pi_{\{i_\textrm{max} , n_\textrm{max}\}}}{\partial \mu_c} & \frac{\partial \Pi_{\{i_\textrm{max} , n_\textrm{max}}\}}{\partial \mu  \Lc^2} \\
\end{array}
\right] \,. 
\end{equation}

Moreover, in order to guarantee the physical admissibility,  additional constraints have to be satisfied which ensure no negative energy terms in the strain energy. The tensor $\Ce$ must be positive definite which requires that $ \langle \Cmicro \bx, \bx \rangle > \langle \Cmacro \bx, \bx \rangle , \forall \bx \in \mathbb{R}^3 $ while the scalar quantity $\mu \Lc^2$ must be strictly positive and $\mu_c$ non-negative. These side constraints lead to 

\begin{equation}
\label{Eq:inequlity}
\begin{aligned}
 \quad  (\mu_\textrm{micro})_\textrm{old} + \Delta \mu_\textrm{micro}   &> \mu_\textrm{macro}  \,, \\
  \quad (\mu^*_\textrm{micro})_\textrm{old} +  \Delta \mu^*_\textrm{micro} &> \mu^*_\textrm{macro}\,, \\
 (\lambda_\textrm{micro})_\textrm{old} + \Delta \lambda_\textrm{micro}  +  (\mu_\textrm{micro})_\textrm{old} +  \Delta \mu_\textrm{micro} &> \lambda_\textrm{macro} + \mu_\textrm{macro}  \,, \\  
   (\mu_c)_\textrm{old} +    \Delta \mu_c & \ge  0 \,.  \\
   (\mu \Lc^2)_\textrm{old} +    \Delta  \mu \Lc^2  & >  0 \,.
 \end{aligned}
\end{equation}
 
In order to guarantee not breaking the previous  constraints, we consider a line search procedure where the step size $\beta$ ensures that we stay in the allowed range of the parameters such that 

\begin{equation}
\left[ \begin{array}{c}
\mu_\textrm{micro} \\ 
\mu^*_\textrm{micro} \\
\lambda_\textrm{micro} \\  
\mu_c \\
\mu  \Lc^2
 \end{array} \right]_\textrm{new} =  \left[ \begin{array}{c}
\mu_\textrm{micro} \\ 
\mu^*_\textrm{micro} \\ 
\lambda_\textrm{micro} \\ 
\mu_c \\
\mu  \Lc^2
 \end{array} \right] + \beta
 \left[
\begin{array}{c}
\Delta \mu_\textrm{micro} \\
\Delta \mu^*_\textrm{micro} \\
\Delta \lambda_\textrm{micro} \\
\Delta \mu_c \\
\Delta \mu \Lc^2  
\end{array}
\right] \,,  \quad \beta \in [0, \beta_\textrm{max}] \, .  
\end{equation}

The maximum step size is defined by solving the inequalities in Equation (\ref{Eq:inequlity}) 

\begin{equation}
\beta_\textrm{max} = \min 
\begin{cases}
 1, \\ 
  h(-\dfrac{(\mu_\textrm{micro})_\textrm{old} - \mu_\textrm{macro}}{\Delta \mu_\textrm{micro}}), \\ 
  h(-\dfrac{(\mu_\textrm{micro}^*)_\textrm{old} - \mu_\textrm{macro}^*}{\Delta \mu_\textrm{micro}^*}), \\ 
h( -  \dfrac{(\lambda_\textrm{micro} + \mu_\textrm{micro})_\textrm{old} - (\lambda_\textrm{macro} + \mu_\textrm{macro})}{\Delta \lambda_\textrm{micro} + \Delta \mu_\textrm{micro}}) \\  
   h(-\dfrac{(\mu_c)_\textrm{old}}{\Delta \mu_c}), \\
   h(-\dfrac{(\mu \Lc^2)_\textrm{old}}{\Delta \mu \Lc^2}) \,,
\end{cases}
\end{equation}

with the definition of the function $h$ as 

\[
h(x) = 
\begin{cases}
x, & x > 0, \\
1, & x \le 0.
\end{cases}
\]

Once the maximum step size is known, we perform a line search  procedure along the direction $\Delta \BLambda$ to determine an optimal step size, i.e. 

\begin{equation}
\begin{aligned}
r^2 = \min_{\beta}  \sum_{n=1}  \sum_{i=1} \lvert \lvert  \Pi^\textrm{het}_{\{i , n\}} -\Pi_{\{i , n\}} ( & \mu_\textrm{micro} +  \beta \, \Delta\Lambda_1 ,  \mu^*_\textrm{micro}   + \beta  \, \Delta\Lambda_2, \\ & \lambda_\textrm{micro} + \beta \, \Delta\Lambda_3, \mu_c + \beta \, \Delta\Lambda_4  , \mu \Lc^2 +   \beta \, \Delta\Lambda_5)  \rvert \rvert^2 \,, 
\end{aligned}
\end{equation} 
 
and the iterative procedure is repeated until a converged  set of parameters is reached. Thus, all static parameters in  the strain energy are defined after implementing this procedure. When a parameter is at the boundary of the admissible range and the update direction would drive it outside its allowed range, the maximum step size for that parameter becomes zero. In this case, this parameter is kept constant from the last iteration and excluded from the update in the current line-search iteration, but it is included again in the optimization in the next iteration.   Because the proposed procedure is based on gradient descent, different initial parameters may lead to different converged final parameters associated with a local minimum rather than the global one. Moreover, the optimization is relatively expensive (few hours) because it is coupled to a finite element code which makes it impractical to explore a large number of initial starts.   To address this limitation, we employ an artificial neural network (ANN) as a surrogate model to the finite element simulation in a first stage to predict the initial parameters.

\subsubsection{An artificial neural network as a first predictor \\} 
\label{sec:ann}

\begin{figure}[h!]
\center
	\unitlength=1mm
	\begin{picture}(180,120)
	\put(10,2){\def\svgwidth{14 cm}{\small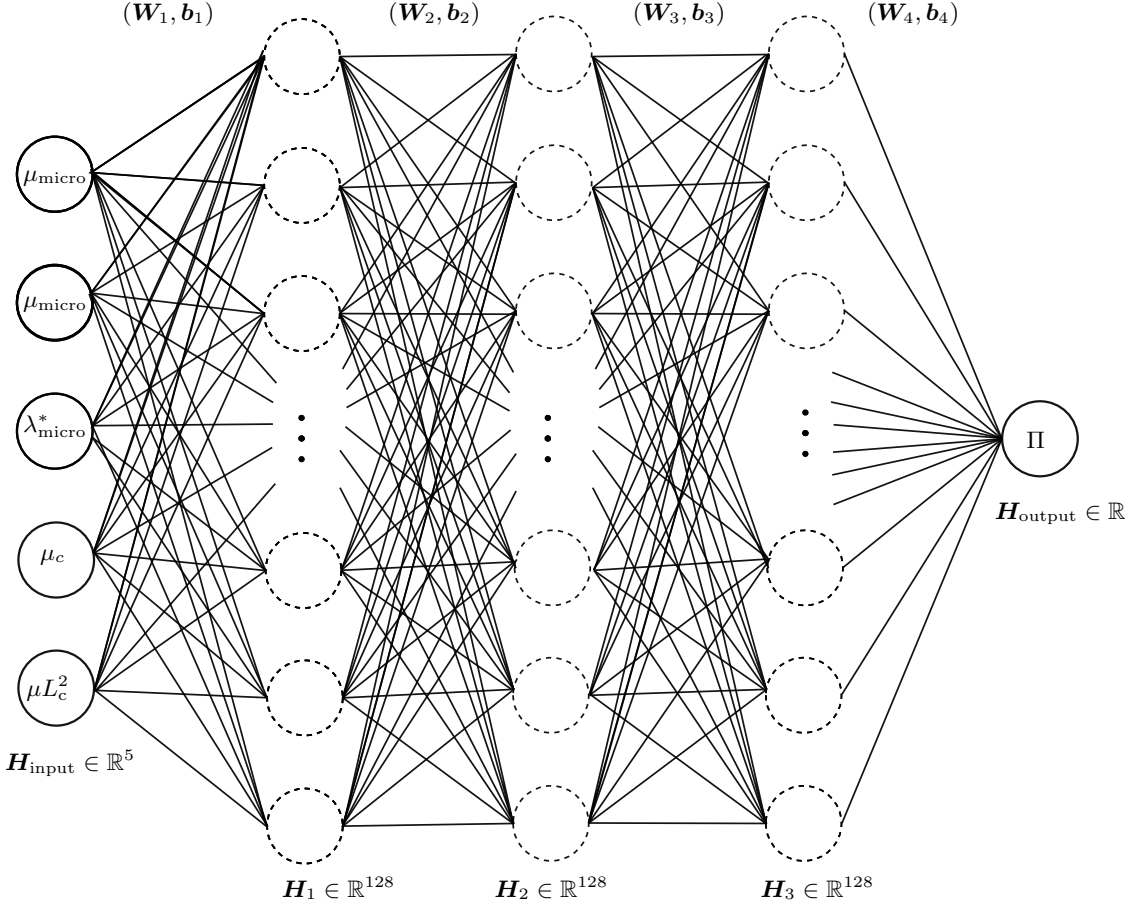}}
	\end{picture}
	\caption{The feed-forward fully-connected neural network as a surrogate model for RMM. The input layer has five features which are the unknown static parameters. Three hidden layers are used, each has 128 neurons and the output is the total energy.  }
	\label{Figure:ann}
\end{figure}

 A feed-forward fully connected neural network is constructed  with three hidden layers shown in Figure \ref{Figure:ann}. We employ a rectified linear unit (ReLU) as an activation function to introduce nonlinearity and enhance learning capability.  Each layer $\bH_{i+1}$ is obtained from the previous layer $\bH_{i}$ as  

\begin{equation}
\bH_{i+1} = \textrm{ReLU}(\bW_{i+1} \bH_{i}  + \bb_{i+1})  \quad \textrm{with} \quad  \textrm{ReLU}(x) = \max{(0,x)} \, , 
\end{equation}
 
where $\bW_{i+1}$ and $\bb_{i+1}$ are the weight matrix and the bias vector for each respective layer $\bH_{i+1}$.
 
The training dataset was generated using the finite element simulation by varying the unknown material parameters within a reasonable and physically motivated range defined as 
$\mu_\textrm{micro} \in [\mu_\textrm{macro},\mu_\textrm{matrix}]$ ,
$\mu_\textrm{micro}^* \in [\mu_\textrm{macro}^*, \mu_\textrm{matrix}^*]$ ,
$\lambda_\textrm{micro} \in [\lambda_\textrm{macro}, \lambda_\textrm{matrix}]$ ,
$\mu_c \in [0,\mu_\textrm{macro}]$ and 
$\mu \Lc^2 \in [0.1 [\textrm{N}],5 \mu_\textrm{macro} l^2]$. 
 Each parameter is divided into equally spaced intervals between its assumed minimum and maximum limits, and all resulting parameter combinations are employed to generate the dataset used for training the surrogate model. However, we observed during the subsequent optimization that one parameter which is $\mu_\textrm{micro}^*$ tends to show higher values than initially  allowed  and therefore we expanded its range to   $\mu_\textrm{micro}^* \in [\mu_\textrm{macro}^*, 10 \mu_\textrm{matrix}^*]$. In total $1,338,444$ datasets were generated with a coarse $5 \times 5$ mesh and used to train the model.  The surrogate neural network achieved a coefficient of determination of $0.999996$ and the mean relative error of less than $1\%$ demonstrating high accuracy.
 
Subsequently, we applied a similar gradient descent procedure to minimize the cost function—defined analogously to that introduced earlier—for the three loading cases and two specimen sizes $n=1,2$. A total of $10^5$ different initial parameter sets were used and the resulting solution with the least cost was chosen as a starting point for the optimization procedure coupled with the finite element code. This prediction reads 

\begin{align*}
\lambda_\textrm{micro}^0 &= 10.36 [\textrm{GPa}]\,,  &\mu_\textrm{micro}^0 &= 7.34[\textrm{GPa}],   \quad \qquad \mu_\textrm{micro}^{*,0} = 58.46[\textrm{GPa}],  \notag \\ 
\mu_c^0 &= 0[\textrm{GPa}], \quad  &(\mu \Lc^2)^0 &= 1256 [\textrm{N}]\,. 
\end{align*}

The production of the training dataset took $40$ hours and training the ANN takes two hours. The search for the best fit using $10^5$ different starting points requires $3.5$ hours. On the other hand, one optimization directly coupled with the finite element calculation takes between an hour and multiple hours depending on the number of required iterations till convergence. Therefore, the effectiveness of this approach is confirmed. The calculations are performed on AMD Ryzen Threadripper PRO 5995WX with 64 cores.    

\subsubsection{Results of the static optimization \\} 

\begin{figure}[h!]
\center
	\unitlength=1mm
    \begin{subfigure}{0.49\textwidth}
	\begin{picture}(95,60)
	\put(0,0){\def\svgwidth{7.5 cm}{\small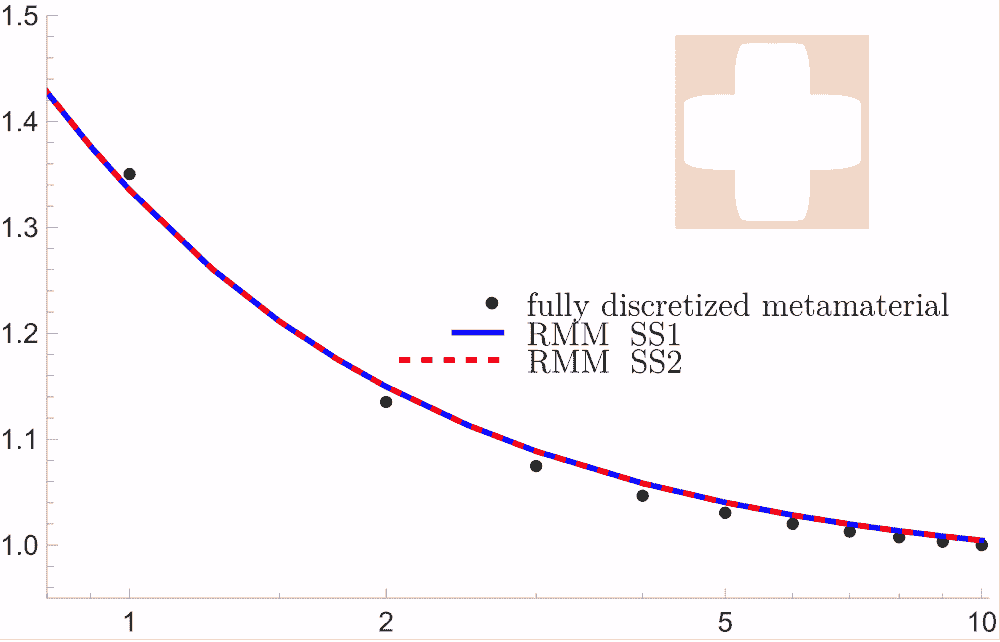}}
	\end{picture}
	\caption{volumetric mode, $\bu = \overline{\bu} = (0.01 x ,  0.01 y)^T$ on $\partial \B$}
    \end{subfigure}
    \begin{subfigure}{0.49\textwidth}
	\begin{picture}(95,60)
	\put(0,0){\def\svgwidth{7.5 cm}{\small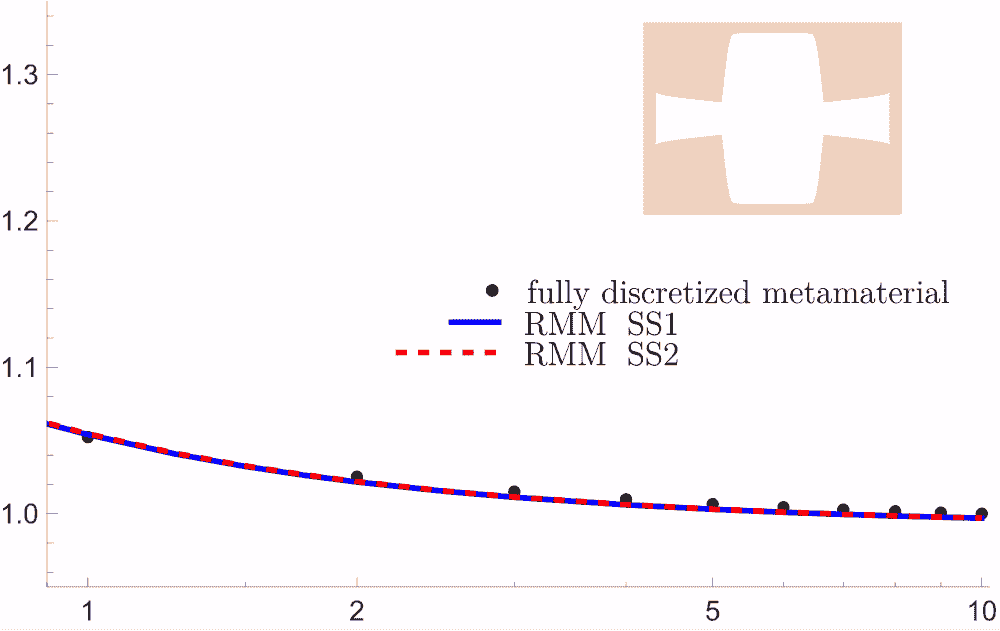}}
	\end{picture}
	\caption{deviatoric mode, $\bu = \overline{\bu} = (0.01 x ,  -0.01 y)^T$ on $\partial \B$}
    \end{subfigure}
    \begin{subfigure}{0.49\textwidth}
	\begin{picture}(90,60)
 	\put(0,0){\def\svgwidth{7.5 cm}{\small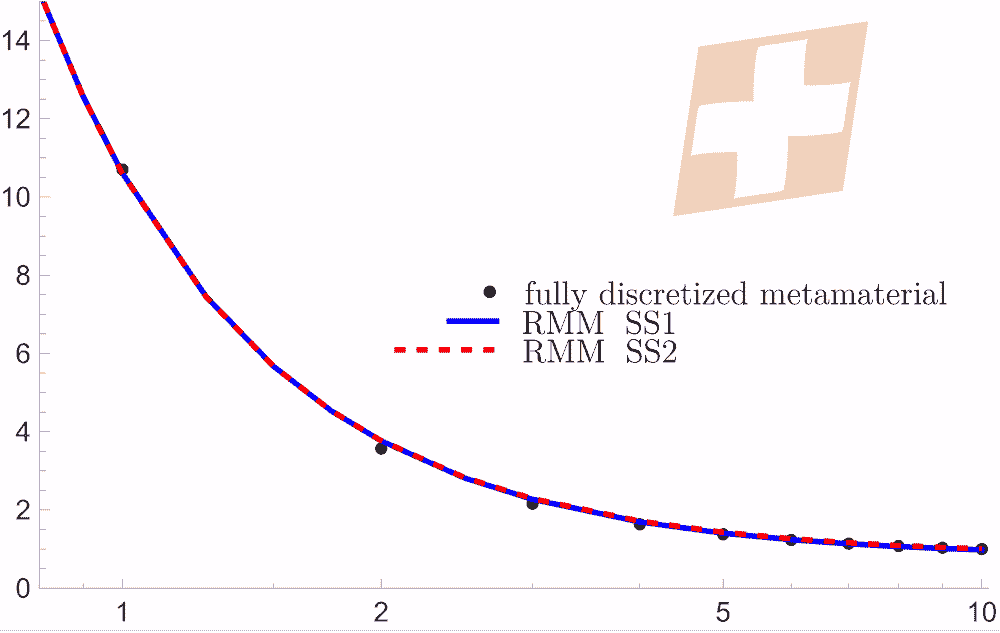}}
	\end{picture}
    \caption{shear mode, $\bu = \overline{\bu} = (0.01 y ,  0.01 x)^T$ on $\partial \B$}
    \end{subfigure}
	\caption{The results of the optimization procedure of the relaxed micromorphic model with the reference microstructure.    }
	\label{Figure:validation}
\end{figure}

The prediction obtained from the ANN surrogate model is taken here as a starting point to run the optimization described in Section \ref{sec:io}. We consider four sizes ($n=1,2,3,4$) in the optimization procedure employing a mesh with $40 \times 40$ elements. We observe that the optimization tends to prefer a negative value for $\mu_c$, however, the optimization procedure constraints $\mu_c$ to be zero. Since the role of $\mu_c$ is crucial for the dispersion curves fitting of the dynamic behavior of the material and cannot vanish, we enforce a small positive value seeking better fitting for the dynamic parameters. The results are summarized  in Table \ref{table:static_optimization}. Indeed introducing a positive value as $\mu_c = 0.1$ GPa does not affect the results of the static fitting (SS2). Figure \ref{Figure:validation} presents the results of the RMM with the optimized static parameters compared with the reference fully detailed microstructures. The relative stiffness is defined as the ratio between the actual stiffness and the macroscopic stiffness obtained by the classical periodic homogenization. Increasing the number of unit-cells $n$ leads to convergence toward a value of 1.  Both static sets (SS1 and SS2) exhibit a high level of agreement. The size-effects in the shear deformation mode are very pronounced  with stiffness ten times stiffer than the macroscopic stiffness.  Volumetric straining demonstrates approximately $40\%$ stiffer response than the macroscopic stiffness whereas deviatoric straining shows minor size-effect less than $5\%$.  At the end of this procedure all the static parameters are univocally identified once for all. The static set of parameters SS2 is the one retained in the following.

\begin{table}[h!]
\centering
\caption{Results of static optimization}
\begin{tabular}{@{} cccccccc @{}}
\toprule
\textbf{static set} & $\mu_\textrm{micro}$ [GPa]& $\mu_\textrm{micro}^*$ [GPa]& $\lambda_\textrm{micro}$ [GPa]& $\mu_c$ [GPa]& $\mu \Lc^2$ [N]& $r^2$ [Joule] & average error ($\%$)\\ 
\midrule
SS1 & $7.48$ & $276.64$ & $11.31$ & $0$ & $1108,1$ & $3.195 \, 10^{-9}$ & $1.69$ \\ 
SS2 & $7.5$ & $356.2$ & $11.41$ & $0.1$ & $1089.6$ & $3.599\, 10^{-9}$ & $1.94$ \\ 
\bottomrule
\end{tabular}
\label{table:static_optimization}
\end{table}

\newpage

\subsection{Identification of the dynamic material parameters; dispersion curves}  
\label{sec:idm}

\begin{figure}[h!]
\center
	\unitlength=1mm
	\begin{picture}(210,80)
 	\put(-5,0){\def\svgwidth{8.5 cm}{\small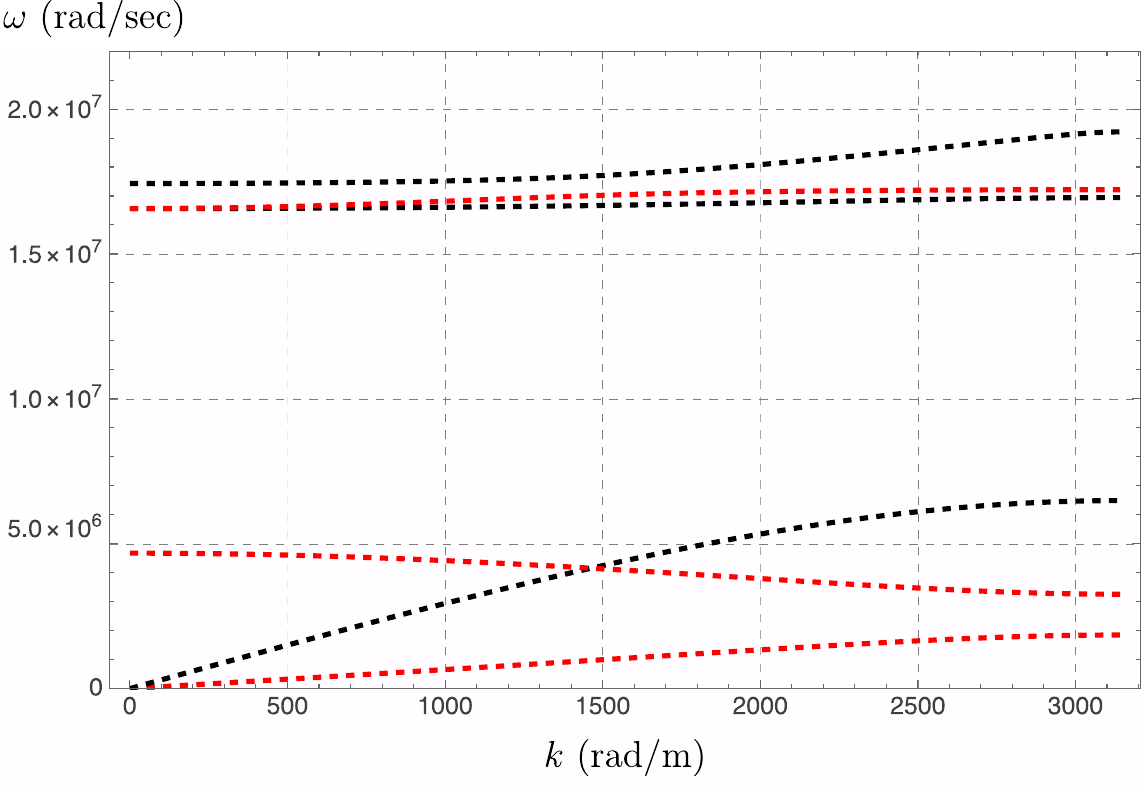}}
 	\put(82,0){\def\svgwidth{8.5 cm}{\small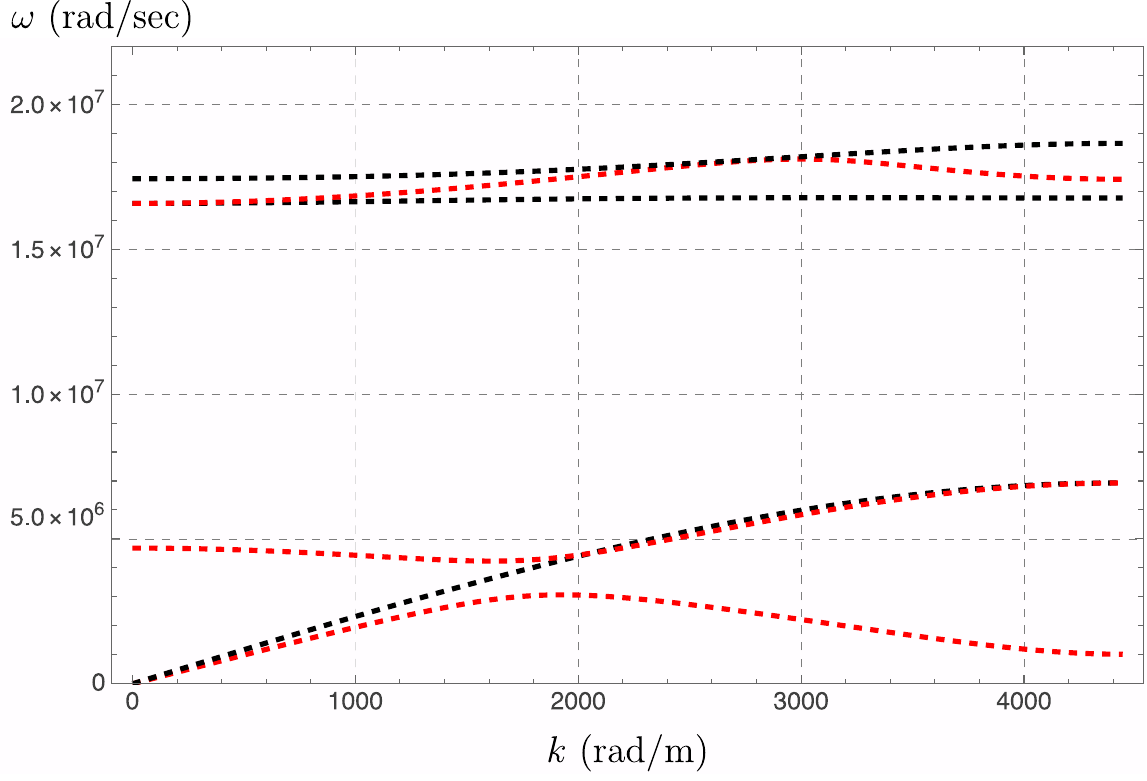}}
 	\end{picture}
    \caption{Dispersion curves of the detailed unit-cell for incidence angles of $0^{\circ}$ and $45^{\circ}$. Red curves correspond to shear branches while black curves correspond to pressure branches. The cut-off frequencies are depicted for shear ($\omega_s^1$,$\omega_s^2$) and for pressure ($\omega_p^1$,$\omega_p^2$). }
    \label{fig:dc_mstd}
\end{figure}

For the parameters appearing in the kinetic energy, we will define them via fitting the analytical dispersion curves of the relaxed micromorphic model with  the ones of the fully discretized microstructure obtained by the Bloch–Floquet analysis. The dispersion curves describe wave propagation of infinite-sized metamaterials.  
 For the heterogeneous microstructure, we considered the first six dispersion curves obtained by COMSOL Multiphysics software for incidence angles of $0^{\circ}$ and $45^{\circ}$; i.e. waves propagate along the horizontal (or vertical) direction and along the unit-cell diagonal direction aligning with the symmetry axes of the unit-cell. We will now discuss the following two cases. In the first case,  we  consider fitting the dispersion curves associated with an incidence angle of $0^\circ$, which coincides  with the direction of propagation of the pressure and shear waves applied to the finite-size metamaterial example in Section \ref{chapter:finitesizes}. In the second case, we will fit the dispersion curves at both incidence angles of $0^{\circ}$ and $45^{\circ}$. Given the targeted class of symmetry (tetragonal), an appropriate procedure must take into account the two directions. However, given the finite-size example targeted in this paper (only normally incident plane wave, see Figures \ref{fig:pre1}-\ref{fig:she2}), the set of parameters fitted only at $0^{\circ}$ can bring some useful insights, also given that the fitting of the relaxed micromorphic model for this unit-cell loses some precision at $0^{\circ}$, if the $45^{\circ}$  is also taken into account. The dispersion curves of the periodic unit-cell with a wavenumber $k$ are shown in Figure \ref{fig:dc_mstd}. The wavevector reads $\bk = (k,0)^T$ at an incidence angle of $0^{\circ}$ and  $\bk = \frac{k}{\sqrt{2}}(1,1)^T$ at an incidence angle of $45^{\circ}$.  Note that the maximal wavenumber value is defined by the periodicity limit which is $ \frac{\pi}{l}$  at an incidence angle of $0^{\circ}$ and  $\frac{\sqrt{2} \pi}{l}$  at an incidence angle of $45^{\circ}$ where $l$ is the size of the unit-cell.
 
\begin{figure}[h!]
\center
	\unitlength=1mm
	\begin{picture}(120,80)
 	\put(0,0){\def\svgwidth{12 cm}{\small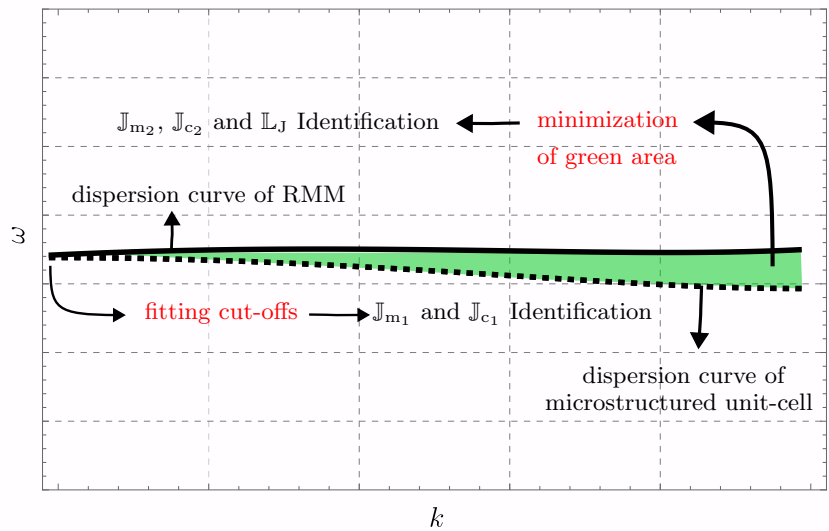}}
	\end{picture}
    \caption{ Depiction of the procedure for determining the dynamic parameters of the relaxed micromorphic model based on fitting dispersion curves.  }
     \label{fig:dc_proc}
\end{figure}

The identification procedure of the dynamic parameters is separated into two stages as depicted in Figure  \ref{fig:dc_proc}. For the case $k=0$, the solution is isotropic (the cutoffs are independent of the propagation's direction) and depends only on the already defined  static parameters besides the coefficients of the tensors $\Jm$ and $\Jc$  associated with $\sym \dot\Bdis$ and $\skew \dot\Bdis$, respectively.  Therefore, fitting the cut-off frequencies is sufficient to identify these coefficients.  The coefficients of the remaining tensors  $\Te$, $\Tc$ and $\LJ$ associated with $\sym \grad \dot\bu$, $\skew \grad \dot\bu$ and $\Curl \dot{\Bdis}$, respectively, are identified by minimizing the difference along varying wavenumber $k$.  For this, we will distinguish two cases; in the first case we minimize the difference of the dispersion curves with an incidence angle $0^{\circ}$ which coincides with the direction of propagation waves impacting the  finite-size specimen. In the second case, we will fit the dispersion curves for both directions, with incidence angles $0^{\circ}$ and $45^{\circ}$, which define the symmetry boundaries of the irreducible Brillouin zone.  

\newpage

\subsubsection{Dispersion curves of RMM \\} 

 We enforce a plane-wave ansatz (out of plane components are neglected) for both displacement and micro-distortion fields as 
 
 \begin{equation}
 u_i (\bx,t) = \psi_i e^{<\bk , \bx> - \omega t} \quad  \textrm{and} \quad \dis_{ij} (\bx,t) = \psi_{ij} e^{<\bk,  \bx> - \omega t}  \quad  \textrm{for} \quad  i,j=1,2 \,. 
 \end{equation}

which we insert in the strong forms in Equations (\ref{eq:equiMic1}) and (\ref{eq:equiMic2}). Here, $\bk = k(\cos \alpha, \sin \alpha)$ is the wavevector, $\alpha$ is the incidence angle of the incoming waves, $k$ is angular wavenumber,  $\omega$ is the angular frequency and $\psi_i,\psi_{ij}$ are the amplitudes of the considered components. The resulting linear system of equations has the form 

\begin{equation}
\bA \BPsi = 0  \quad \textrm{with} \quad  \bA  \in \CC^{6\times6}  \quad \textrm{and} \quad \BPsi  \in \RR^{6} 
\end{equation}

with  $\BPsi = (\psi_1,\psi_2,\psi_{11},\psi_{12},\psi_{21},\psi_{22})^T$. The non-trivial solution of the previous system of equations leads to $\det \bA = 0$ whose roots are in the so-called dispersion relations $\omega = \omega (\bk, \alpha)$. More details can be found in \cite{VosRizNefMad:MaL:2023}.

\subsubsection{Identification of  $\Jm$ and $\Jc$: the cut-off frequencies  \\} 

For the case $k=0$, the dispersion relations simplify into a six-order equation in $\omega^2$ which is depending on the known static parameters as well as on the coefficients of the tensors  $\Jm$ and $\Jc$. The values of the dispersion curves at $k=0$ are known as cut-off frequencies. The presence of the curvature terms ($\curl \Bdis$ and $\curl \ddot\Bdis$) does not change the values of the cut-offs.  The non-negative roots of the equation are 

\begin{equation}
\label{eq:co}
\begin{aligned}
&\omega_{1,2} = 0 \,, \quad \omega_3 = \sqrt{\frac{\mu_c}{\rho \Lc^2 M_{c_1}}}  \,, \quad \omega_4 = \sqrt{\frac{\mu_{e}^* + \mu_\textrm{micro}^* }{\rho \Lc^2 M_{m_1}^*}} \,, \\  &\omega_5 = \sqrt{\frac{\mu_{e} + \mu_\textrm{micro} }{\rho \Lc^2 M_{m_1}}} \,,  \quad 
\omega_6 = \sqrt{\frac{\lambda_{e} + \mu_{e} + \lambda_\textrm{micro}  + \mu_\textrm{micro} }{\rho \Lc^2 (\Lambda_{m_1} + M_{m_1})}}\,.
\end{aligned}
\end{equation}

The unknown coefficients  of the tensors  $\Jm$ and $\Jc$ are then calculated by assigning the numerical values of the cut-off frequencies obtained from the reference microstructure to the expressions (\ref{eq:co}) of the RMM: $\omega_3 = \bar\omega_s^1$, $\omega_4 = \bar\omega_{s}^2$,  $\omega_5 = \bar\omega_{p}^1$ and $\omega_6=\bar\omega_p^2$ leading to

 \begin{equation}
 \begin{aligned}
 &\rho \Lc^2 M_{c_1} &=& \frac{\mu_c}{(\bar\omega_s^1)^2}, &\quad \rho \Lc^2 M_{m_1}^* &= \frac{\mu_{e}^* + \mu_\textrm{micro}^*}{(\bar\omega_{s}^2)^2},  \\  
  &\rho \Lc^2 M_{m_1} &=& \frac{\mu_{e} + \mu_\textrm{micro} }{(\bar\omega_{p}^1)^2},   
  &\rho \Lc^2 (\Lambda_{m_1} + M_{m_1}) &= \frac{\lambda_{e} + \mu_{e} + \lambda_\textrm{micro}  + \mu_\textrm{micro}}{(\bar\omega_p^2)^2} \,, 
 \end{aligned}
 \end{equation}
 
 where the results are shown in Table \ref{tab:cutoffs}.

\begin{table}[h!]
\centering
\caption{Results of cut-off frequencies fitting}
\begin{tabular}{@{} cccc @{}}
\toprule
 $\rho \Lc^2 M_{c_1}$ [kg  N/m$^3$] & $\rho \Lc^2 M_{m_1}^*$  [kg  N/m$^3$] & $\rho \Lc^2 M_{m_1}$  [kg  N/m$^3$] & $\rho \Lc^2 \Lambda_{m_1}$  [kg  N/m$^3$] \\
\midrule
 $4.567 \cdot 10^{-6}$ & $1.3 \cdot 10^{-3}$ & $1.287 \cdot 10^{-4}$ & $-2.427  \cdot 10^{-5}$ \\
\bottomrule
\end{tabular}
\label{tab:cutoffs}
\end{table}

\subsubsection{Identification of  $\Te$, $\Tc$ and $\LJ$: fitting the dispersion curves on one direction \\} 
\label{iden:incedent0}

 In the following, we fit the dispersion curves obtained for an incidence angle $\alpha = 0^{\circ}$  and wavevector $\bk = (k,0)^T$.  For the case when the reference system is aligned with the direction of wave's propagation and in turn with the symmetry axes of the unit-cell, the determinant is reduced to $\det \bA^{{0}^\circ} = \det \bA_P^{{0}^{\circ}} \, \det \bA_S^{{0}^{\circ}} $ where the  determinant $\det \bA_P^{{0}^{\circ}}$ is associated with the pressure waves and the determinant $\det \bA_S^{{0}^\circ}$ is associated with the shear waves. This allows to uncouple the solution to pressure and shear waves. The solution of  $\det \bA_P^{{0}^\circ} = 0$ (roots of a third-order polynomial in $\omega^2$) leads to dispersion curves associated to pressure as

\begin{equation}
 \omega^{{0}^\circ,p}_i =  \omega^{{0}^\circ,p}_i (k, \overbrace{\mu_e,\lambda_e,\mu_\textrm{micro},  \lambda_\textrm{micro},  \mu \Lc^2}^\text{known  static parameters},  \overbrace{ \underbrace{\rho \Lc^2 \Lambda_{m_1}, \rho \Lc^2 \ M_{m_1}}_\textrm{known from cut-offs}, \underbrace{ \rho \Lc^2 \ \Lambda_{m_2}, \rho \Lc^2 \ M_{m_2}, \rho \Lc^4 \ M}_\textrm{yet to be identified }}^\textrm{dynamic parameters})\,, \quad i={1,2,3} 
\end{equation}

and  the solution of  $\det \bA_S^{{0}^\circ} = 0$ (roots of a third-order polynomial in $\omega^2$) leads to dispersion curves associated to shear  as 

\begin{equation}
 \omega^{{0}^\circ,s}_i =  \omega^{{0}^\circ,s}_i (k, \overbrace{\mu_e^*, \mu_\textrm{micro}^*, \mu_c,  \mu \Lc^2}^\text{known static parameters},  \overbrace{\underbrace{\rho \Lc^2 \ M_{c_1}, \rho \Lc^2  M_{m_1}^*}_\text{known from cut-offs}, \underbrace{\rho \Lc^2 \ M_{c_2}, \rho \Lc^2  M_{m_2}^*, \rho \Lc^4 \ M}_\textrm{yet to be identified}}^\textrm{dynamic parameters})\,, \quad i={1,2,3} \,.
\end{equation}

The remaining parameters are identified through a least-squares error minimization between the analytical dispersion curves of the RMM and the reference curves obtained from the Bloch–Floquet analysis of the microstructured periodic unit-cell (see Figure \ref{fig:dc_proc}). Two cases are considered: in the first case, the curvature terms ($\curl\Bdis$ and $\curl \ddot\Bdis$) are neglected; in the second case, they are included.  The optimization is defined as

\begin{equation}
\label{eq:optone}
\sum_{i=1}^6  \int_k [ r_i (\omega_i^{{0}^\circ} (k)-\omega_i^{{0}^\circ,\textrm{mstd}}(k))]^2 \textrm{d}k \,  \rightarrow \textrm{min} \,, 
\end{equation}

where $r_i$ is the weighting coefficient associated with the $i-$th dispersion branch $\omega_i$.  The integral is approximated by the sum of discrete values at chosen wavenumbers $k_j$ as 

\begin{equation}
\sum_{i=1}^6  \sum_{k_j} [r_i (\omega_i^{{0}^\circ} (k_j)-\omega_i^{{{0}^\circ},\textrm{mstd}}(k_j))]^2 \textrm{d}k \,  \rightarrow \textrm{min} \,. 
\end{equation}

The discrete angular wavenumbers are chosen as  $k_j = \frac{\pi}{l} \cdot (\frac{1}{5},\frac{2}{5},\frac{3}{5},\frac{4}{5},1)$ and  we set the weighting coefficients $r_i=2$ for the acoustic branches and $r_i=1$ for the optic branches giving greater emphasis to fitting the acoustic curves.

The only terms appearing in both pressure and shear dispersion relations are the ones associated with the curvature. For  $ \rho \Lc^4 M= 0$ and $\mu \Lc^2 = 0$, the dispersion relations of shear and pressure waves are decoupled and the optimization is split into two optimization problems; the first optimization for shear dispersion curves leading to the identification of 
$ \rho \Lc^2 \ M_{m_2}$ and $ \rho \Lc^2 \ \Lambda_{m_2}$  and the second optimization for pressure dispersion curves leading to the identification of $\rho \Lc^2 \ M_{c_2}$ and $\rho \Lc^2  M_{m_2}^*$.  Taking the curvature into consideration leads to the coupling of shear and pressure dispersion relations where all the remaining unknown dynamic parameters have to be identified simultaneously. The optimization is implemented using the Interior Point OPTimizer (IPOPT) package in Mathematica  (via the function IPOPTMinimize).   The results are summarized in Table  \ref{tab:fitting1}.

\begin{table}[h!]
\centering
\caption{Results of dispersion relations fitting (on one direction)}
\begin{tabular}{@{} lccccc @{}}
\toprule
\textbf{set} & $ \rho \Lc^2 \ M_{m_2}$& $\rho \Lc^2 \ \Lambda_{m_2}$  & $\rho \Lc^2 \ M_{c_2}$  & $\rho \Lc^2  M_{m_2}^*$  &  $\rho \Lc^4  M$ \\
& [kg  N/m$^3$]  & [kg  N/m$^3$]  & [kg  N/m$^3$]  & [kg  N/m$^3$]  & [kg  N/m] \\
\midrule
without $\Curl \Bdis$ & $2.777 \cdot 10^{-5}$ & $2.773 \cdot 10^{-5}$ & $2.021 \cdot 10^{-5}$ & $2.026 \cdot 10^{-5}$ & 0 \\
with  $\Curl \Bdis$ & $2.71 \cdot 10^{-5}$ & $2.71 \cdot 10^{-5}$  & $2.036 \cdot 10^{-5}$ & $2.036 \cdot 10^{-5}$ & $1.545 \cdot 10^{-10}$\\ 
\bottomrule
\end{tabular}
\label{tab:fitting1}
\end{table}

\begin{figure}[htbp]
  \centering
  \begin{subfigure}{0.49\textwidth}
  \includegraphics[width=\textwidth]{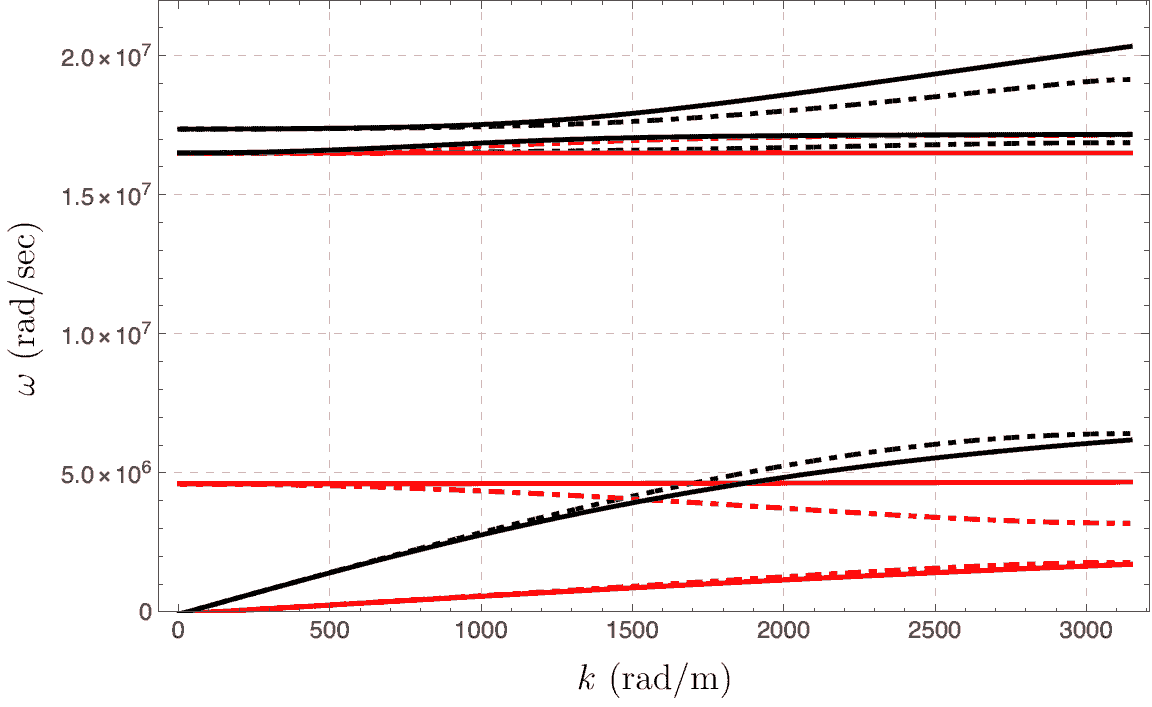}
    \caption{without  $\Curl \Bdis$, $\alpha = 0^\circ$.}
  \end{subfigure}
  \begin{subfigure}{0.49\textwidth}
  \includegraphics[width=\textwidth]{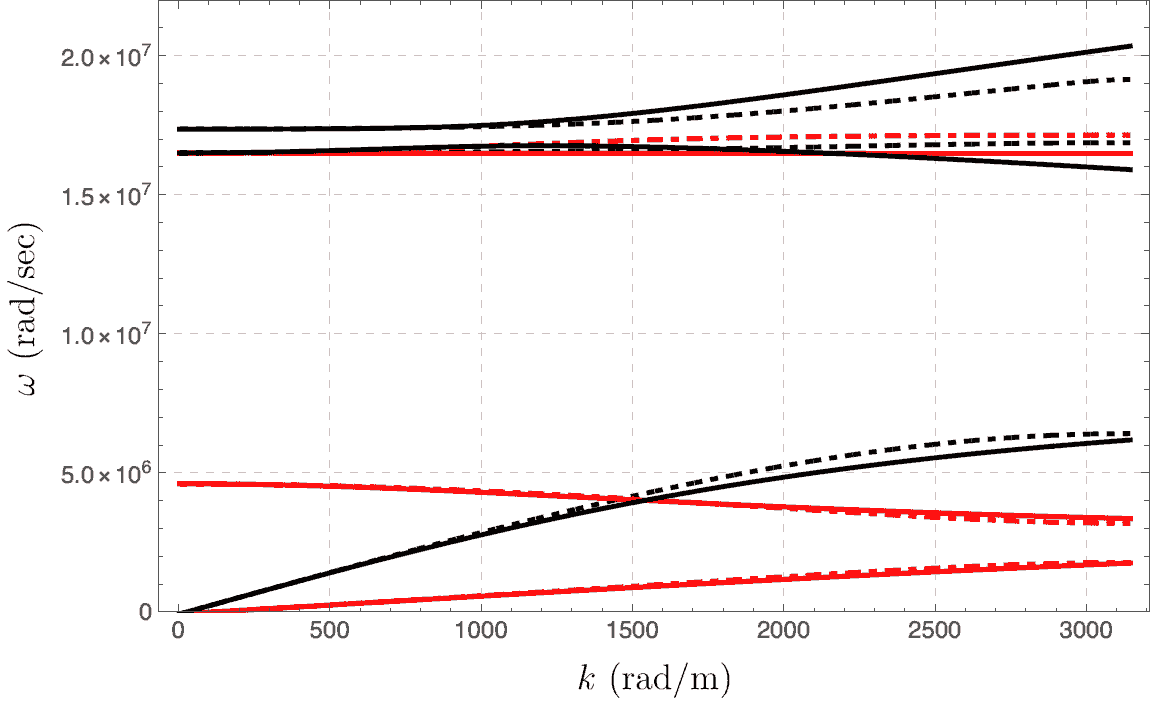}
  \caption{with   $\Curl \Bdis$, $\alpha = 0^\circ$.}
  \end{subfigure}
  \begin{subfigure}{0.49\textwidth}
  \includegraphics[width=\textwidth]{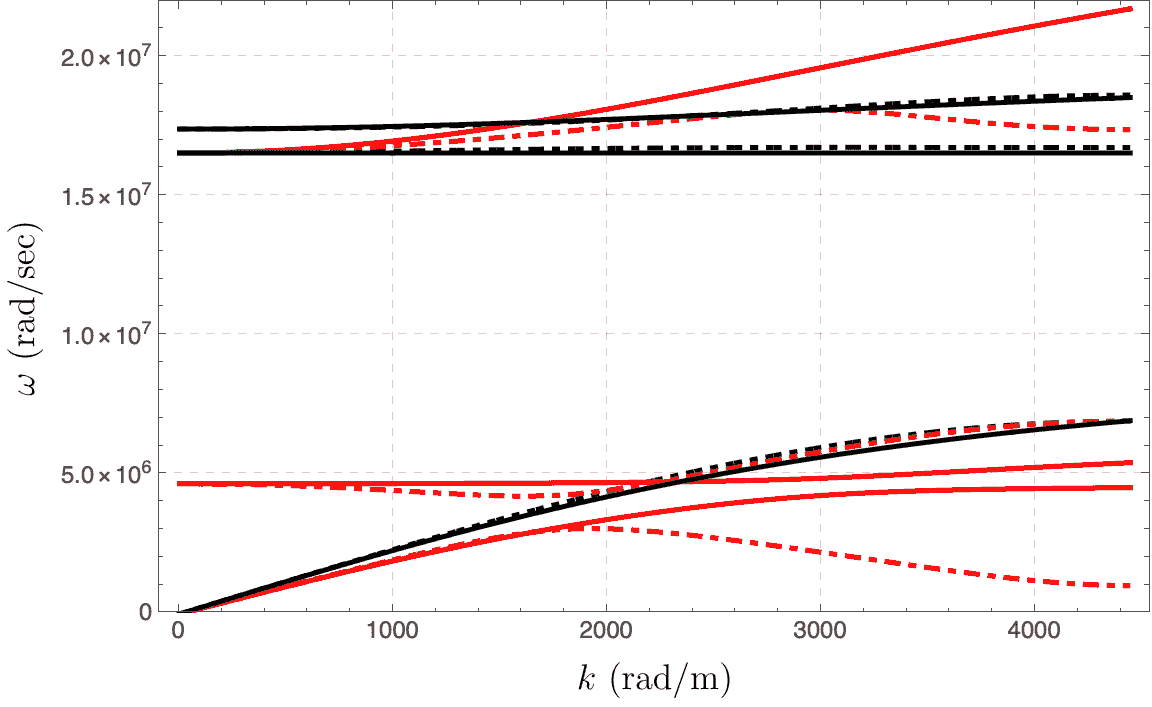}
    \caption{without  $\Curl \Bdis$, $\alpha = 45^\circ$.}
  \end{subfigure}
  \begin{subfigure}{0.49\textwidth}
  \includegraphics[width=\textwidth]{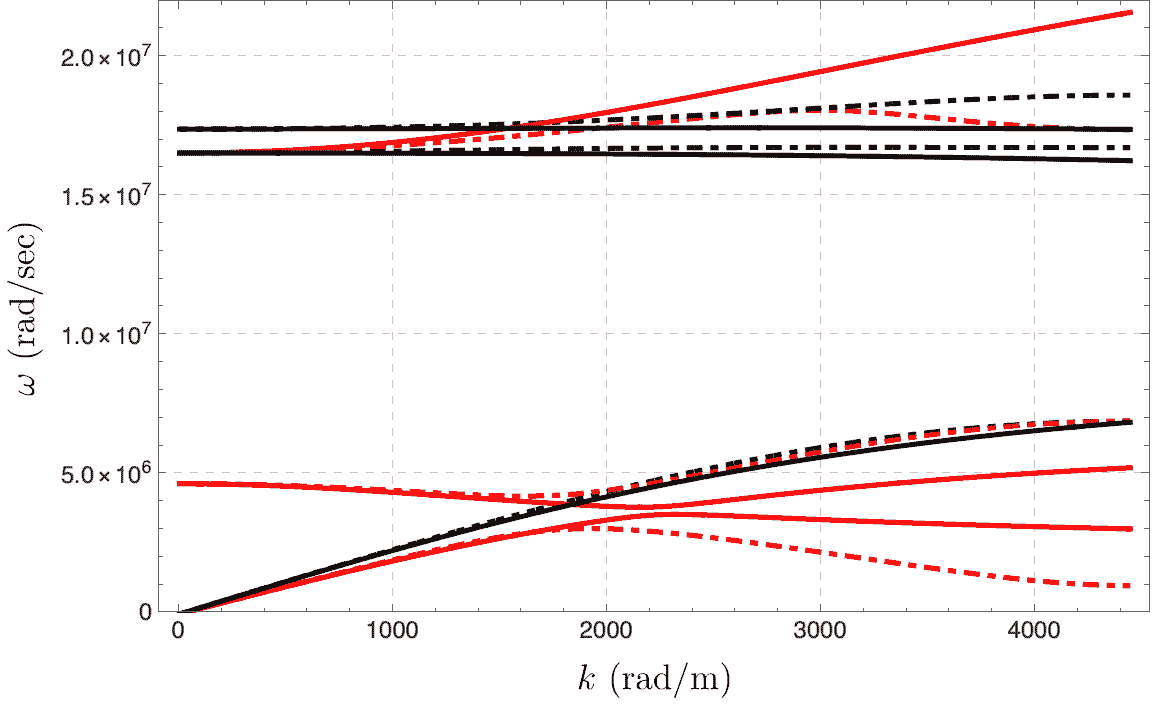}
  \caption{with  $\Curl \Bdis$,  $\alpha = 45^\circ$. }
  \end{subfigure}
  \caption{Dispersion curves for an incidence angle $\alpha = 0$ which are used for the fitting and for an incidence angle $\alpha = 45^\circ$ which are not used for the fitting. The continuous lines are for the RMM while dashed lines are for microstructured reference solution of Bloch-Floquet analysis. The red curves are for shear waves while black curves are for pressure waves.   }
 \label{fig:doptimization}
\end{figure}

In Figure \ref{fig:doptimization}, we plot the dispersion curves of the RMM beside the reference dispersion curves of Bloch–Floquet analysis conducted on the heterogeneous periodic unit-cell. Figure \ref{fig:doptimization} shows the dispersion curves for an incidence angle $0^{\circ}$ which are used for the fitting  and the dispersion curves at an incidence angle $45^{\circ}$ which we do not consider during the optimization procedure.  We notice that adding the $\Curl \Bdis$  enhances the results significantly especially because of the possibility of catching the behavior of modes with negative group velocity. The difference is most pronounced on the shear branch with a cut-off frequency at nearly $467 \cdot 10^6$ rad/sec (with negative group velocity). The shear branch decreases with increasing wavenumber and this only possible for the RMM when the curvature is considered. 

\subsubsection{Identification of  $\Te$, $\Tc$ and $\LJ$; fitting the dispersion curves on two directions \\} 
\label{iden:incedent045}

 We considered here the fitting on two directions. These directions correspond to an incidence angle  $\alpha = 0^\circ$ with a wavevector $\bk = (k,0)^T$ and  an incidence angle  $\alpha = 45^\circ$ with a  wavevector $\bk = \frac{1}{\sqrt{2}} (k,k)^T$ where $k$ is the wavenumber.  These two directions present the symmetry directions bounding the irreducible Brillouin zone.  Given the class of symmetry, a good fitting in these two directions implies a good fitting in the intermediate directions.  
 
 The dispersion relations at an incidence angle $\alpha = 0^\circ$ are shown in Section \ref{iden:incedent0}. Since the unit-cell has tetragonal symmetry,  the determinant of the acoustic  tensor for an incidence angle $\alpha = 45^\circ$ is reduced  as well to $\det \bA^{{45}^\circ} = \det \bA_P^{{45}^\circ} \, \det \bA_S^{{45}^\circ} $ where the  determinant $\det \bA_P^{{45}^\circ}$ is associated with the pressure waves and the determinant $\det \bA_P^{{45}^\circ}$ is associated with the shear waves.   The solution of  $\det \bA_P^{{45}^\circ} = 0$ leads to dispersion curves associated to pressure  as

\begin{equation}
\begin{aligned}
 \omega^{{45}^\circ,p}_i =  \omega^{{45}^\circ,p}_i ( & k, \overbrace{\mu_e,\mu_e^*,\lambda_e,\mu_\textrm{micro},\mu_\textrm{micro}^*,  \lambda_\textrm{micro},  \mu \Lc^2}^\text{known static parameters}, \\ &  \overbrace{ \underbrace{\rho \Lc^2 \Lambda_{m_1}, \rho \Lc^2 \ M_{m_1},  \rho \Lc^2 \ M_{m_1}^*}_\textrm{known from cut-offs}, \underbrace{ \rho \Lc^2 \ \Lambda_{m_2}, \rho \Lc^2 \ M_{m_2}, \rho \Lc^2 \ M_{m_2}^*, \rho \Lc^4 \ M}_\textrm{yet to be defined }}^\textrm{dynamic parameters})\,, \quad i={1,2,3} 
 \end{aligned}
\end{equation}

and  the solution of  $\det \bA_S^{{45}^\circ} = 0$ leads to dispersion curves associated to shear  as 

\begin{equation}
 \omega^{{45}^\circ,s}_i =  \omega^{{45}^\circ,s}_i (k, \overbrace{\mu_e, \mu_\textrm{micro}, \mu_c,  \mu \Lc^2}^\text{known static parameters},  \overbrace{\underbrace{\rho \Lc^2 \ M_{c_1}, \rho \Lc^2  M_{m_1}}_\text{known}, \underbrace{\rho \Lc^2 \ M_{c_2}, \rho \Lc^2  M_{m_2}, \rho \Lc^4 \ M}_\textrm{yet to be defined}}^\textrm{dynamic parameters})\,, \quad i={1,2,3} \,.
\end{equation}

\begin{figure}[htbp]
  \centering
  \begin{subfigure}{0.49\textwidth}
  \includegraphics[width=\textwidth]{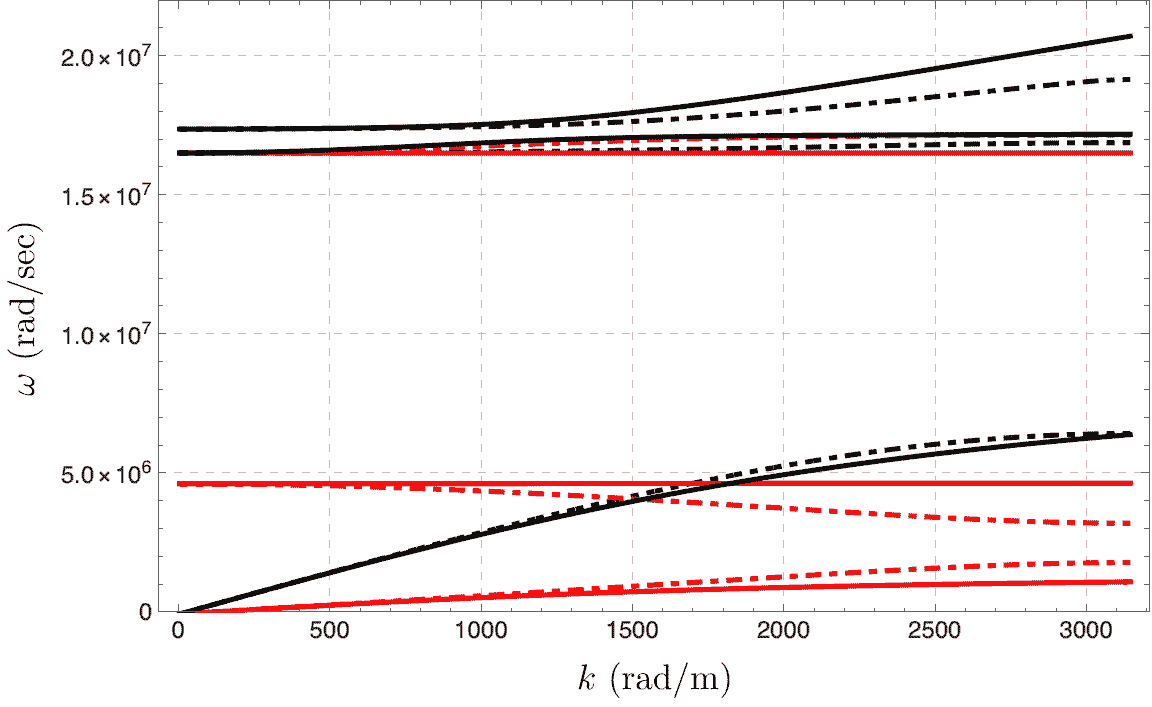}
    \caption{without  $\Curl \Bdis$, $\alpha = 0^\circ$.}
  \end{subfigure}
  \begin{subfigure}{0.49\textwidth}
  \includegraphics[width=\textwidth]{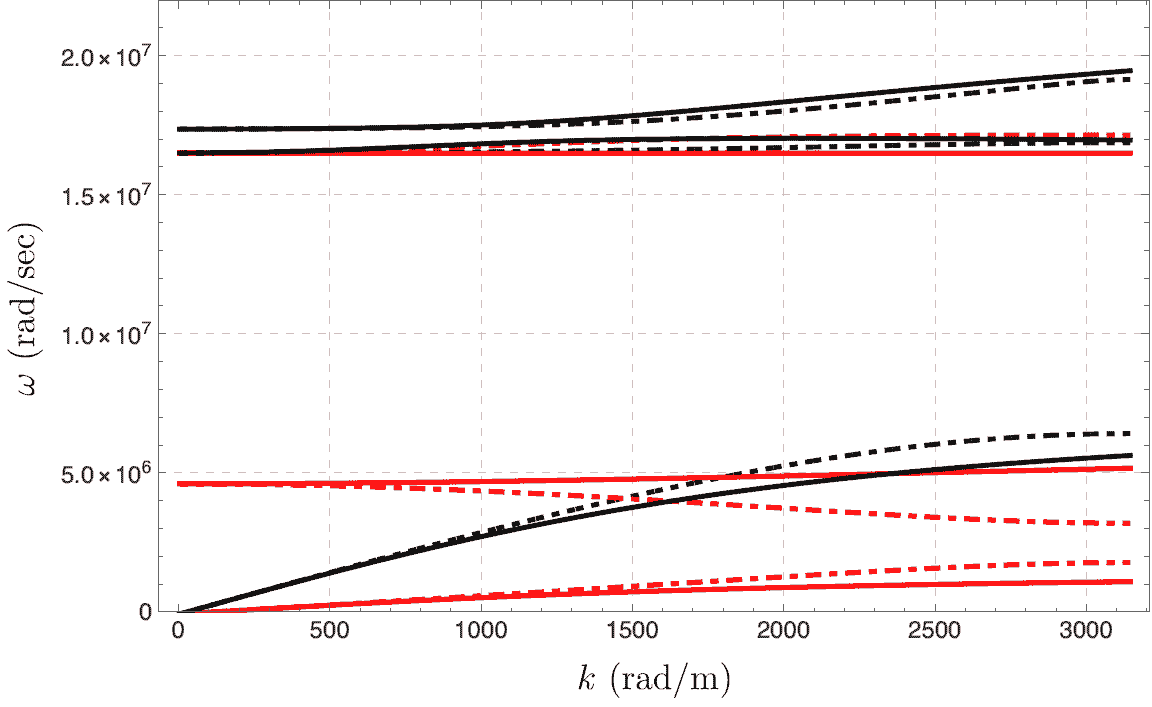}
  \caption{with   $\Curl \Bdis$, $\alpha = 0^\circ$.}
  \end{subfigure}
  \begin{subfigure}{0.49\textwidth}
  \includegraphics[width=\textwidth]{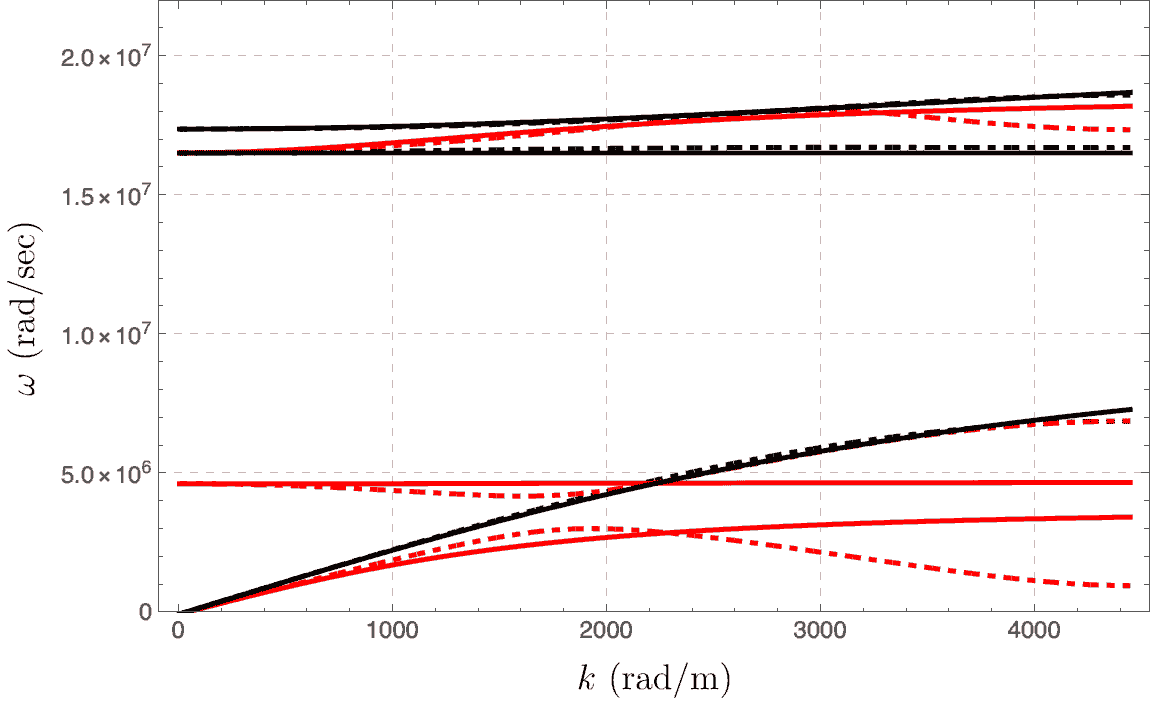}
    \caption{without  $\Curl \Bdis$, $\alpha = 45^\circ$.}
  \end{subfigure}
  \begin{subfigure}{0.49\textwidth}
  \includegraphics[width=\textwidth]{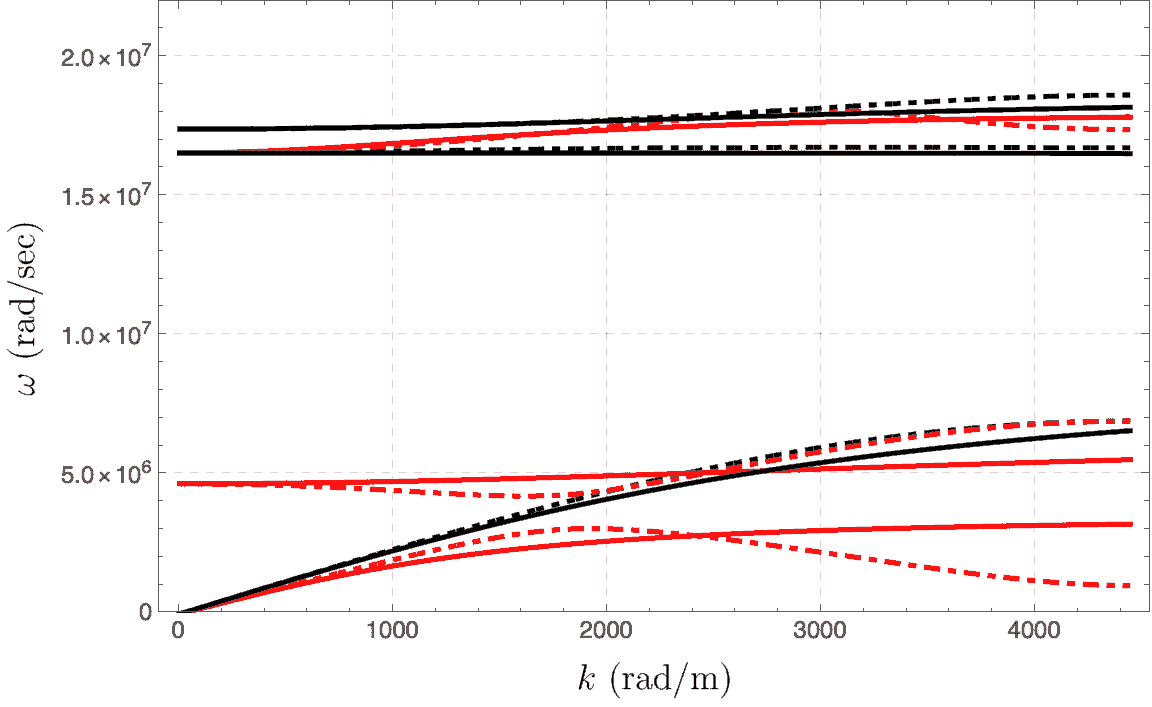}
  \caption{with  $\Curl \Bdis$.  $\alpha = 45^\circ$. }
  \end{subfigure}
  \caption{Dispersion curves at incidence angles of $0^\circ$  and $45^\circ$ which are used for the fitting procedure. The continuous lines are for the RMM while dashed lines are for microstructured reference solution of Bloch-Floquet analysis. The red curves are for shear waves while black curves are for pressure waves.   }
 \label{fig:doptimization2}
\end{figure}

The optimization, incorporating two incidence angles as in  analogy to Equation (\ref{eq:optone}), is defined as 

\begin{equation}
\sum_{a=0^\circ,45^\circ} \sum_{i=1}^6  \int_k [ r_i (\omega_i^{\alpha} (k)-\omega_i^{\alpha,\textrm{mstd}}(k))]^2 \textrm{d}k \,  \rightarrow \textrm{min} \,, 
\end{equation}

which is approximated by the sum of discrete values at chosen wavenumbers $k_j$ as

\begin{equation}
\sum_{a=0^\circ,45^\circ} \sum_{i=1}^6  \sum_{k_j} [r_i (\omega_i^{{0}^\circ} (k_j)-\omega_i^{{{0}^\circ},\textrm{mstd}}(k_j))]^2 \textrm{d}k \,  \rightarrow \textrm{min} \,. 
\end{equation}

where the discrete values are taken as $k_j = (\frac{\pi}{l}) \cdot (\frac{1}{5},\frac{2}{5},\frac{3}{5},\frac{4}{5},1)^T$ at an incidence  angle of $0^{\circ}$ and as $k_j = (\frac{\sqrt{2} \pi}{l}) \cdot (\frac{1}{5},\frac{2}{5},\frac{3}{5},\frac{4}{5},1)^T$ at an incidence angle of $45^{\circ}$.  Similar to the the fitting on one direction,  we set the weighting coefficients to 2 for the acoustic curves and 1 for the optic curves.  The results of the optimization procedure are shown in Figure \ref{fig:doptimization2} and the resulting parameters are summarized in Table \ref{tab:fitting2}. It is worth mentioning that the uniqueness of the outcome of the fitting procedure introduced in Sections \ref{iden:incedent0} and \ref{iden:incedent045} is not guaranteed, see \cite{FelJenPatMad:2025:nli}, and future work will explore the introduction of side constraints that may be derived from finite-size metamaterials with different cuts.

\begin{table}[h!]
\centering
\caption{Results of dispersion relations fitting (on two directions)}
\begin{tabular}{@{} lccccc @{}}
\toprule
\textbf{set} & $ \rho \Lc^2 \ M_{m_2}$& $\rho \Lc^2 \ \Lambda_{m_2}$  & $\rho \Lc^2 \ M_{c_2}$  & $\rho \Lc^2  M_{m_2}^*$  &  $\rho \Lc^4  M$ \\
& [kg  N/m$^3$]  & [kg  N/m$^3$]  & [kg  N/m$^3$]  & [kg  N/m$^3$]  & [kg  N/m] \\
\midrule
without $\Curl \Bdis$ & $3.214 \cdot 10^{-5}$ & $-2.685 \cdot 10^{-6}$ & $2.86 \cdot 10^{-4}$ & $2.766 \cdot 10^{-5}$ & 0 \\
with  $\Curl \Bdis$ & $1.553 \cdot 10^{-5}$ & $-1.537 \cdot 10^{-5}$  & $2.473 \cdot 10^{-4}$ & $9.356 \cdot 10^{-5}$ & $2.533 \cdot 10^{-11}$\\ 
\bottomrule
\end{tabular}
\label{tab:fitting2}
\end{table}

It can be noticed that the complete fitting involving both directions at $0^\circ$  and $45^\circ$ lowers the precision of the RMM with respect to the fitting at $0^\circ$ , especially with respect to the branches with negative group velocity. This is mainly due to the shape of the shear acoustic curve at $45^\circ$ which is passing from positive to negative group velocity while increasing $k$. To improve the RMM's ability to catch such more complex behavior, extra constitutive parameters should be added to the model which is out of the scope of this paper. The introduced fitting procedure will be used in future work to fit the RMM's parameters on unit-cells that show monotonically increasing acoustic curves. In these cases the optimization potential of the proposed procedure will be further unveiled.

\newpage

\section{Finite-size metamaterial}
\label{chapter:finitesizes}

\begin{figure}[h!]
\center
	\unitlength=1mm
	\begin{picture}(150,125)
	\put(20,5){\def\svgwidth{11 cm}{\small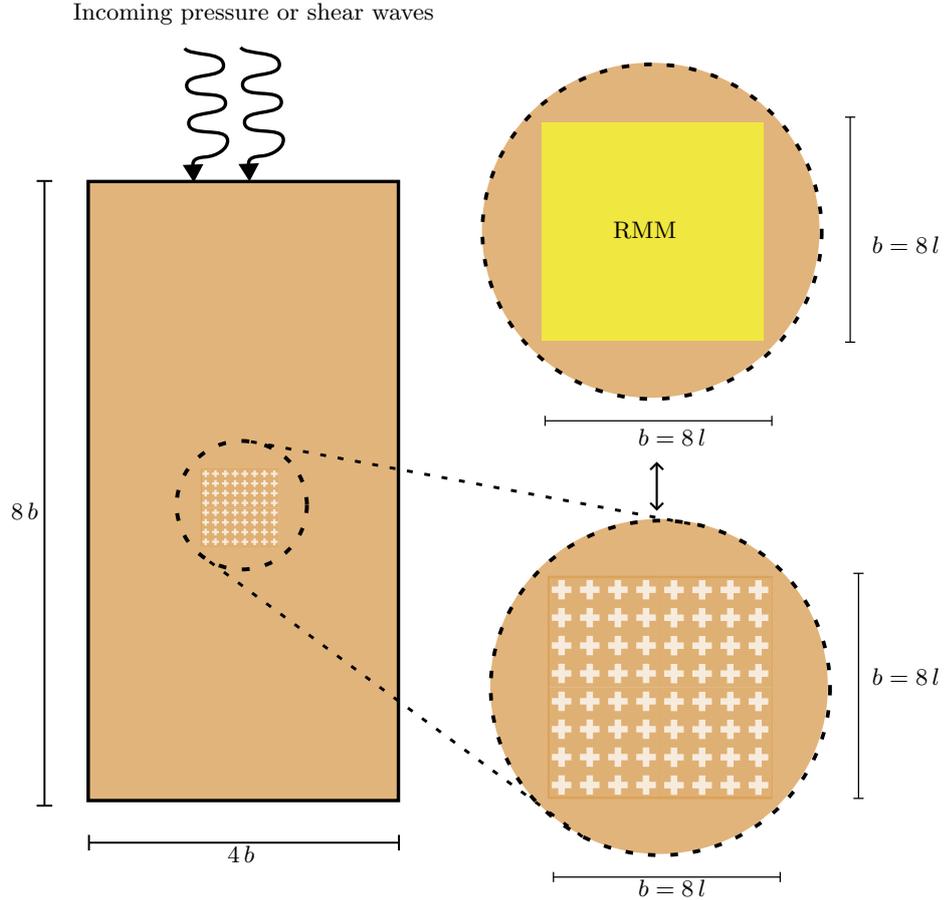}}
	\end{picture}
	\caption{The geometry of the finite-size example.    }
	\label{Figure:finite_size}
\end{figure}

We consider a finite-size metamaterial consisting of $8\times8$ unit-cells. The metamaterial specimen is embedded in a homogeneous matrix which is modeled as a Cauchy continuum with the same material parameters of the base material. The total size of the surrounding Cauchy continuum is $4 b \times 8 b$ where $b$ is the total size of the metamaterial block $b = 8l$ and the metamaterial block is located at the center of the domain. An illustration of the numerical setup is shown in Figure \ref{Figure:finite_size}.

We induce pressure and shear incoming waves with different frequencies and an amplitude equals to $u_0 = 10^{-12} $ m which travel from the upper edge to the lower edge and hit the metamaterial block which cause different reflection/transmission patterns.  We analyze this example numerically by performing a frequency-domain analysis in COMSOL Multiphysics. The results of the reference microstructured solution will be compared with the solution of the RMM. For the RMM simulations, we examine the dynamical parameters obtained via fitting in one propagation direction and in two propagation directions. Since the focus is on interaction between the size-effects and dynamical behavior of finite-sized metamaterials,  we consider the following two cases

$\bullet$ {\bf Reduced relaxed micromorphic model (RRMM)} without curvature terms (i.e. no $\Curl \Bdis$ or $\Curl \dot \Bdis$) which leads to an internal variable model that does not account for size-effects, however, still captures the band-gap. The behavior of the static part of the reduced relaxed micromorphic model is close to the soft periodic homogenization limit. Even if, in theory, $\Cmicro$ cannot be calculated from the static fitting procedure presented in Section \ref{sec:se} because the curvature term is not presented. However, it is sensible to ensure that the static parameters are fitted by taking into account the curvature terms in order to predict the size-effect behavior under static testing.  The static parameters are also kept for the limiting case when no curvature is present. The remaining dynamic parameters are refitted by the dynamic fitting in Section \ref{sec:idm} for the case without curvature. When considering the finite-size specimen in the experiment of Figure \ref{Figure:finite_size}, no boundary condition on $\Bdis$ can be considered at the RRMM/Cauchy interfaces as far as no curvature is retained in the RRMM.

$\bullet$ {\bf Relaxed micromorphic model (RMM)} with curvature terms (both $\Curl \Bdis$ or $\Curl \dot \Bdis$ are active). Here the stiffness of the static part accounts for the hardening caused by the size-effects.  Moreover, we enforce the consistent boundary condition ($\Bdis \cdot \Btau = \nabla \bu \cdot \Btau $) on  the interface between the RMM and surrounding Cauchy continuum.

\begin{figure}[!ht]
    \centering
    \begin{subfigure}{0.7\textwidth}
        \includegraphics[width=\textwidth]{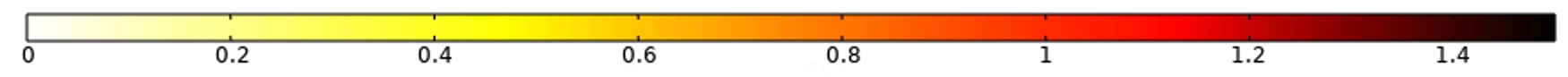}
        \put(0,10){\textbf{\large ${\lvert u \lvert}/{u_0}$}}
    \end{subfigure}

	\begin{subfigure}{0.0\textwidth}    
    \begin{picture}(0,0)
        \put(-6,20){\rotatebox{90}{\textbf{$12.56 \cdot 10^{5}$ rad/sec}}}
    \end{picture}
	\end{subfigure}
    \begin{subfigure}{0.14\textwidth}
        \includegraphics[width=\textwidth]{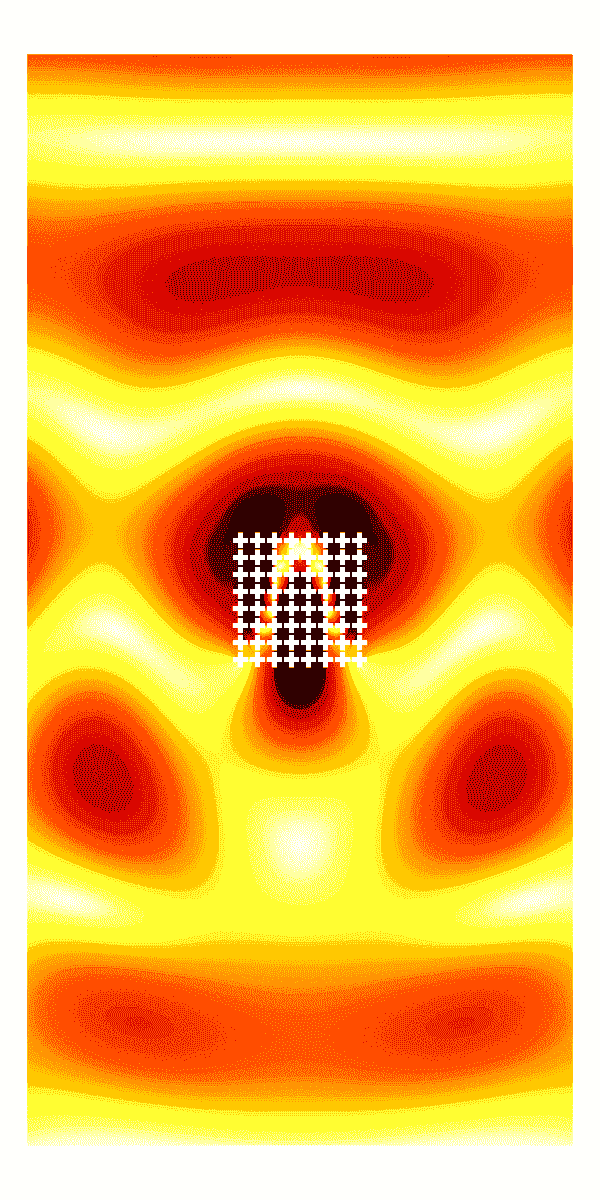}
    \end{subfigure}
    \begin{subfigure}{0.14\textwidth}
        \includegraphics[width=\textwidth]{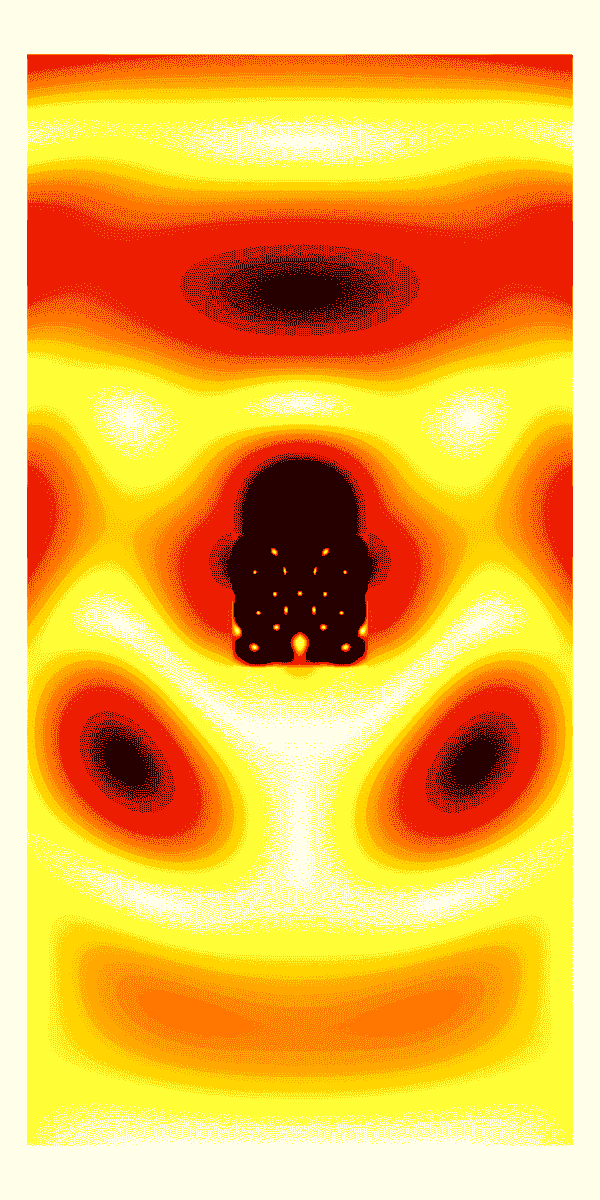}
    \end{subfigure}
    \begin{subfigure}{0.14\textwidth}
        \includegraphics[width=\textwidth]{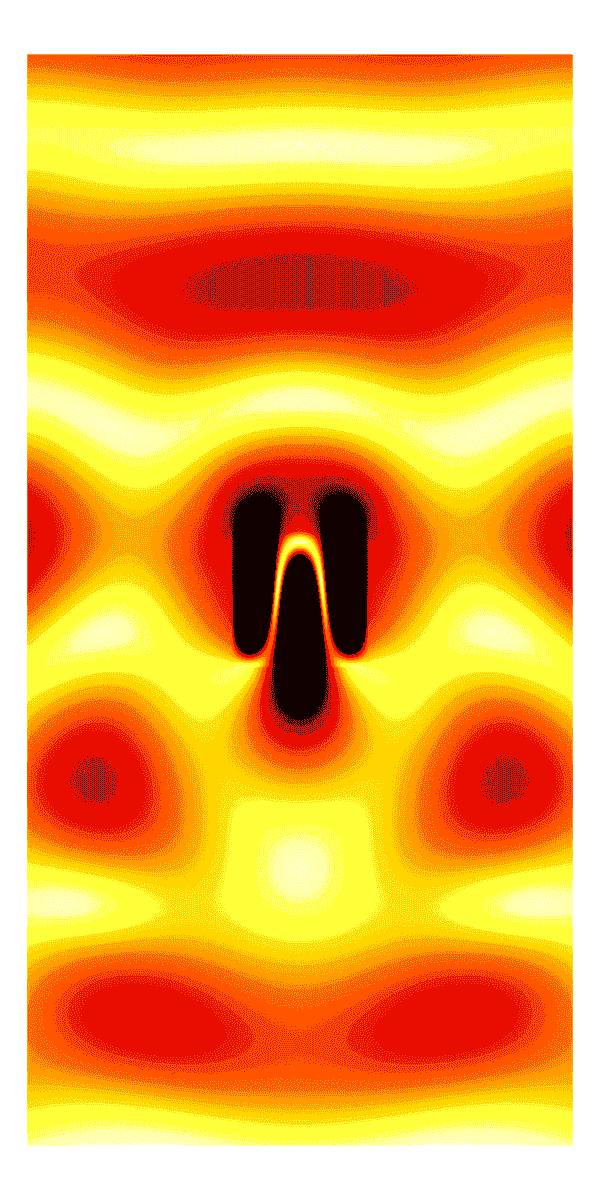}
    \end{subfigure}
    \begin{subfigure}{0.14\textwidth}
        \includegraphics[width=\textwidth]{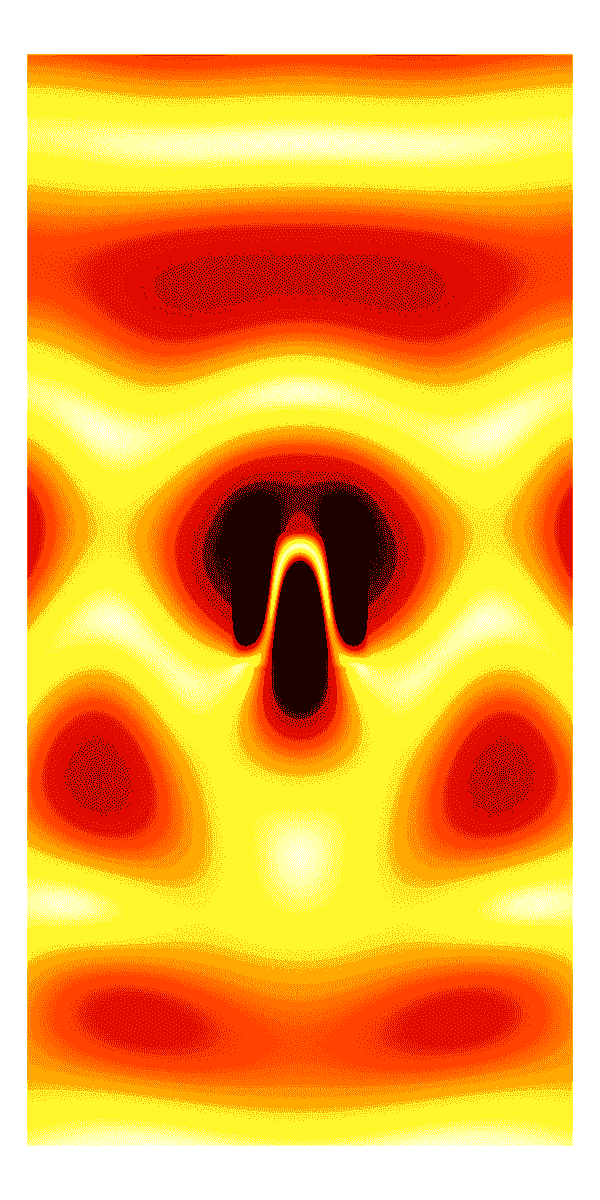}
    \end{subfigure}
    \begin{subfigure}{0.14\textwidth}
        \includegraphics[width=\textwidth]{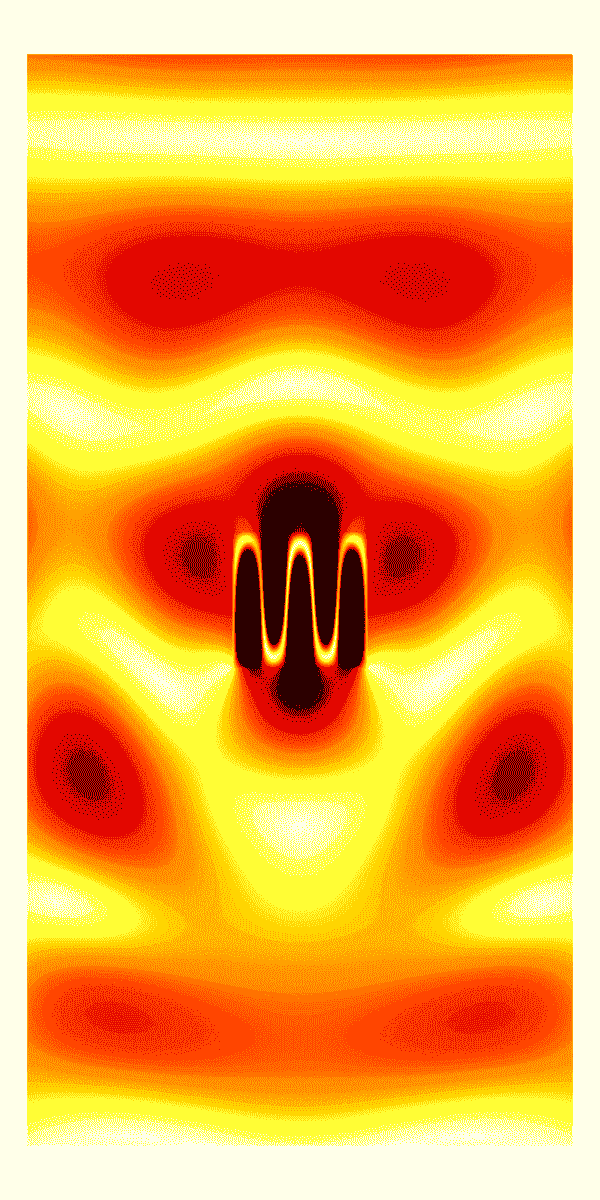}
    \end{subfigure}
    \begin{subfigure}{0.14\textwidth}
        \includegraphics[width=\textwidth]{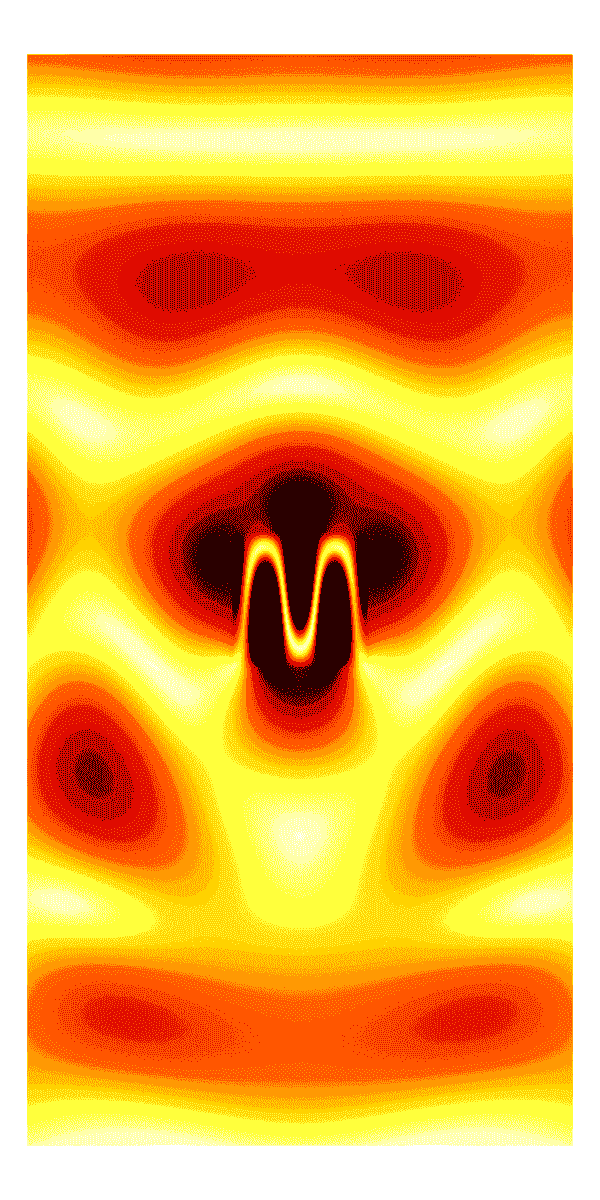}
    \end{subfigure}

    \centering
	\begin{subfigure}{0.0\textwidth}
    \begin{picture}(0,0)
        \put(-6,20){\rotatebox{90}{\textbf{$25.13 \cdot 10^{5}$ rad/sec}}}
    \end{picture}
	\end{subfigure}
    \begin{subfigure}{0.14\textwidth}
        \includegraphics[width=\textwidth]{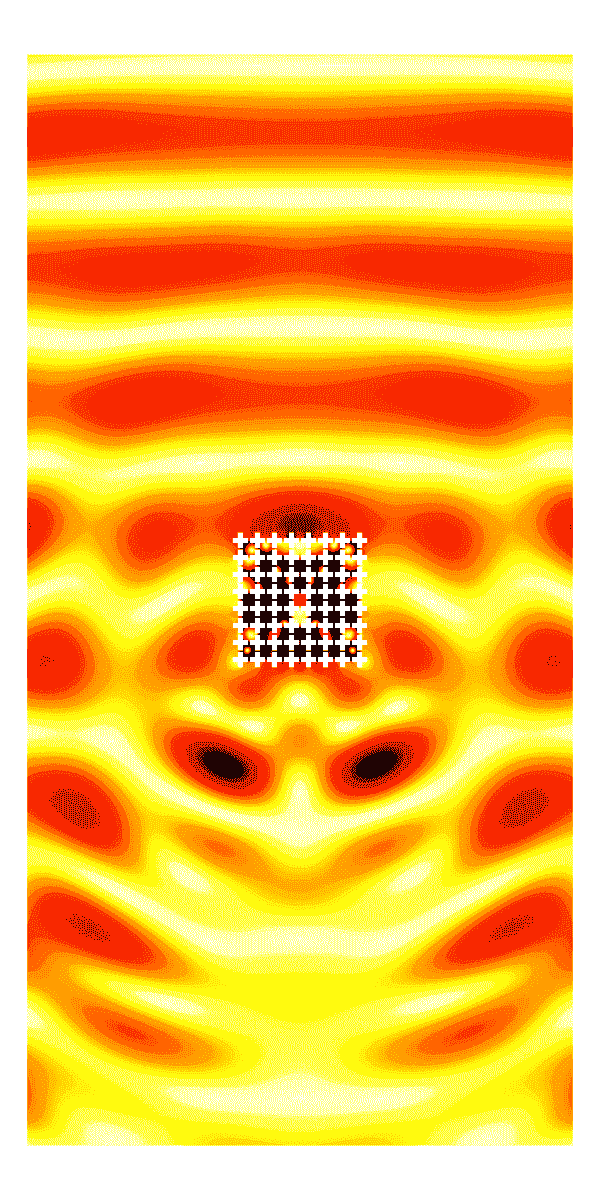}
    \end{subfigure}
    \begin{subfigure}{0.14\textwidth}
        \includegraphics[width=\textwidth]{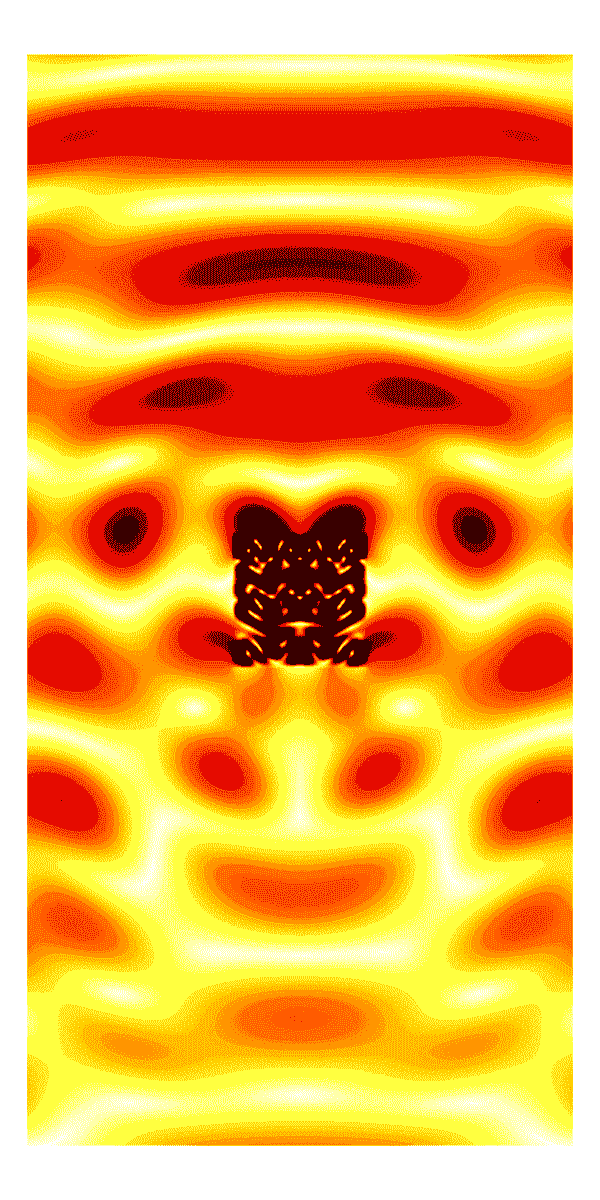}
    \end{subfigure}
    \begin{subfigure}{0.14\textwidth}
        \includegraphics[width=\textwidth]{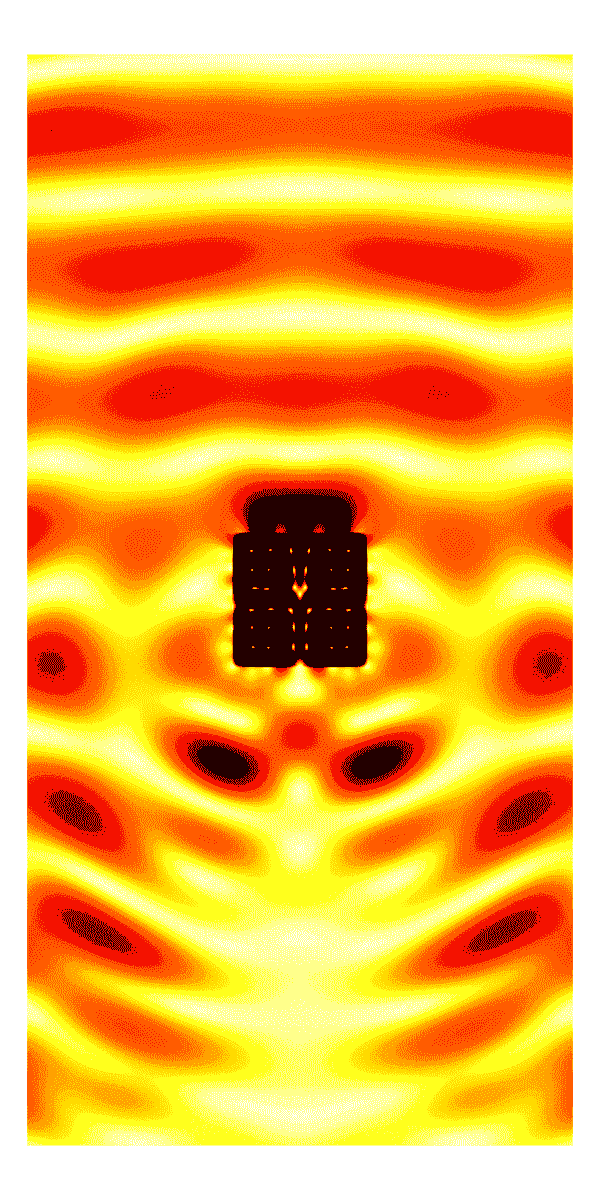}
    \end{subfigure}
    \begin{subfigure}{0.14\textwidth}
        \includegraphics[width=\textwidth]{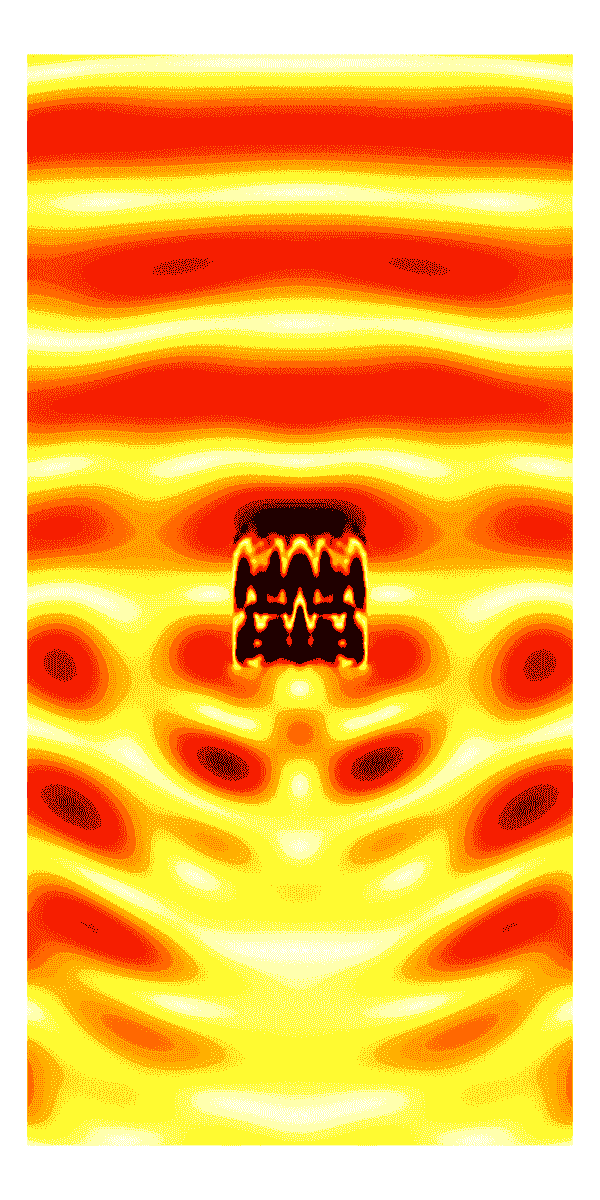}
    \end{subfigure}
    \begin{subfigure}{0.14\textwidth}
        \includegraphics[width=\textwidth]{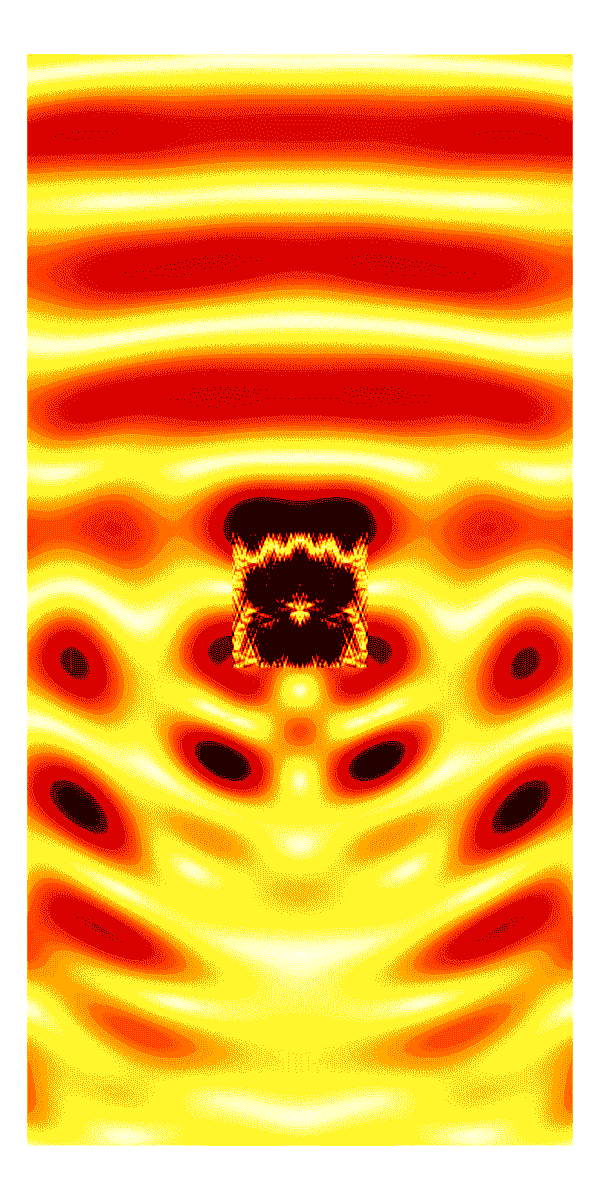}
    \end{subfigure}
    \begin{subfigure}{0.14\textwidth}
        \includegraphics[width=\textwidth]{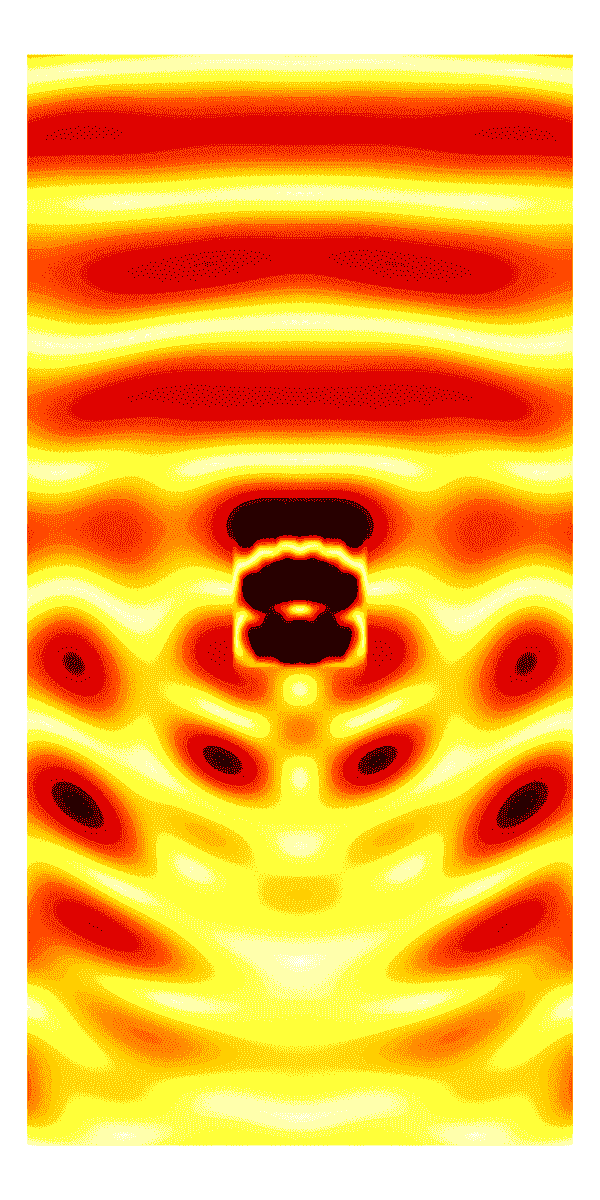}
    \end{subfigure}
    
	\begin{subfigure}{0.0\textwidth}
    \begin{picture}(0,0)
        \put(-6,20){\rotatebox{90}{\textbf{$37.7 \cdot 10^{5}$ rad/sec}}}
    \end{picture}
	\end{subfigure}
    \begin{subfigure}{0.14\textwidth}
        \includegraphics[width=\textwidth]{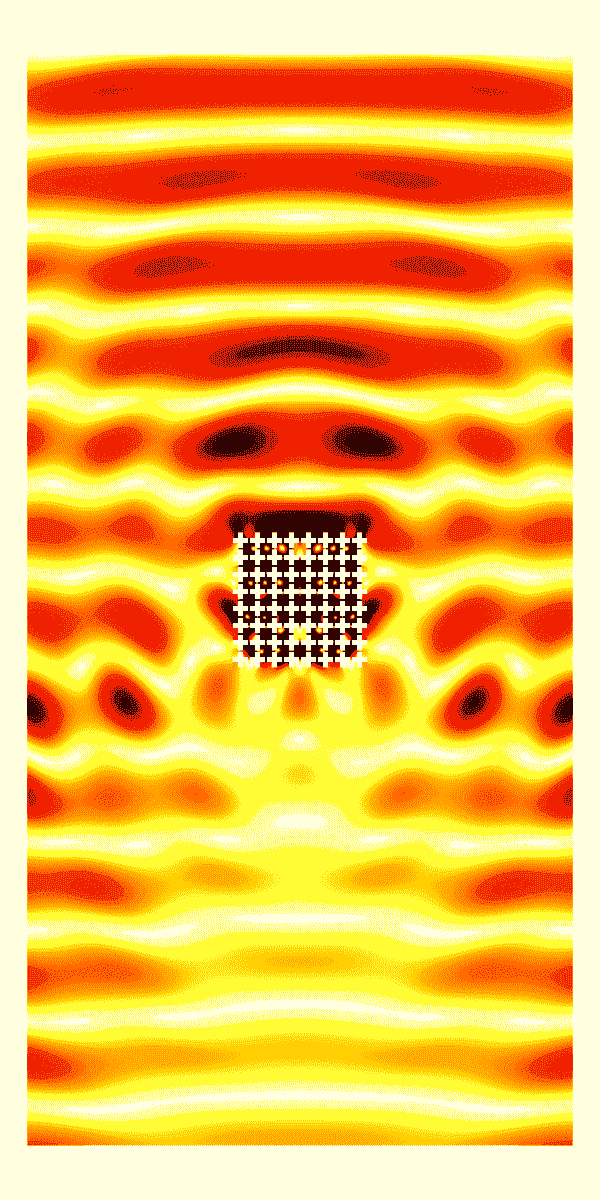}
    \end{subfigure}
    \begin{subfigure}{0.14\textwidth}
        \includegraphics[width=\textwidth]{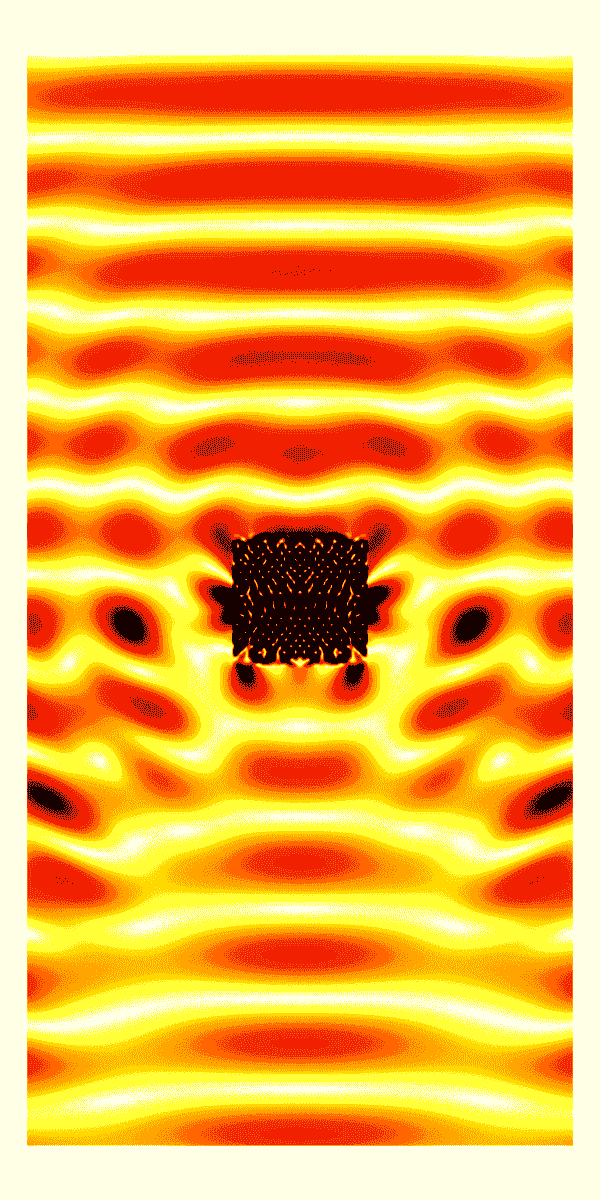}
    \end{subfigure}
    \begin{subfigure}{0.14\textwidth}
        \includegraphics[width=\textwidth]{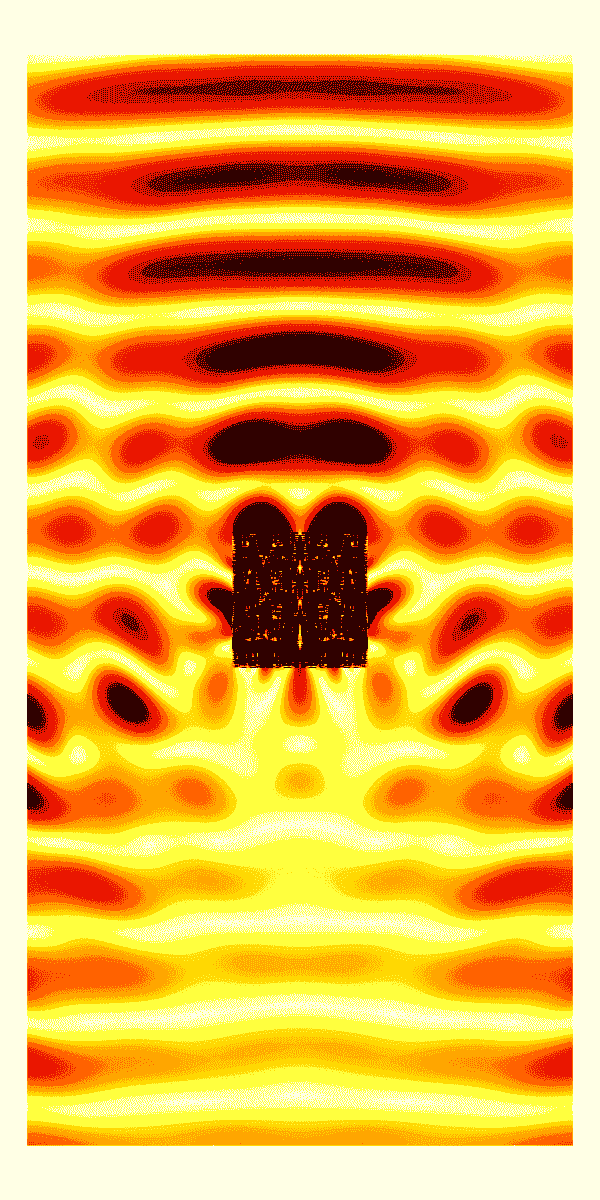}
    \end{subfigure}
    \begin{subfigure}{0.14\textwidth}
        \includegraphics[width=\textwidth]{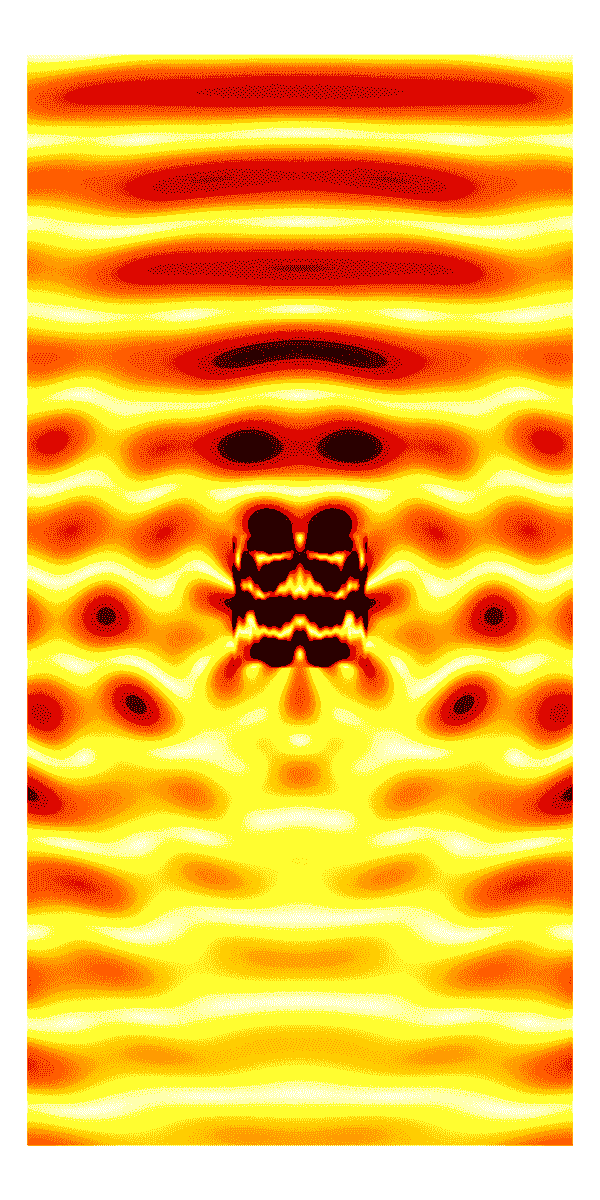}
    \end{subfigure}
    \begin{subfigure}{0.14\textwidth}
        \includegraphics[width=\textwidth]{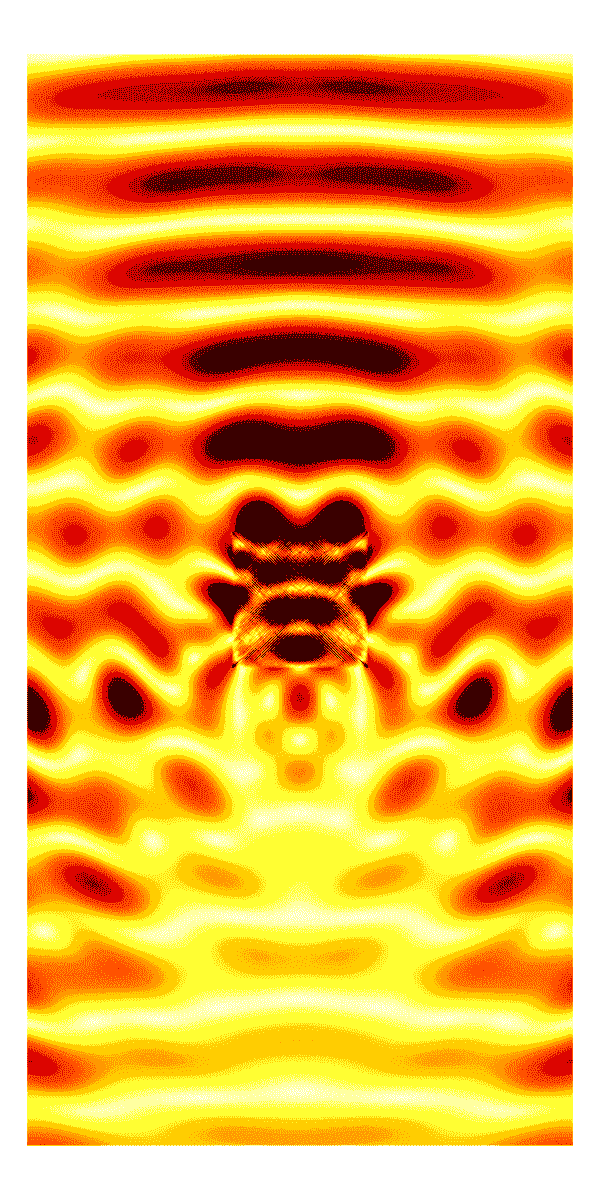}
    \end{subfigure}
    \begin{subfigure}{0.14\textwidth}
        \includegraphics[width=\textwidth]{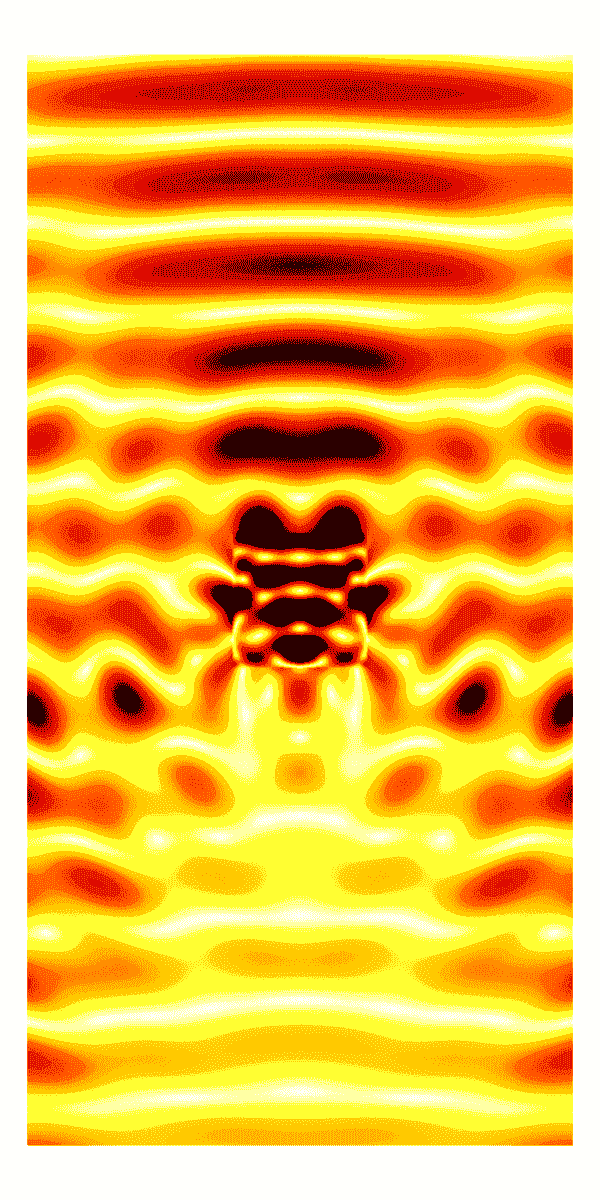}
    \end{subfigure}

    	\begin{subfigure}{0.0\textwidth}
    \begin{picture}(0,0)
        \put(-6,20){\rotatebox{90}{\textbf{$50.27 \cdot 10^{5}$ rad/sec}}}
    \end{picture}
	\end{subfigure}
    \begin{subfigure}{0.14\textwidth}
        \includegraphics[width=\textwidth]{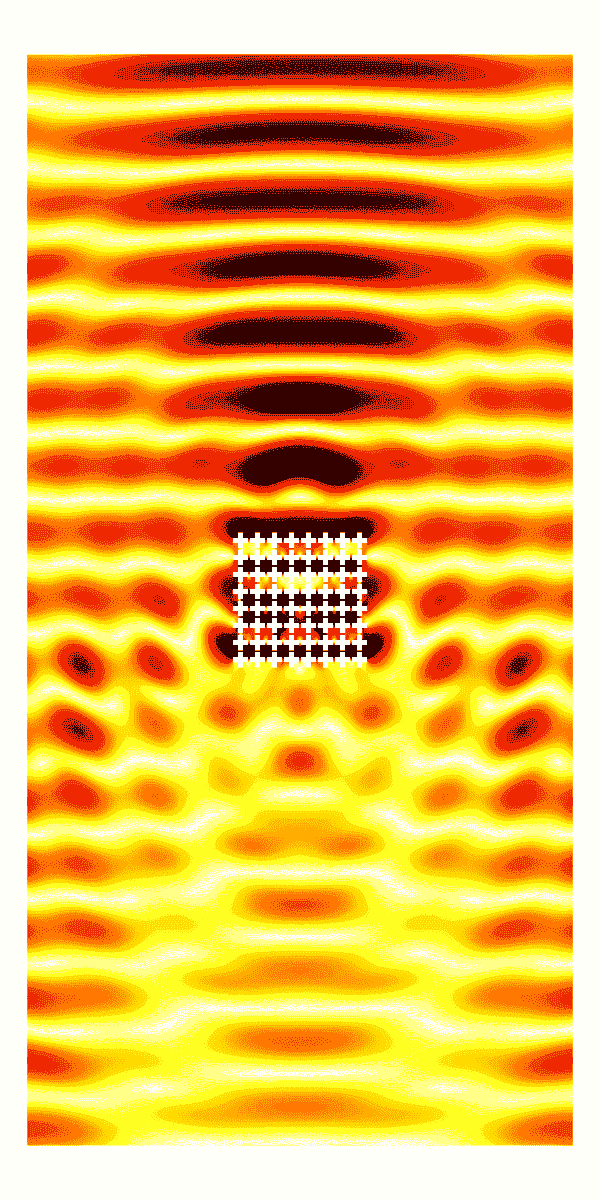}
    \end{subfigure}
    \begin{subfigure}{0.14\textwidth}
        \includegraphics[width=\textwidth]{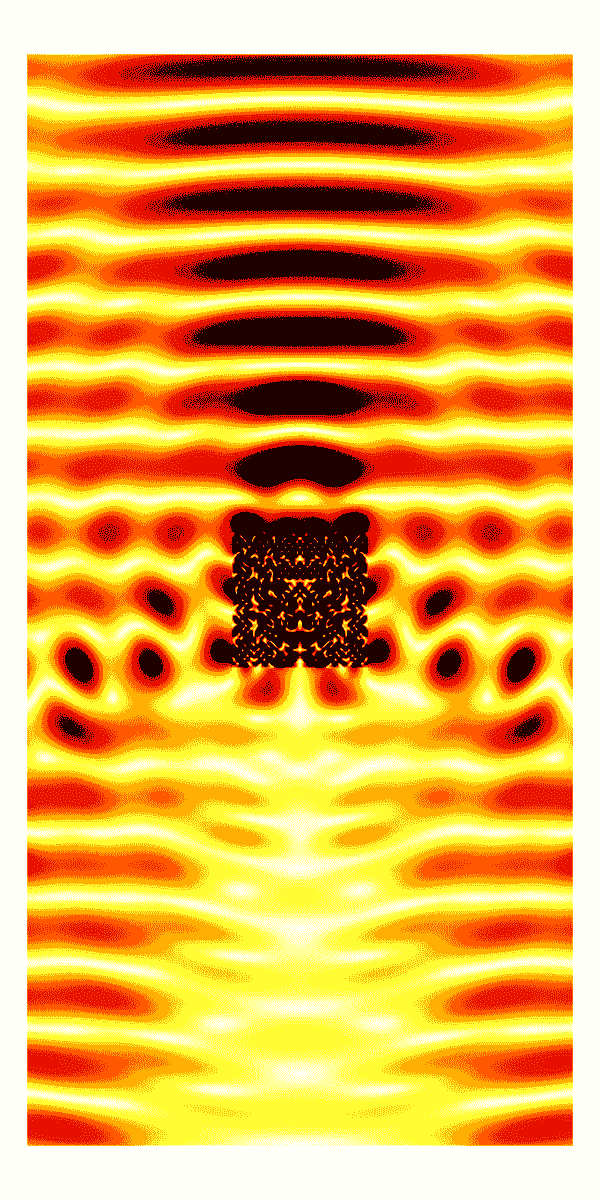}
    \end{subfigure}
    \begin{subfigure}{0.14\textwidth}
        \includegraphics[width=\textwidth]{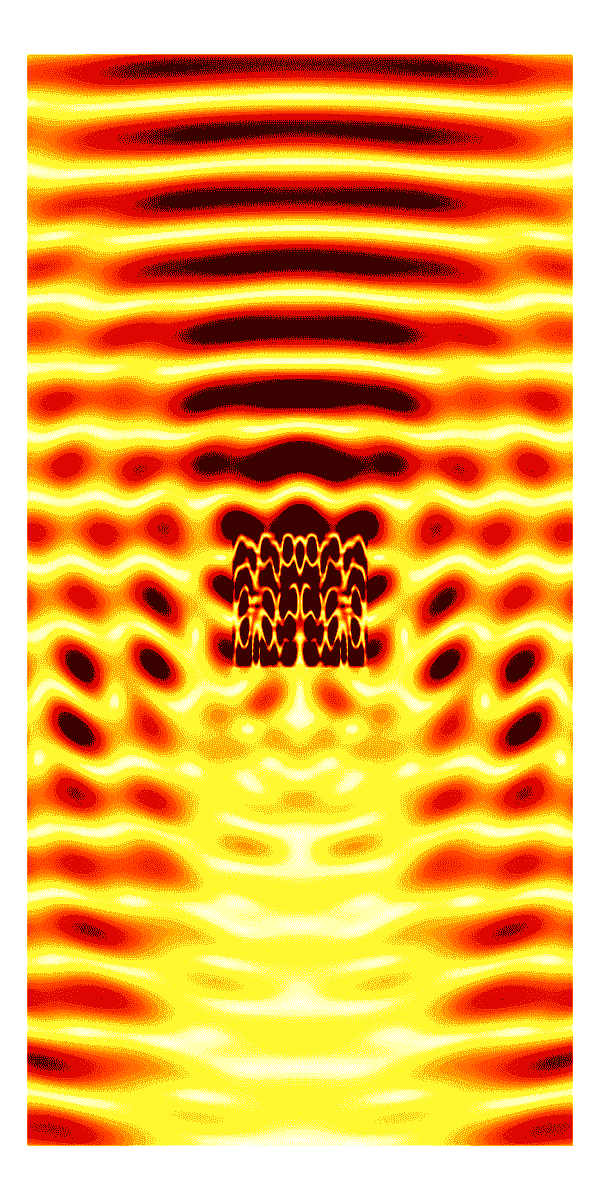}
    \end{subfigure}
    \begin{subfigure}{0.14\textwidth}
        \includegraphics[width=\textwidth]{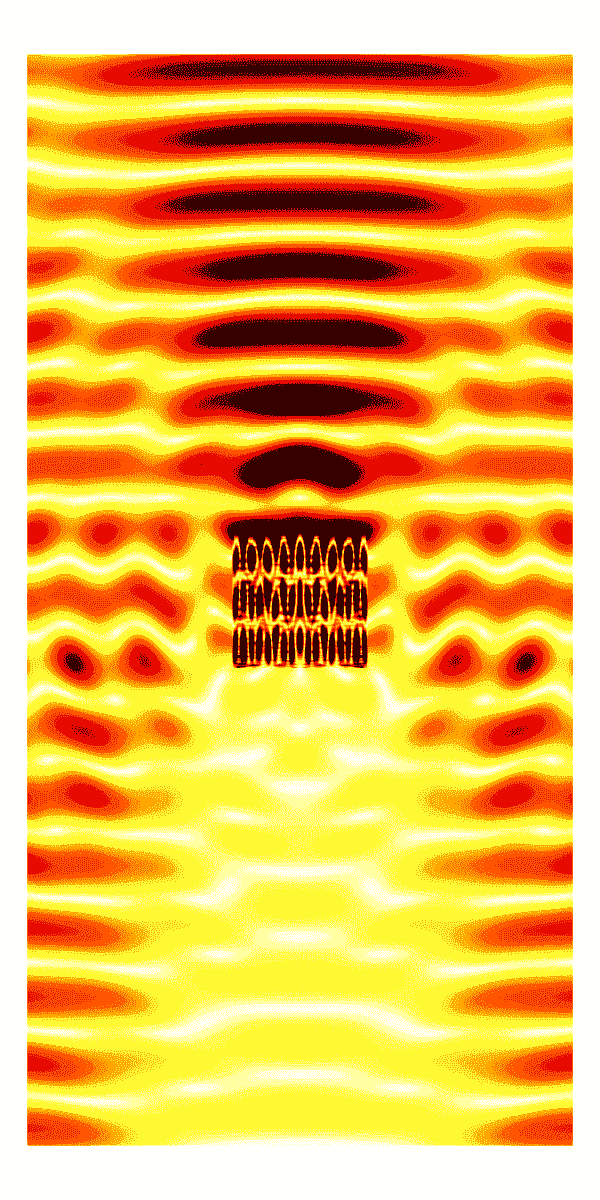}
    \end{subfigure}
    \begin{subfigure}{0.14\textwidth}
        \includegraphics[width=\textwidth]{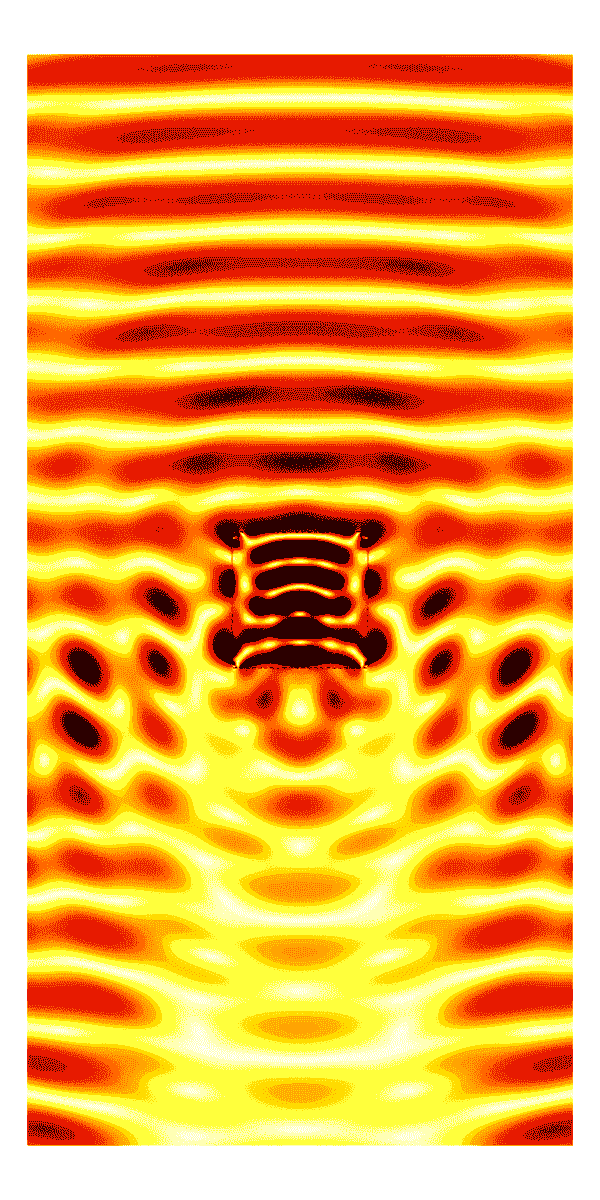}
    \end{subfigure}
    \begin{subfigure}{0.14\textwidth}
        \includegraphics[width=\textwidth]{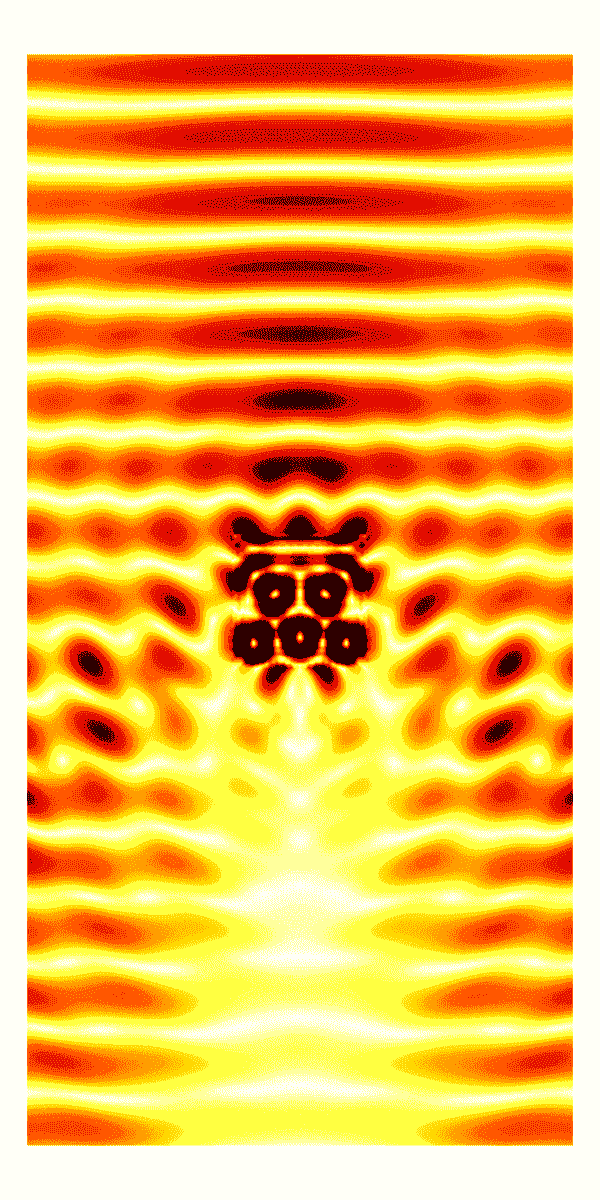}
    \end{subfigure}
        
	\begin{subfigure}{0.0\textwidth}
    \begin{picture}(0,0)
        \put(-6,40){\rotatebox{90}{\textbf{$62.83 \cdot 10^{5}$ rad/sec}}}
    \end{picture}
	\end{subfigure}
    \begin{subfigure}{0.14\textwidth}
        \includegraphics[width=\textwidth]{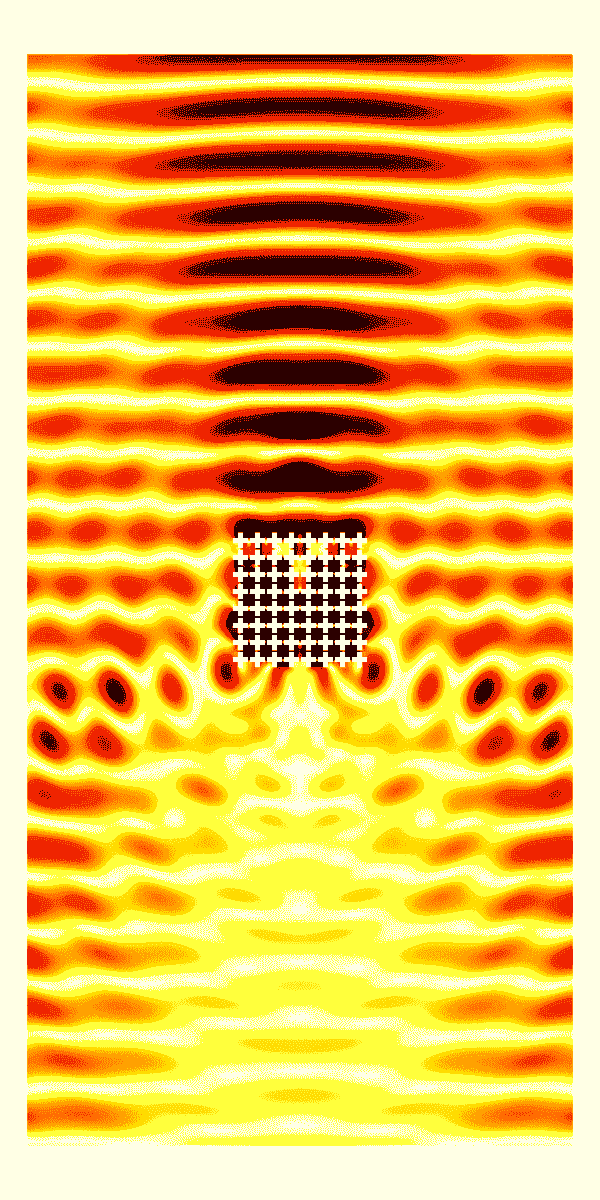}
    \caption*{\footnotesize microstructured}
    \end{subfigure}
    \begin{subfigure}{0.14\textwidth}
        \includegraphics[width=\textwidth]{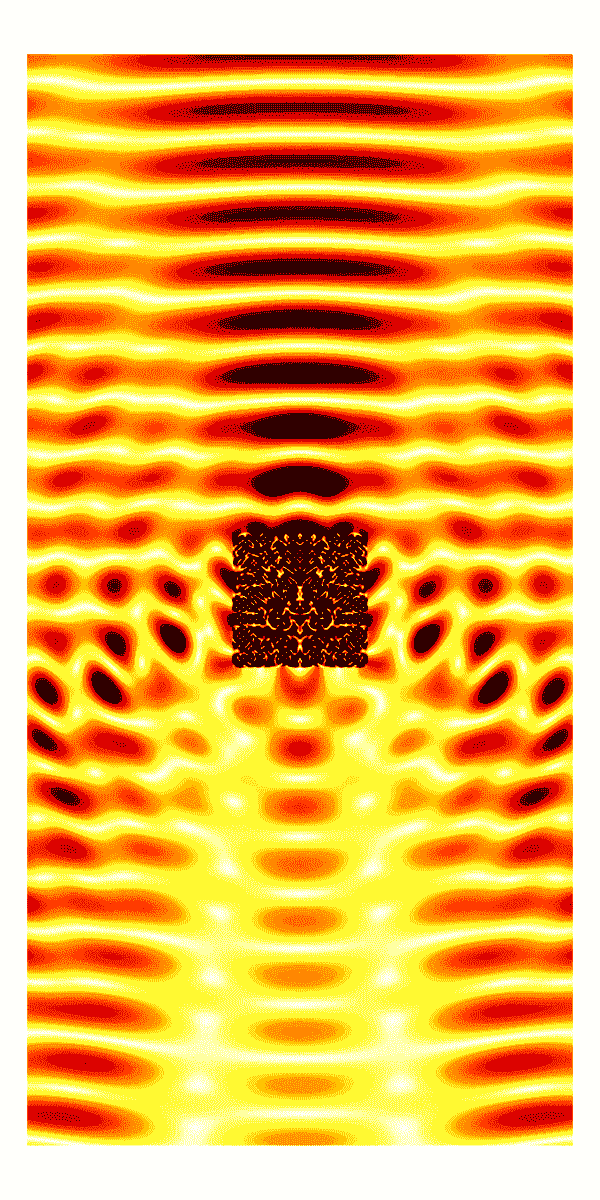}
   \caption*{\footnotesize macro-Cauchy}
    \end{subfigure}
    \begin{subfigure}{0.14\textwidth}
        \includegraphics[width=\textwidth]{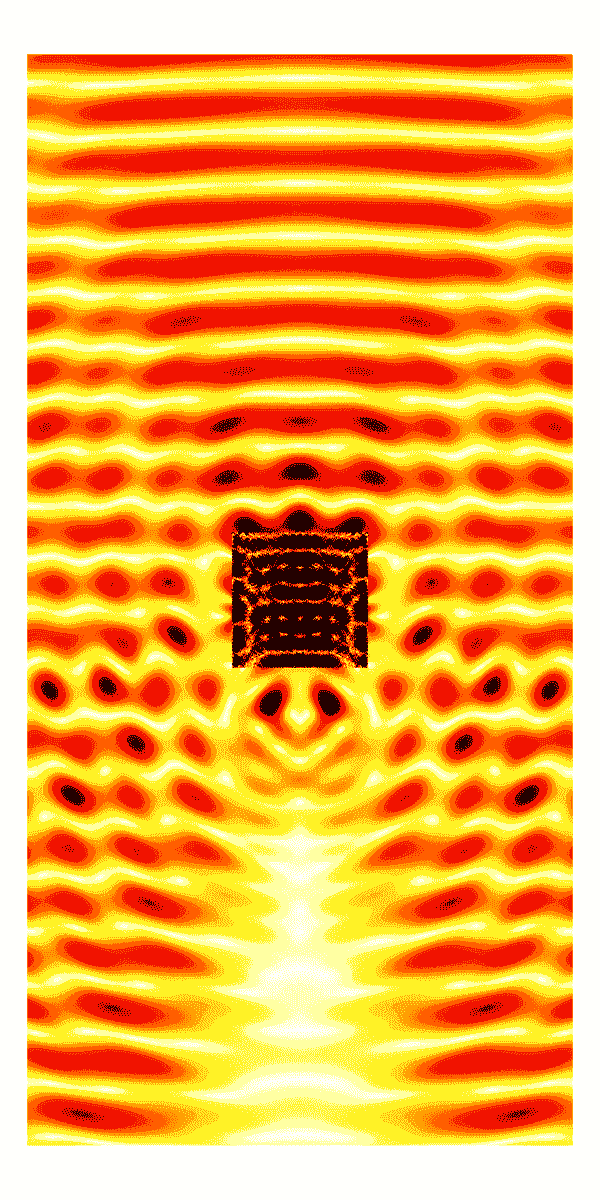}
    \caption*{\footnotesize RRMM(1)}
    \end{subfigure}
    \begin{subfigure}{0.14\textwidth}
        \includegraphics[width=\textwidth]{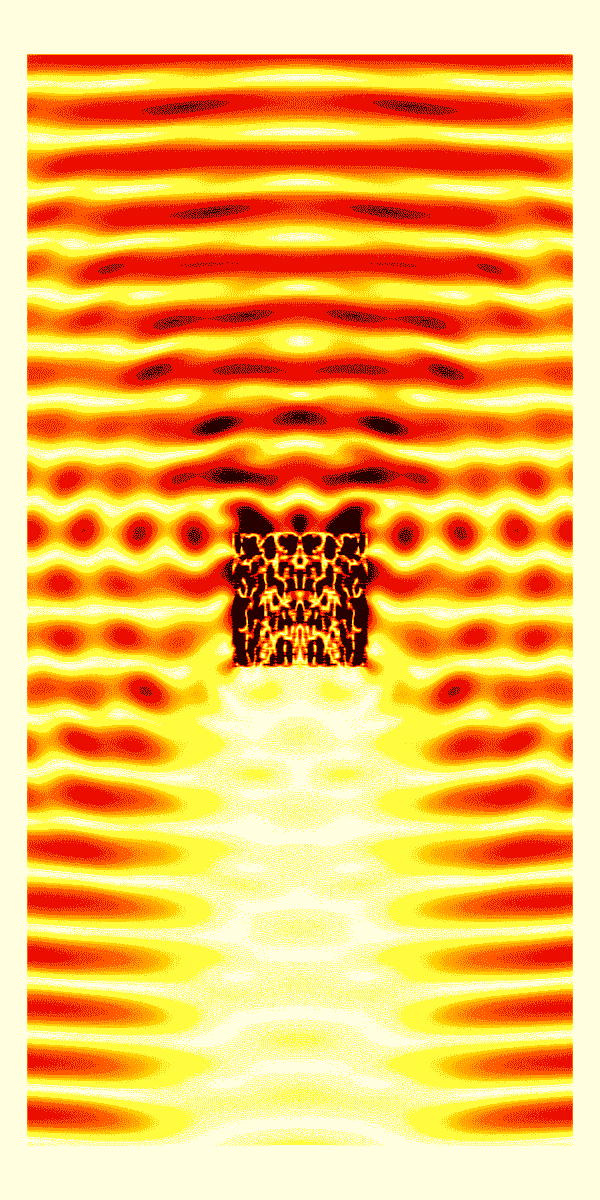}
    \caption*{\footnotesize RMM(1)}  
    \end{subfigure}
    \begin{subfigure}{0.14\textwidth}
        \includegraphics[width=\textwidth]{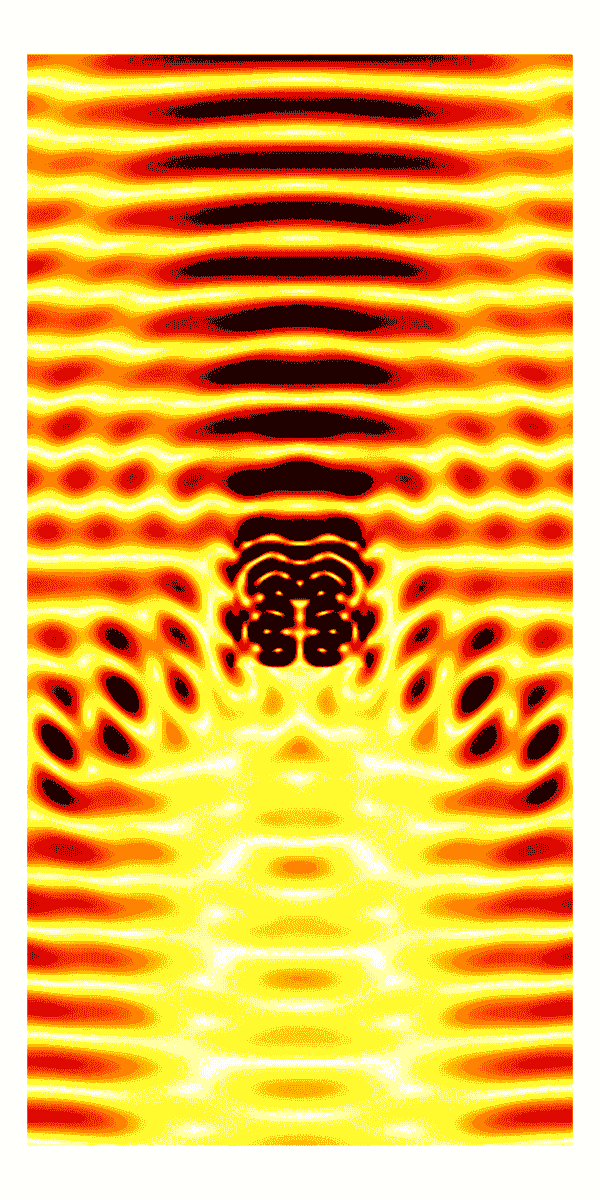}
    \caption*{\footnotesize RRMM(2)}  
    \end{subfigure}
    \begin{subfigure}{0.14\textwidth}
        \includegraphics[width=\textwidth]{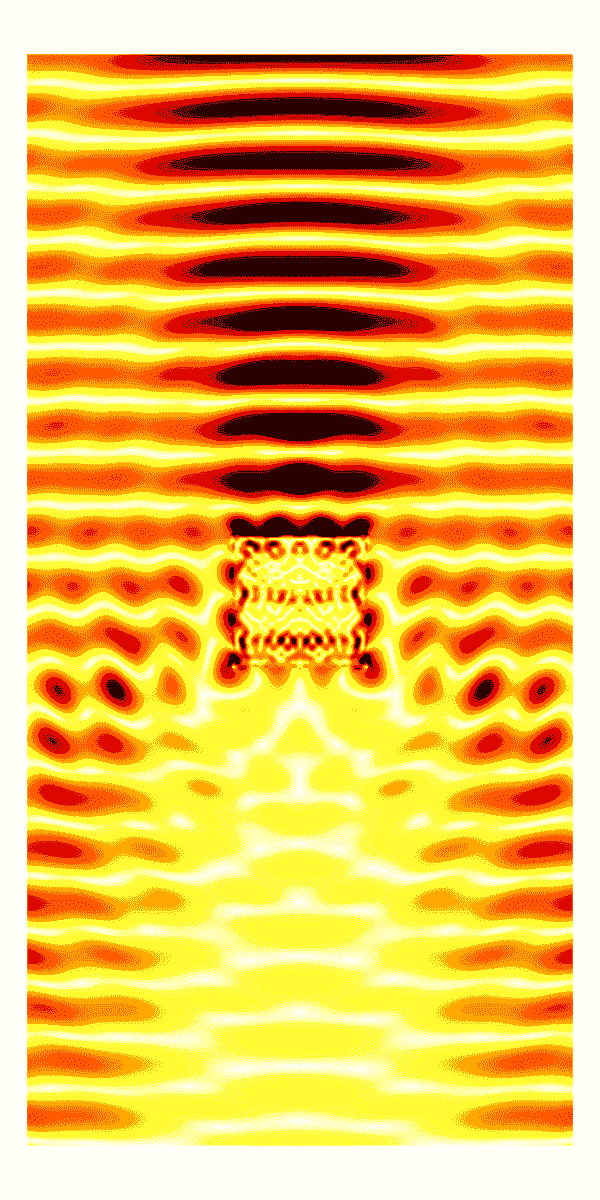}
    \caption*{\footnotesize RMM(2)}  
    \end{subfigure}
\caption{Results of finite-size scattering pattern for an incident pressure wave (lower frequency range up to the lower band-gap limit) with the material parameters obtained by fitting the dispersion curves in one direction at $0^\circ$ (marked as 1) and in two directions at $0^\circ$  and $45^\circ$  (marked as 2).}
\label{fig:pre1}
\end{figure}

\begin{figure}[!ht]
    \centering
    \begin{subfigure}{0.7\textwidth}
        \includegraphics[width=\textwidth]{figures/legend.jpg}
        \put(0,10){\textbf{\large ${\lvert u \lvert}/{u_0}$}}
    \end{subfigure}
    
	\begin{subfigure}{0.0\textwidth}
    \begin{picture}(0,0)
        \put(-6,30){\rotatebox{90}{\textbf{$75.4 \cdot 10^{5}$ rad/sec}}}
    \end{picture}
	\end{subfigure}
    \begin{subfigure}{0.14\textwidth}
        \includegraphics[width=\textwidth]{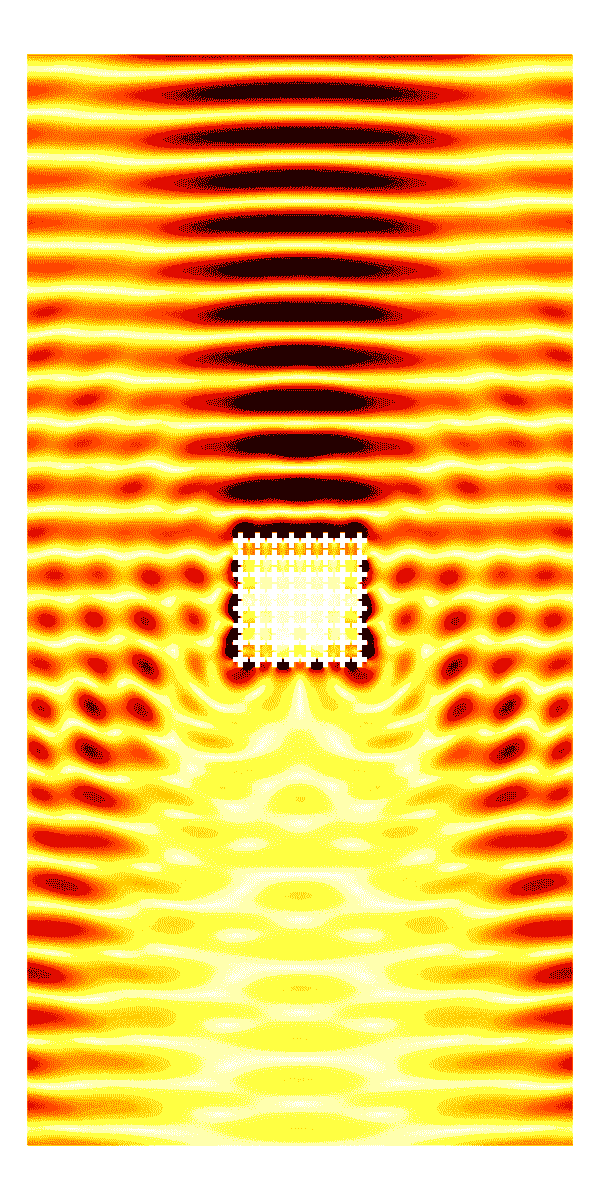}
    \end{subfigure}
    \begin{subfigure}{0.14\textwidth}
        \includegraphics[width=\textwidth]{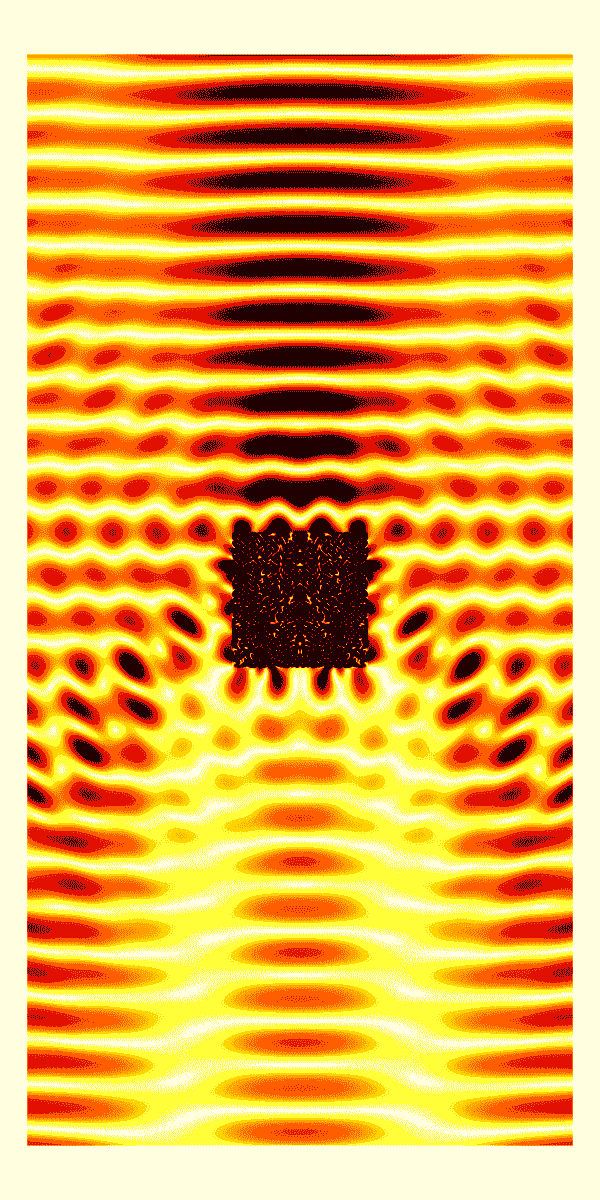}
    \end{subfigure}
    \begin{subfigure}{0.14\textwidth}
        \includegraphics[width=\textwidth]{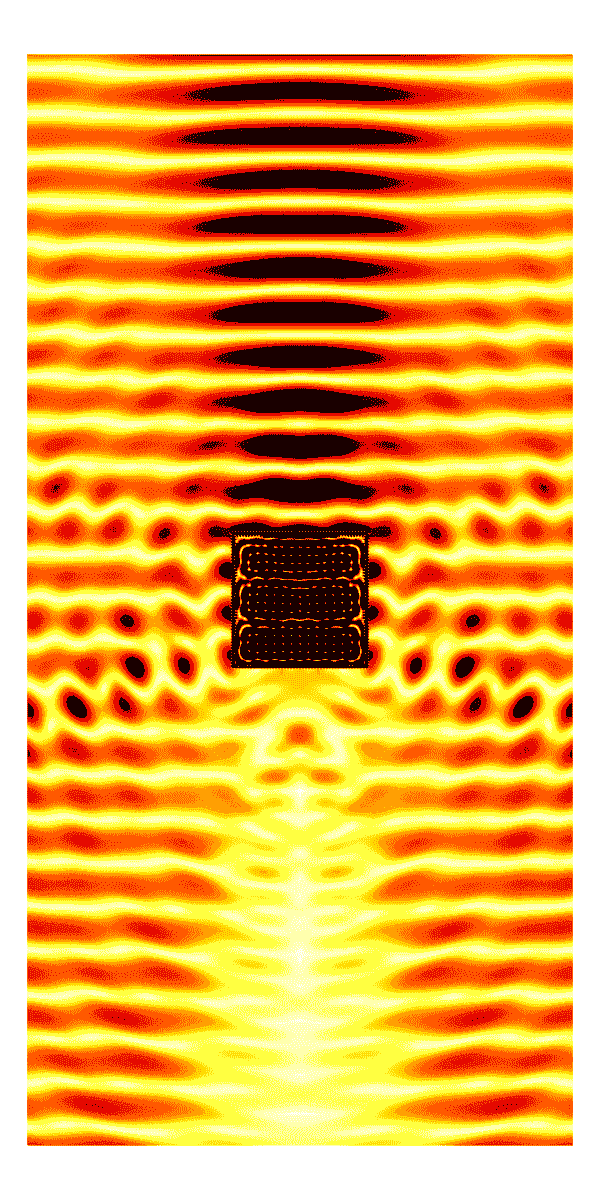}
    \end{subfigure}
    \begin{subfigure}{0.14\textwidth}
        \includegraphics[width=\textwidth]{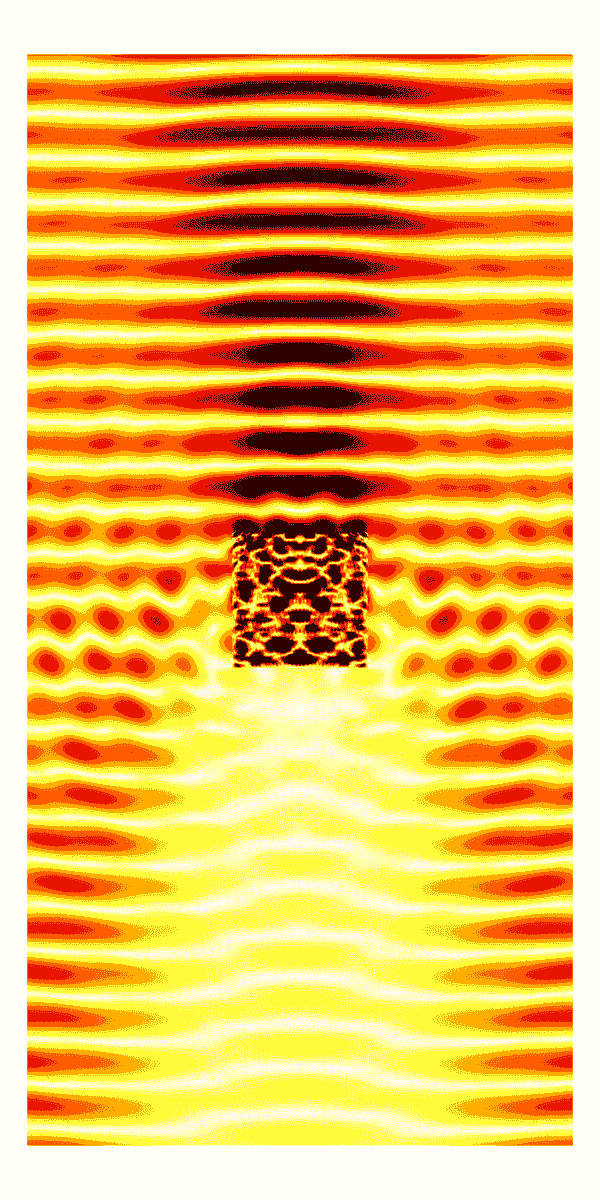}
    \end{subfigure}
    \begin{subfigure}{0.14\textwidth}
        \includegraphics[width=\textwidth]{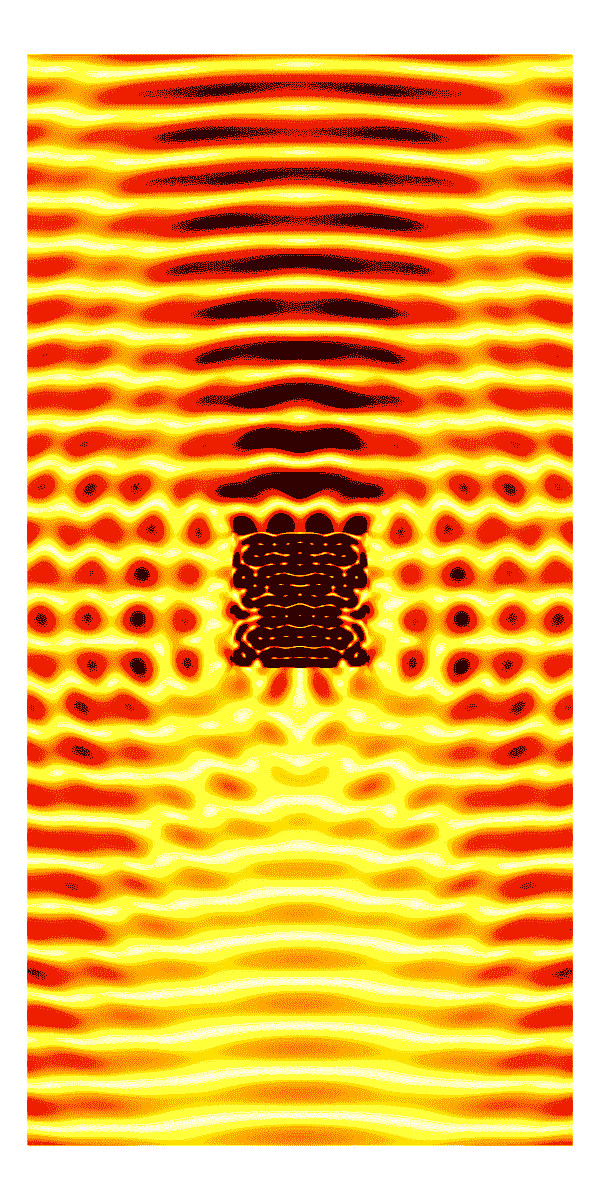}
    \end{subfigure}
    \begin{subfigure}{0.14\textwidth}
        \includegraphics[width=\textwidth]{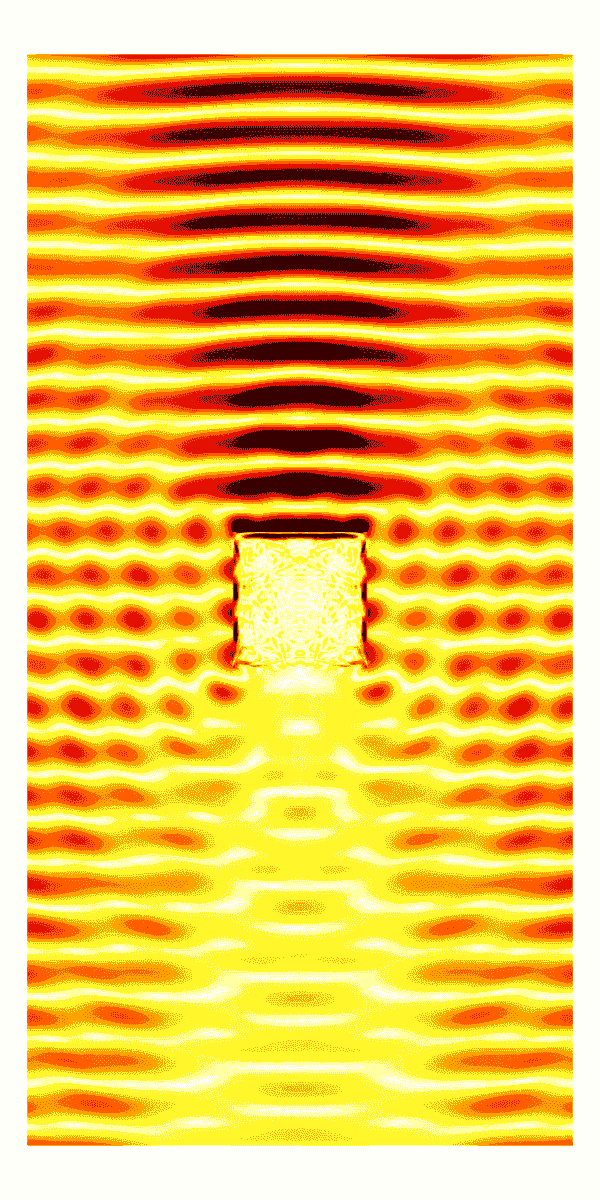}
    \end{subfigure}

	\begin{subfigure}{0.0\textwidth}
    \begin{picture}(0,0)
        \put(-6,30){\rotatebox{90}{\textbf{$87.96 \cdot 10^{5}$ rad/sec}}}
    \end{picture}
	\end{subfigure}
    \begin{subfigure}{0.14\textwidth}
        \includegraphics[width=\textwidth]{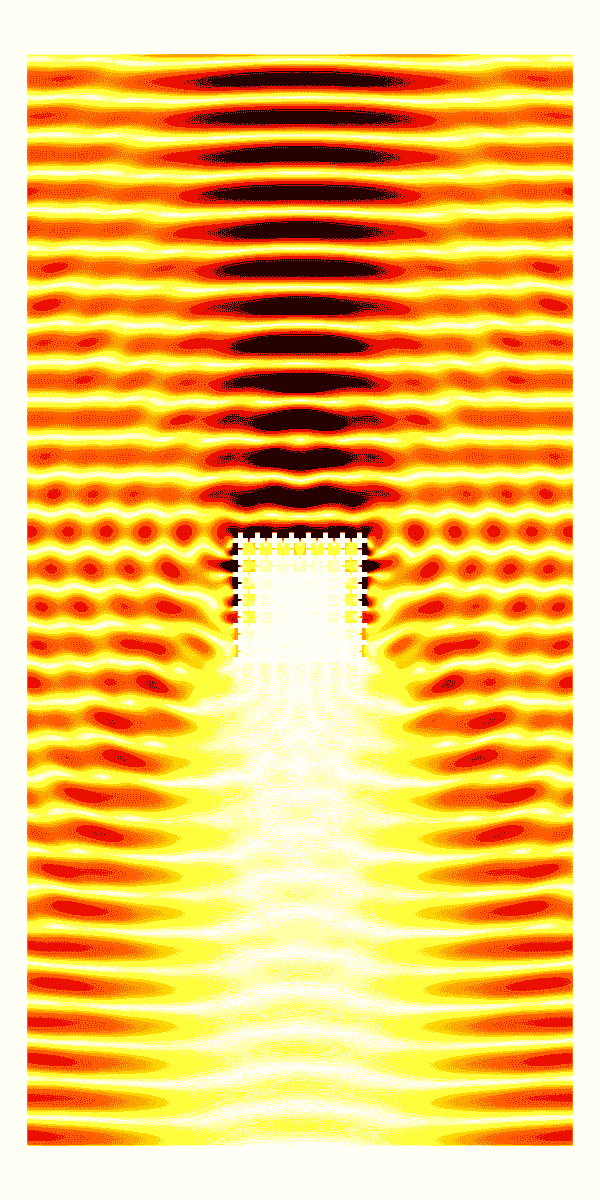}
    \end{subfigure}
    \begin{subfigure}{0.14\textwidth}
        \includegraphics[width=\textwidth]{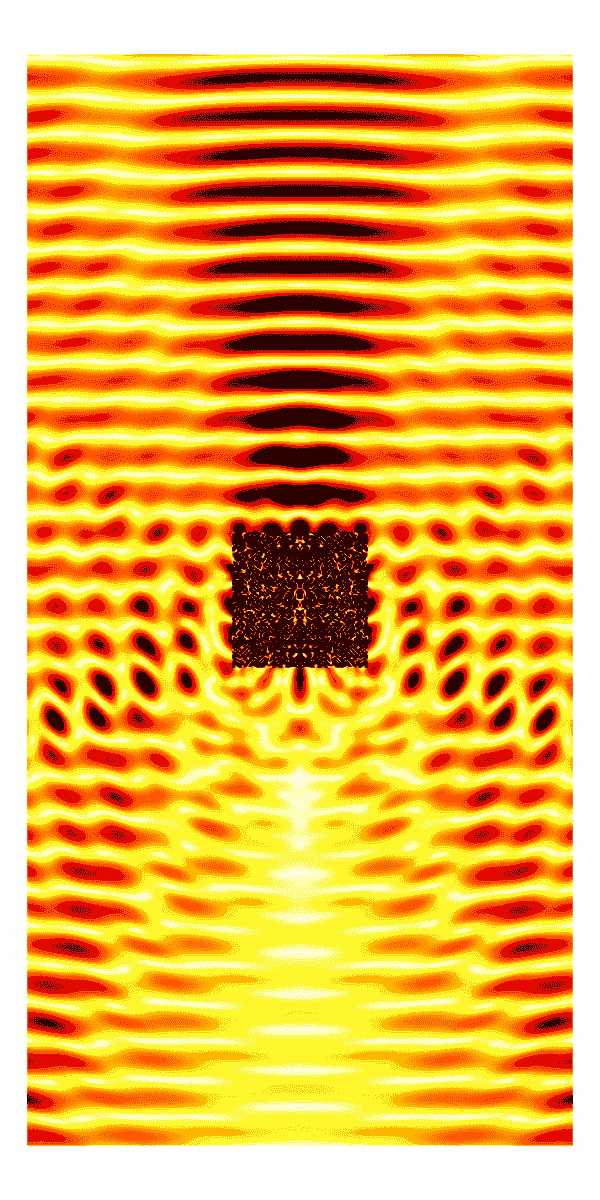}
    \end{subfigure}
    \begin{subfigure}{0.14\textwidth}
        \includegraphics[width=\textwidth]{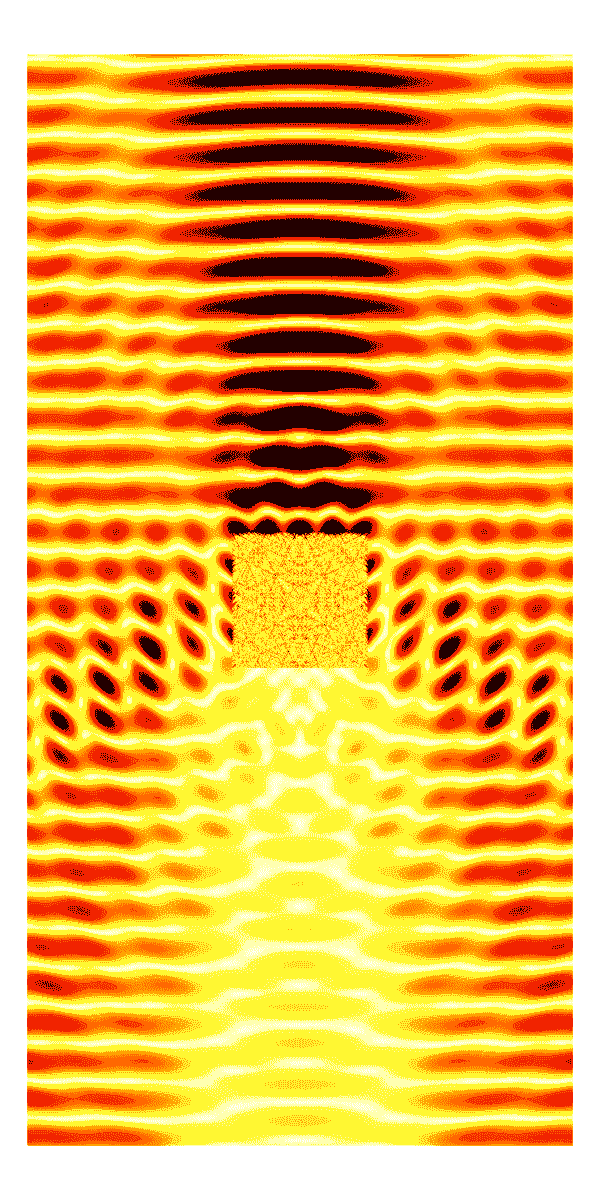}
    \end{subfigure}
    \begin{subfigure}{0.14\textwidth}
        \includegraphics[width=\textwidth]{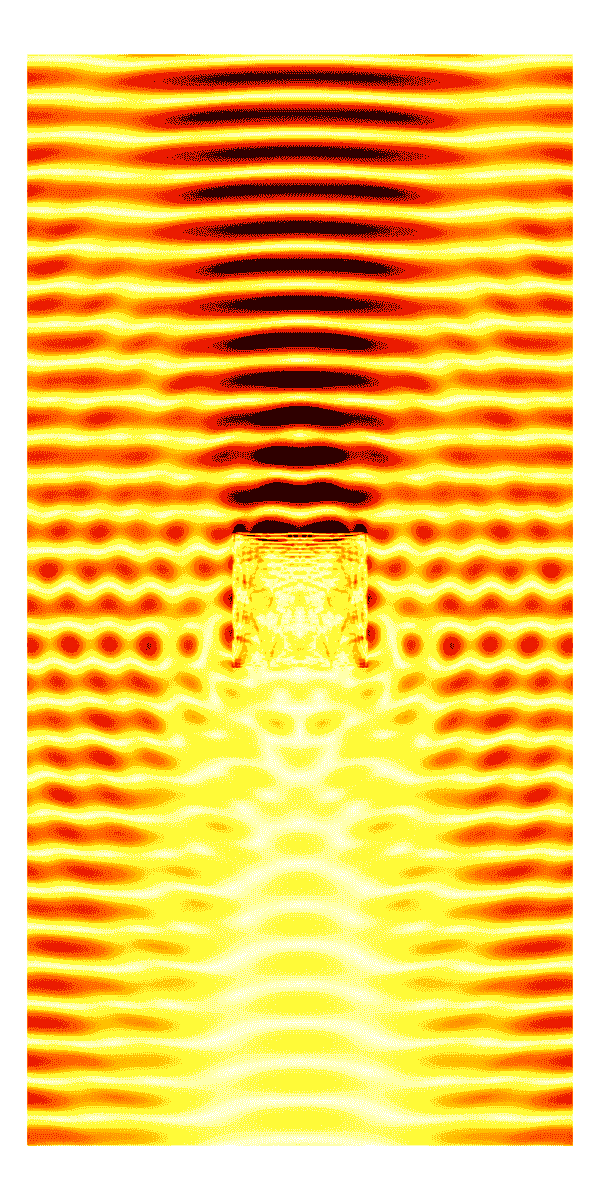}
    \end{subfigure}
    \begin{subfigure}{0.14\textwidth}
        \includegraphics[width=\textwidth]{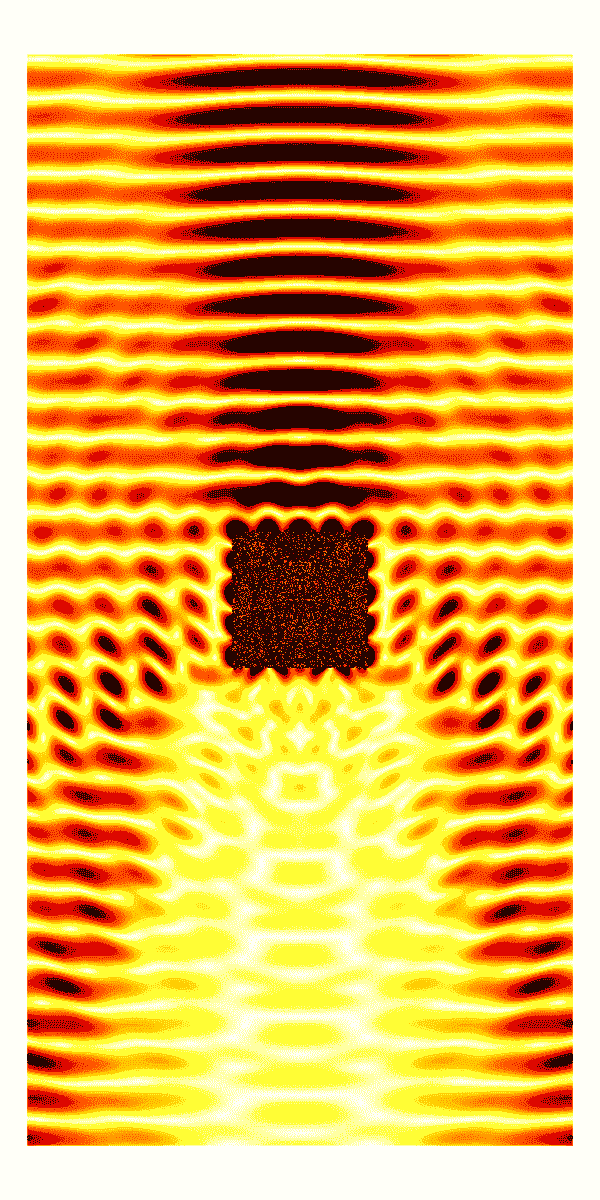}
    \end{subfigure}
    \begin{subfigure}{0.14\textwidth}
        \includegraphics[width=\textwidth]{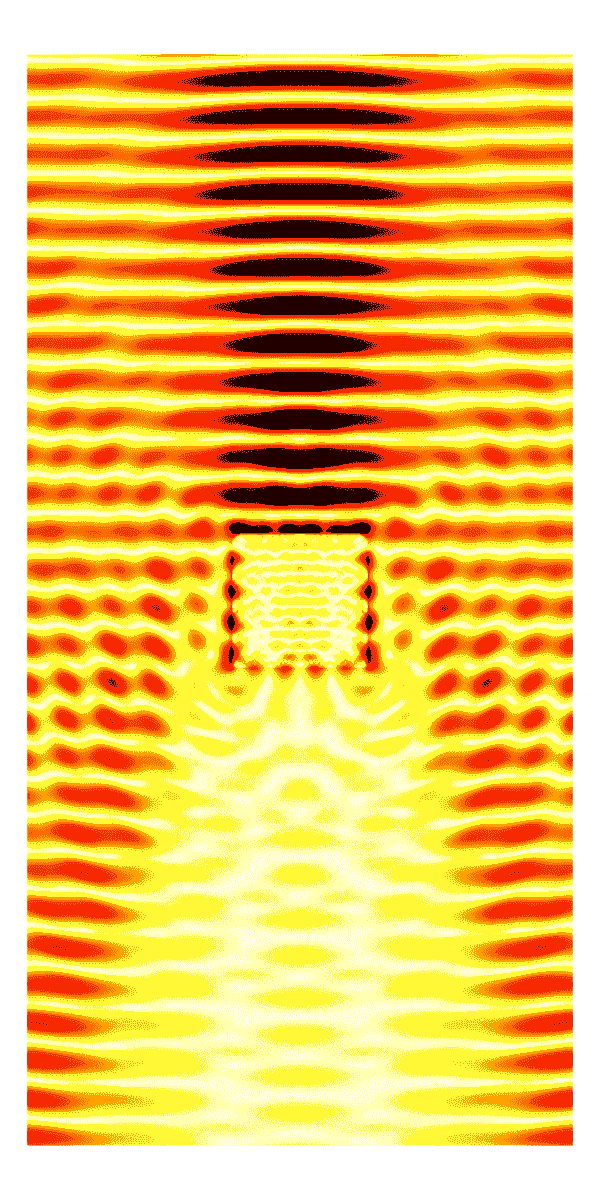}
    \end{subfigure}

	\begin{subfigure}{0.0\textwidth}
    \begin{picture}(0,0)
        \put(-6,30){\rotatebox{90}{\textbf{$10.05 \cdot 10^{6}$ rad/sec}}}
    \end{picture}
	\end{subfigure}
    \begin{subfigure}{0.14\textwidth}
        \includegraphics[width=\textwidth]{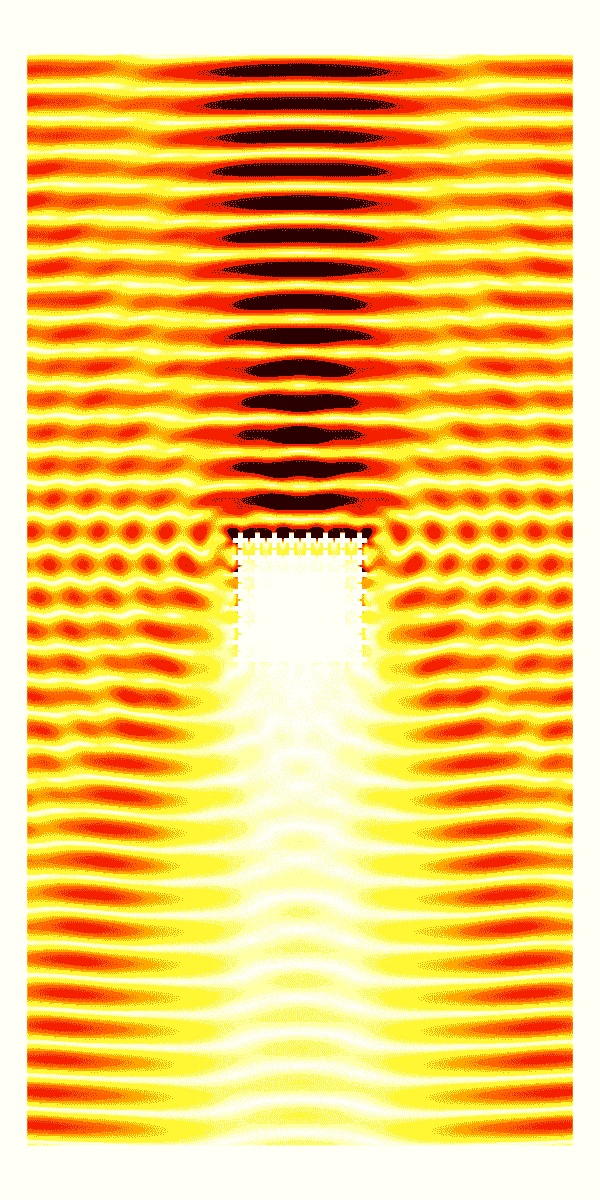}
    \end{subfigure}
    \begin{subfigure}{0.14\textwidth}
        \includegraphics[width=\textwidth]{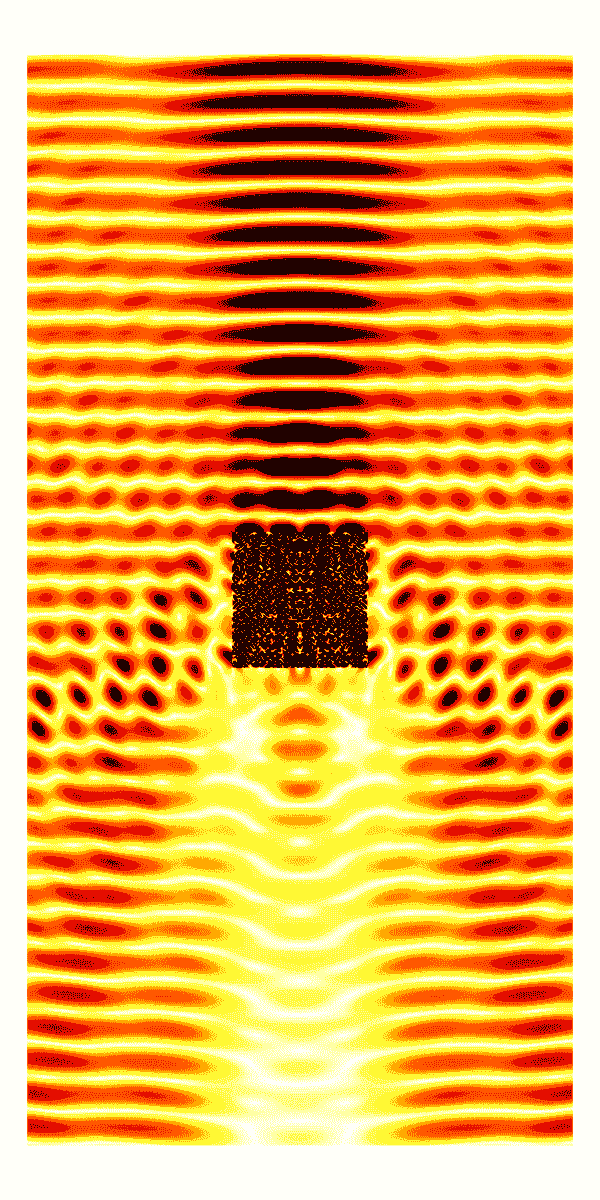}
    \end{subfigure}
    \begin{subfigure}{0.14\textwidth}
        \includegraphics[width=\textwidth]{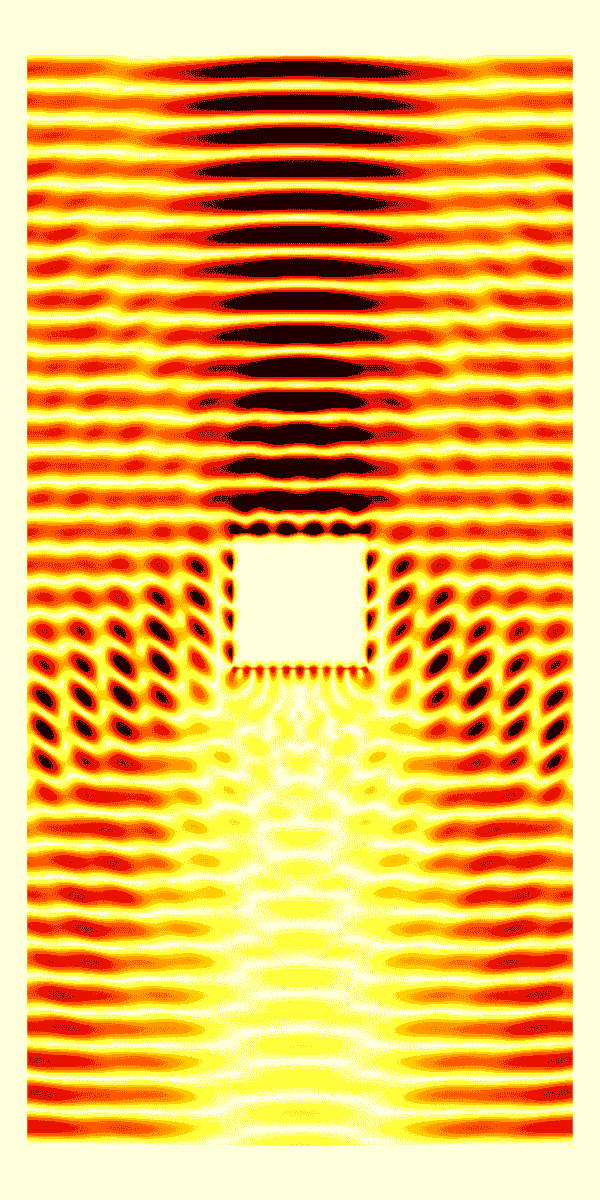}
    \end{subfigure}
    \begin{subfigure}{0.14\textwidth}
        \includegraphics[width=\textwidth]{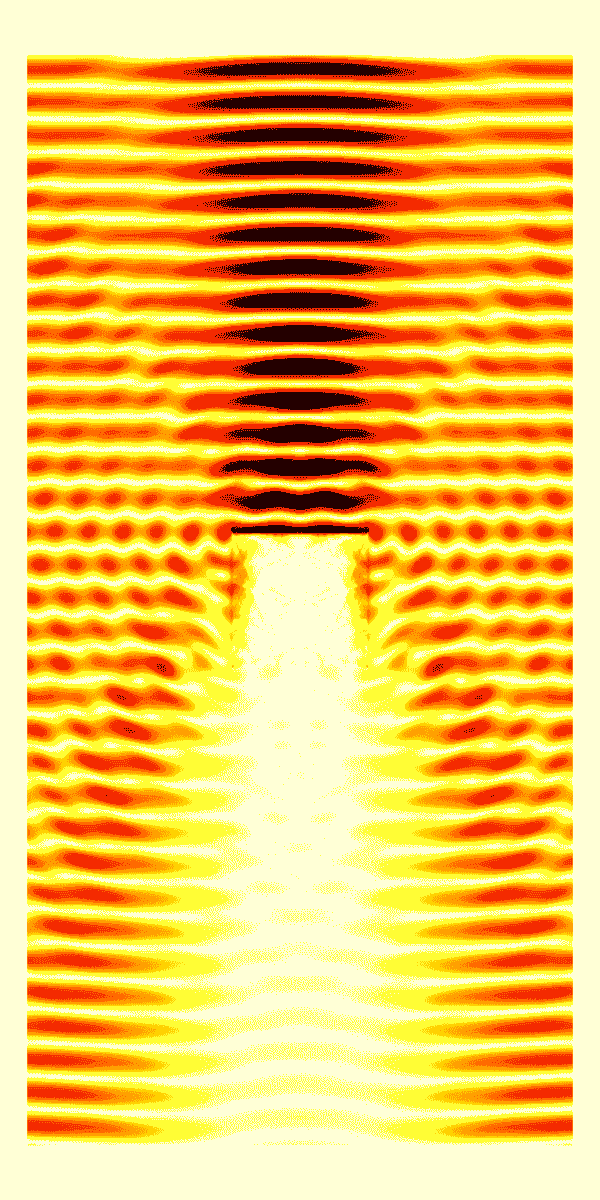}
    \end{subfigure}
    \begin{subfigure}{0.14\textwidth}
        \includegraphics[width=\textwidth]{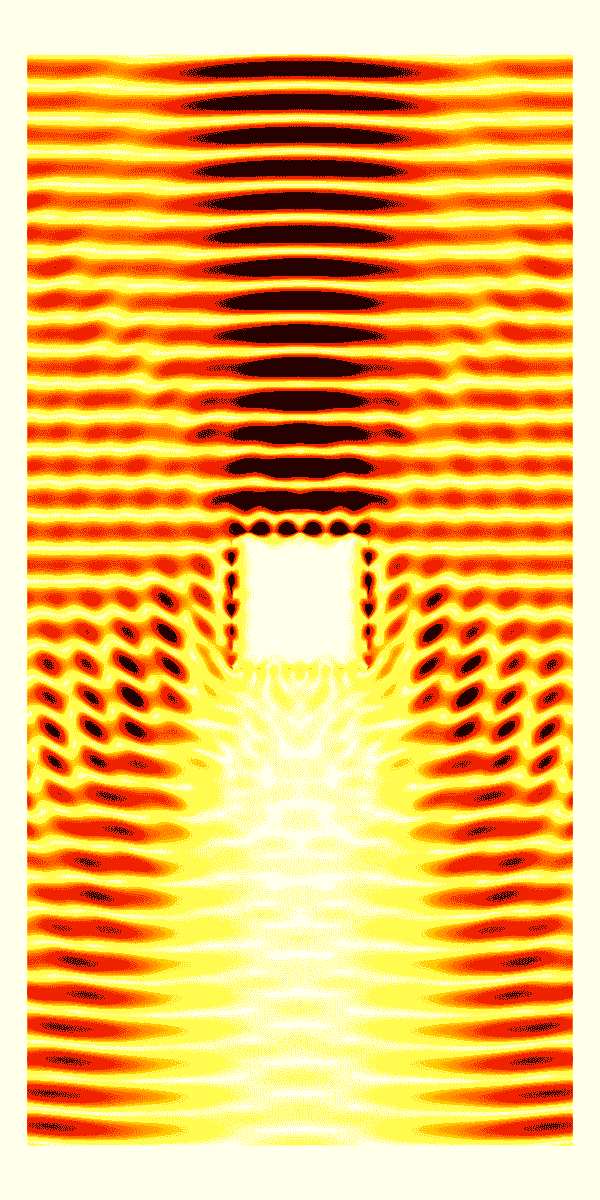}
    \end{subfigure}
    \begin{subfigure}{0.14\textwidth}
        \includegraphics[width=\textwidth]{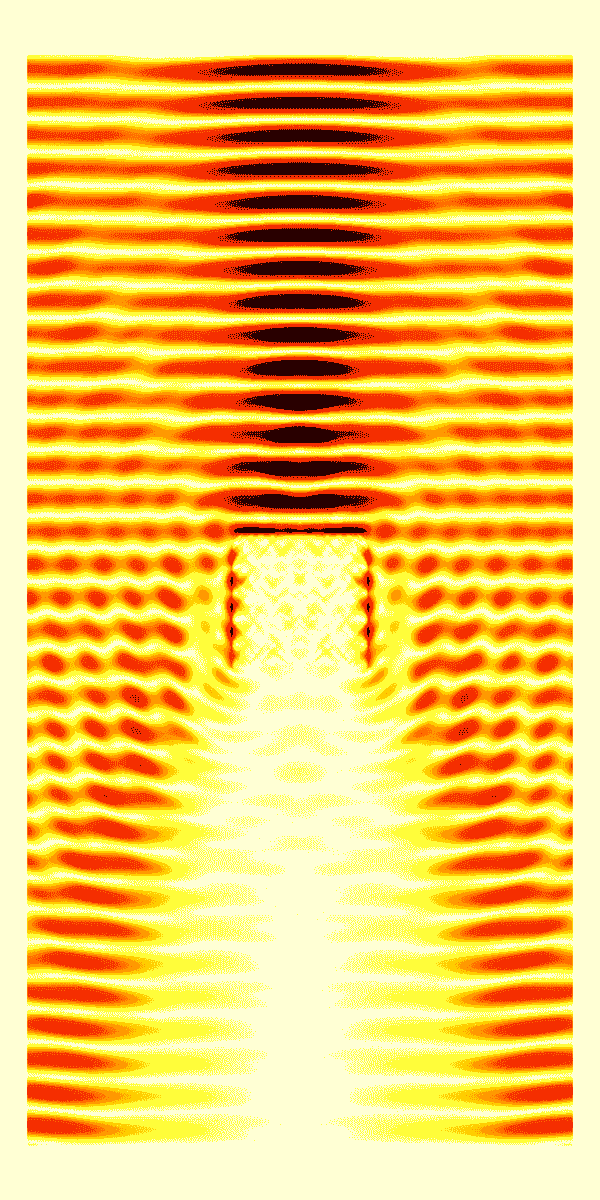}
    \end{subfigure}

	\begin{subfigure}{0.0\textwidth}
    \begin{picture}(0,0)
        \put(-6,30){\rotatebox{90}{\textbf{$11.31 \cdot 10^{6}$ rad/sec}}}
    \end{picture}
	\end{subfigure}
    \begin{subfigure}{0.14\textwidth}
        \includegraphics[width=\textwidth]{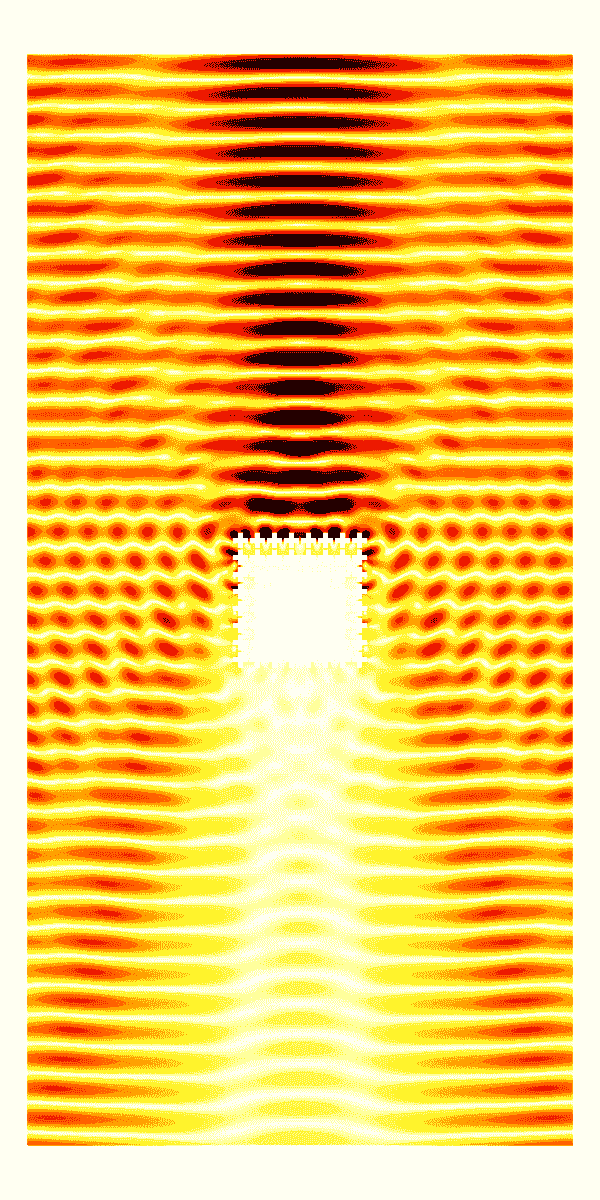}
    \end{subfigure}
    \begin{subfigure}{0.14\textwidth}
        \includegraphics[width=\textwidth]{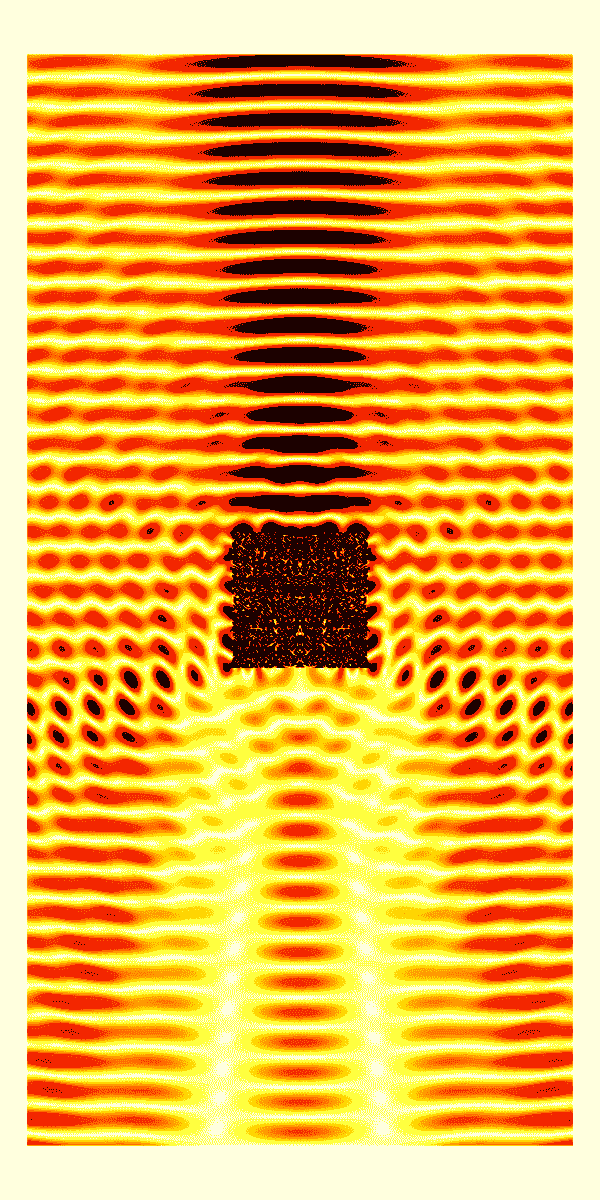}
    \end{subfigure}
    \begin{subfigure}{0.14\textwidth}
        \includegraphics[width=\textwidth]{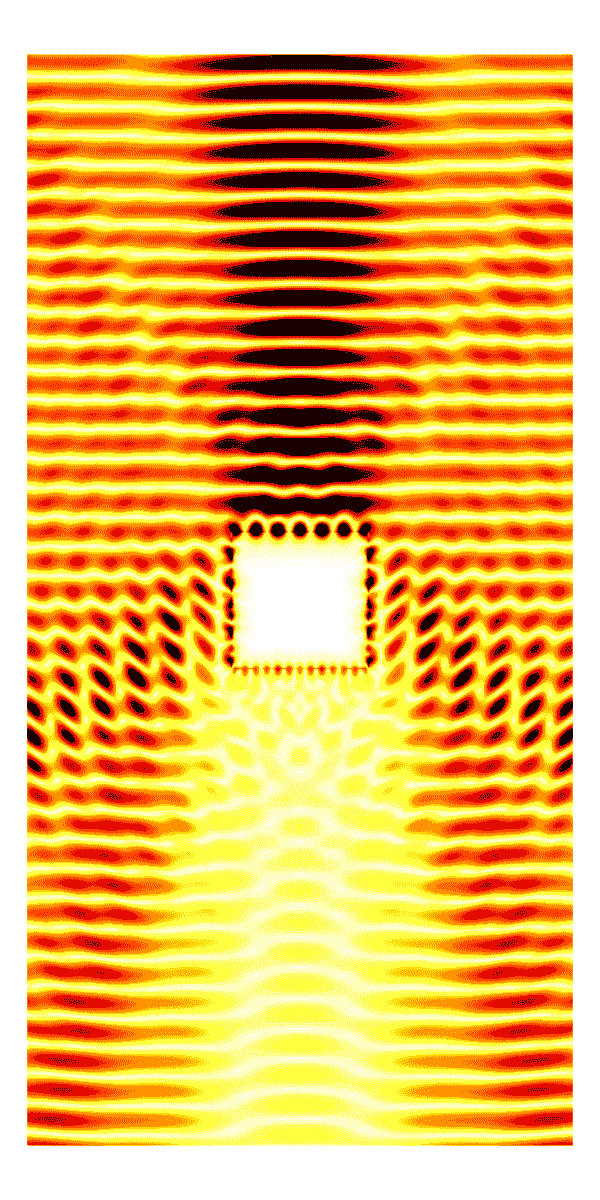}
    \end{subfigure}
    \begin{subfigure}{0.14\textwidth}
        \includegraphics[width=\textwidth]{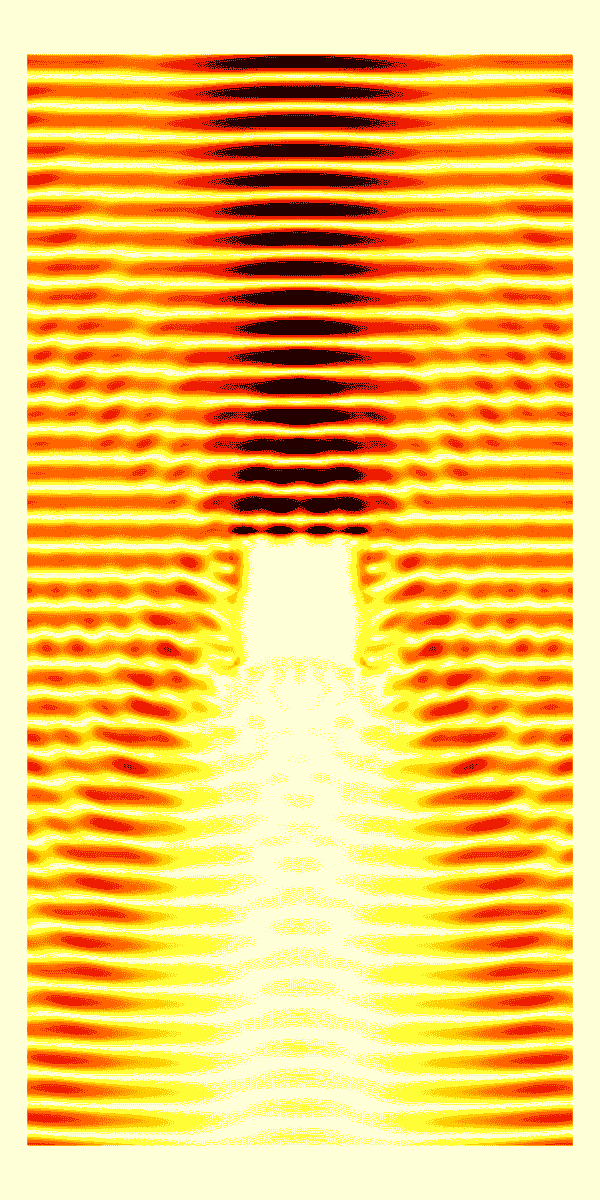}
    \end{subfigure}
    \begin{subfigure}{0.14\textwidth}
        \includegraphics[width=\textwidth]{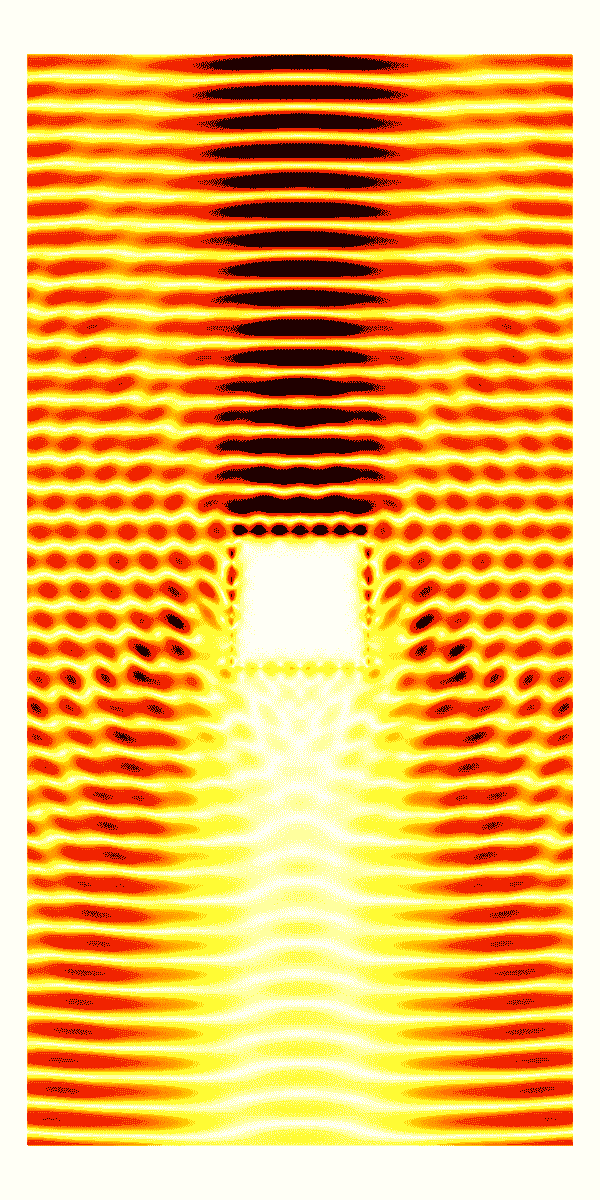}
    \end{subfigure}
    \begin{subfigure}{0.14\textwidth}
        \includegraphics[width=\textwidth]{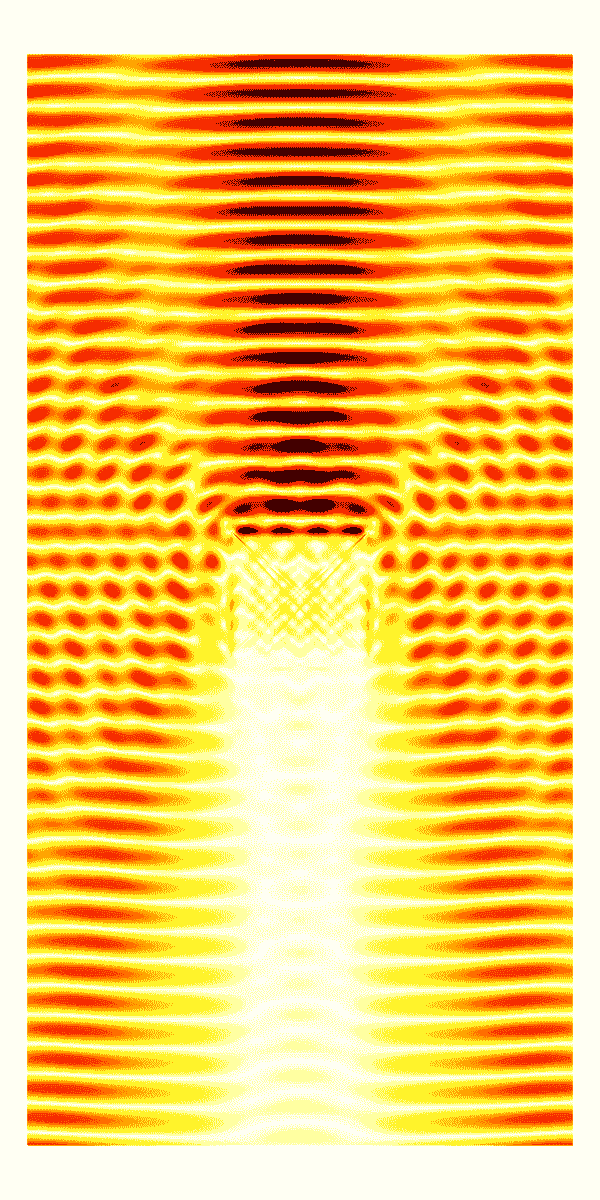}
    \end{subfigure}

	\begin{subfigure}{0.0\textwidth}
    \begin{picture}(0,0)
        \put(-6,40){\rotatebox{90}{\textbf{$12.56 \cdot 10^{6}$ rad/sec}}}
    \end{picture}
	\end{subfigure}
    \begin{subfigure}{0.14\textwidth}
        \includegraphics[width=\textwidth]{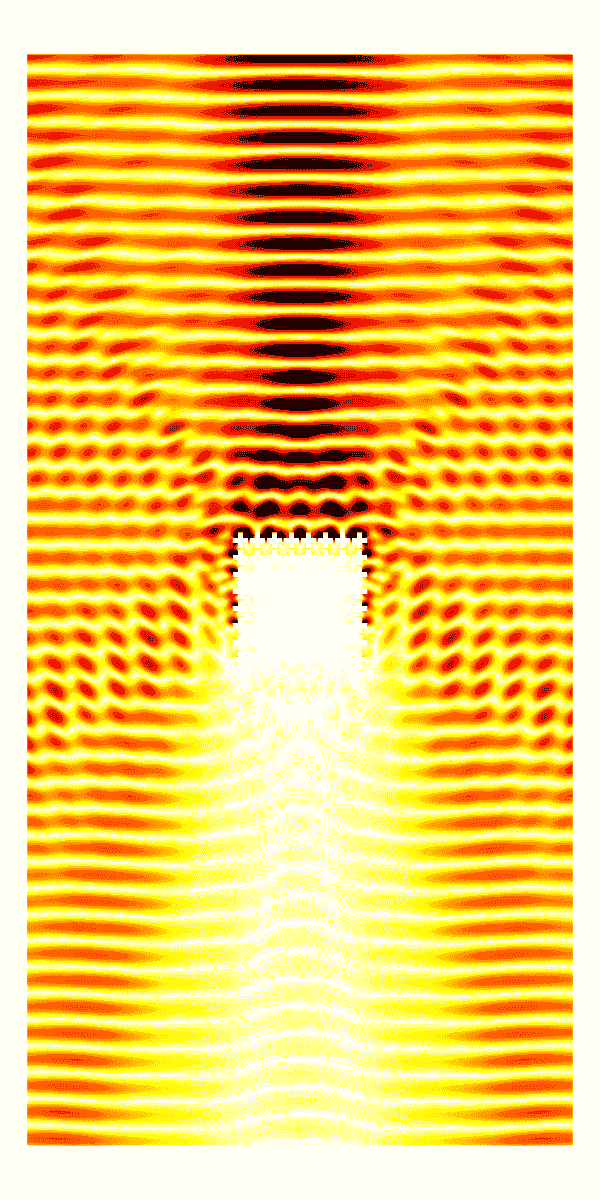}
    \caption*{\footnotesize microstructured}
    \end{subfigure}
    \begin{subfigure}{0.14\textwidth}
        \includegraphics[width=\textwidth]{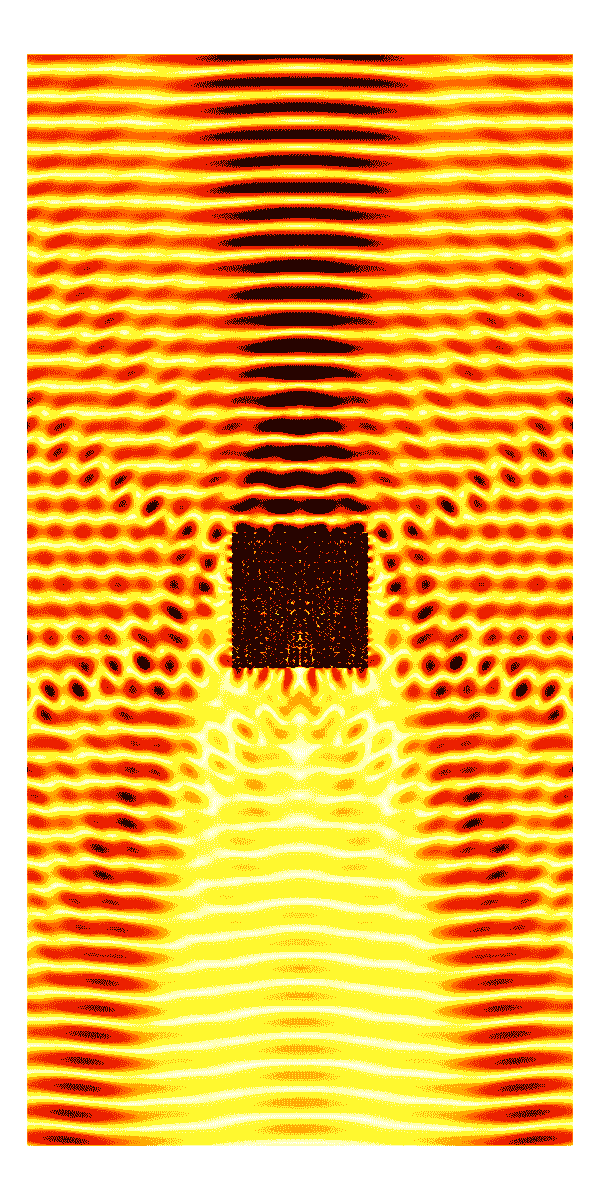}
   \caption*{\footnotesize macro-Cauchy}
    \end{subfigure}
    \begin{subfigure}{0.14\textwidth}
        \includegraphics[width=\textwidth]{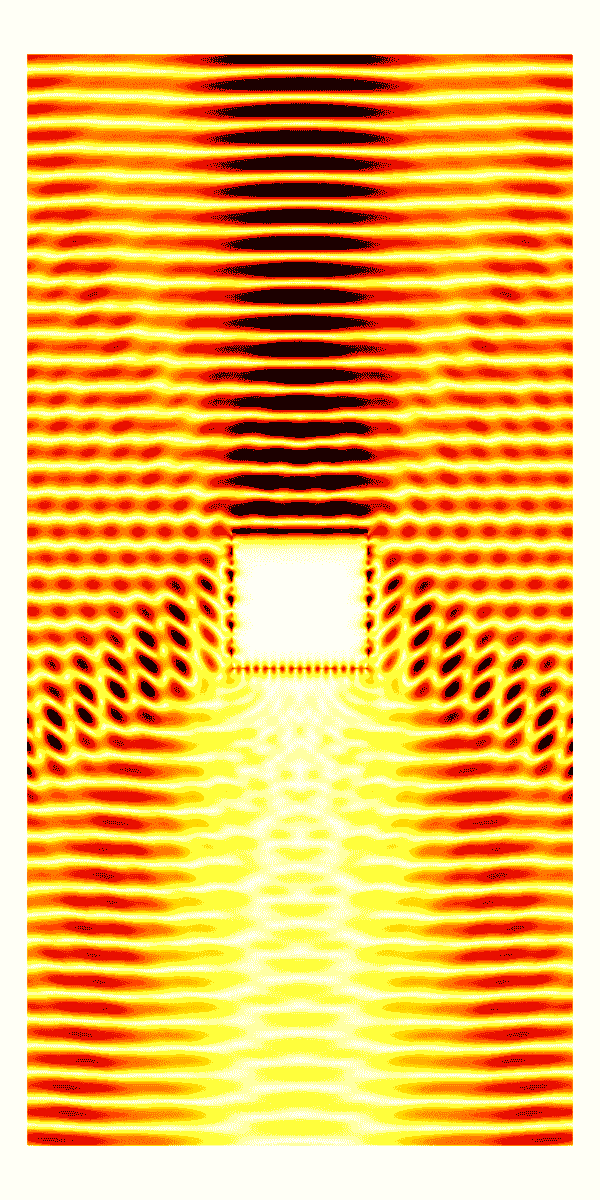}
    \caption*{\footnotesize RRMM(1)}
    \end{subfigure}
    \begin{subfigure}{0.14\textwidth}
        \includegraphics[width=\textwidth]{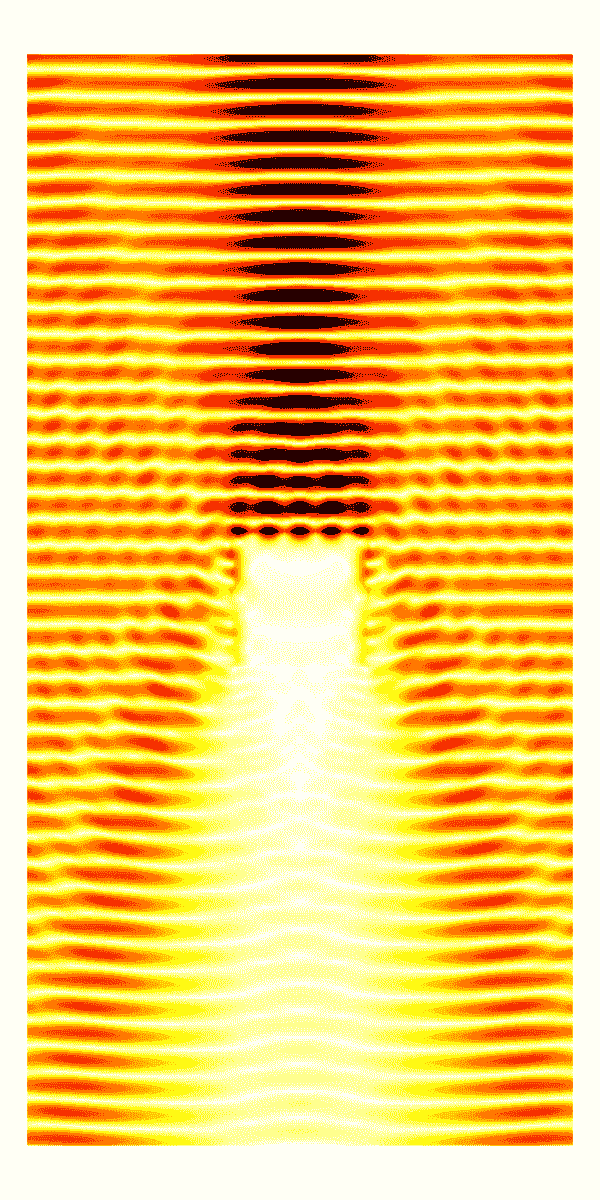}
    \caption*{\footnotesize RMM(1)}  
    \end{subfigure}
    \begin{subfigure}{0.14\textwidth}
        \includegraphics[width=\textwidth]{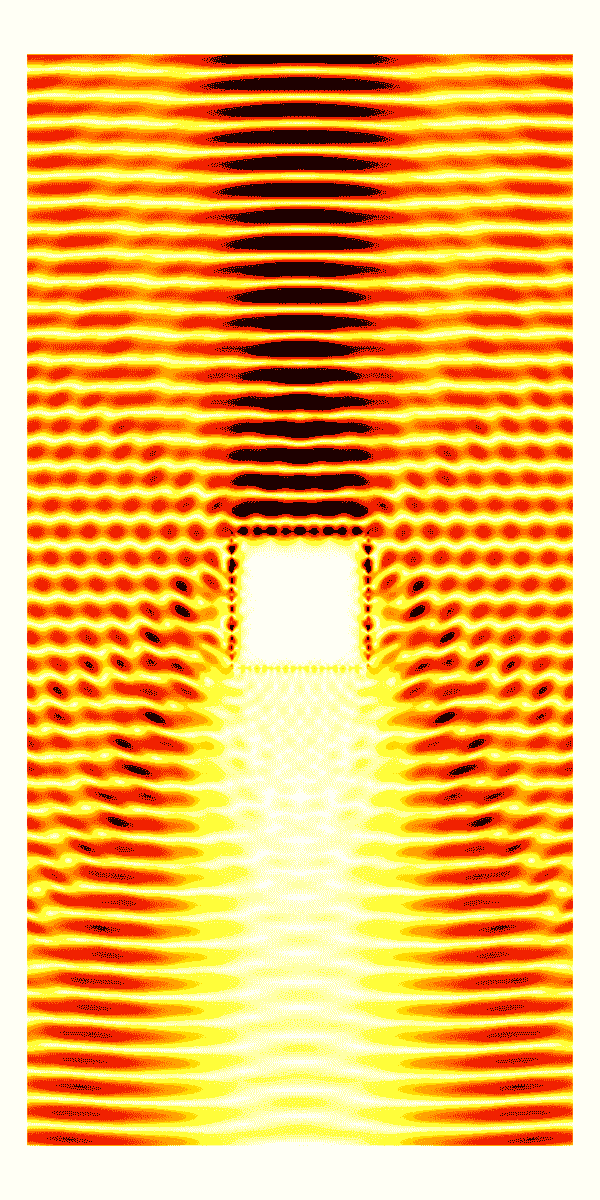}
     \caption*{\footnotesize RRMM(2)}  
    \end{subfigure}
    \begin{subfigure}{0.14\textwidth}
        \includegraphics[width=\textwidth]{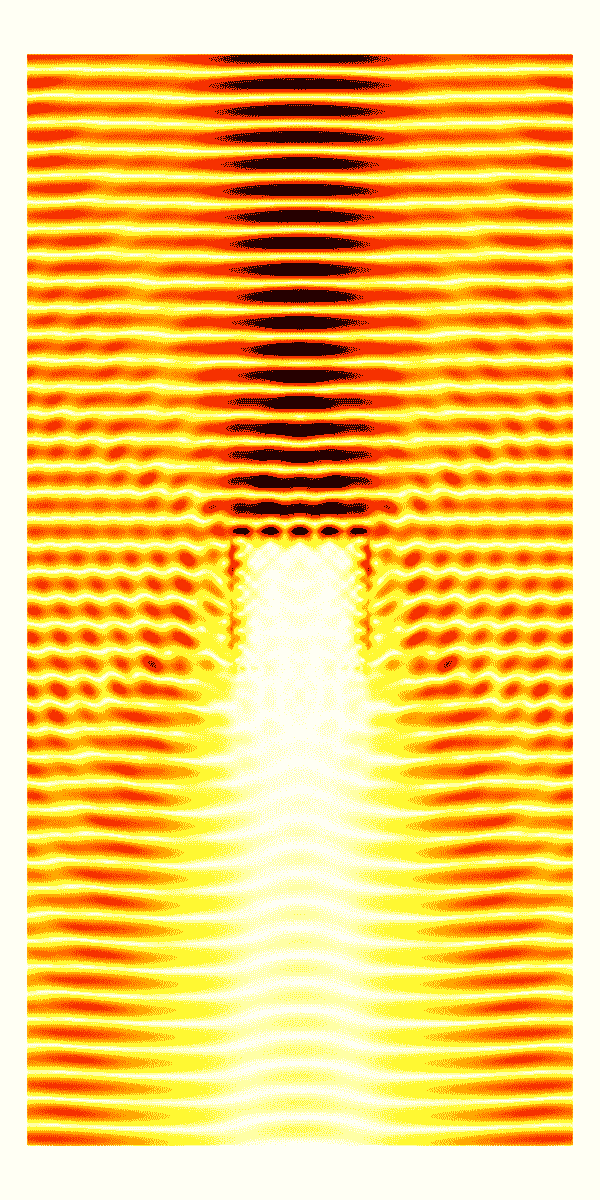}
     \caption*{\footnotesize RMM(2)}  
    \end{subfigure}
        
\caption{Results of finite-size scattering pattern for an incident pressure wave (band-gap frequency range) with the material parameters obtained by fitting the dispersion curves in one direction at $0^\circ$ (marked as 1) and in two directions at $0^\circ$  and $45^\circ$  (marked as 2).}
\label{fig:pre2}
\end{figure}

\begin{figure}[!ht]
    \centering
    \begin{subfigure}{0.7\textwidth}
        \includegraphics[width=\textwidth]{figures/legend.jpg}
        \put(0,10){\textbf{\large ${\lvert u \lvert}/{u_0}$}}
    \end{subfigure}
    
	\begin{subfigure}{0.0\textwidth}
    \begin{picture}(0,0)
        \put(-6,30){\rotatebox{90}{\textbf{$12.56 \cdot 10^{5}$ rad/sec}}}
    \end{picture}
	\end{subfigure}
    \begin{subfigure}{0.14\textwidth}
        \includegraphics[width=\textwidth]{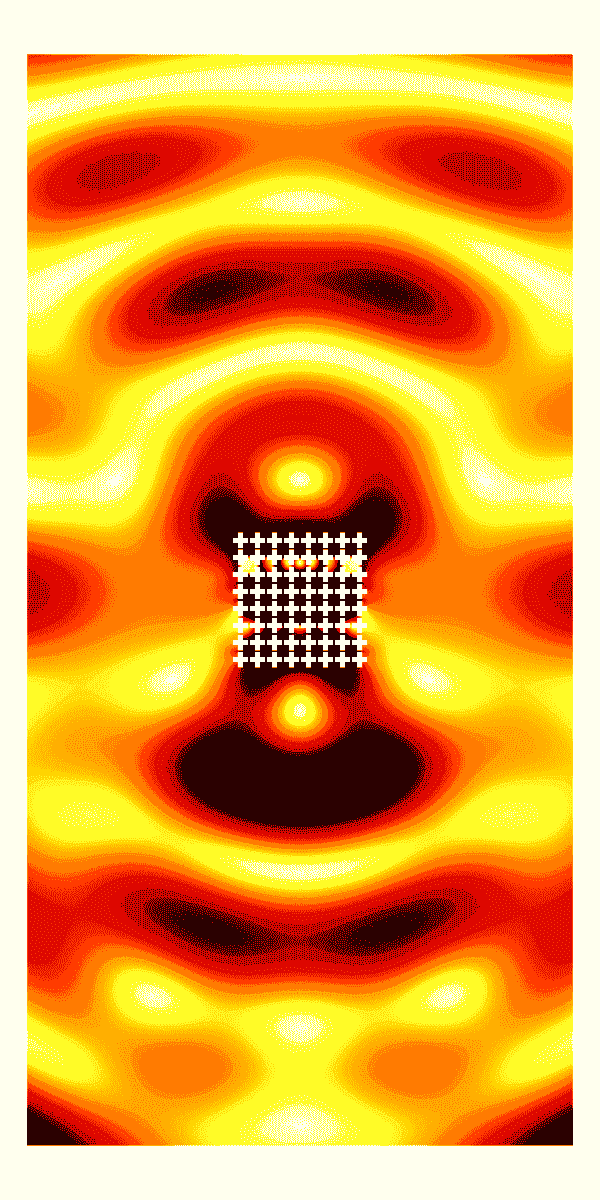}
    \end{subfigure}
    \begin{subfigure}{0.14\textwidth}
        \includegraphics[width=\textwidth]{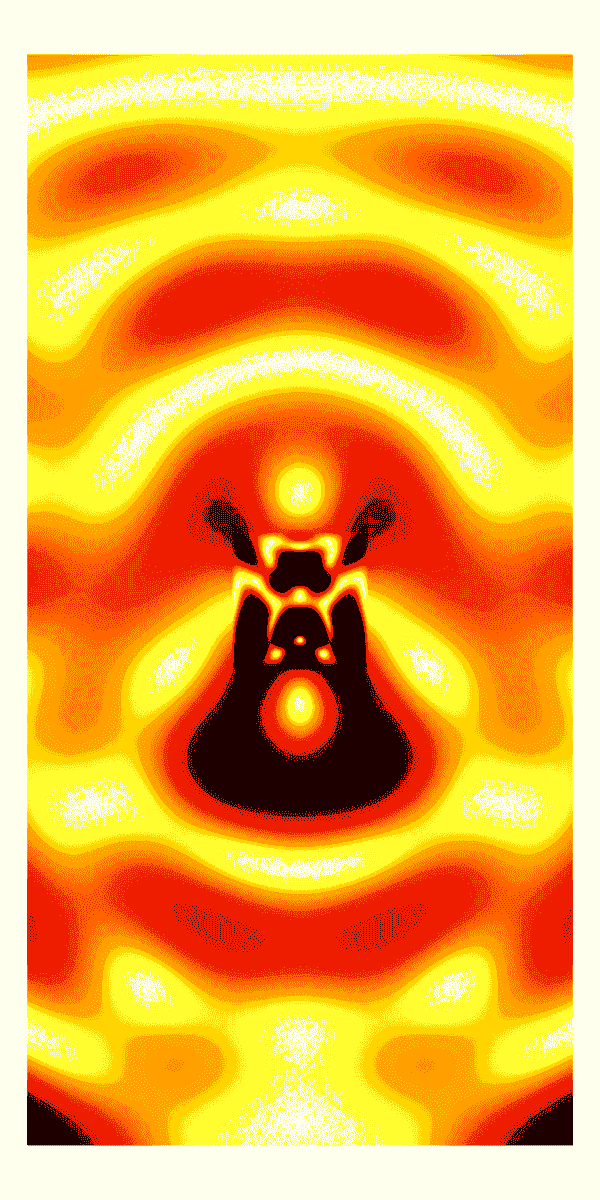}
    \end{subfigure}
    \begin{subfigure}{0.14\textwidth}
        \includegraphics[width=\textwidth]{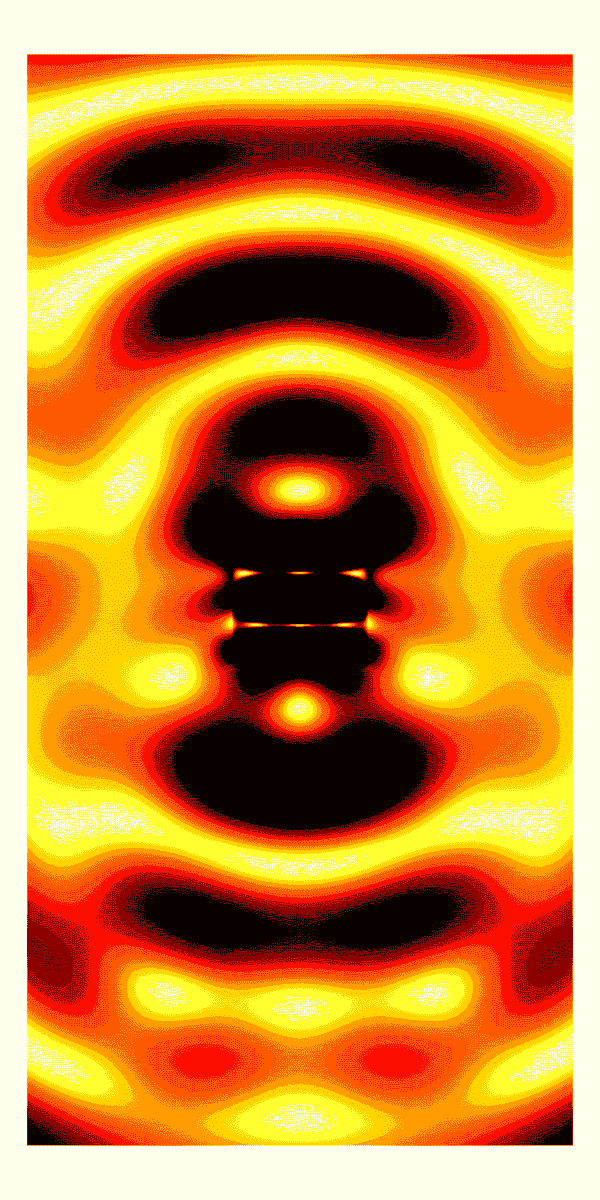}
    \end{subfigure}
    \begin{subfigure}{0.14\textwidth}
        \includegraphics[width=\textwidth]{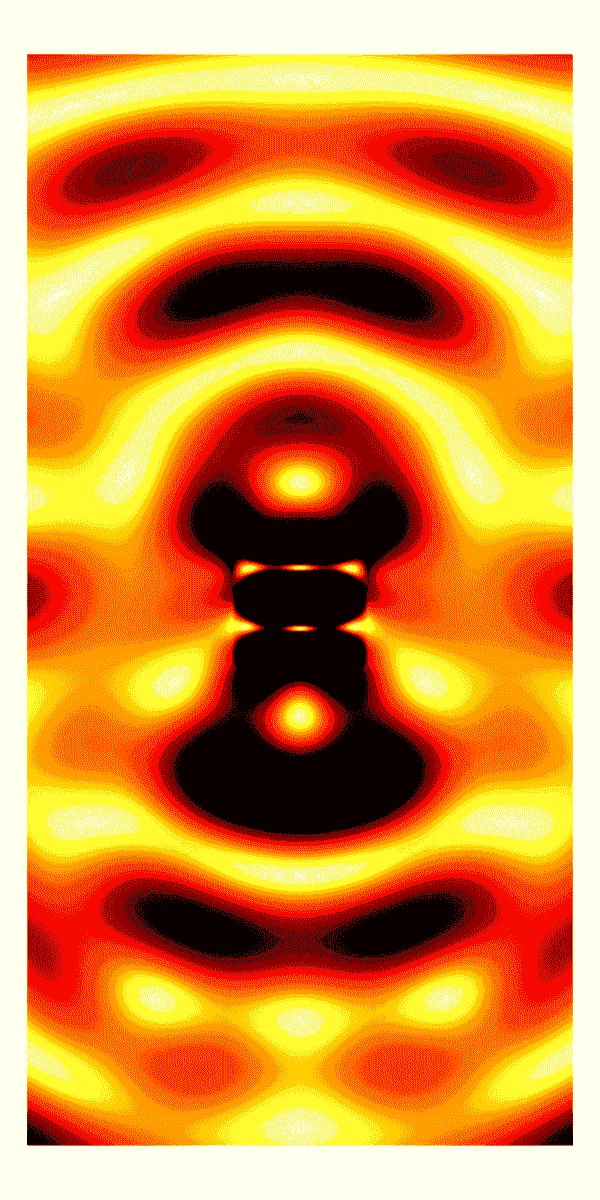}
    \end{subfigure}
    \begin{subfigure}{0.14\textwidth}
        \includegraphics[width=\textwidth]{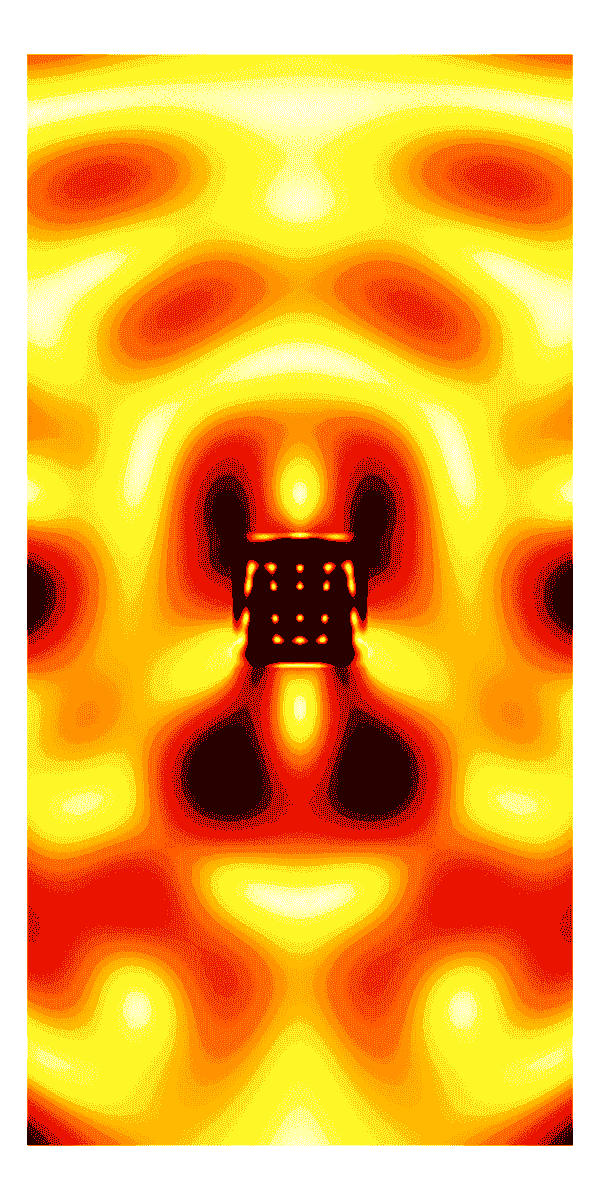}
    \end{subfigure}
    \begin{subfigure}{0.14\textwidth}
        \includegraphics[width=\textwidth]{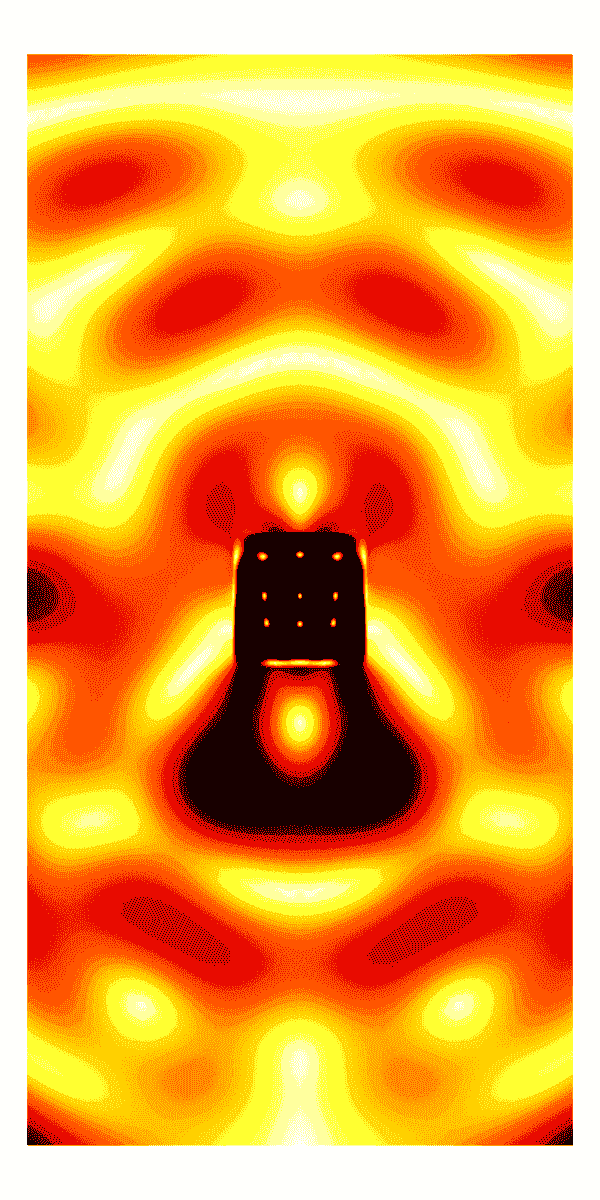}
    \end{subfigure}

    \centering
	\begin{subfigure}{0.0\textwidth}
    \begin{picture}(0,0)
        \put(-6,30){\rotatebox{90}{\textbf{$25.13 \cdot 10^{5}$ rad/sec}}}
    \end{picture}
	\end{subfigure}
    \begin{subfigure}{0.14\textwidth}
        \includegraphics[width=\textwidth]{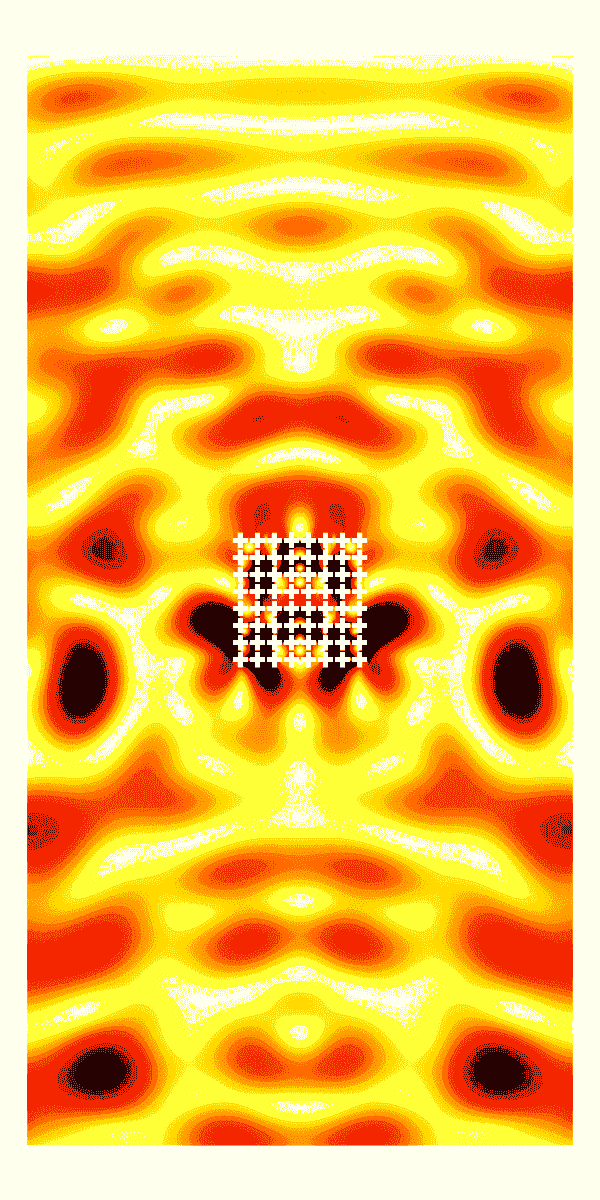}
    \end{subfigure}
    \begin{subfigure}{0.14\textwidth}
        \includegraphics[width=\textwidth]{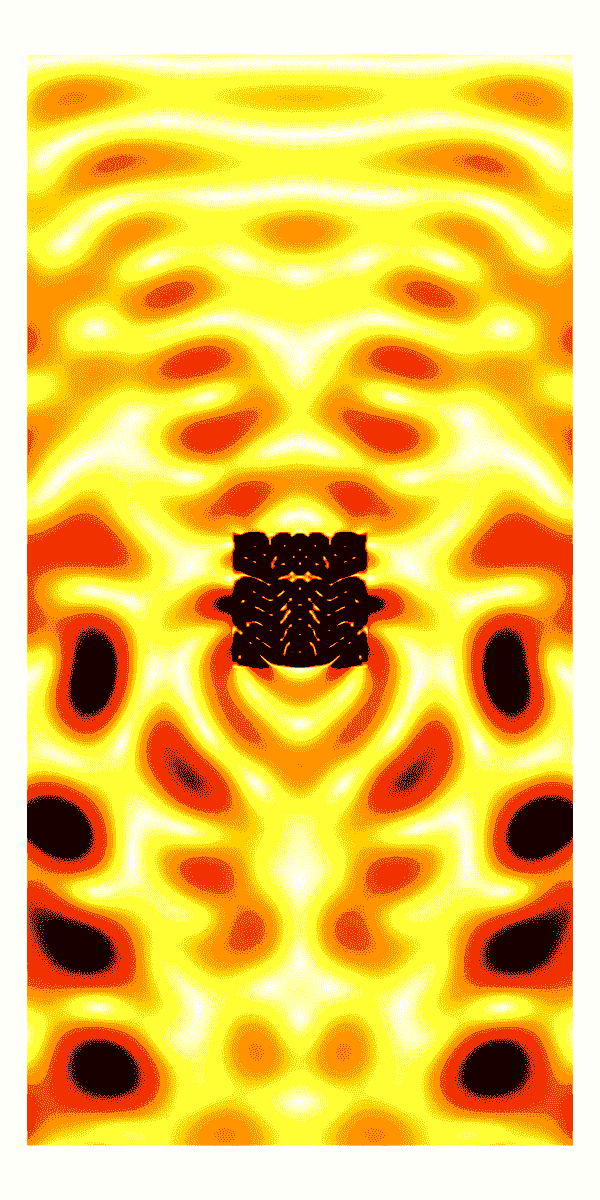}
    \end{subfigure}
    \begin{subfigure}{0.14\textwidth}
        \includegraphics[width=\textwidth]{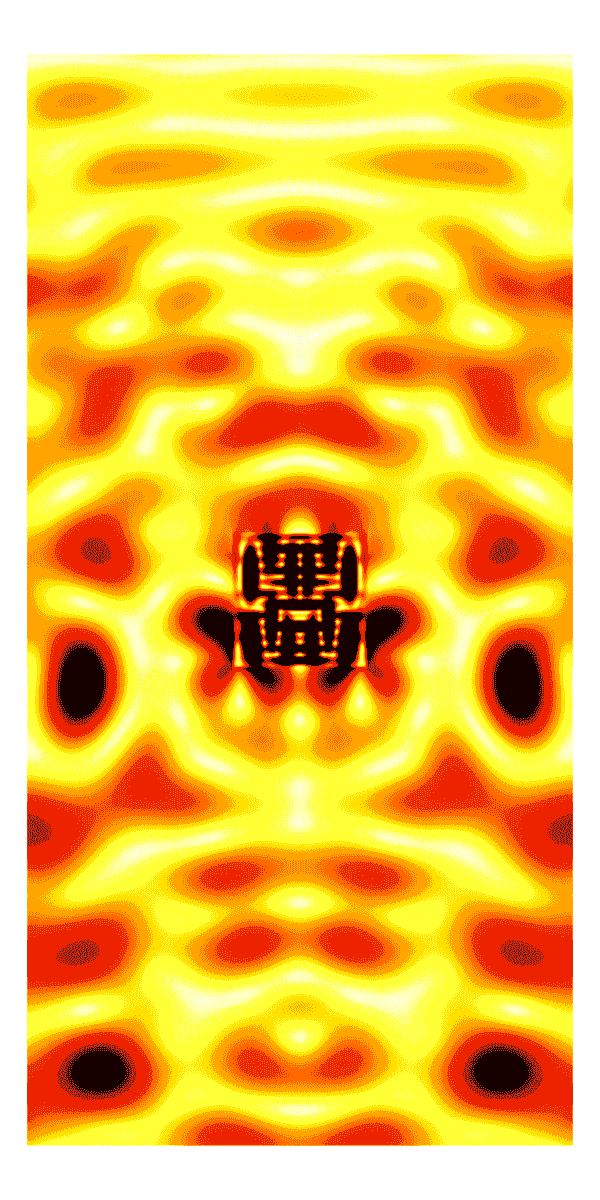}
    \end{subfigure}
    \begin{subfigure}{0.14\textwidth}
        \includegraphics[width=\textwidth]{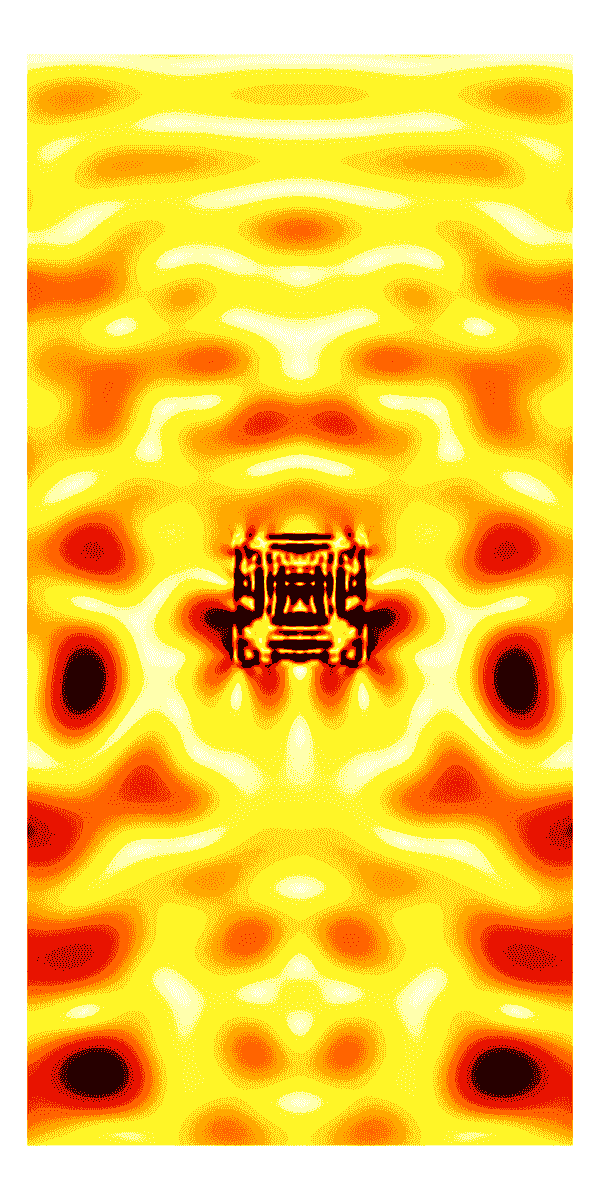}
    \end{subfigure}
    \begin{subfigure}{0.14\textwidth}
        \includegraphics[width=\textwidth]{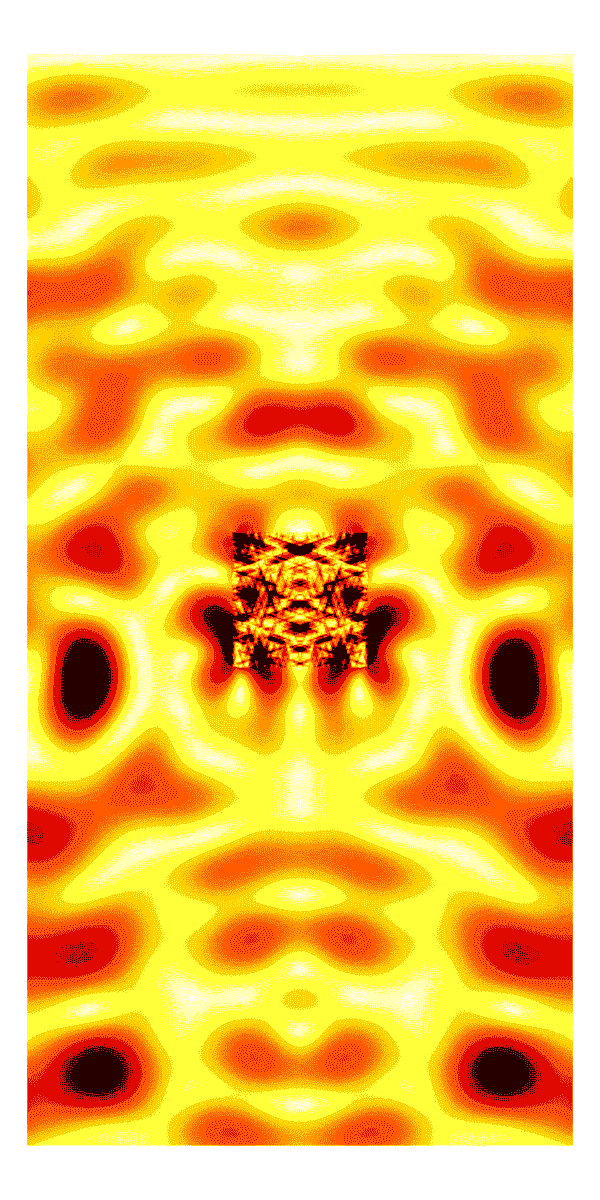}
    \end{subfigure}
    \begin{subfigure}{0.14\textwidth}
        \includegraphics[width=\textwidth]{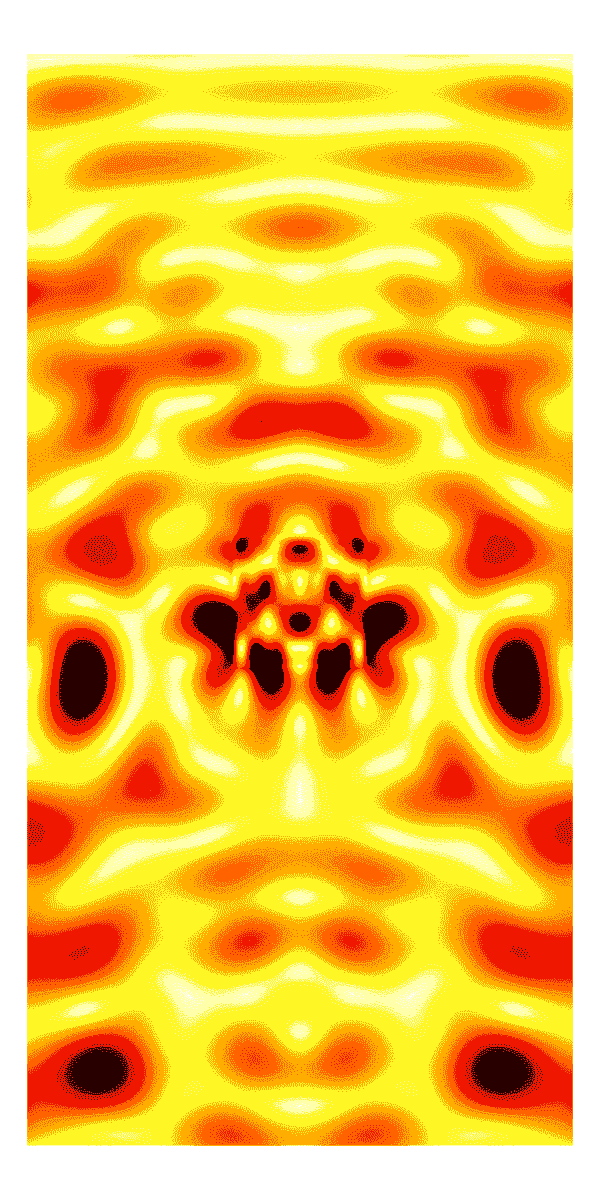}
    \end{subfigure}
    
	\begin{subfigure}{0.0\textwidth}
    \begin{picture}(0,0)
        \put(-6,30){\rotatebox{90}{\textbf{$37.7 \cdot 10^{5}$ rad/sec}}}
    \end{picture}
	\end{subfigure}
    \begin{subfigure}{0.14\textwidth}
        \includegraphics[width=\textwidth]{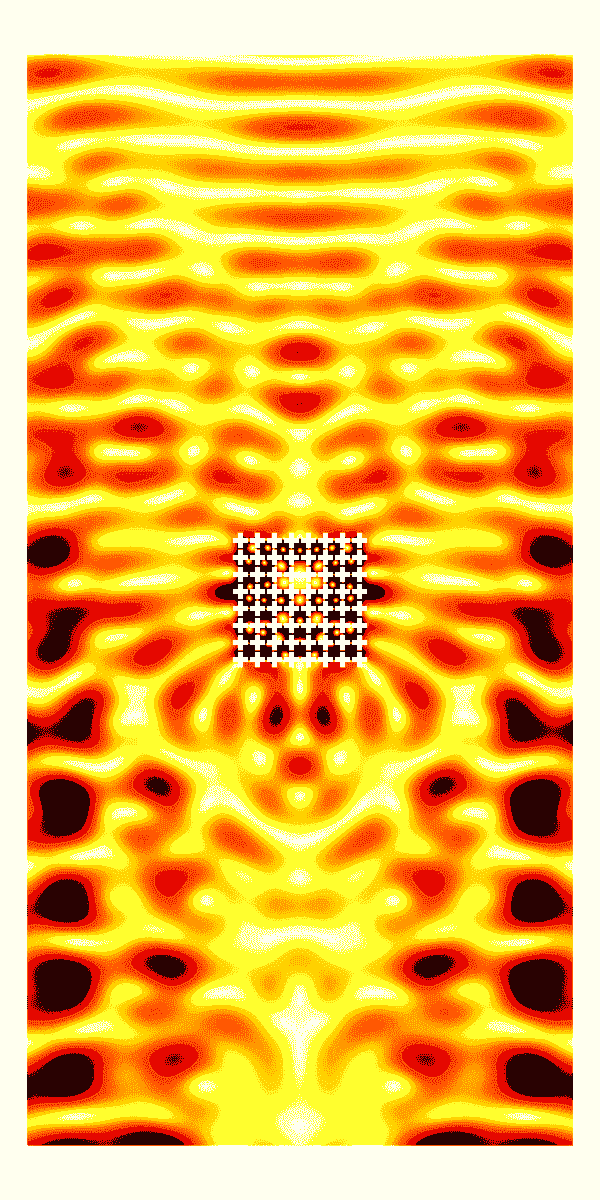}
    \end{subfigure}
    \begin{subfigure}{0.14\textwidth}
        \includegraphics[width=\textwidth]{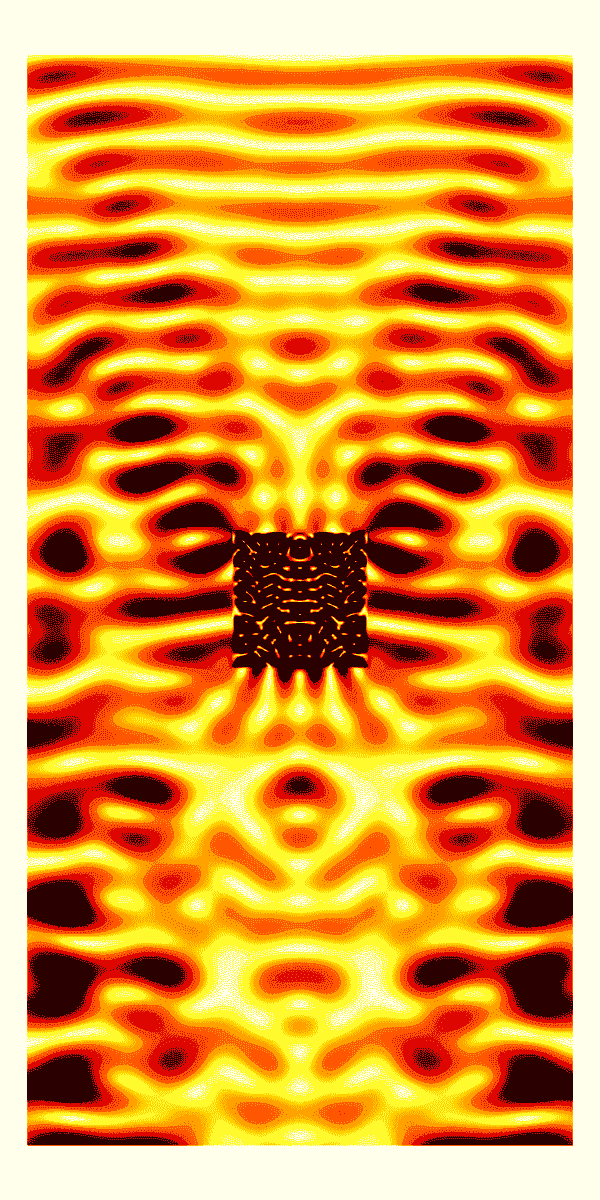}
    \end{subfigure}
    \begin{subfigure}{0.14\textwidth}
        \includegraphics[width=\textwidth]{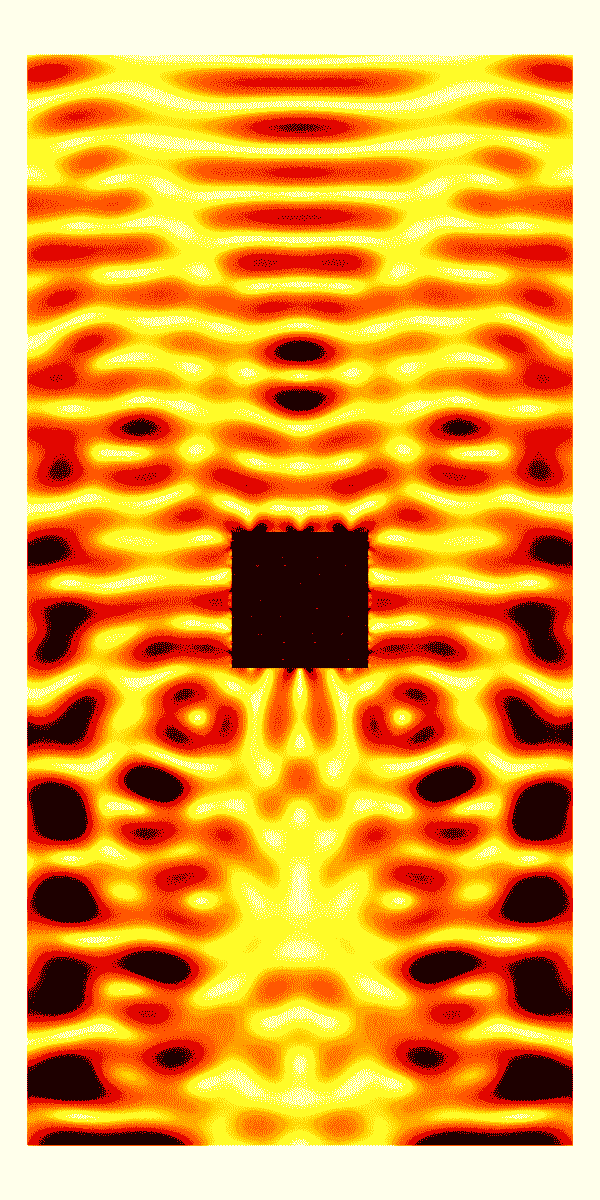}
    \end{subfigure}
    \begin{subfigure}{0.14\textwidth}
        \includegraphics[width=\textwidth]{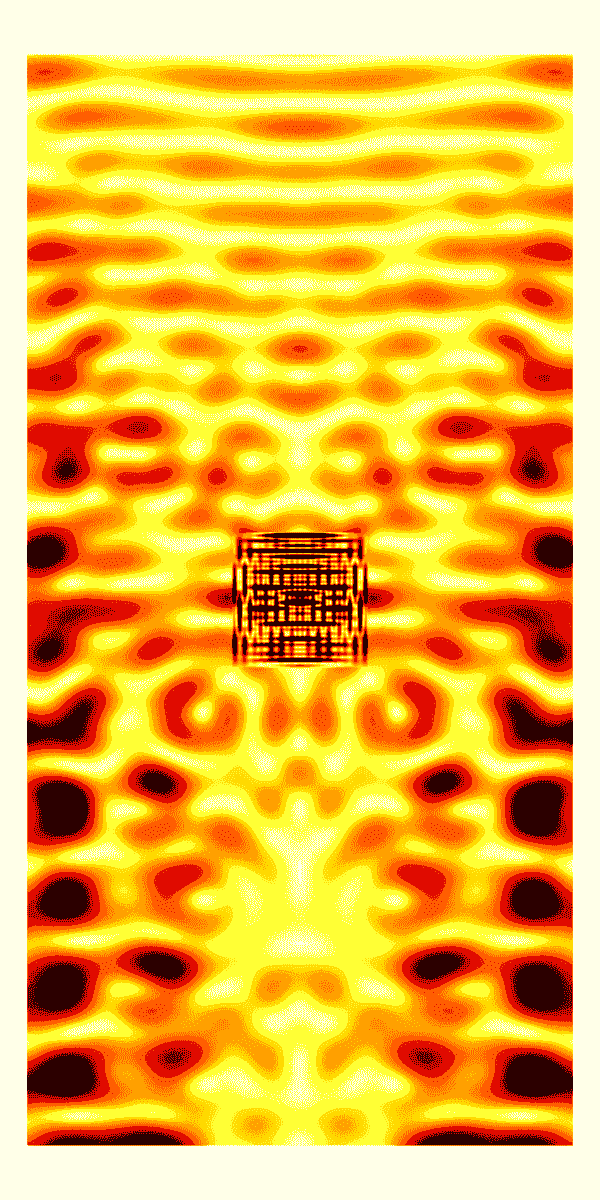}
    \end{subfigure}
    \begin{subfigure}{0.14\textwidth}
        \includegraphics[width=\textwidth]{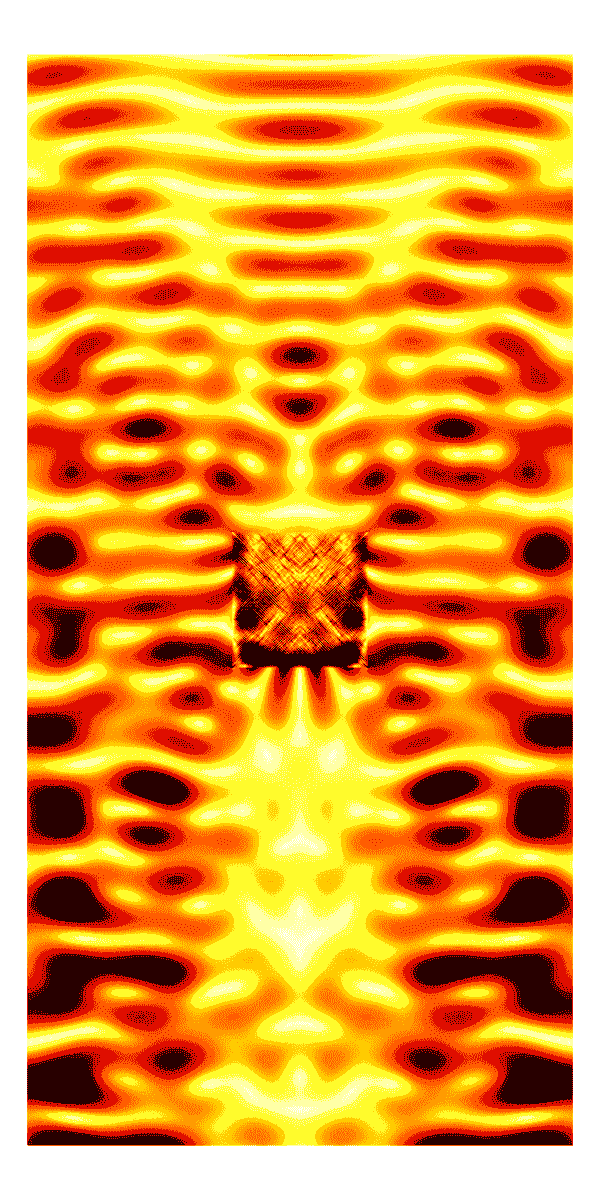}
    \end{subfigure}
    \begin{subfigure}{0.14\textwidth}
        \includegraphics[width=\textwidth]{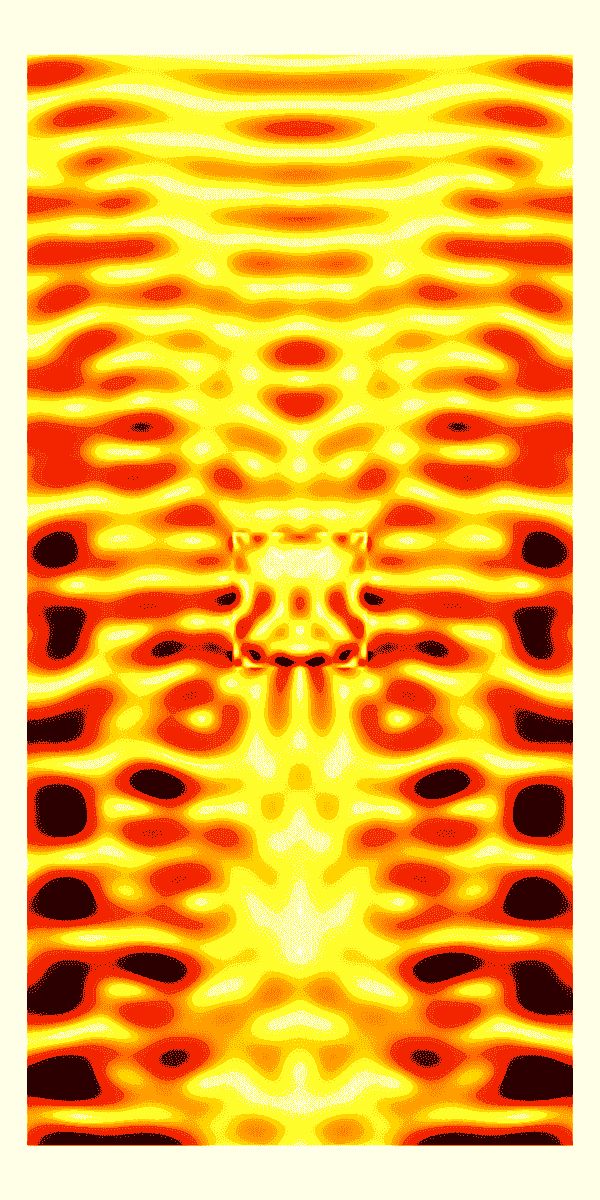}
    \end{subfigure}

    	\begin{subfigure}{0.0\textwidth}
    \begin{picture}(0,0)
        \put(-6,30){\rotatebox{90}{\textbf{$50.27 \cdot 10^{5}$ rad/sec}}}
    \end{picture}
	\end{subfigure}
    \begin{subfigure}{0.14\textwidth}
        \includegraphics[width=\textwidth]{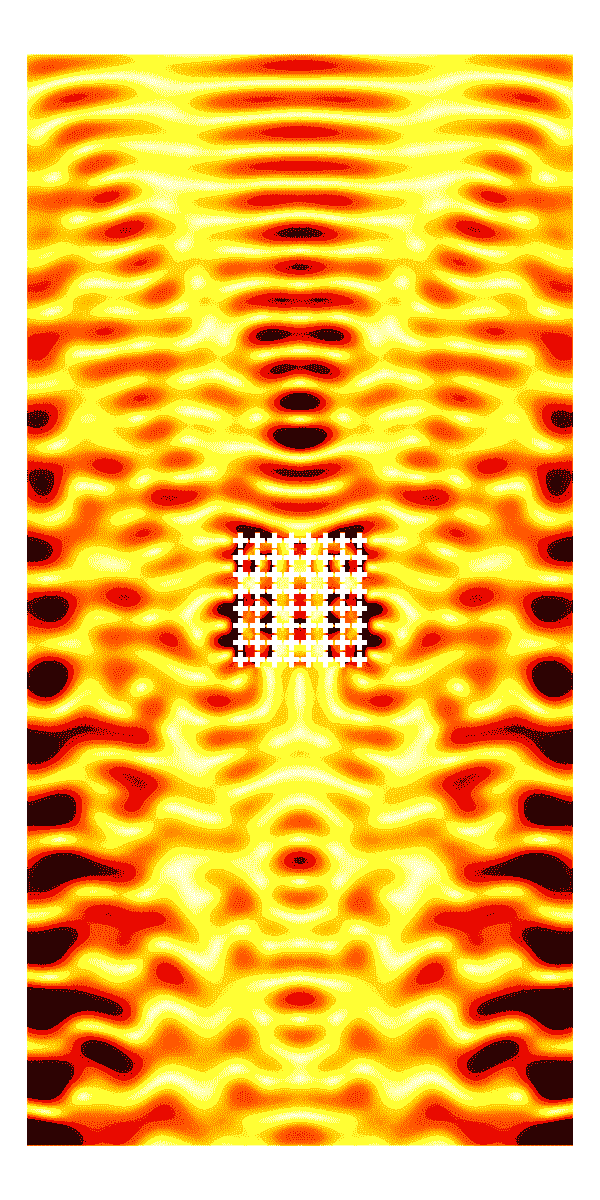}
    \end{subfigure}
    \begin{subfigure}{0.14\textwidth}
        \includegraphics[width=\textwidth]{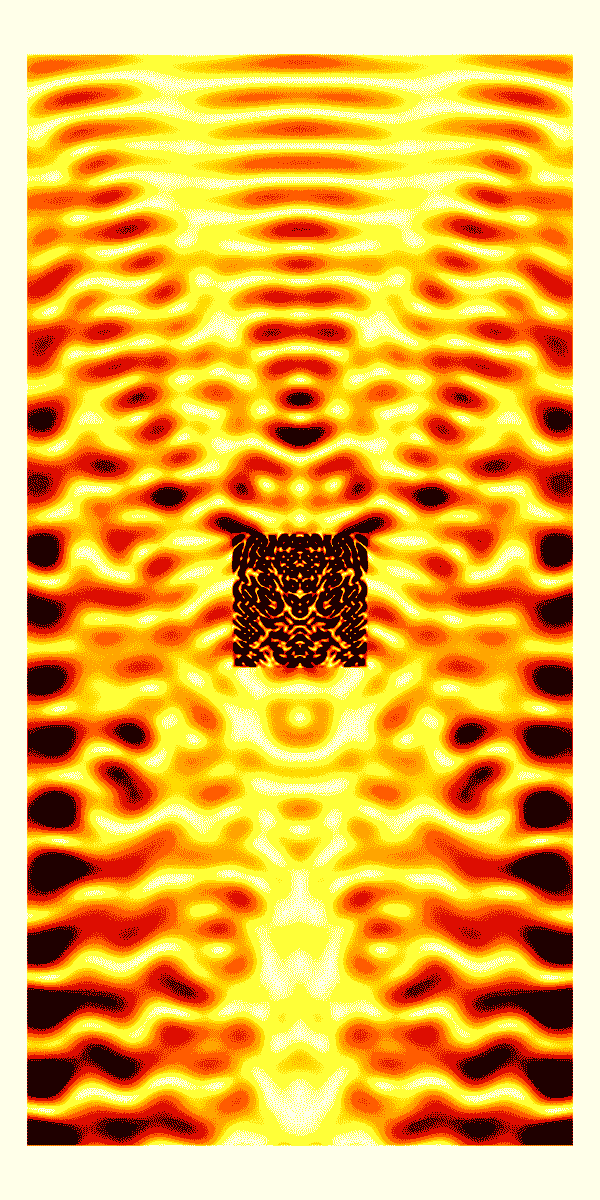}
    \end{subfigure}
    \begin{subfigure}{0.14\textwidth}
        \includegraphics[width=\textwidth]{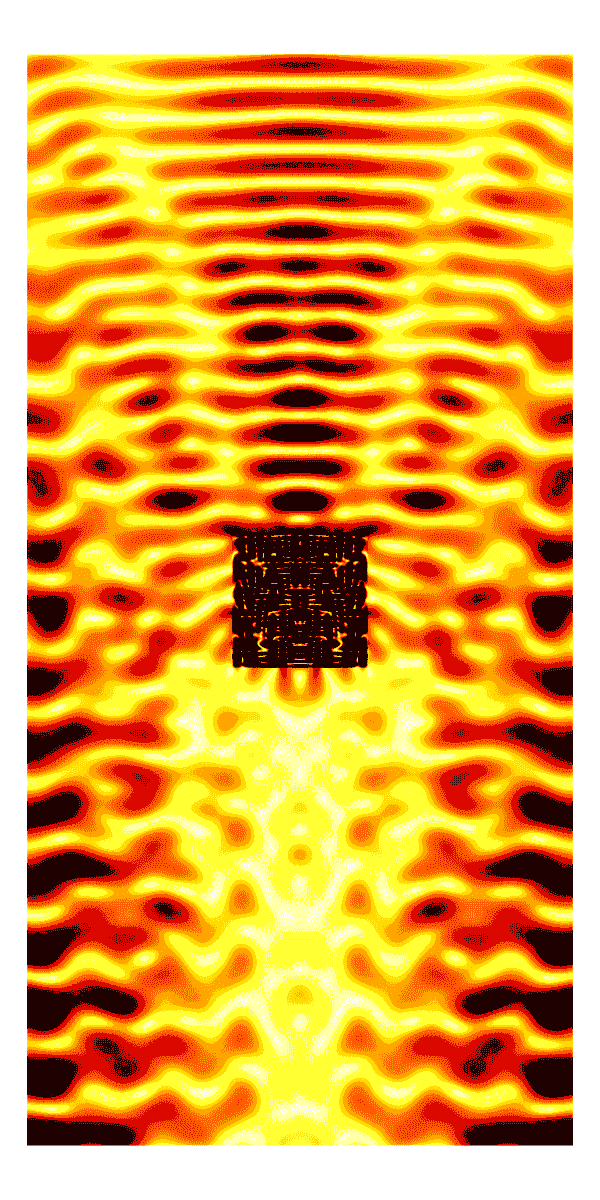}
    \end{subfigure}
    \begin{subfigure}{0.14\textwidth}
        \includegraphics[width=\textwidth]{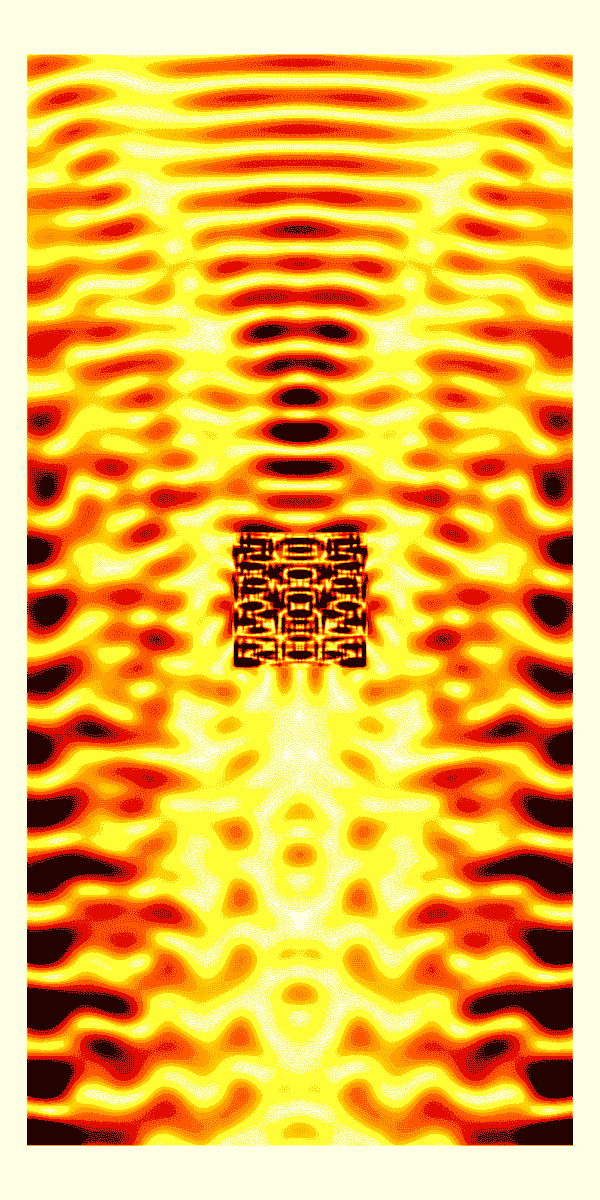}
    \end{subfigure}
    \begin{subfigure}{0.14\textwidth}
        \includegraphics[width=\textwidth]{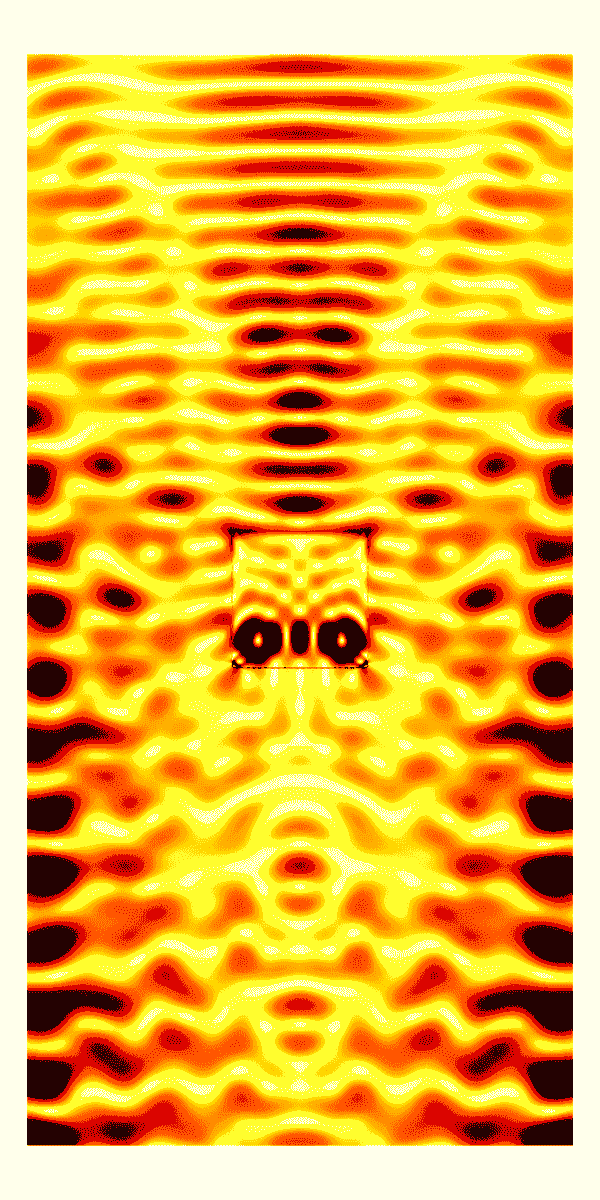}
    \end{subfigure}
    \begin{subfigure}{0.14\textwidth}
        \includegraphics[width=\textwidth]{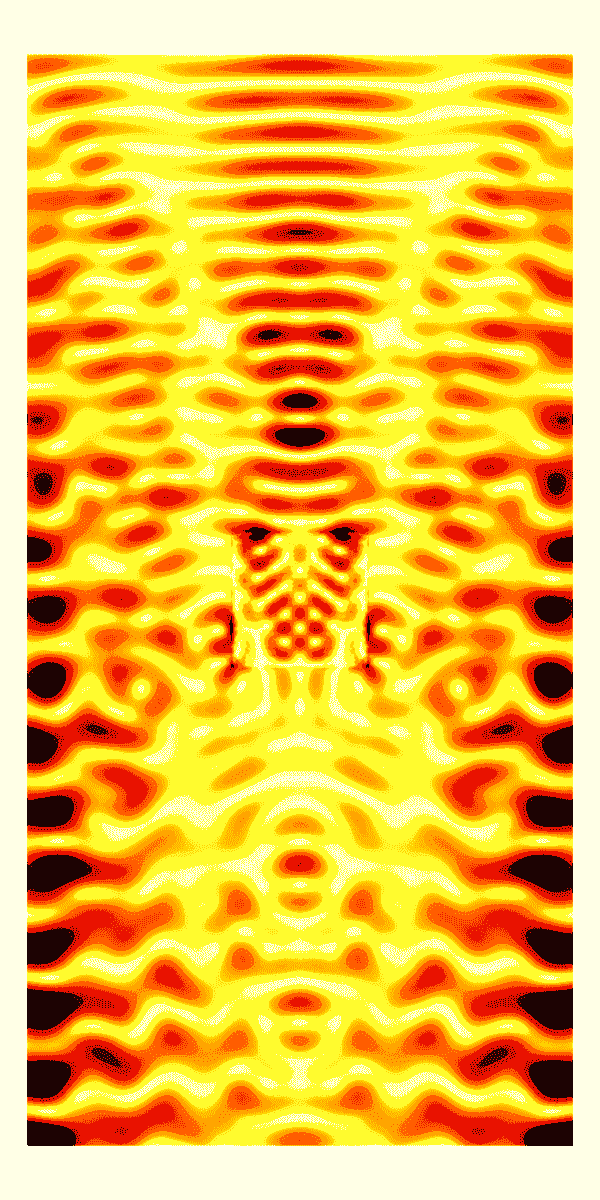}
    \end{subfigure}
    
	\begin{subfigure}{0.0\textwidth}
    \begin{picture}(0,0)
        \put(-6,40){\rotatebox{90}{\textbf{$62.83 \cdot 10^{5}$ rad/sec}}}
    \end{picture}
	\end{subfigure}
    \begin{subfigure}{0.14\textwidth}
        \includegraphics[width=\textwidth]{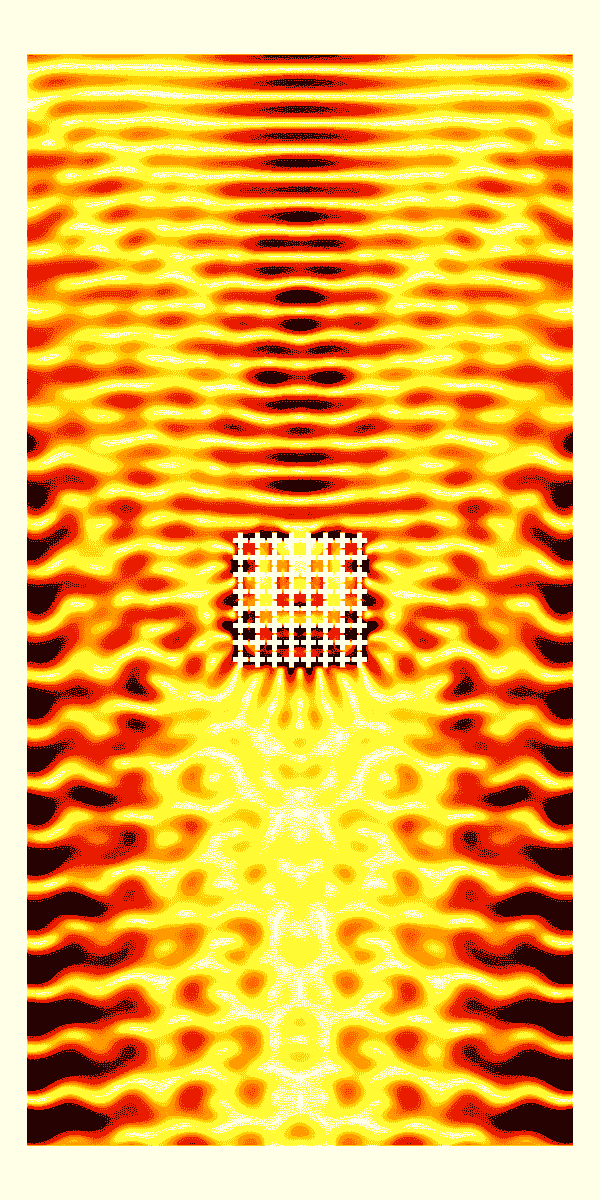}
    \caption*{\footnotesize microstructured}
    \end{subfigure}
    \begin{subfigure}{0.14\textwidth}
        \includegraphics[width=\textwidth]{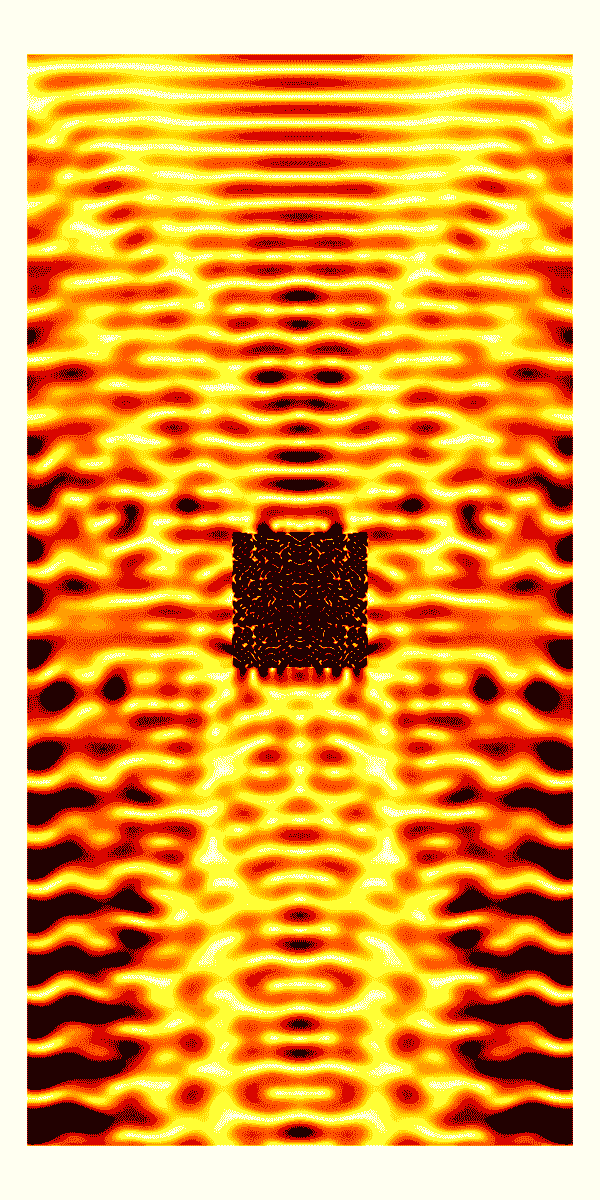}
   \caption*{\footnotesize macro-Cauchy}
    \end{subfigure}
    \begin{subfigure}{0.14\textwidth}
        \includegraphics[width=\textwidth]{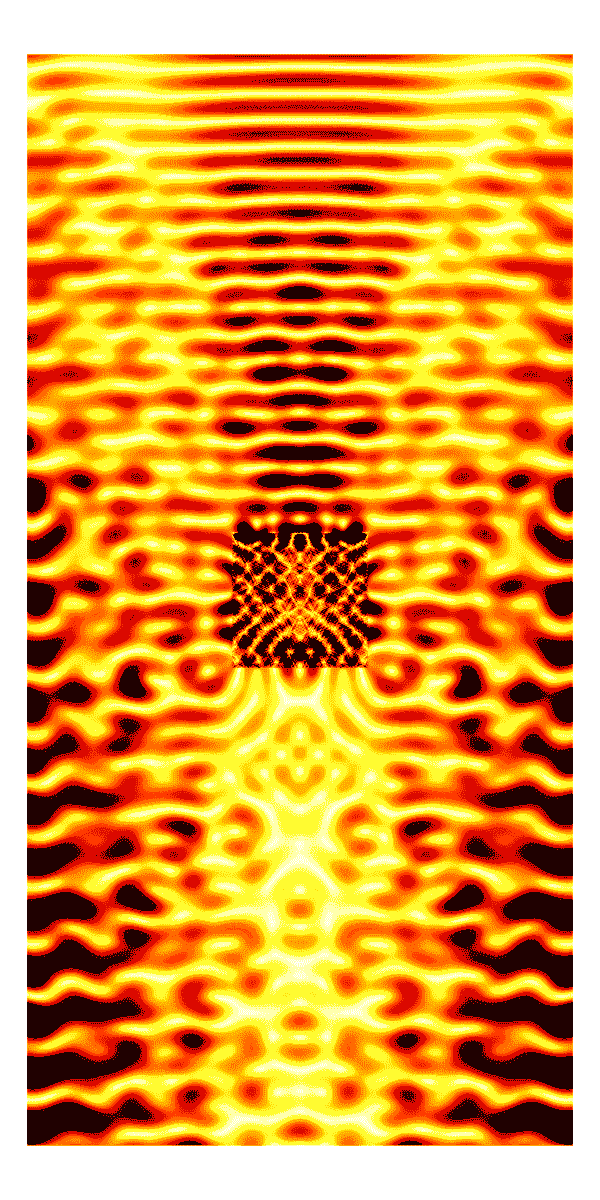}
    \caption*{\footnotesize RRMM(1)}
    \end{subfigure}
    \begin{subfigure}{0.14\textwidth}
        \includegraphics[width=\textwidth]{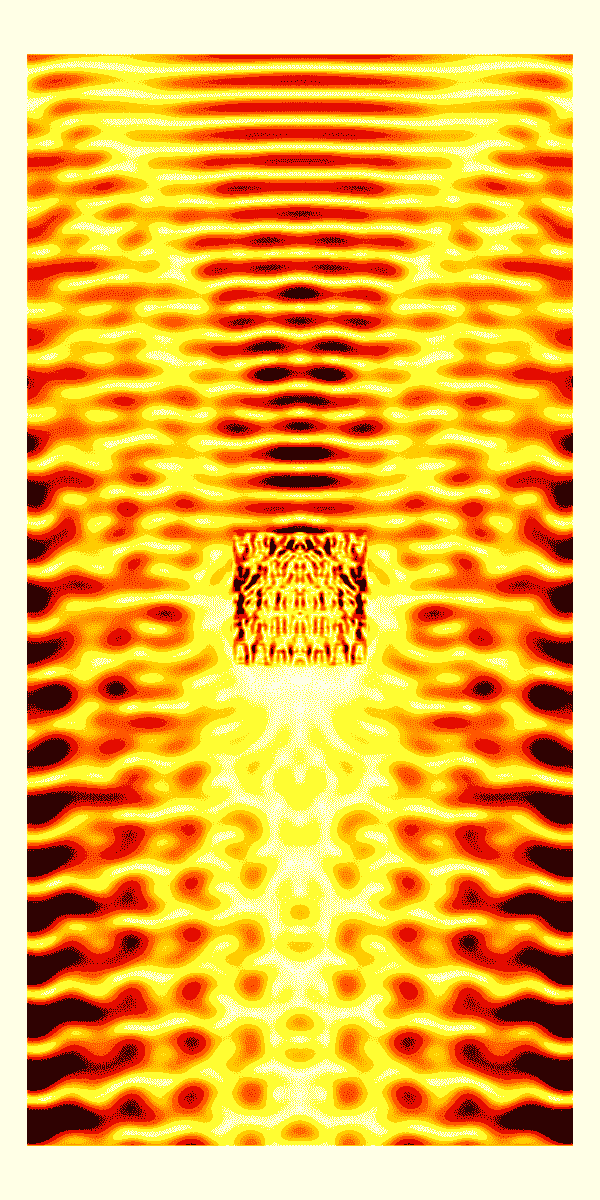}
    \caption*{\footnotesize RMM(1)}  
    \end{subfigure}
    \begin{subfigure}{0.14\textwidth}
        \includegraphics[width=\textwidth]{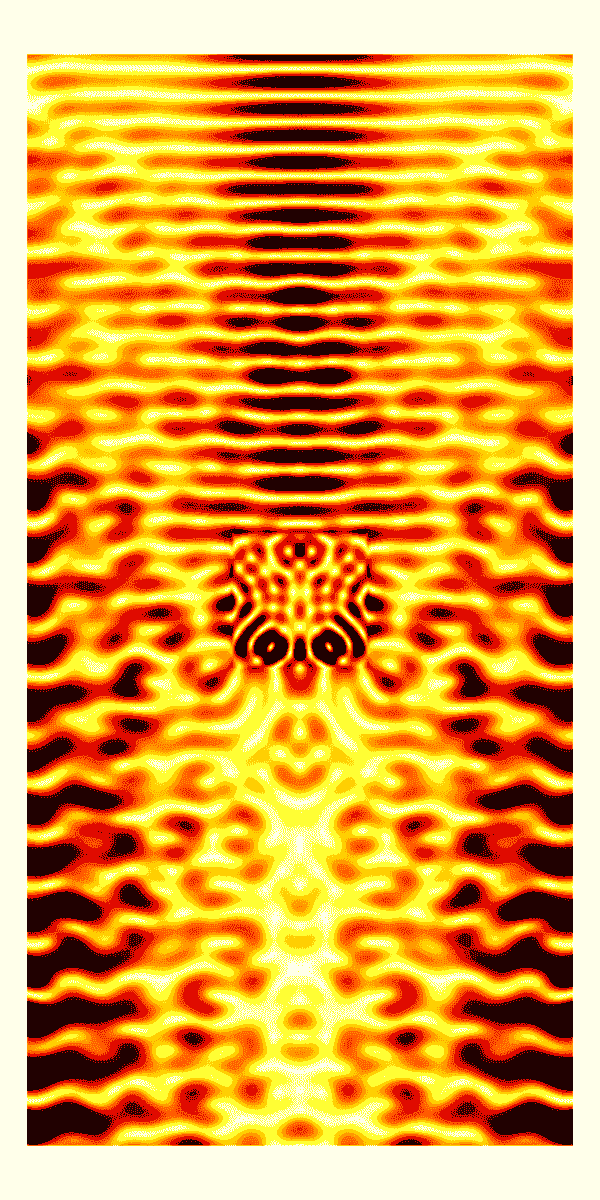}
    \caption*{\footnotesize RRMM(2)}
    \end{subfigure}
    \begin{subfigure}{0.14\textwidth}
        \includegraphics[width=\textwidth]{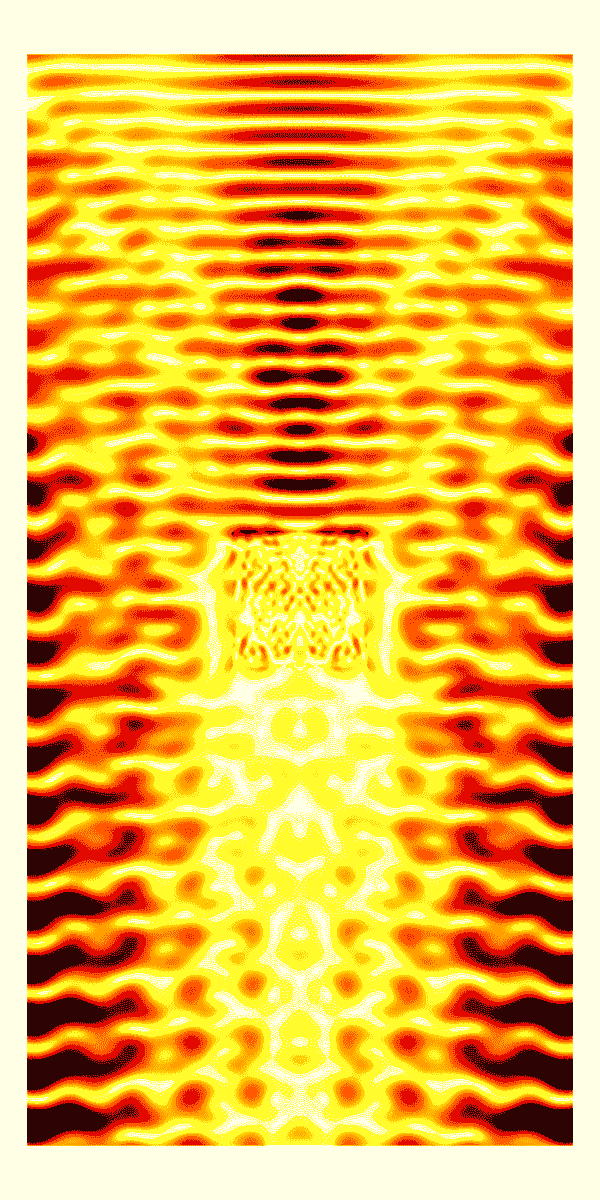}
    \caption*{\footnotesize RMM(2)}  
    \end{subfigure}
\caption{Results of finite-size scattering pattern for an incident shear wave (lower frequency range up to the lower band-gap limit) with the material parameters obtained by fitting the dispersion curves in one direction at $0^\circ$ (marked as 1) and in two directions at $0^\circ$  and $45^\circ$  (marked as 2).}
\label{fig:she1}
\end{figure}

\begin{figure}[!ht]
    \centering
    \begin{subfigure}{0.7\textwidth}
        \includegraphics[width=\textwidth]{figures/legend.jpg}
        \put(0,10){\textbf{\large ${\lvert u \lvert}/{u_0}$}}
    \end{subfigure}
    
	\begin{subfigure}{0.0\textwidth}
    \begin{picture}(0,0)
        \put(-6,30){\rotatebox{90}{\textbf{$75.4 \cdot 10^{5}$ rad/sec}}}
    \end{picture}
	\end{subfigure}
    \begin{subfigure}{0.14\textwidth}
        \includegraphics[width=\textwidth]{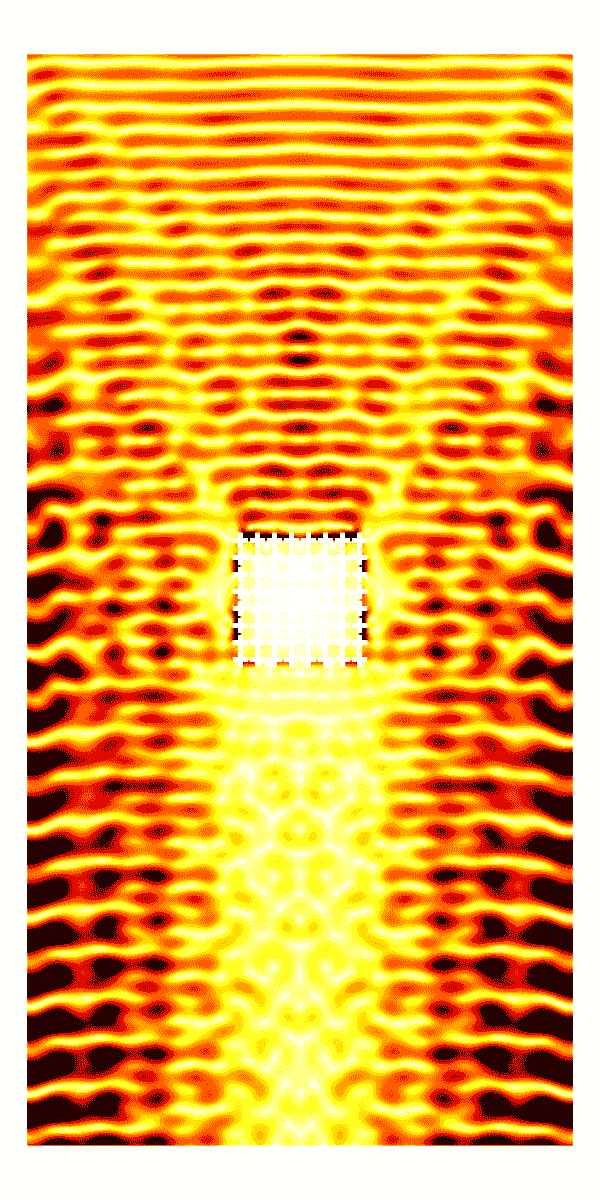}
    \end{subfigure}
    \begin{subfigure}{0.14\textwidth}
        \includegraphics[width=\textwidth]{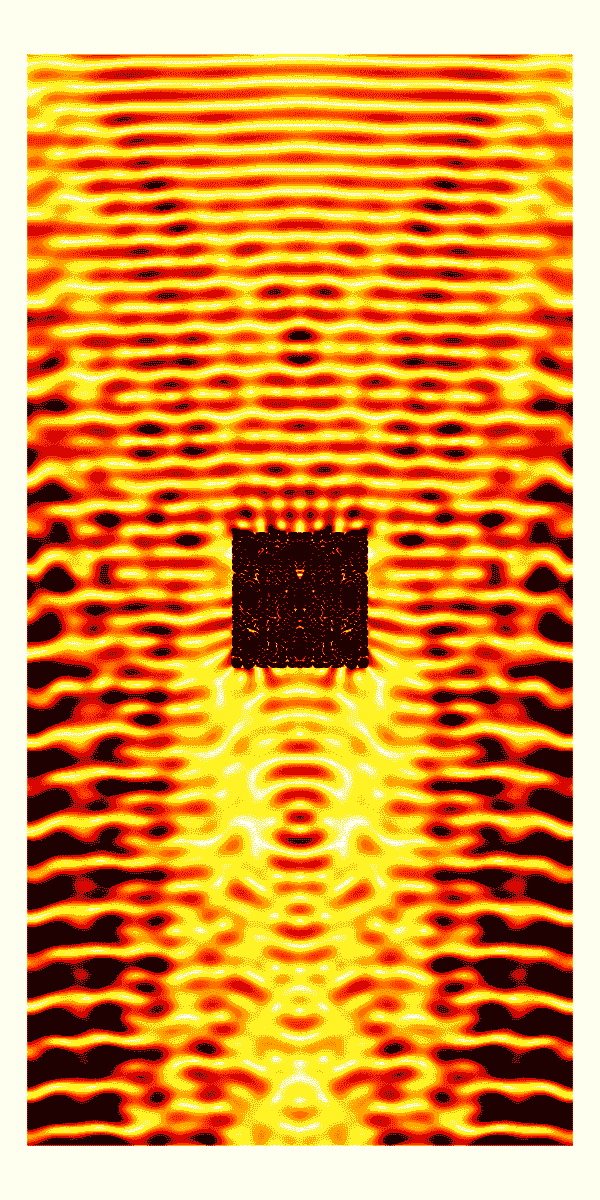}
    \end{subfigure}
    \begin{subfigure}{0.14\textwidth}
        \includegraphics[width=\textwidth]{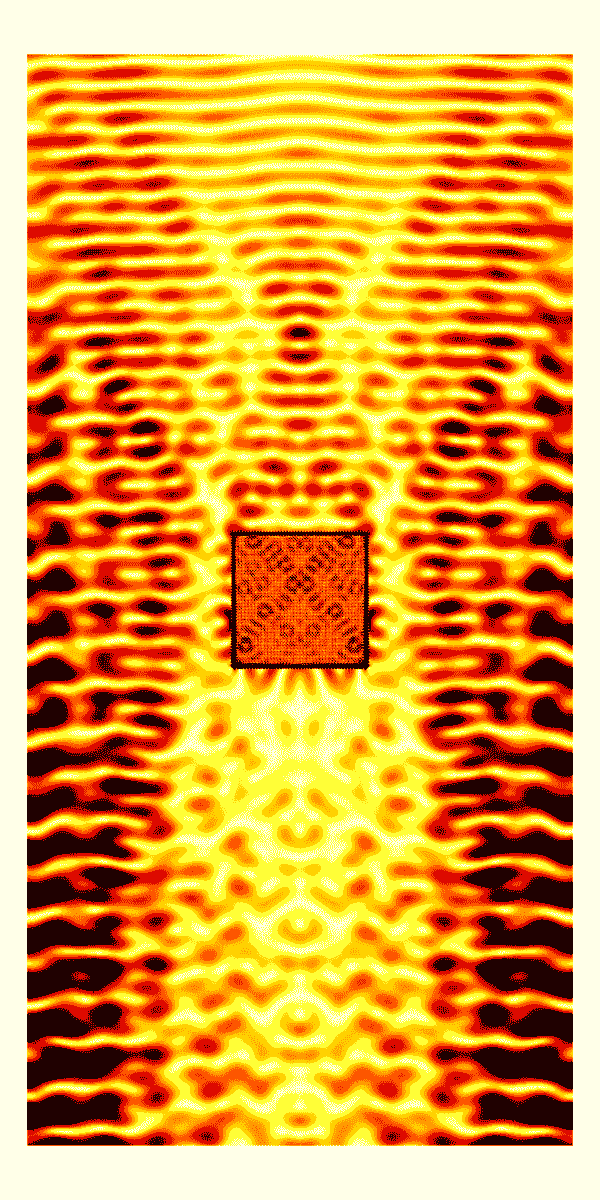}
    \end{subfigure}
    \begin{subfigure}{0.14\textwidth}
        \includegraphics[width=\textwidth]{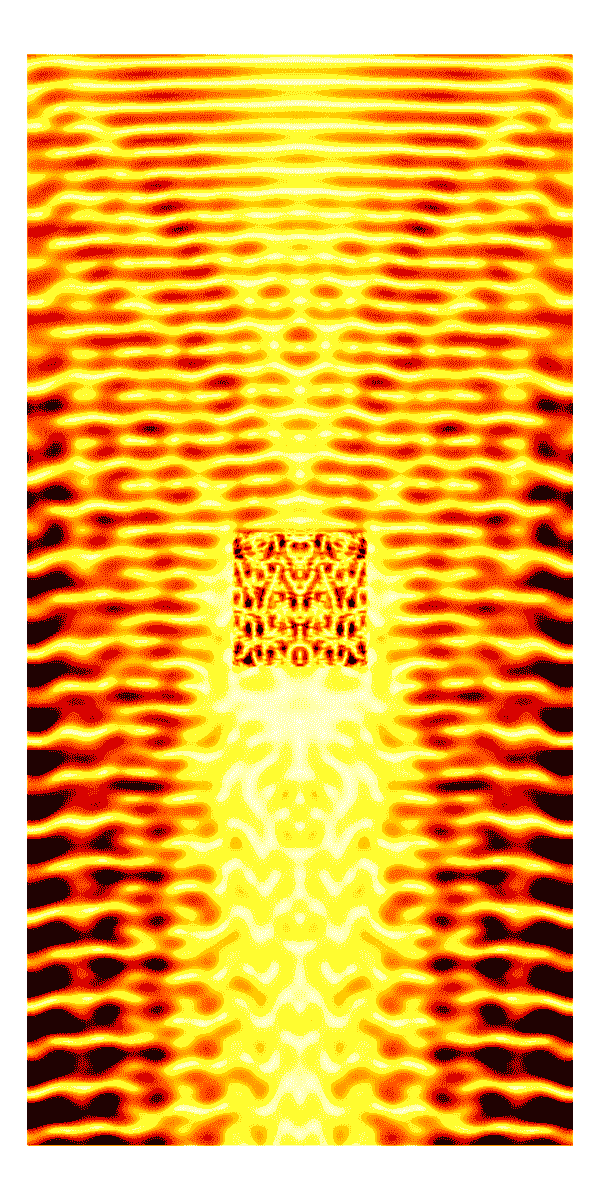}
    \end{subfigure}
    \begin{subfigure}{0.14\textwidth}
        \includegraphics[width=\textwidth]{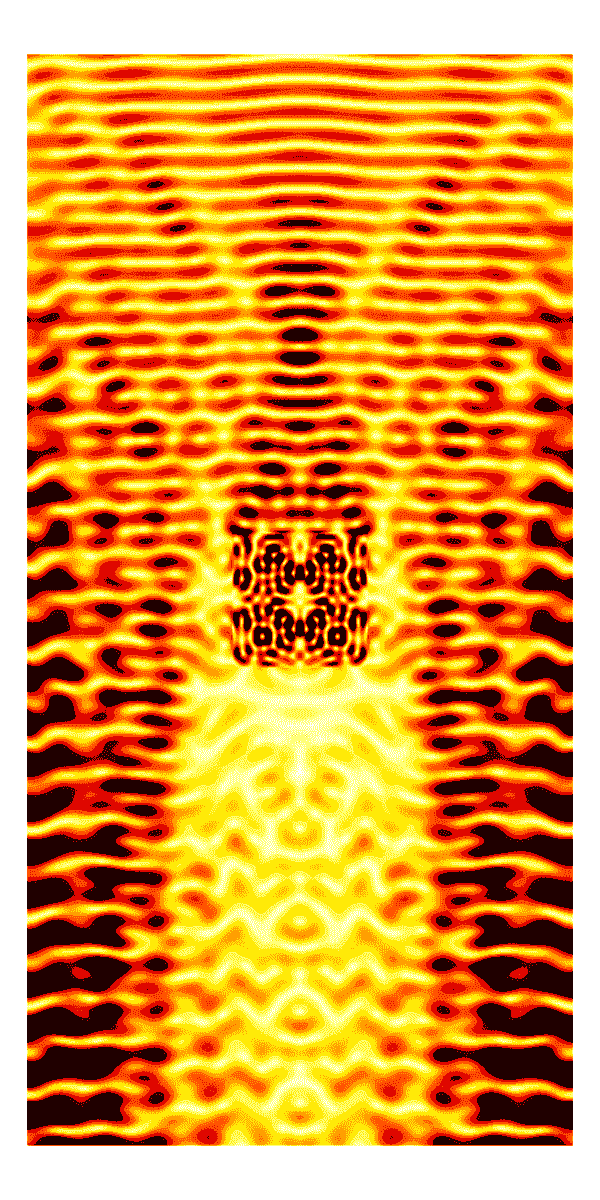}
    \end{subfigure}
    \begin{subfigure}{0.14\textwidth}
        \includegraphics[width=\textwidth]{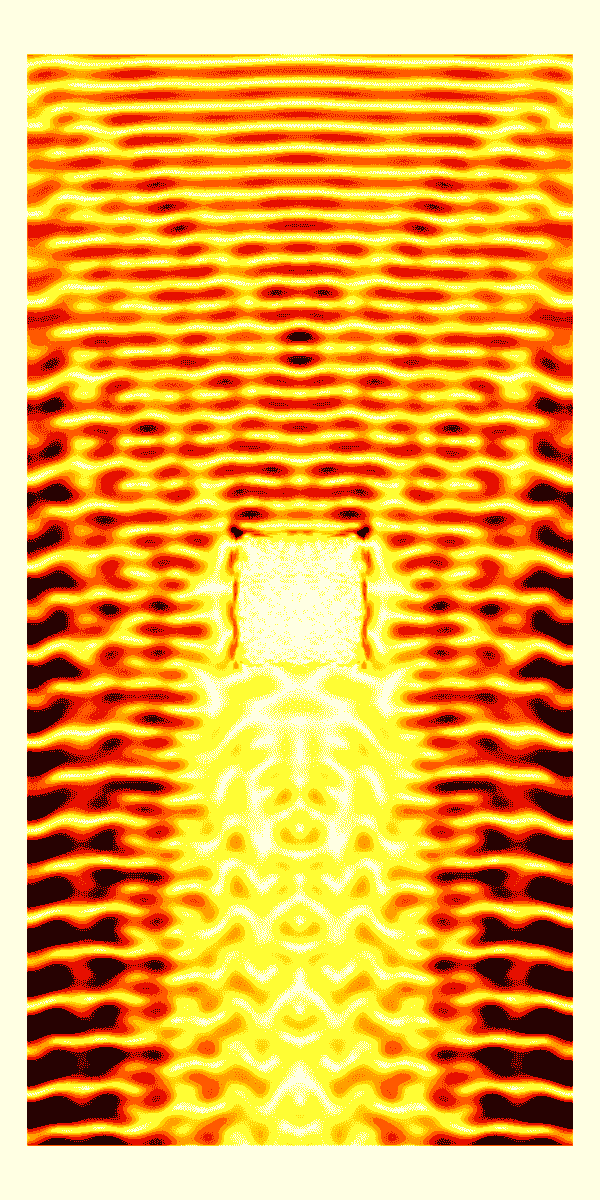}
    \end{subfigure}

	\begin{subfigure}{0.0\textwidth}
    \begin{picture}(0,0)
        \put(-6,30){\rotatebox{90}{\textbf{$87.96 \cdot 10^{5}$ rad/sec}}}
    \end{picture}
	\end{subfigure}
    \begin{subfigure}{0.14\textwidth}
        \includegraphics[width=\textwidth]{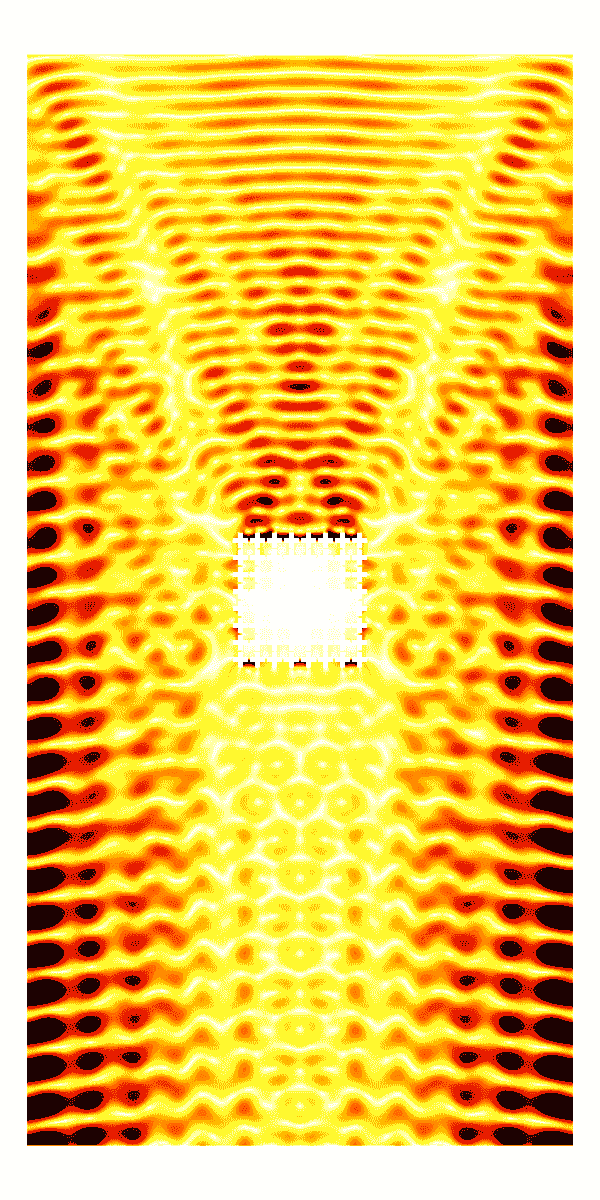}
    \end{subfigure}
    \begin{subfigure}{0.14\textwidth}
        \includegraphics[width=\textwidth]{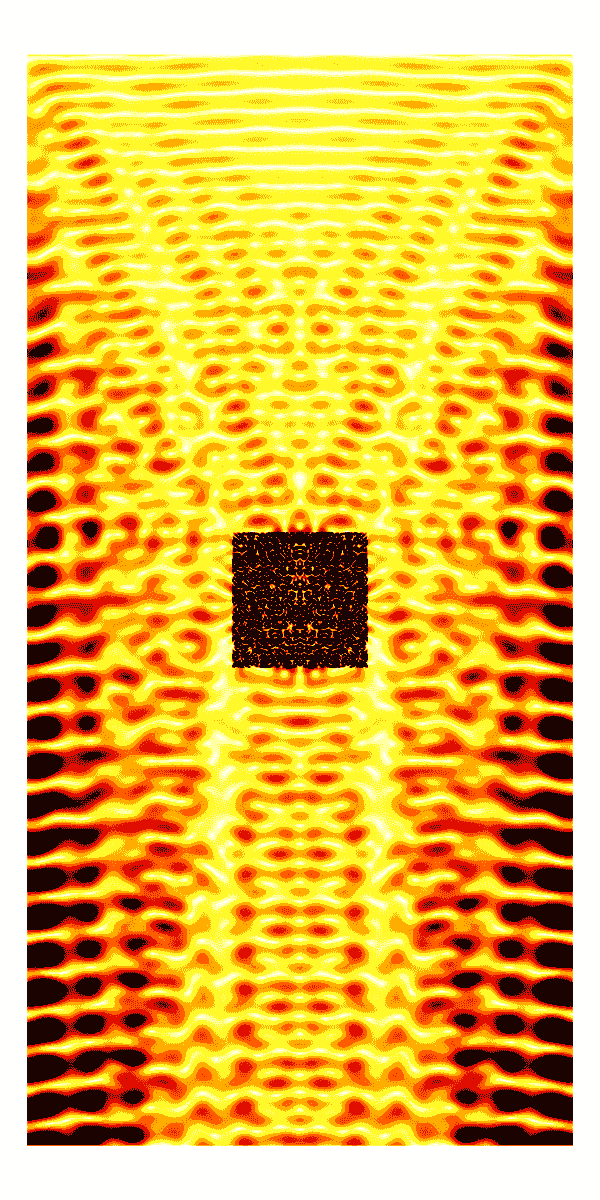}
    \end{subfigure}
    \begin{subfigure}{0.14\textwidth}
        \includegraphics[width=\textwidth]{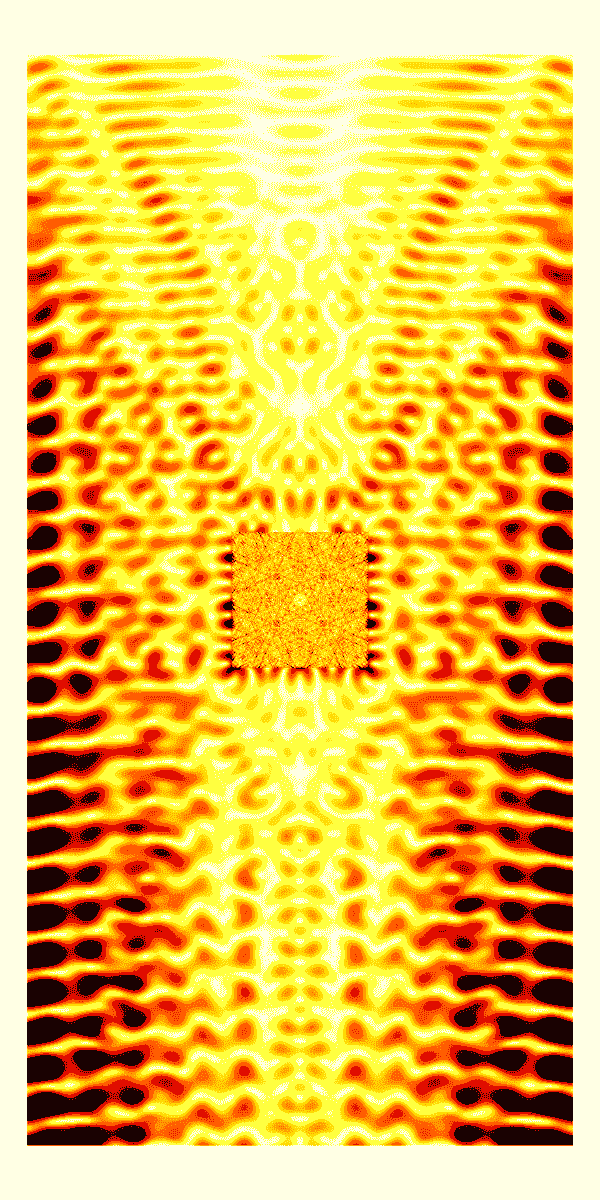}
    \end{subfigure}
    \begin{subfigure}{0.14\textwidth}
        \includegraphics[width=\textwidth]{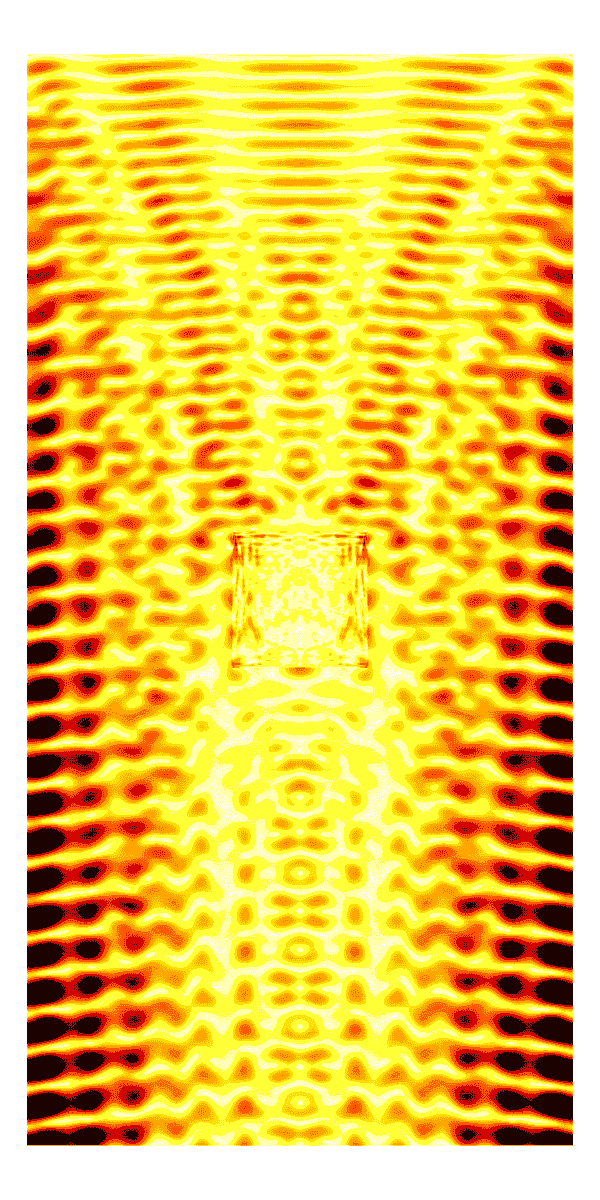}
    \end{subfigure}
    \begin{subfigure}{0.14\textwidth}
        \includegraphics[width=\textwidth]{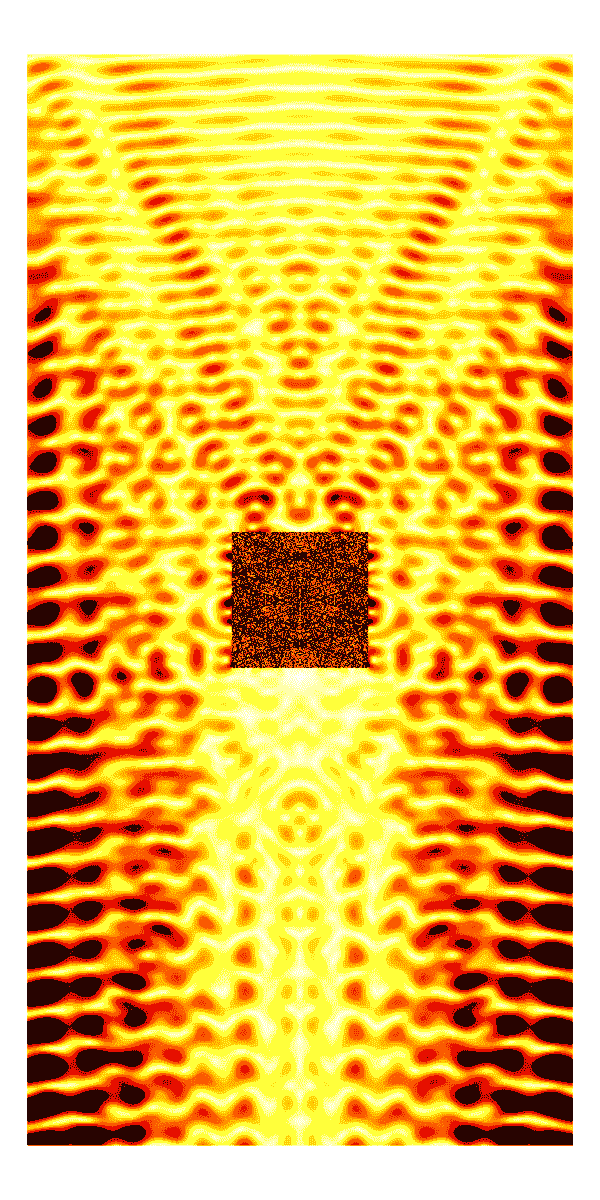}
    \end{subfigure}
    \begin{subfigure}{0.14\textwidth}
        \includegraphics[width=\textwidth]{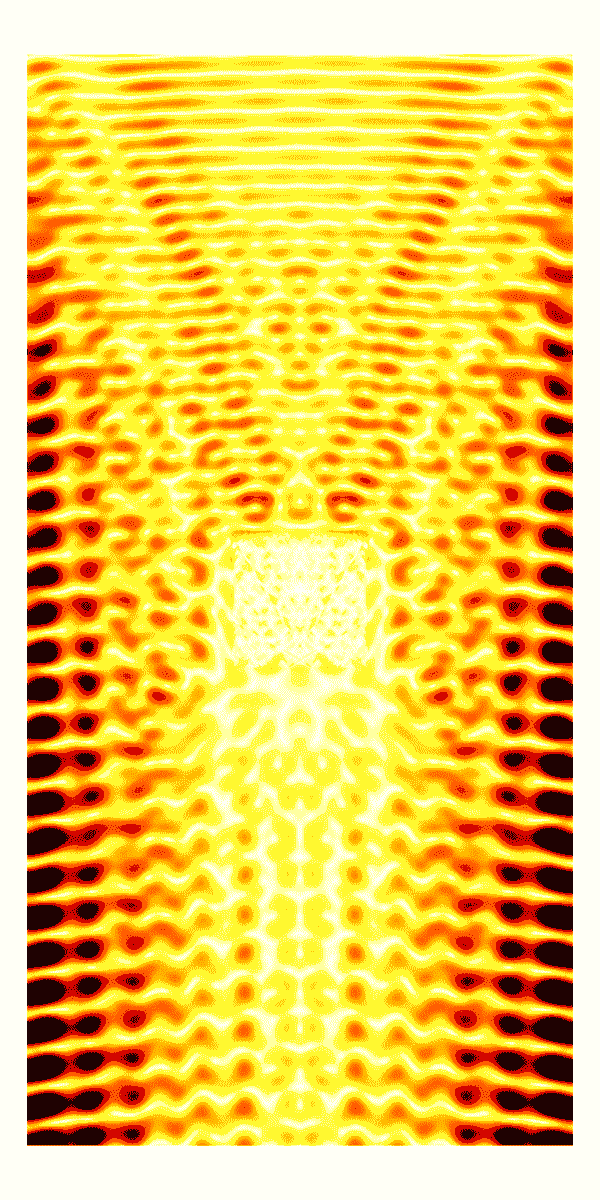}
    \end{subfigure}

	\begin{subfigure}{0.0\textwidth}
    \begin{picture}(0,0)
        \put(-6,30){\rotatebox{90}{\textbf{$10.05 \cdot 10^{6}$ rad/sec}}}
    \end{picture}
	\end{subfigure}
    \begin{subfigure}{0.14\textwidth}
        \includegraphics[width=\textwidth]{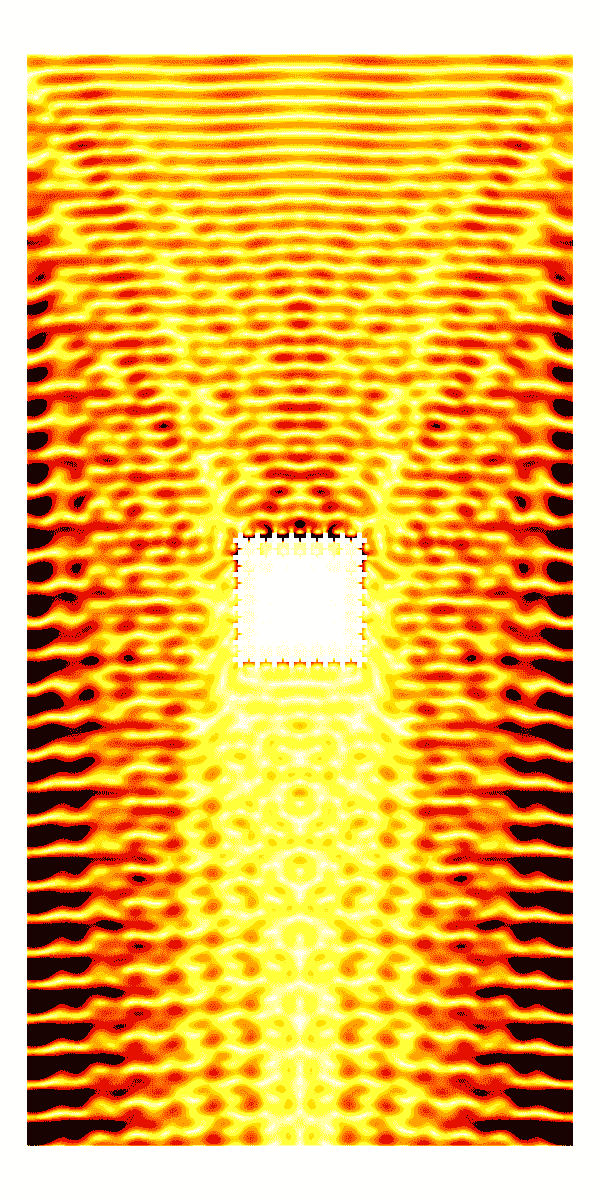}
    \end{subfigure}
    \begin{subfigure}{0.14\textwidth}
        \includegraphics[width=\textwidth]{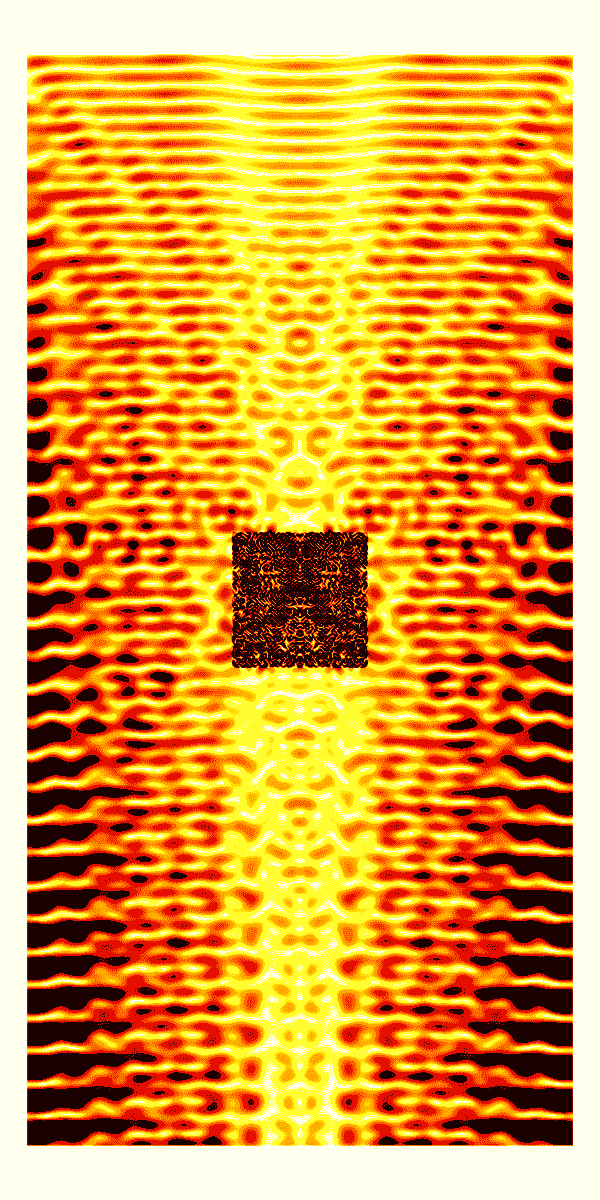}
    \end{subfigure}
    \begin{subfigure}{0.14\textwidth}
        \includegraphics[width=\textwidth]{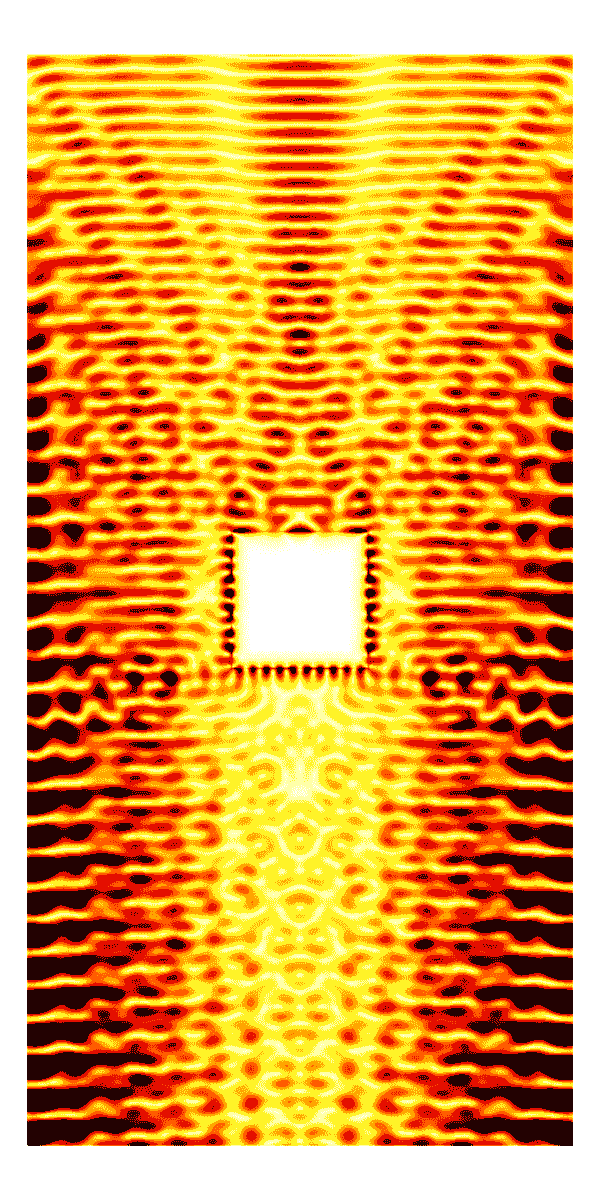}
    \end{subfigure}
    \begin{subfigure}{0.14\textwidth}
        \includegraphics[width=\textwidth]{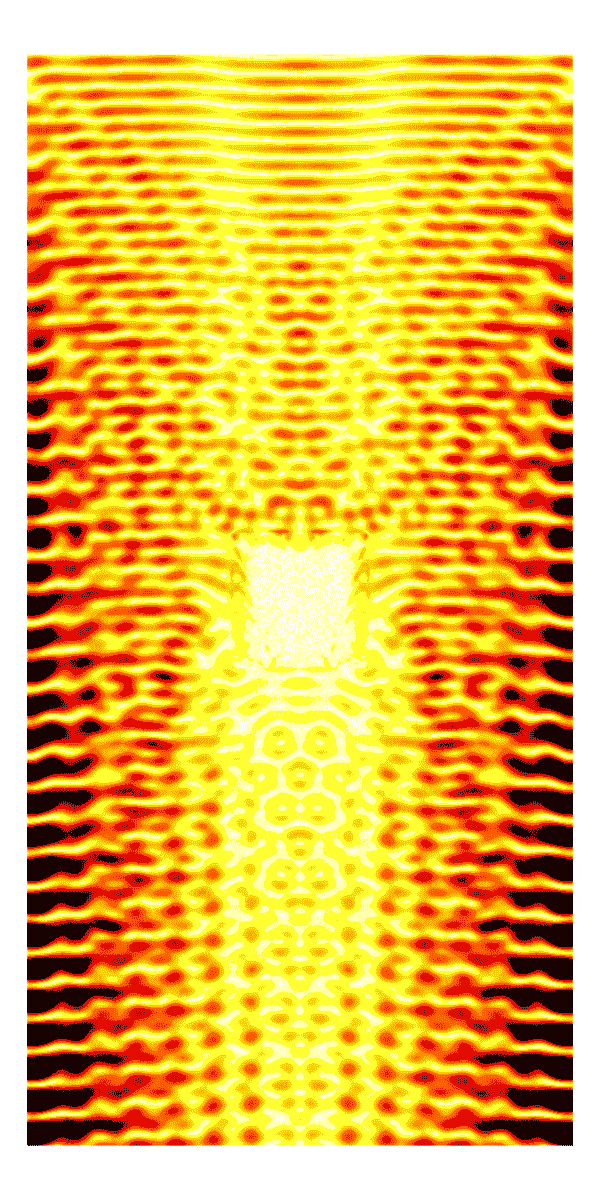}
    \end{subfigure}
    \begin{subfigure}{0.14\textwidth}
        \includegraphics[width=\textwidth]{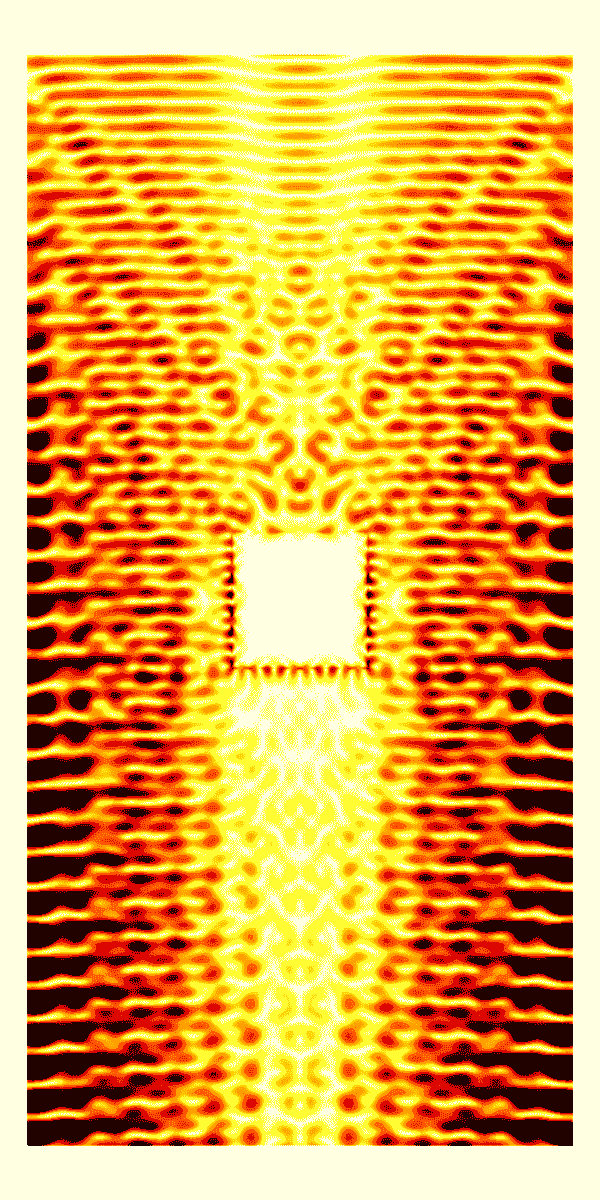}
    \end{subfigure}
    \begin{subfigure}{0.14\textwidth}
        \includegraphics[width=\textwidth]{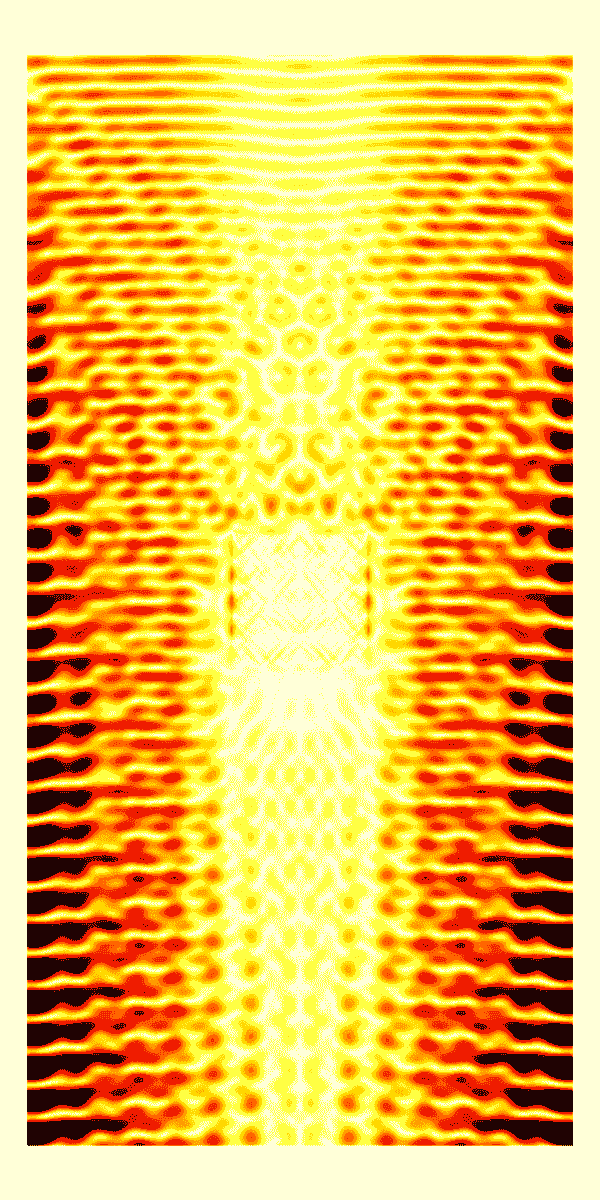}
    \end{subfigure}
    
	\begin{subfigure}{0.0\textwidth}
    \begin{picture}(0,0)
        \put(-6,30){\rotatebox{90}{\textbf{$11.31 \cdot 10^{6}$ rad/sec}}}
    \end{picture}
	\end{subfigure}
    \begin{subfigure}{0.14\textwidth}
        \includegraphics[width=\textwidth]{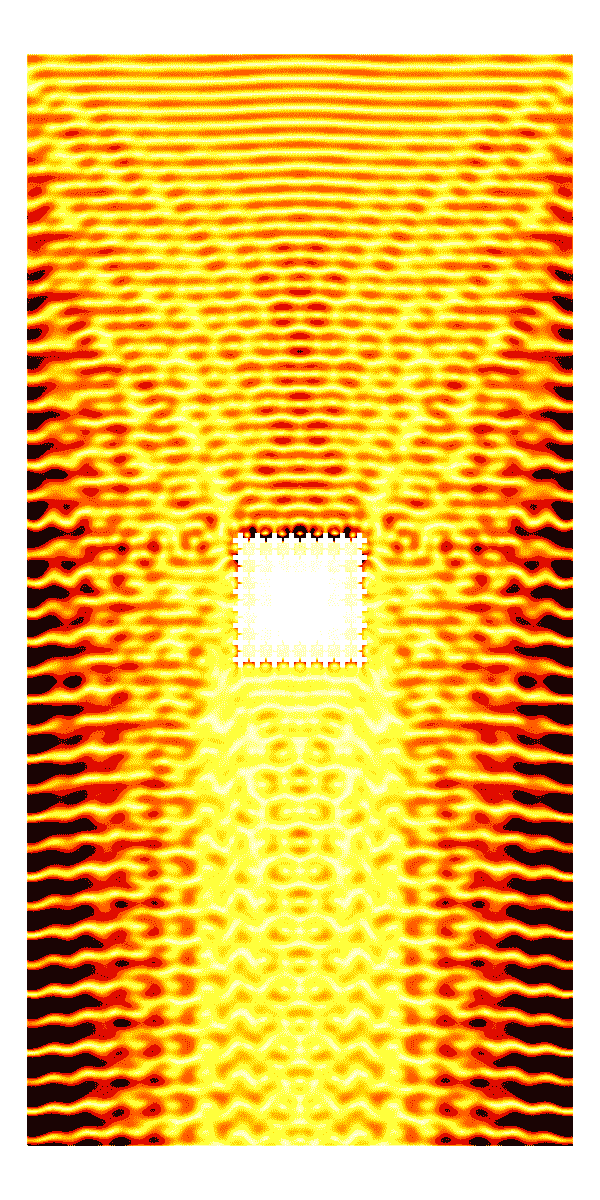}
    \end{subfigure}
    \begin{subfigure}{0.14\textwidth}
        \includegraphics[width=\textwidth]{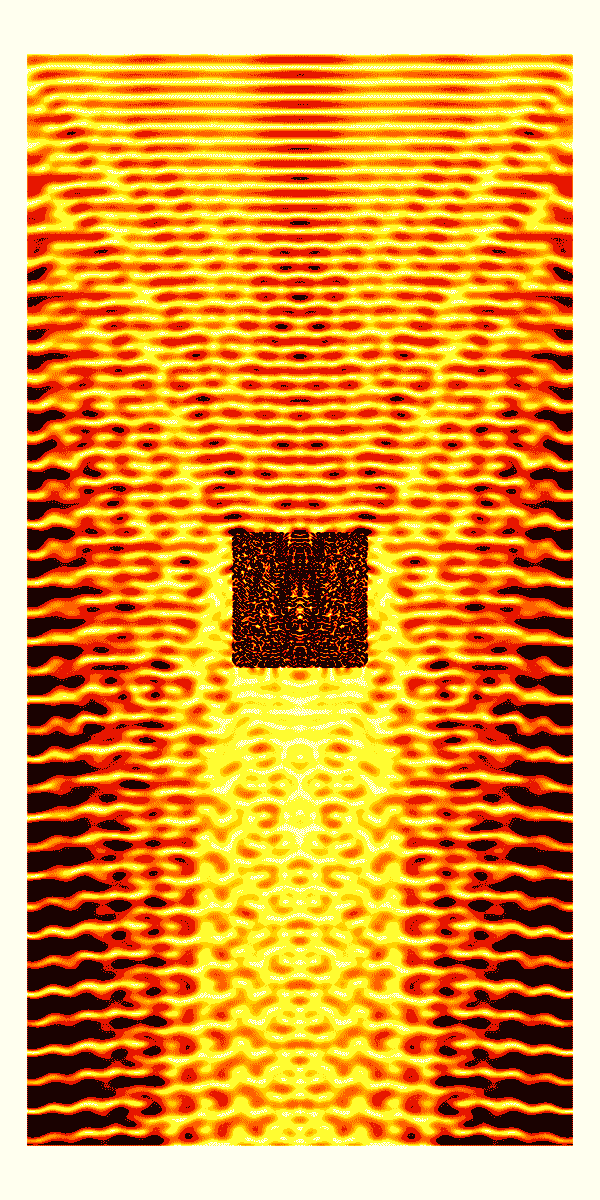}
    \end{subfigure}
    \begin{subfigure}{0.14\textwidth}
        \includegraphics[width=\textwidth]{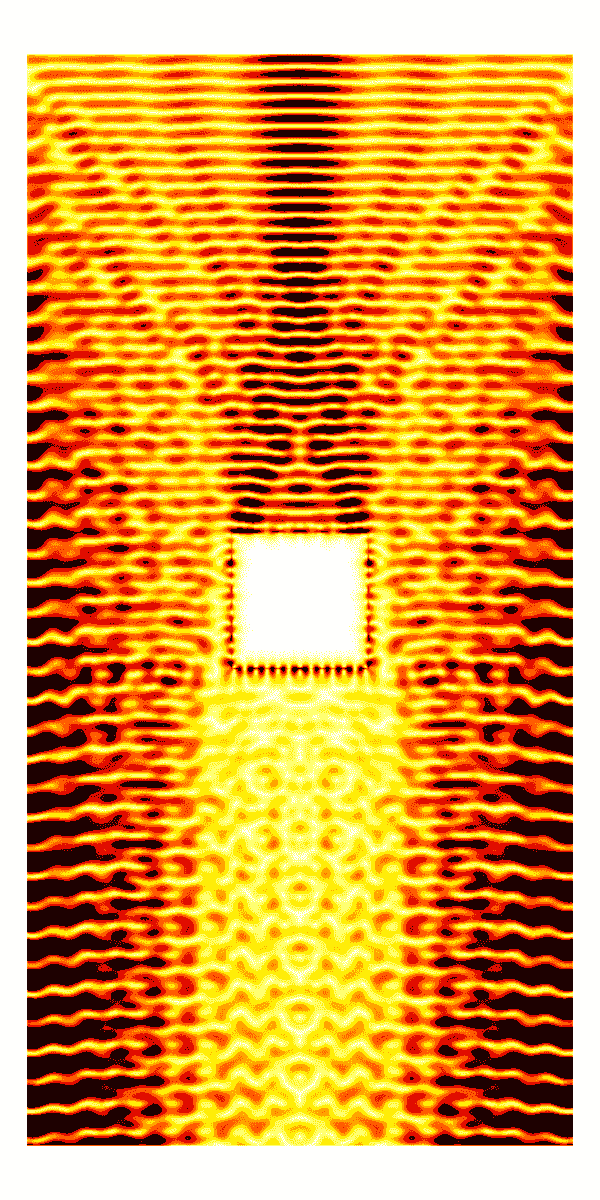}
    \end{subfigure}
    \begin{subfigure}{0.14\textwidth}
        \includegraphics[width=\textwidth]{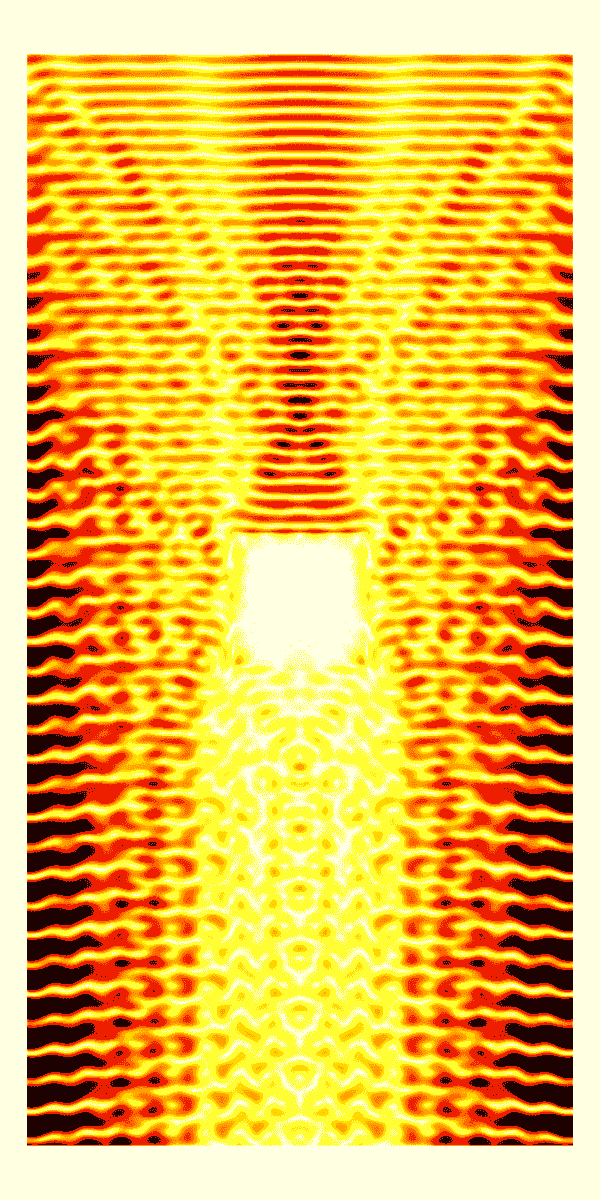}
    \end{subfigure}
    \begin{subfigure}{0.14\textwidth}
        \includegraphics[width=\textwidth]{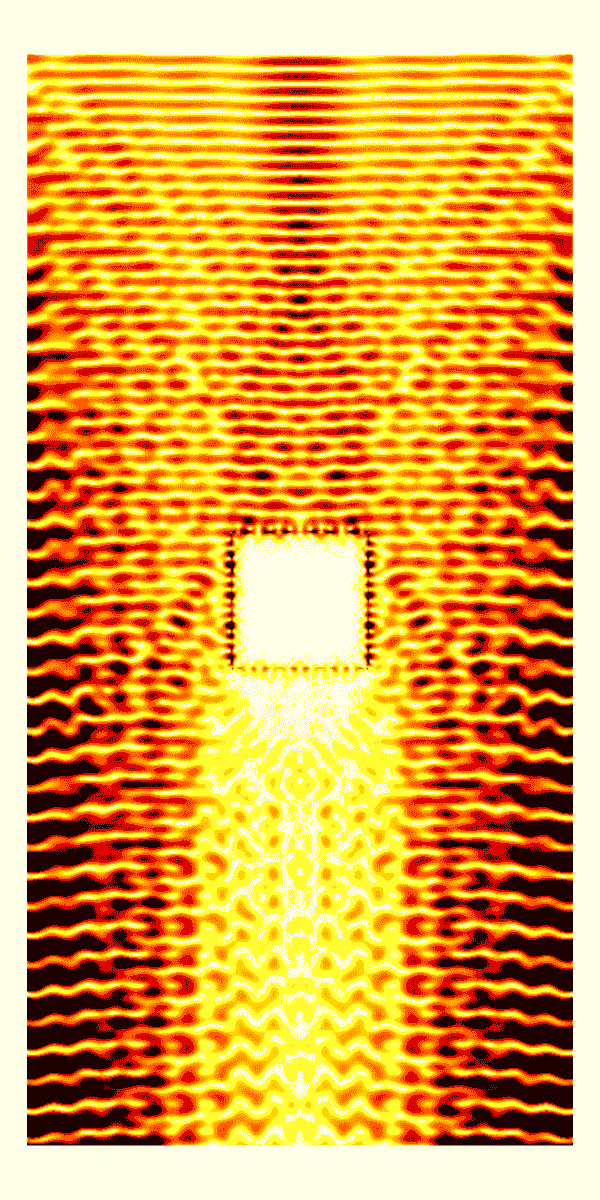}
    \end{subfigure}
    \begin{subfigure}{0.14\textwidth}
        \includegraphics[width=\textwidth]{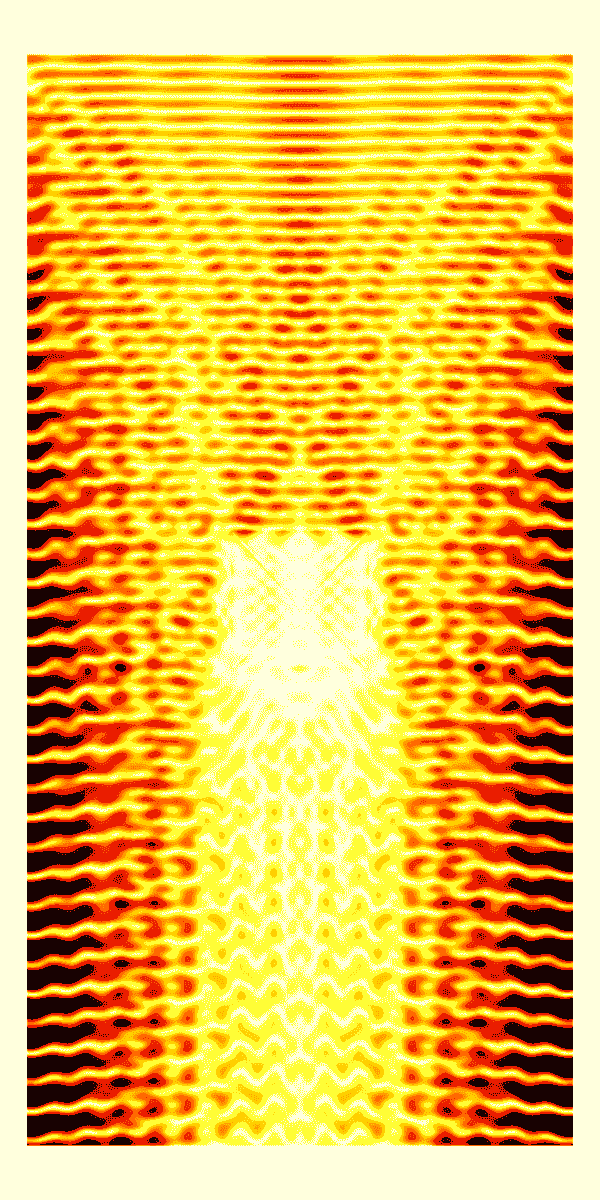}
    \end{subfigure}    
   
	\begin{subfigure}{0.0\textwidth}
    \begin{picture}(0,0)
        \put(-6,40){\rotatebox{90}{\textbf{$12.56 \cdot 10^{6}$ rad/sec}}}
    \end{picture}
	\end{subfigure}
    \begin{subfigure}{0.14\textwidth}
        \includegraphics[width=\textwidth]{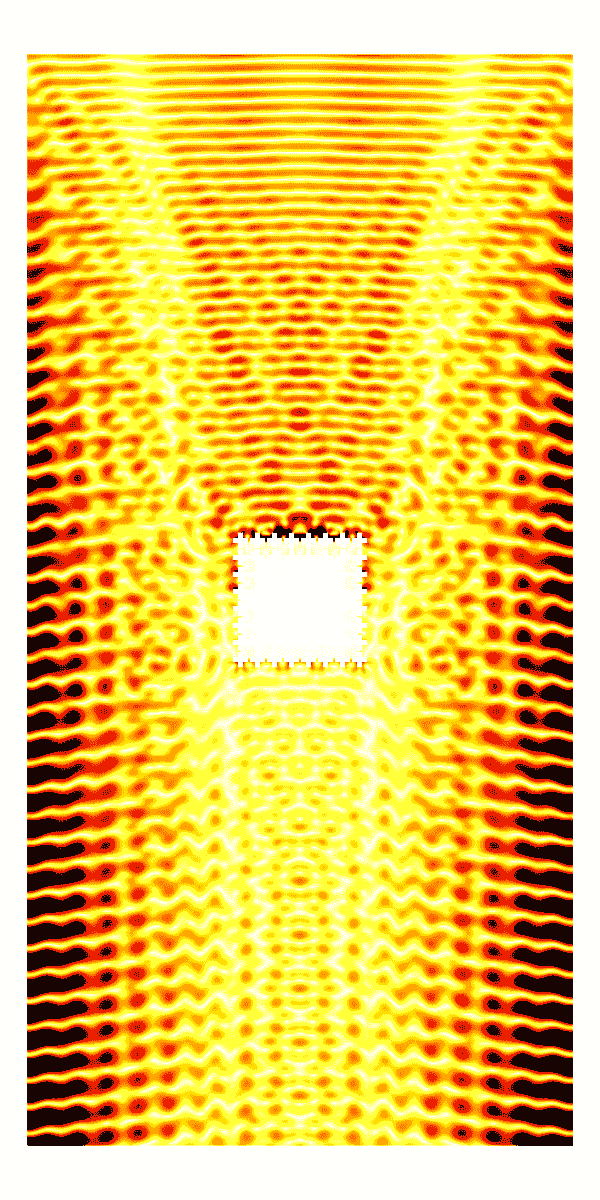}
    \caption*{\footnotesize microstructured}
    \end{subfigure}
    \begin{subfigure}{0.14\textwidth}
        \includegraphics[width=\textwidth]{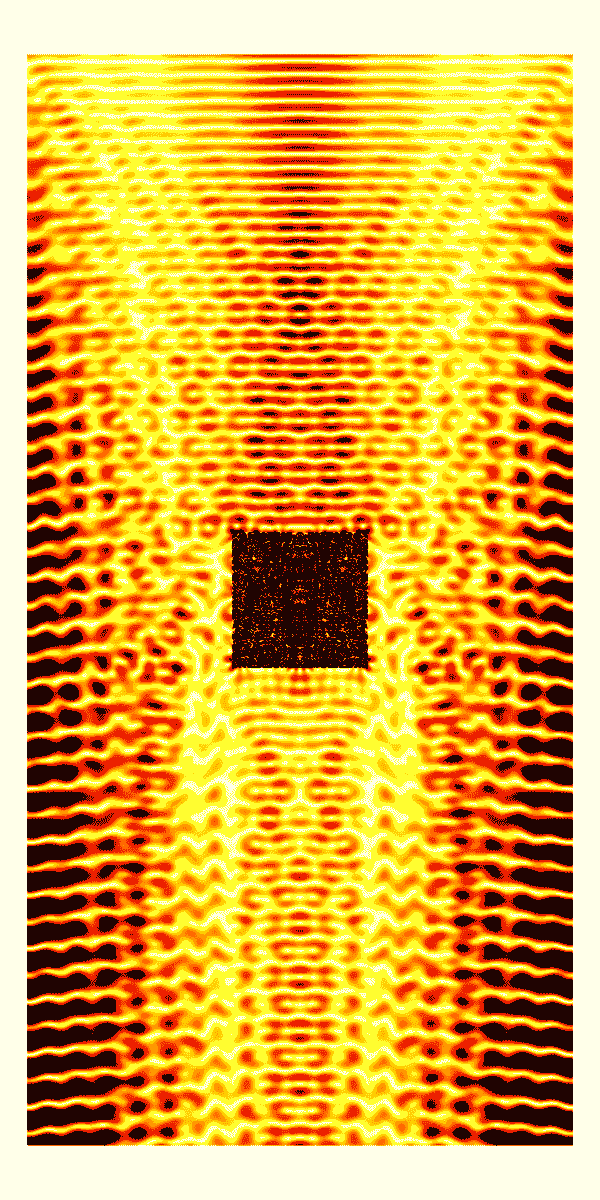}
   \caption*{\footnotesize macro-Cauchy}
    \end{subfigure}
    \begin{subfigure}{0.14\textwidth}
        \includegraphics[width=\textwidth]{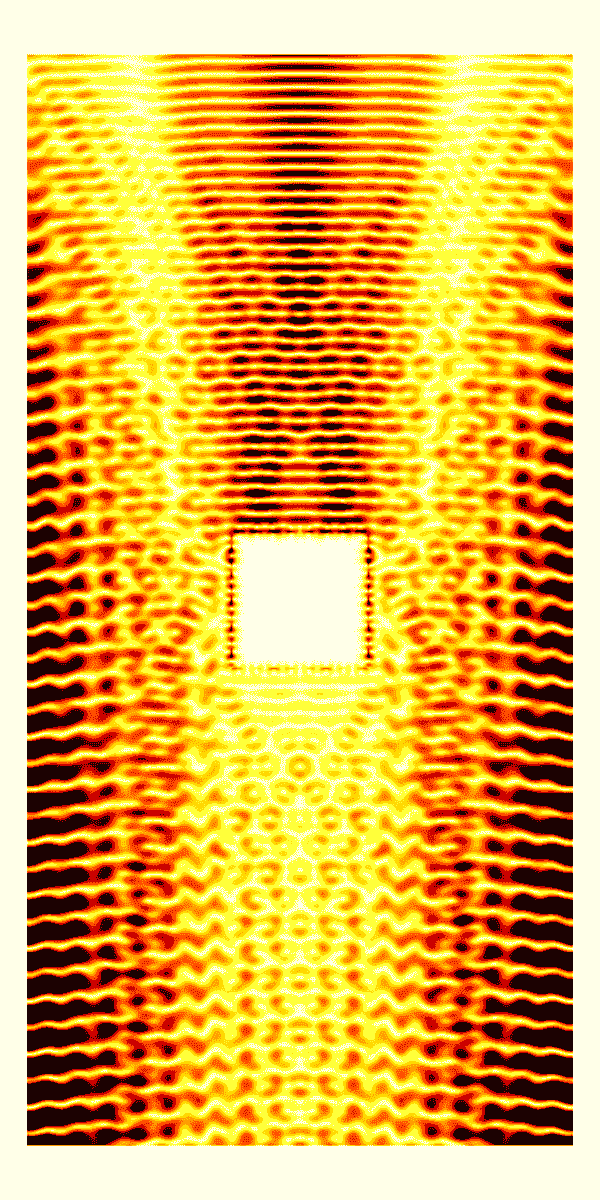}
    \caption*{\footnotesize RRMM(1)}
    \end{subfigure}
    \begin{subfigure}{0.14\textwidth}
        \includegraphics[width=\textwidth]{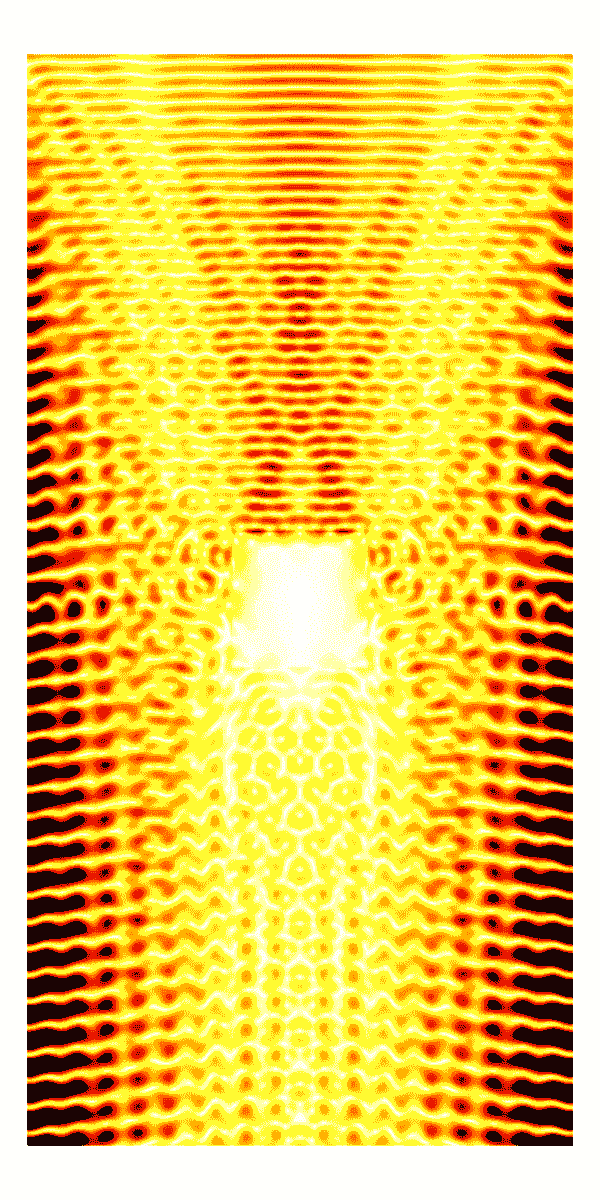}
    \caption*{\footnotesize RMM(1)}  
    \end{subfigure}
    \begin{subfigure}{0.14\textwidth}
        \includegraphics[width=\textwidth]{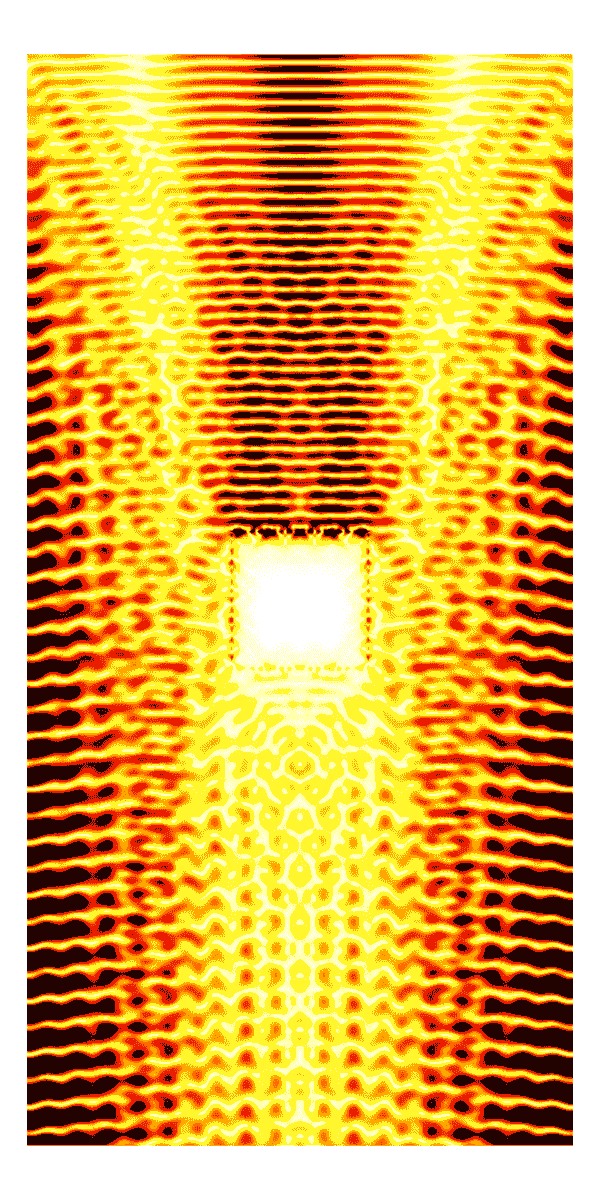}
    \caption*{\footnotesize RRMM(2)}
    \end{subfigure}
    \begin{subfigure}{0.14\textwidth}
        \includegraphics[width=\textwidth]{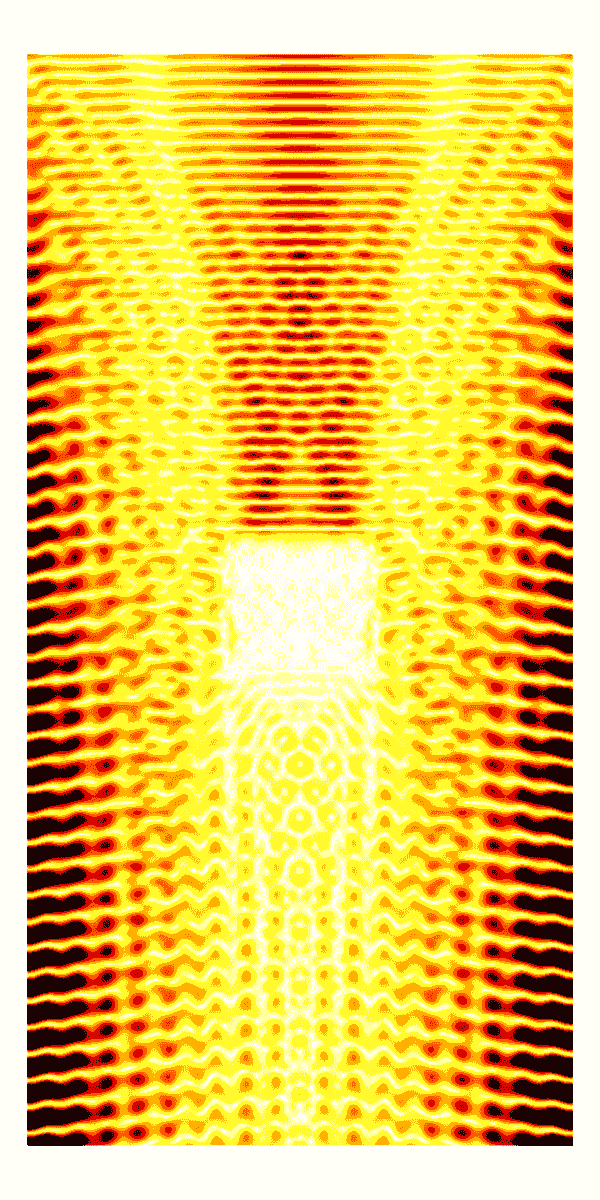}
    \caption*{\footnotesize RMM(2)}  
    \end{subfigure}
\caption{Results of finite-size scattering pattern for an incident shear wave (band-gap frequency range) with the material parameters obtained by fitting the dispersion curves in one direction at $0^\circ$ (marked as 1) and in two directions at $0^\circ$  and $45^\circ$  (marked as 2).}
\label{fig:she2}
\end{figure}

Different excitation frequencies are examined starting from a low frequency and gradually increasing until the middle of the band-gap is reached. The normalized displacement  ${\lvert u \lvert}/{u_0}$ contours of the reference full discretized microstructure, the macroscopic Cauchy continuum, the RMM and the RRMM solutions are shown in Figures \ref{fig:pre1} and  \ref{fig:pre2}  for pressure waves and  in Figures \ref{fig:she1} and  \ref{fig:she2}  for shear waves.

 \subsection{Discussion of the results}

Overall, the relaxed micromorphic model reproduces the pattern observed in the reference solution with a fully discretized  microstructure for most of the investigated frequencies in a highly satisfactory manner, for both shear and pressure waves. The classical first-order homogenization model (macro-Cauchy) fails to reproduce the observed response of the microstructured reference solution especially for higher frequencies. The capability of the relaxed micromorphic model is clearly demonstrated  qualitatively and quantitatively for both formulations with and without curvature.  

At low frequencies up to $25.13 \cdot 10^{5}$ rad/sec, the results of the RMM and the RRMM with the parameters fitted in one propagation direction show better agreement with the reference solution than those based on the parameters fitted in two propagation directions. This holds for both shear and pressure waves and can be attributed  to the poorer fitting of the acoustic branches in the two-directions fitting optimization compared to the one-direction fitting procedure.  We observe at a frequency of $50.27 \cdot 10^{5}$ rad/sec for pressure waves that the RMM and the RRMM with fitted parameters in two propagation directions behave poorly due to the imperfect fitting together with the fact that the coupling between pressure and shear waves in not activated at those low frequencies. On the other hand, we notice at a frequency of $62.83 \cdot 10^{5}$ rad/sec for pressure waves that the RMM and the RRMM with parameters fitted in one propagation direction lead to worse agreement because at this frequency the wavelength becomes comparable to the unit-cell's size thus triggering multiple reflection/transmission patterns that do not allow to ensure the separation between shear and pressure modes.

 At frequencies from $62.83 \cdot 10^{5}$ rad/sec to $12.56 \cdot 10^{6}$ rad/sec for pressure waves, the RMM formulation with curvature terms ($\Curl \Bdis$ and $\Curl \ddot\Bdis $) and parameters fitted in two directions shows the best agreement with the reference microstructured solution. For shear waves at frequencies from $37.7 \cdot 10^{5}$ rad/sec to $75.4 \cdot 10^{5}$ rad/sec and from $11.31 \cdot 10^{6}$ rad/sec to $12.56 \cdot 10^{6}$ rad/sec, the RMM formulation with the curvature terms and parameters fitted in two directions also shows the best agreement with reference solution.  For frequencies from $87.96 \cdot 10^{5}$ rad/sec to $10.05 \cdot 10^{6}$ rad/sec  with shear waves, the pattern obtained in the reference solution cannot be produced by either the RMM and the RRMM. This behavior may be related to the less accurate fitting near the lower edge of the band-gap.

It is worth mentioning that the  dispersion curves fitting procedure is carried out over the irreducible Brillouin zone which is bounded by the periodicity limits defined by wavenumbers of $\frac{\pi}{l}$ and  $\frac{\sqrt{2} \pi}{l}$ at incidence angles of $0^\circ$ and $45^\circ$, respectively. While this applies for the microstructured unit-cell, the dispersion curves of the RMM/RRMM, as a homogeneous continuum, do not exhibit such periodicity limits and the analytical dispersion curves of the RMM/RRMM interfere the band-gap for wavenumbers beyond the irreducible Brillouin zone leading to some non-vanishing responses for RMM/RRMM inside the band-gap reflected in some spurious pattern inside RMM/RRMM domain for frequencies from $75.4 \cdot 10^{5}$ rad/sec to $87.96 \cdot 10^{5}$ rad/sec for both shear and pressure waves. Whereas the fitting procedure should consider that the asymptotes of the analytical dispersion curves of the RMM/RRMM should not intersect the band-gap is an open question for future works. 

For all tested frequencies and for both pressure and shear waves, accounting for the curvature in the RMM has a major positive effect and leads to better agreement with the reference solution compared to the formulation without curvature (RRMM) when the parameters are  fitted using the same procedure (i.e. in one or two propagation directions). The main difference between the two formulations is that the amplitudes of the reflected and scattered waves are higher in the RRMM than in the reference microstructured solution whereas the RMM generally does not exhibit this behavior. A physical interpretation is that activating the curvature term damps this excess energy.  The strain energy in the RRMM operates near the soft macroscopic bound. In contrast, the RMM exhibits a stiffer response where the elasticity tensor $\Cmicro$ contributes directly to the overall static energy. In other words, indeed accounting for the static size-effects through incorporating the curvature leads to a stiffer response and a better agreement with the reference microstructured solution. This is more pronounced on high frequencies when the wavelength is comparable to the unit-cell size. 

We can summarize that considering the curvature terms in the RMM positively affects the model's predictability especially at higher frequencies where the wavelength becomes comparable to the unit-cell's size when the scale-separation hypothesis is broken. We can point out that such differences at high frequencies can be related to the higher interfacial stiffness associated with the generated tractions at the specimen's interfaces. Accounting for these enhanced tractions at the interfaces improves the quality of the solution. Moreover, the constitutive forces of such tractions is intrinsically related to the bulk behavior and cannot account for independent micro-inertial effects occurring at the interfaces. This also implies that the RMM (or any generalized continuum)  cannot account for different metamaterial's cuts at the considered interfaces. To address these issues, we will develop in future work a generalization of the concept of "interface forces" introduced in \cite{DemVosMad:2025:eif,LeoFelGiaJen:2024:esf} to account for such interfacial effects associated with different cuts of the metamaterial. 

\section{Conclusions} 
\label{sec:con}

We introduced the relaxed micromorphic model together with its associated variational formulation and the resulting strong form. A two-stage optimization procedure for the fitting of the RMM elastic coefficients is presented. The static parameters are identified via size-effect fitting for linear deformation modes. For this purpose a least-squares minimization procedure is developed, which employs a cost function built on the total energy  difference with a line search algorithm, while  a feed-forward fully connected neural network is used as an accelerator to reduce the computation time by providing an initial prediction. Once the static parameters are fitted, the dynamic parameters are then determined via a least squares fitting procedure of the dispersion curves either using one or two propagation directions giving rise to two sets of parameters. For each of the two sets of parameters, two cases were considered: with and without curvature terms. We studied a finite-size example and examined the results of fitting in one and  in two directions with and without curvature terms. The results obtained when accounting for curvature show better agreement with the reference microstructured solution. This indicates that the static size-effect phenomenon interacts with the dynamic response of mechanical metamaterials especially when the wavelength becomes comparable to the unit-cell size (i.e. dispersive and band-gap frequency ranges). This interaction is reflected in a stiffer response and increased damping induced by size-effects.  This can also be interpreted that including the curvature modifies the interface's effective properties giving rise to more realistic reflection/transmission patterns. Future work will focus on two points: the first point is working with metamaterials for which the RMM produces a better fitting, and the second point is further exploring the interfacial effects associated with different metamaterial's cuts and how to reproduce these different responses via enhanced interfacial traction or inertia seeking to improve the qualitative and quantitative model's performance.

\section*{Acknowledgement} 
Mohammad Sarhil thanks Robert Martin from the University of Duisburg-Essen for the helpful support related to the machine learning. 
Mohammad Sarhil, Leonardo Andres Perez Ramirez, Max Jendrik Voss and Angela Madeo acknowledge support from the European Commission through the funding of the ERC Consolidator Grant META-LEGO, no. 101001759. 

\FloatBarrier

\bibliographystyle{abbrvnat}
{\footnotesize
\bibliography{references.bib}
}

\end{document}